\documentclass{amsart}
\usepackage[english]{babel}
\usepackage{graphicx}
\usepackage{amssymb}
\usepackage{amsmath}
\usepackage[foot]{amsaddr}
\usepackage{amsfonts}
\usepackage{bbm}
\usepackage{xcolor}
\usepackage{float}
\usepackage{algorithm}
\usepackage{algpseudocode}

\usepackage{amsmath}
\usepackage{amssymb}
\usepackage{mathrsfs}
\usepackage{color}
\usepackage{graphicx}
\usepackage[percent]{overpic}
\usepackage{url}
\usepackage{enumerate}
\usepackage{overpic}
\usepackage{amsthm}
\usepackage{moreverb}
\usepackage[pdftex,colorlinks,bookmarksopen,bookmarksnumbered,citecolor=red,urlcolor=red]{hyperref}

\def\f{{\bm f}}
\def\F{{\bm F}}

\def\n{{\bm n}}
\def\u{{\bm u}}

\def\x{{\bm x}}

\def\0{\boldsymbol{0}}

\def\dt{\partial_t}

\def\cl {\nonumber \\}
\def\el {\nonumber }

\newtheorem{rem}{Remark}[section]

\newcommand{\bm}[1]{\mbox{\boldmath{$#1$}}}
\def\div{\nabla\cdot}

\newcommand{\rev}[2][red]{{\textcolor{#1}{#2}}}

\usepackage{tikz}
\usepackage[top=1in, bottom=1.25in, left=1.25in, right=1.25in]{geometry}

\graphicspath{{images/}}

\begin{document}

\title[]{A POD-Galerkin reduced order model for Navier-Stokes equations in stream function-vorticity formulation} % in a Finite Volume environment

\author{Michele Girfoglio$^1$, Annalisa Quaini$^2$ and Gianluigi Rozza$^1$}
\address{$^1$ mathLab, Mathematics Area, SISSA, via Bonomea 265, I-34136 Trieste, Italy}
\address{$^2$ Department of Mathematics, University of Houston, Houston TX 77204, USA}

\begin{abstract}
We develop a Proper Orthogonal Decomposition (POD)-Galerkin based Reduced Order Model (ROM) for the efficient numerical simulation 
of the parametric Navier-Stokes equations in the stream function-vorticity formulation. Unlike previous works, we choose different reduced coefficients for the vorticity and stream function fields. In addition, for parametric studies we use a global POD basis space obtained from a database of time dependent full order snapshots related to sample points in the parameter space. %We also provide a validation of the FOM solver against the solutions provided by the Navier Stoks equations in velocity-pressure formulation.  
We test the performance of our ROM strategy with the vortex merger benchmark. Accuracy and efficiency are assessed for
both time reconstruction and physical parametrization. % rWe consider two test cases: a channel flow over a step at $Re = 600$ and a lid driven cavity flow at $Re = 1000$ % implementation of the Leray
%model that combines a three-step algorithm called Evolve-Filter-Relax (EFR) with a
%computationally efficient finite volume method. The main novelty of our ROM lies in
%We test the performance of our ROM
%approach on two benchmark problems: 2D and 3D unsteady flow past a cylinder at
%Reynolds number 0 ≤ Re ≤ 100. The accuracy of the ROM is assessed against results
%obtained with the full order model for velocity, pressure fields and time evolution of the
%aerodynamics coefficients.

%we build a computational pipeline for the assessment of the Navier-Stokes equations in a stream function-vorticity formulation in a Finite Volume environment. We address the matter  both at full order and reduced order framework. The high fidelity solver is validated against the solutions provided by the Navier-Stokes equations in velocity-pressure formulation. On the other hand, the accuracy of the reduced order model is assessed with respect to the results obtained with the full order model. We consider two test cases: a channel flow over a step at $Re = 600$ and a lid driven cavity flow at $Re = 1000$.
\end{abstract}

\maketitle


\section{Introduction}\label{sec:intro}

The formulation of the 2D Navier--Stokes equations in 
terms of stream function and vorticity represents an
attractive alternative to the model in primitive variables for two main reasons: (i) there are only two scalar unknowns and (ii) the divergence free constraint for the velocity 
is automatically satisfied by the definition of the stream function. Computational studies on this formulation can be found
in, e.g., \cite{Behr1991,Lequeurre2020,Minev2015,Sousa2005,Tezduyar1990}. %\michele{Ovviamente Anna questi riferimenti sono lontani dall'essere esaustivi, giusto per mettere qualcosa su psi omega a livello FOM, anzi forse ne sono anche troppi, ne taglierei qualcuno in caso!} \anna{Ne ho tolto uno :) Usano tutti elementi finiti?}

While there exists an abundance of literature on Reduced Order Models (ROMs) for the Navier--Stokes equations formulated in primitive variables starting from different Full Order Methods (FOMs), e.g. Finite Element methods \cite{Veroy2005, Burkardt2006, Ballarin2014} or Finite Volume methods \cite{Lorenzi2016, Stabile2018, Girfoglio_JCP}, 
%\michele{qui ho evitato di citare tutti i nostri lavori perchè intenderei inserire qui soltanto i rif. principali, piu' di ampio respiro; tutti gli altri nostri lavori li ho inseriti nelle conclusioni quando menziono il filtro non lineare}.
it is only relatively recently
that ROMs have been applied to the stream function-vorticity formulation \cite{Ahmed2019,Ahmed2020,Dumon2010,Pawar2020-2,Pawar2020,Pawar2021,San2015}.

%\michele{qui invece ho cercato di essere esaustivo nelle citazioni perche' importanti per enfatizzare le novita', che poi sono più o meno tutte del gruppo di Iliescu}.
%\anna{Si', assicuriamoci che siano tutti. Anche per queste referenze e' utile riportare se usano tutti FE, perche' questo sarebbe un ulteriore elemento di distinzione per noi.} \michele{Ho visto che quelli che fanno PGD, riferimento 9, usano FV...tutti gli altri sono Iliescu style quindi suppongo FE...a questo eviterei di menzionare FV come elemento di novità, sei d'accordo? Per quanto riguarda invece l'accoppiamento QGE + filtro siamo i primi a usare FV e in quel caso ovviamente lo evidenzieremo :)}
In this paper, we develop a POD–Galerkin ROM for the stream function-vorticity formulation.  The main building blocks of our approach are:
\begin{itemize}
\item[-] the collection of a database of simulations using a computationally efficient finite volume method; 
\item[-] the extraction of the most energetic modes representing the system dynamics through Proper Orthogonal Decomposition (POD);
\item[-] a Galerkin projection on the space spanned by these most energetic modes for the computation of stream function and vorticity reduced coefficients.
\end{itemize}

Two are the main novelties of our approach. First,  
unlike previous works \cite{Ahmed2020,Ahmed2019,Pawar2021,Pawar2020-2,Pawar2020,San2015} we consider different coefficients for the approximation of the vorticity and stream function fields. This choice leads to two important consequences: (i)  the stream function basis functions do not depend on the particular vorticity basis functions, 
but are instead computed directly from the stream function high-fidelity solutions during the offline phase; (ii) the reduced spaces for the stream function and vorticity can have different dimensions. The second novelty pertains the parametric study with respect to two crucial model parameters (Reynolds number and strength of the forcing term) for which we use a global POD basis space computed by time dependent FOM snapshots associated to sample points in the parameter space. This is a difference with respect to \cite{Pawar2020, Pawar2020-2, Ahmed2020, Ahmed2019, Pawar2021}, where a POD basis is computed for each parameter in the training set and the the basis functions for new parameter
values are found via interpolation of the basis functions associated to the training set.

The work in this paper represents an intermediate step towards the development of new FOM and ROM approaches for the quasi-geostrophic equations that are usually written in terms of stream function and (potential) vorticity. See \cite{QGE_review} for a recent review.  

All the FOM simulations presented in this work have been performed with OpenFOAM\textsuperscript{\textregistered} \cite{Weller1998}, an open source Finite Volume C++ library widely used by commercial and academic organizations. 
For the Navier--Stokes equations in primitive variables,
OpenFOAM features several partitioned algorithms (PISO \cite{PISO}, PIMPLE \cite{PIMPLE} and SIMPLE \cite{SIMPLE})
based on the Chorin-Temam projection scheme \cite{Temam2001}. 
To the best of our knowledge, no solver for the stream function-vorticity formulation has been shared with the large 
OpenFOAM\textsuperscript{\textregistered} community. Thus, we have
implemented such solver at the FOM level. %\anna{Michele, tu contribuisci a OPENFOAM? Intendo se il tuo codice per la function-vorticity formulation puo' essere sottomesso per essere incluso in OPENFOAM.} 
The ROM computations have been carried out
with ITHACA-FV \cite{Stabile2018}, an in-house implementation of several ROM techniques within OpenFOAM\textsuperscript{\textregistered}.
%An important outcome of this work is that the code created for it is incorporated in an open-source library \cite{} and therefore is readily shared with the community.
%The use of a ROM differential filter, i.e. a ROM spatial filter that uses an
%explicit lengthscale (i.e., the filtering radius). This approach has not been
%explored in depth by the ROM stabilization and closure community.
%- The computation of the pressure field at ROM level (through the pressure

%We remark that this work is an intermediate step towards future works focused on the development of ROMs for Quasi-Geostrophic Equations (QGE) (see \cite{QGE_review} for a recent review) that are usually expressed in terms of stream function and (potential) vorticity.

The rest of this paper is organized as follows. In Sec.~\ref{sec:FOM}, we describe the full order model and the numerical method we use for its time and space discretization. Sec.~\ref{sec:ROM} presents the reduced order model. The numerical experiments are reported in Sec.~\ref{sec:num_res}. Finally, conclusions and future perspectives are provided in Sec.~\ref{sec:conc}.

\section{The Full Order Model}\label{sec:FOM}
%...In Sec 3.1 we introduce the procedure we use to construct
%a POD-Galerkin ROM and 
%We are going to 
%In Sec. \ref{sec:uvp}, we briefly introduce the Navier-Stokes equations in stream function-vorticity formulation.present the strategy we choose for pres-
%sure stabilization at reduced order level. Finally, Sec. 3.3 describes the lifting
%function method we apply to enforce non-homogeneous Dirichlet boundary con-
%ditions for the velocity field at the reduced order level. The FOM computations
%are carried out using Open[47], an in-house open source C++ library.
%In Sec. \ref{sec:psizeta}
%\subsection{The Navier-Stokes equations in velocity-pressure formulation}\label{sec:uvp}
\subsection{The Navier-Stokes equations in stream function-vorticity formulation}\label{sec:psizeta}
We consider the motion of a two-dimensional incompressible, viscous fluid in a fixed domain $\Omega \subset \mathbb{R}^2$ 
over a time interval of interest $(t_0, T)$. The flow is described by the incompressible Navier-Stokes equations:
\begin{align}
\dt \u + \div \left(\u \otimes \u\right) - \dfrac{1}{Re} \Delta \u + \nabla p & = \f \quad \mbox{ in }\Omega \times (t_0,T),\label{eq:ns-mom}\\
\div \u & = 0\quad\, \mbox{ in }\Omega \times(t_0,T),\label{eq:ns-mass}
\end{align}
where eq.~(\ref{eq:ns-mom}) states the conservation of linear momentum and eq.~(\ref{eq:ns-mass}) represents the conservation of mass. 
Here, $\u(x,y,t) = (u(x,y,t), v(x,y,t),0)$ is the fluid velocity, $\dt$ denotes the time derivative, $p(x,y,t)$ is the pressure and $Re$ is the Reynolds number. In \eqref{eq:ns-mom}, we take into account possible body forces $\f (x,y,t)$. We focus on forcing terms that can be expressed as product of two functions:
one function that depends only on space and the other that depends only on time, i.e.~$\f (x,y,t)= f_2(t) \f_1 (x,y)$. See Remark \ref{rem:rem1} for more details about this choice. 
%\anna{Michele, sia la funzione in spazio che quella in tempo sono vettoriali? In caso, tra le due non ci dovrebbe essere un dot product, altrimenti la forzante diventa uno scalare.} \michele{Credo che la scelta più consistente sia considerare quella time dependent scalare e vettoriale quella space dependent, che ne dici? Di fatto, siccome lavoriamo direttamente su $\psi$-$\omega$ il problema a livello pratico non si pone mai perchè ragioniamo direttamente in termini scalari}
%\anna{Ho messo la funzione in tempo scalare perche'sono d'accordo che sia il modo piu' semplice di risolvere la questione.} \michele{Perfetto!}
%, $\u_{in}$, $\u_w$ and $\u_0$ are 
%given.

Let $\partial_x$ and $\partial_y$ denote the derivative with respect to the $x$ and $y$ spatial coordinate, respectively. 
By applying the curl operator $\nabla \times$ to eq.~(\ref{eq:ns-mom}), we obtain the governing equation for the vorticity field $\bm{\omega}(x,y,t) = \nabla \times \u = (0, 0, \omega) = (0, 0, \partial_x v - \partial_y u)$
\begin{align*}
\dt \omega+ \div \left(\u \omega \right) - \dfrac{1}{Re} \Delta \omega & = F \quad \mbox{ in }\Omega \times (t_0,T),
\end{align*}
where $\F = (0, 0, F) = \nabla \times \f$. %denotes, with an abuse of notation, the curl of body forces (if any). 
The incompressibility constraint (\ref{eq:ns-mass}) leads to the introduction of the stream function $\bm{\psi}(x,y,t) = (0,0,\psi)$ such that
$\u = \nabla \times \bm{\psi}$, or, equivalently, $(\partial_y \psi, -\partial_x \psi) = (u, v)$ .
The stream function $\psi$ and vorticity $\omega$ are linked by a Poisson equation
%u = \dfrac{\partial \psi}{\partial y}, \quad v = -\dfrac{\partial \psi}{\partial x}. \label{eq:ns-mom-upsi}
%\end{align}
%Equation (\ref{eq:ns-mom-upsi}) can also expressed in the following way
%\begin{align}
%\u = \nabla \times \bm{psi} \label{eq:ns-mom-upsi2}
\begin{align*}
-\Delta \psi = \omega  \quad \mbox{ in }\Omega \times (t_0,T). 
\end{align*}

To close problem \eqref{eq:ns-mom}-\eqref{eq:ns-mass}, we need to provide initial data $\u(x,y,t_0) = \u_0$ and enforce proper boundary conditions. In this work, we consider the following slip condition on the entire boundary:
%slip boundary conditions for $\u$ and homogeneous Neumann boundary conditions for $p$:
\begin{align}
\u \cdot \n & = 0 \quad\ \mbox{and} \quad %\left(\nabla \u \cdot \n\right)\cdot \bm{t} = 0
\partial_n \left(\u \cdot \bm{t} \right) = 0\quad \mbox{ on } \partial\Omega \times(t_0,T)\label{eq:bc-ns-d}, %\\
%\nabla p \cdot \n & = 0\quad\ \mbox{ on } \partial\Omega \times(t_0,T).\label{eq:bc-ns-d2}%\\
%\left(\dfrac{1}{Re} \nabla \u - p\right) \n& = \bm{0} \quad \quad\mbox{ on } \partial\Omega_{out} \times(t_0,T),\label{eq:bc-ns-n}
\end{align}
where $\n$ is the outward unit normal and $\bm{t}$ the unit tangent vector to $\partial\Omega$.
In terms of $\psi$ and $\omega$, we consider initial data $\omega(x,y,t_0) = \omega_0$ and express \eqref{eq:bc-ns-d} as $\psi = \partial_n \omega = 0$. 
%\anna{E con l'eq.~\eqref{eq:bc-ns-d2} cosa si fa?} \textcolor{red}{guarda di fatto non la uso, probabilmente è collegato al fatto che il contributo di pressione scompare nella formulazione $\psi$-$\omega$} \anna{Allora mi chiedo se ci sia bisogno di riportarla. Di fatto \eqref{eq:bc-ns-d}
%e' sufficiente come BC per \eqref{eq:ns-mom}-\eqref{eq:ns-mass}.} 
%\michele{hai ragionissima! La BC su p ci serve solo numericamente quando sostituiamo la continuità con l'eq. di Poisson per la pressione (che magari riporterò per completessa nel paragrafo dei risultati, quando confrontiamo le due formulazione di NSE). Metto a posto e commento la nostra conversazione! :D}
%Concerning %Concerning the boundary conditions, for NSE, we impose slip boundary conditionsin this work we impose consider no slip boundary condition
%for the inflow, one could assume stream function and vorticity to be know and set $(\partial_y \psi, -\partial_x \psi) = (u_{in}, v_{in})$ and $\zeta = \partial_x v_{in} - \partial_y u_{in}$, respectively, as well as for the outflow, normal derivative vanished, i.e. $\partial_n \zeta = \partial_n \psi = 0$ METTERE REFERENZE. On the solid walls the stream function value is set as an arbitrary constant, $\psi = c$, whilst we enforce the kinematic condition for the vorticity, SCRIVERE FORMULA E REFEREZE (see, e.g., ).

Summarizing, the Navier-Stokes equations in stream function-vorticity formulation, which represent our full order model, are given by
\begin{align}
\dt \omega + \div \left(\left(\nabla \times \bm{\psi}\right) \omega \right) - \dfrac{1}{Re} \Delta \omega & = F \quad \mbox{ in }\Omega \times (t_0,T),\label{eq:psizeta-mom}\\
-\Delta \psi & = \omega \quad \mbox{ in }\Omega \times (t_0,T),\label{eq:psizeta-mass}
\end{align}
endowed with boundary conditions 
\begin{align}
\psi & = 0\quad\ \mbox{ on } \partial\Omega \times(t_0,T),\label{eq:bc-pz-d1}\\
\partial_n \omega & = 0 \quad\ \mbox{ on } \partial\Omega \times(t_0,T),\label{eq:bc-pz-d22}
%\partial_n \psi  = 0, \quad \partial_n \zeta & = 0 \quad \quad\mbox{ on } \partial\Omega_{out} \times(t_0,T).\label{eq:bc-pz-n}
\end{align}
and initial data $\omega(x,y,t_0) = \omega_0$.
%\anna{Forse dobbiamo dire qualcosa su \eqref{eq:bc-pz-d1}
%perche' sembra che venga fuori dal niente.}
%\textcolor{red}{Possiamo dire che la 7) in buona sostanza ci dice che il bordo del dominio è una linea di corrente}

% and  The situation is more delicate with vorticity boundary conditions
%on Γw and several suggestions can be found in the literature. One common choice, see, e.g.,
%[13, 24, 40, 42, 43, 44], is the kinematic condition% on Γin, and
%letting the normal vorticity derivative vanish is a reasonable outflow boundary condition [13,
%26]: (∇w)n = 0 on Γout. The situation is more delicate with vorticity boundary conditions
%on Γw and several suggestions can be found in the literature. One common choice, see, e.g.,
%[13, 24, 40, 42, 43, 44], is the kinematic condition
 %where without inow or outow in the cavity ows we consider, the stream function is zeroon the boundary,|@=0. The boundary conditions (2) translate into boundary conditions forthe stream functio
%In order to characterize the flow regime under consideration, we define the Reynolds number as
%\begin{equation}\label{eq:re}
%Re = \frac{U L}{\nu},
%\end{equation}
%where $U$ and $L$ are characteristic macroscopic velocity and length, respectively. For an internal flow in a cylindrical pipe, $U$ is the mean sectional velocity and $L$ is the diameter. %For large Reynolds numbers, inertial forces are dominant over viscous forces and vice versa. 

\subsection{Time and space discretization}
Let us start with the time discretization of the FOM (\ref{eq:psizeta-mom})-(\ref{eq:psizeta-mass}). Let $\Delta t \in \mathbb{R}$, $t^n = t_0 + n \Delta t$, with $n = 0, ..., N_T$ and $T = t_0 + N_T \Delta t$. We denote by $y^n$ the approximation of a generic quantity $y$ at the time $t^n$. %In the following we will denote by $\Omega$ the domain of the equations.
Problem \eqref{eq:psizeta-mom}-\eqref{eq:psizeta-mass} discretized in time by a Backward Differentiation Formula of order 1 (BDF1) reads: given $\omega^0= \omega_0$, for $n \geq 0$ find the solution $(\psi^{n+1}, \omega^{n+1})$ of system:
\begin{align}
\frac{1}{\Delta t}\, \omega^{n+1} + \div \left(\left(\nabla \times \bm{\psi}^{n+1}\right)\omega^{n+1}\right) - \dfrac{1}{Re}\Delta\omega^{n+1} & = b^{n+1},\label{eq:disc_filter_ns-1} \\
\nabla \psi^{n+1} &= \omega^{n+1},\label{eq:disc_filter_ns-2}
\end{align}
%with boundary condition
%\begin{align}
%\psi^{n+1} & = 0 \quad\ \mbox{ on } \partial\Omega \times(t_0,T).\label{eq:bc-pz-d-disc}
%(\partial_y \psi^{n+1}, -\partial_x \psi^{n+1}) & = (u^{n+1}_{in}, v^{n+1}_{in}) \quad\ \mbox{ on } \partial\Omega_{in} \times(t_0,T),\label{eq:bc-pz-d-disc}\\
%\psi^{n+1} & =  c \quad\ \mbox{ on } \partial\Omeg16a_{w} \times(t_0,T),\label{eq:bc-pz-d2-disc}\\
%\partial_n \psi^{n+1} & = 0 \quad \quad\mbox{ on } \partial\Omega_{out} \times(t_0,T).\label{eq:bc-pz-n}
%\end{align}
%In Eq.  (\ref{eq:disc_filter_ns-1}) 
where $b^{n+1} = F^{n+1} + \omega^n/\Delta t$. 
In order to contain the computational cost required
to approximate the solution to problem \eqref{eq:disc_filter_ns-1}-\eqref{eq:disc_filter_ns-2}, we opt for a segregated algorithm. %by solving equation \eqref{eq:disc_filter_ns-1} first, replacing the advection
%field $\nabla \times \bm{\psi}^{n+1}$ with a suitable extrapolation $\nabla \times \bm{\psi}^{*}$, and then solving the equation \eqref{eq:disc_filter_ns-1}. 
Given the vorticity $\omega^{n}$, at $t^{n+1}$ such algorithm requires to:

\begin{itemize}
    \item [i)] Find the vorticity $\omega^{n+1}$ such that
    \begin{align}
\frac{1}{\Delta t}\, \omega^{n+1} + \div \left(\left(\nabla \times \bm{\psi}^{*}\right)\omega^{n+1}\right) - \dfrac{1}{Re}\Delta\omega^{n+1} & = b^{n+1},\label{eq:disc_filter_ns-1-bis}
     \end{align}
where $\bm{\psi}^{n+1}$ in \eqref{eq:disc_filter_ns-1} is replaced by an extrapolation $\bm{\psi}^{*}$. Since we are using BDF1, we set
$\bm{\psi}^{*} =\bm{\psi}^{n}$.
%\begin{rem}\label{rem1}
%We adopt a first order extrapolation for the convective velocity although a BDF2 scheme is used for the time discretization of problem \eqref{eq:psizeta-mom}. This is what the
%NSE solvers in OpenFOAM do, so discretization \eqref{eq:disc_filter_ns-1-bis} would make it a fair comparison
%between our ?? \anna{(manca qualcosa qui)} and NSE algorithms.
%\end{rem}
    \item [ii)] Find $\psi^{n+1}$ such that
    \begin{align}
\nabla \psi^{n+1} = -\omega^{n+1}.\label{eq:disc_filter_ns-2-bis}
\end{align}
\end{itemize}

%where $\b^{n+1} = \f^{n+1} + (4\u^n - \u^{n-1})/(2\Delta t)$. 
%Obviously, other \rev{time} discretization schemes are possible. However, for clarity of exposition we will restrict the description of the approach to the case of BDF2.
For the space discretization of problem  (\ref{eq:disc_filter_ns-1-bis})-(\ref{eq:disc_filter_ns-2-bis}), we adopt a Finite Volume (FV) approximation that is derived directly from the integral form of the governing equations. For this purpose, 
we partition the computational domain $\Omega$ into cells or control volumes $\Omega_i$,
with $i = 1, \dots, N_{c}$, where $N_{c}$ is the total number of cells in the mesh. 
%When we set $\chi = 1$, we will refer to the algorithm as EF, instead of EFR, since there is no relaxation in practice.
The integral form of eq. \eqref{eq:disc_filter_ns-1-bis} for each volume $\Omega_i$ is given by:

\begin{align}\label{eq:zetaFV}
\frac{1}{\Delta t}\, \int_{\Omega_i} \omega^{n+1} d\Omega &+ \int_{\Omega_i} \div \left(\left(\nabla \times \bm{\psi}^{n}\right) \omega^{n+1}\right)  d\Omega \cl
&- \dfrac{1}{Re} \int_{\Omega_i} \Delta\omega^{n+1} d\Omega 
= \int_{\Omega_i}b^{n+1} d\Omega.
\end{align}
By applying the Gauss-divergence theorem, eq.~\eqref{eq:zetaFV} becomes:
\begin{align}\label{eq:zetFV2}
\frac{1}{\Delta t}\, \int_{\Omega_i} \omega^{n+1} d\Omega &+  \int_{\partial \Omega_i} \left(\left(\nabla \times \bm{\psi}^{n}\right) \omega^{n+1}\right) \cdot d\textbf{A} \cl
&- \dfrac{1}{Re}\int_{\partial \Omega_i} \nabla\omega^{n+1} \cdot d\textbf{A}  = \int_{\Omega_i}b^{n+1} d\Omega, 
\end{align}
where $\textbf{A}$ is the surface vector associated with the boundary 
of $\Omega_i$.

%\textcolor{red}{Ho introdotto il numero di volumi finiti $N$ in cui discretizziamo il dominio computazionale} 
Let  \textbf{A}$_j$ be the surface vector of each face of the control volume, 
with $j = 1, \dots, M$. 
Each term in eq.~\eqref{eq:zetFV2} is approximated as follows: 
\begin{itemize}
%\anna{Michele: usiamo $i$ come subindex per
%il control volume e anche per le facce di ciascun volume, non dovremmo usare un indice diverso?} \textcolor{red}{Hai ragione e scusa per la negligenza! Ho modificato usando $j$ come indice di conteggio per le facce di ciascun volumetto di controllo, ho lasciato $i$ come indice di conteggio per $\partial \Omega_i$ che rappresenta l'intera frontiera associata al volumetto di controllo $\Omega_i$ e ho tolto gli indici dai differenziali $d\textbf{A}$ e $d\Omega$}
\item[-] \textit{Convective term}: 
\begin{align}
\int_{\partial \Omega_i} \left(\left(\nabla \times \bm{\psi}^{n}\right) \omega^{n+1}\right) \cdot d\textbf{A} &\approx \sum_j^{} \left(\left(\nabla \times \bm{\psi}_j^{n}\right) \omega^{n+1}_j\right) \cdot \textbf{A}_j = \sum_j^{} \varphi_j^n \omega^{n+1}_j, \label{eq:conv} \\ 
\varphi_j^n &= \left(\nabla \times \bm{\psi}_j^{n}\right) \cdot \textbf{A}_j.  \label{eq:phin}
\end{align} 
In \eqref{eq:conv}, $\nabla \times \bm{\psi}_j^{n}$ is the extrapolated convective velocity
and $\omega^{n+1}_j$ is the vorticity, both relative to the centroid of each control volume face. In \eqref{eq:phin}, $\varphi_j^n$ is the convective flux associated to $\nabla \times \bm{\psi}^{n}$ through face $j$ of the control volume. In OpenFOAM\textsuperscript{\textregistered} solvers, 
the convective flux at the cell faces is typically a linear interpolation of the values from the adjacent cells. %As for the face 
%value $\omega^{n+1}_j$, % \anna{(Non capisco, $\varphi_j^n$ e' il convective flux di cui parliamo nella frase precedente ma qui c'e' $\omega^{n+1}_j$)} \michele{(Si, abbiamo parlato di come calcoliamo $\varphi_j^n$ e ora lo facciamo anche per $\omega^{n+1}_j$. Non e' chiaro cosi'?)},
We also need to approximate $\omega^{n+1}$ at cell face $j$. Different interpolation methods can be applied, including central, upwind, second order upwind, and blended differencing schemes \cite{jasakphd}. In this work, we use a Central Differencing (CD) scheme. % \michele{CITARE LAX}. 
%\textcolor{red}{Anna, vedi che ho cambiato un po' la descrizione della discretizzazione del termine convettivo. Questo e' il remark sull'interpolazione di cui ti parlavo e inoltre vengono specificati i due schemi convettivi che useremo e il cui confronto diventa cruciale per FDA!}
%As it will seen, the choice of the interpolation schemes could be strongly affect %Both velocities have to be determined from the cell centroid values by making use of appropriate interpolation schemes.
%In this work, we use a limited linear differencing scheme is used for the convective term \cite{OpenFoam2015}.
%\anna{(possiamo dare una referenza?)} 
%\textcolor{red}{(ho messo il riferimento al sito di OpenFOAM (specificando la versione che e' stata usata nel lavoro e l'anno di rilascio); al momento non mi vengono in mente altre referenze!)}. 
%Finally, in eq.~\eqref{eq:conv} $\varphi^*_j = \u^{*}_j \cdot \textbf{A}_j$ is the mass flux related to $\u^{*}$ through each face of the control volume.
\item[-] \textit{Diffusion term}: 
\begin{align}%\label{eq:diff}
\int_{\partial \Omega_i} \nabla\omega^{n+1} \cdot d\textbf{A} \approx \sum_j^{} (\nabla\omega^{n+1})_j \cdot \textbf{A}_j, \el
\end{align} 
where $(\nabla\omega^{n+1})_j$ is the gradient of $\omega^{n+1}$ at face $j$. 
%In the case the value of the gradient ${\nabla \phi}_i$ of a generic conservative variable $\phi$ in control volume $\Omega_i$ is needed, as in eq. \eqref{eq:diff}, this one can be computed dividing the expression in eq. \eqref{eq:grad} by the volume of the cell $\Omega_i$.
Let us briefly explain how $(\nabla\omega^{n+1})_j$ is approximated with
second order accuracy on a structured, orthogonal mesh. Let $P$ and $Q$ be two neighboring control volumes (see Fig.~\ref{fig:gradient_image}).
%\anna{Ho due piccoli problemi: 1. nel testo parliamo di structured orthogonal mesh, mentre la figura e' per una mesh non-ortogonale; 2. se la figura e' identica a una che abbiamo gia' usato potremmo avere problemi di copyright. Per risolvere 2, basta poco, tipo spostare la $j$ dall'angolo ad una posizione piu' centrale.} \michele{Cambiata la figura :)}
\begin{figure}[h!]
\centering
\includegraphics[width=0.5\textwidth]{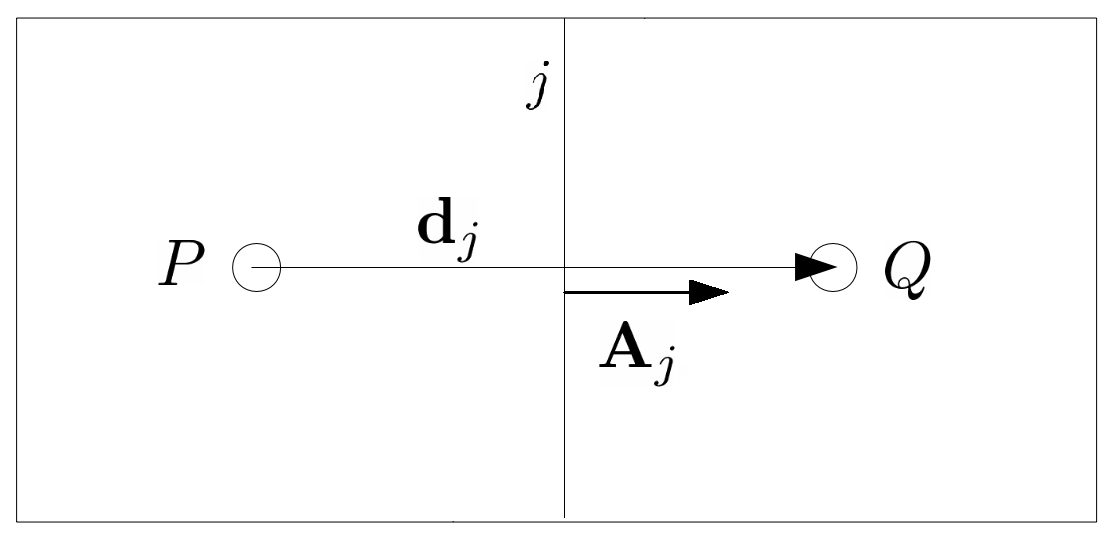}
\caption{Close-up view of two orthogonal control volumes in a 2D configuration.}
\label{fig:gradient_image}
\end{figure}
%as in Fig.~\ref{fig:gradient_image}. 
%\anna{Siccome stiamo discutendo il caso di griglie strutturate, tolto il riferimento alla figura con la griglia non-strutturata
%per evitare confusione.}
The term $(\nabla\omega^{n+1})_j$ is evaluated by subtracting the value of vorticity at the cell centroid on the $P$-side of the face (denoted with $\omega^{n+1}_P$) from the value of vorticity at the centroid on the $Q$-side (denoted with $\omega^{n+1}_Q$) and dividing by the magnitude of the distance vector $\textbf{d}_j$ connecting the two cell centroids:
\begin{align}
(\nabla\omega^{n+1})_j \cdot \textbf{A}_j = \dfrac{\omega^{n+1}_Q - \omega^{n+1}_P}{|\textbf{d}_j|} |\textbf{A}_j|. \el
\end{align} 
%\anna{Michele, e' $\omega^{n+1}_Q - \omega^{n+1}_P$ perche' la normale punta da P a Q, giusto?} \michele{direi di si} \anna{Allora aggiungerei la solita figura che mettiamo per rendere il testo piu' chiaro.} \michele{ok, fatto!!}
%For non-structured, non-orthogonal meshes
%(see Fig.~\ref{fig:gradient_image}), an explicit non-orthogonal correction has to be added to the orthogonal component
%in order to preserve second order accuracy. See \cite{jasakphd} for details.
%\anna{Quello che hai appena spiegato vale per una componente della velocita' e una data direzione, giusto? Forse puo' valere la pena di fare una figura per renderlo piu' immediato, temo
%che a parole risulti un po' confuso. } \textcolor{red}{Immagine messa e inserito qualche dettaglio in piu'!} 
\end{itemize}

Let us denote with $\omega^{n+1}_i$ and $b^{n+1}_i$ the average vorticity 
and source term in control volume $\Omega_i$, respectively.
Moreover, we denote with $\omega^{n+1}_{i,j}$ the vorticity
associated to the centroid of face $j$ normalized by the volume of $\Omega_i$.
Then, the discretized form of eq.~\eqref{eq:zetFV2}, divided by the control volume 
$\Omega_i$, can be written as:
\begin{align}\label{eq:evolveFV-1.1_disc}
\frac{1}{\Delta t}\, \omega^{n+1}_i &+ \sum_j^{} \varphi^n_j \omega^{n+1}_{i,j} - \dfrac{1}{Re} \sum_j^{} (\nabla\omega^{n+1}_i)_j \cdot \textbf{A}_j  = b^{n+1}_i.
\end{align}

Next, we deal with the space approximation of the eq.~\eqref{eq:disc_filter_ns-2-bis}. After using Gauss-divergence theorem, the integral form of eq. \eqref{eq:disc_filter_ns-2-bis} reads:
\begin{align}\label{eq:psiFV}
 - \int_{\partial \Omega_i} \nabla\psi^{n+1} \cdot d\textbf{A}  = \int_{\Omega_i}{\omega}^{n+1} d\Omega.
\end{align}
Once we approximate the integrals and divide by the control volume $\Omega_i$,
eq.~\eqref{eq:psiFV} becomes:
\begin{align}\label{eq:psiFV2}
-\sum_j^{} (\nabla\psi^{n+1}_i)_j \cdot \textbf{A}_j  =\omega^{n+1}_i.
\end{align}
In eq. \eqref{eq:psiFV2}, $(\nabla\psi^{n+1}_i)_j$ is the gradient of $\psi^{n+1}$ at faces $j$ and it is approximated in the same way as $(\nabla\omega^{n+1}_i)_j$.
Finally, the fully discretized form of problem \eqref{eq:disc_filter_ns-1-bis}-\eqref{eq:disc_filter_ns-2-bis} is given by system \eqref{eq:evolveFV-1.1_disc}, \eqref{eq:psiFV2}.

As mentioned in Sec.~\ref{sec:intro}, for the implementation of the numerical scheme described in this section we chose the finite volume C++ library OpenFOAM\textsuperscript{\textregistered} \cite{Weller1998}.

\section{The Reduced Order Model}\label{sec:ROM}
%The Reduced Order Model (ROM) we propose is an extension of the model introduced in \cite{Stabile2017, Stabile2018}. In Sec~\ref{sec:ROM_1} we introduce the procedure we use to construct a POD-Galerkin ROM
%and in Sec.~\ref{sec:ROM_2} we present the strategy we choose for pressure stabilization at reduced order level. 
%enforce non-homogeneous Dirichlet boundary conditions for the velocity field at the reduced order level. 
%The ROM computations are carried out using ITHACA-FV \cite{RoSta17}, 
%an in-house open source C++ library. %the main notions of
%projection-based ROMs are recalled. Subsection 3.2 introduces the POD technique and the general procedure
%used to construct a POD-Galerkin ROM. Finally
%subsection 3.4 outlines the treatment of non-homogeneous boundary conditions at the reduced order level.

%\subsection{A POD-Galerkin projection method}\label{sec:ROM_1}

The main idea of reduced order modeling for parametrized PDEs is the assumption that 
solutions live in a low dimensional manifold. Thus, any solution can be approximated as a 
linear combination of a reduced number of global basis functions. 

We approximate vorticity field $\omega$ and stream function $\psi$ as linear combinations of the dominant modes (basis functions), which are assumed to be dependent on space variables only,
multiplied by scalar coefficients that depend on time and/or parameters. %For simplicity, we will consider $Re$ as the only parameter. 
We arrange all the parameters the problem depends upon in a vector $\bm{\pi}$ that belongs to a $d$-dimensional parameter space $\mathcal{P}$ in $\mathbb{R}^d$,
where $d$ is the number of parameters.
Thus, we have:
%The modes, $\bm{\varphi}_i(\bm{x})$, $\bm{\overline{\varphi}}_i(\bm{x})$, $\psi_i(\bm{x})$ and $\overline{\psi}_i(\bm{x})$, , while the coefficients, $\beta_i(\bm{\pi}, t)$, $\overline{\beta_i}(\bm{\pi}, t)$, $\gamma_i(\bm{\pi}, t)$ and $\overline{\gamma_i}(\bm{\pi}, t)$, are allowed to have temporal  dependency:
%\textcolor{blue}{$i$ lo abbiamo gia' usato per i volumi di controllo ma non credo che riadoperare lo stesso indice qui crei confusione. Che ne pensi?}
\begin{align}
\omega \approx \omega_r = \sum_{i=1}^{N_{\omega}^r} \beta_i(\bm{\pi}, t) {\varphi}_i(\bm{x}), \quad 
\psi \approx \psi_r = \sum_{i=1}^{N_{\psi}^r} \gamma_i(\bm{\pi}, t) \xi_i(\bm{x}). \label{eq:ROM_1} %\\
%\u \approx \u_r = \sum_{i=1}^{N_u^r} \overline{\beta_i}(\bm{\pi}, t) \bm{\overline{\varphi}}_i(\bm{x}), \quad
%\overline{q} \approx \overline{q}_r = \sum_{i=1}^{N_{\overline{q}}^r} \overline{\gamma_i}(\bm{\pi}, t) \overline{\psi}_i(\bm{x}). \label{eq:ROM_2}
\end{align}
In \eqref{eq:ROM_1}, $N_{\Phi}^r$, $\Phi = \omega, \psi$, denotes the cardinality of a reduced basis for the space field $\Phi$ belongs to. We remark that we consider different coefficients for the approximation of the vorticity $\omega$ and stream function $\psi$ fields, unlike previous works \cite{Ahmed2019,Ahmed2020,Pawar2020-2,Pawar2020,Pawar2021,San2015}. %\michele{verificare di aver citato tutto!}.
This choice will be justified numerically in Sec.~\ref{sec:time}.

\begin{rem}\label{rem:rem1}
As mentioned earlier, we only consider a body force given by the product between a space dependent function and a time dependent function. For the stream function-vorticity formulation, this means:
\begin{equation*}
F(x,y,z) = F_2(t) F_1(x,y).
\end{equation*}
Thanks to this assumption, the forcing term is already expressed in the form of \eqref{eq:ROM_1} and does not require further treatment. %On the other hand, when time and space dependencies of $F$ are not split, one could find the basis. %We are releasing this constraint as it leads to a substantial simplification. Moreover, 
%\michele{da completare...}
%of indicator function a defined in (13). For this reason, we do not need to compute a reduced
%order approximation of ṽ.
\end{rem}

Using \eqref{eq:ROM_1} to approximate $\omega^{n+1}$ and $\psi^{n+1}$ in \eqref{eq:disc_filter_ns-1-bis}-\eqref{eq:disc_filter_ns-2-bis},
we obtain
\begin{align}
\frac{1}{\Delta t}\, \omega_r^{n+1} + \div \left(\left(\nabla \times \bm{\psi}_r^{*}\right)\omega_r^{n+1}\right) - \dfrac{1}{Re}\Delta\omega_r^{n+1} & = b_r^{n+1}, \label{eq:red-1.1} \\
\nabla \psi_r^{n+1} &= -\omega_r^{n+1}, \label{eq:red-1.2}
\end{align}
where $b_r^{n+1} = F^{n+1} + \omega_r^n/\Delta t$ and we set $\bm{\psi}_r^{*} = \bm{\psi}_r^{n}$. %$b_r^{n+1} = (4\omega_r^n - \omega_r^{n-1})/(2\Delta t)$.

In the literature, one can find several techniques to generate the reduced basis spaces, e.g.,~Proper Orthogonal Decomposition (POD), the Proper Generalized Decomposition and the Reduced Basis with a greedy sampling strategy.
See, e.g., \cite{ModelOrderReduction,ChinestaEnc2017,Chinesta2011, Dumon20111387,Kalashnikova_ROMcomprohtua,quarteroniRB2016,Rozza2008}. 
We generate the reduced basis spaces with the method of snapshots. Next, we briefly describe how this method works.

Let $\mathcal{K}= \{\bm{\pi}^1, \dots, \bm{\pi}^{N_k}\}$ be a finite dimensional training set of samples chosen inside the parameter space $\mathcal{P}$.
We solve the FOM described in Sec.~\ref{sec:FOM}  
for each $\bm{\pi}^k \in \mathcal{K}$ and for each time instant $t^j \in \{t^1, \dots, t^{N_t}\} \subset (t_0, T]$. %and $t^k \in \{t^1, \dots, t^{N_t}\} \subset (t_0, T]$. %  \subset \mathcal{P}$. The total number of snapshots $N_s$ is given by $N_s = N_k \cdot N_t$. 
The snapshots matrices are obtained from the full-order snapshots:  %$\bm{\mathcal{S}}_\u$, $\bm{\mathcal{S}}_q$ and $\bm{\mathcal{S}}_\overline{q}$, are given by $N_s$ full-order snapshots:
\begin{align}\label{eq:space}
\bm{\mathcal{S}}_{{{\Phi}}} = [{{\Phi}}(\bm{\pi}^1, t^1), \dots, {{\Phi}}(\bm{\pi}^{N_k}, t^{N_t})] \in \mathbb{R}^{N_{\Phi}^h \times N_s} \quad
\text{for} \quad {{\Phi}} = \{\omega_h, \psi_h\},
\end{align}
where $N_s = N_t \cdot N_k$ is the total number of the snapshots, $N_{\Phi}^h$
is the dimension of the space $\Phi$ belong to in the FOM, and the subscript $h$ indicates a solution computed with the FOM. %Note that ${\Phi}$ could be either a scalar or
%a vector field. 
The POD problem consists in finding, for each value of the dimension of the POD space $N_{POD} = 1, \dots, N_s$, the scalar coefficients $a_1^1, \dots, a_1^{N_s}, \dots, a_{N_s}^1, \dots, a_{N_s}^{N_s}$ and functions ${\zeta}_1, \dots, {\zeta}_{N_s}$, that minimize the error between the snapshots and their projection onto the POD basis. In the $L^2$-norm, we have
\begin{align}
E_{N_{POD}} = \text{arg min} \sum_{i=1}^{N_s} ||{{\Phi}_i} - \sum_{k=1}^{N_{POD}} a_i^k {\zeta}_k || \quad \forall N_{POD} = 1, \dots, N_s    \cl
\text{with} \quad ({\zeta}_i, {\zeta}_j)_{L_2(\Omega)} = \delta_{i,j} \quad \forall i,j = 1, \dots, N_s. \label{eq:min_prob}
\end{align}

It can be shown \cite{Kunisch2002492} that problem~\eqref{eq:min_prob} is equivalent to the following eigenvalue problem
\begin{align}
\bm{\mathcal{C}}^{{\Phi}} \bm{Q}^{{\Phi}} &= \bm{Q}^{{\Phi}} \bm{\Lambda}^{{\Phi}}, \label{eq:eigen_prob} \\
\mathcal{C}_{ij}^\Phi &= ({\Phi}_i, {\Phi}_j)_{L_2(\Omega)} \quad \text{for} \quad i,j = 1, \dots, N_s,
\end{align}
where $\bm{\mathcal{C}}^{{\Phi}}$ is the correlation matrix computed from the snapshot matrix $\bm{\mathcal{S}}_{{{\Phi}}}$, $\bm{Q}^{{\Phi}}$ is the matrix of eigenvectors and $\bm{\Lambda}^{{\Phi}}$ is a diagonal matrix whose diagonal entries are the
eigenvalues of $\bm{\mathcal{C}}^{{\Phi}}$. 
%\anna{Forse dire ``vector of eigenvalues'' e' un po' fuorviante, $\bm{\Lambda}^{{\Phi}}$ e' una matrice con gli autovalori sulla diagonale.}16
Then, the basis functions are obtained as follows:
\begin{align}\label{eq:basis_func}
{\zeta}_i = \dfrac{1}{N_s \Lambda_i^\Phi} \sum_{j=1}^{N_s} {\Phi}_j Q_{ij}^\Phi.
\end{align}
The POD modes resulting from the aforementioned methodology are:
%\anna{Non capisco questa notazione: e' una matrice o uno spazio?}
\begin{align}\label{eq:spaces}
L_\Phi = [{\zeta}_1, \dots, {\zeta}_{N_\Phi^r}] \in \mathbb{R}^{N_\Phi^h \times N_\Phi^r},
\end{align}
where the values of $N_\Phi^r < N_s$ are chosen according to the eigenvalue decay of the vectors of eigenvalues \bm{\Lambda}.
Then, the reduced order model can be obtained through a Galerkin projection of the governing equations onto the POD spaces. 

In order to write the algebraic system associated with the reduced problem \rev{\eqref{eq:red-1.1}-\eqref{eq:red-1.2}}, we
introduce the following matrices:
\begin{align}
&M_{r_{ij}} = ({\varphi}_i, {\varphi}_j)_{L_2(\Omega)}, \quad \widetilde{M}_{r_{ij}} = ({\xi}_i, {\varphi}_j)_{L_2(\Omega)}, \quad A_{r_{ij}} = ({\varphi}_i, \Delta {\varphi}_j)_{L_2(\Omega)}, \label{eq:matrices_evolve1} \\
&B_{r_{ij}} = ({\xi}_i, \Delta \xi_j)_{L_2(\Omega)}, \quad G_{r_{ijk}} = (\varphi_i, \div \left(\left(\nabla \times \xi_j\right) \varphi_k\right))_{L_2(\Omega)}, \label{eq:matrices_evolve2} \\
&H_{r_{ij}} = ({\varphi}_i, F_1)_{L_2(\Omega)},
\label{eq:matrices_evolve3} %\\
\end{align}
where ${\varphi}_i$ and $\xi_i$ are the basis functions in \eqref{eq:ROM_1}. 
At time $t^{n+1}$, the reduced algebraic system 
for \rev{\eqref{eq:red-1.1}-\eqref{eq:red-1.2}} reads: given $\bm{\beta}^{n}$ and $\bm{\gamma}^{n}$ %\michele{(e anche gamma no? o non lo hai citato per qualche motivo particolare?)}, 
find vectors $\bm{\beta}^{n+1}$ and $\bm{\gamma}^{n+1}$ containing the values of coefficients $\beta_i$ and $\gamma_i$ in \eqref{eq:ROM_1} at time $t^{n+1}$ such that
\begin{align}
&\bm{M}_r\left(\dfrac{\bm{\beta}^{n+1} - \bm{\beta}^{n}}{\Delta t}\right) + \left(\bm{\gamma}^n\right)^T\bm{G}_r\bm{\beta}^{n+1} - \dfrac{1}{Re} \bm{A}_r \bm{\beta}^{n+1} = \bm{H}_rF_2^{n+1}, \label{eq:reduced_1} \\
& \bm{B}_r \bm{\gamma}^{n+1} + \widetilde{\bm{M}_r} \bm{\beta}^{n+1} = 0. \label{eq:reduced_2}
%\cl
%&\quad \quad  + \bm{B}_r \bm{\gamma}^{n+1} = \dfrac{\rho}{\Delta t} \widetilde{\bm{M}}_r \left(2\overline{\bm{\beta}}^{n} - \dfrac{1}{2}\overline{\bm{\beta}}^{n-1}\right), \label{eq:reduced_1} \\
%&\bm{P}_r \bm{\beta}^{n+1} = 0, \label{eq:reduced_2}
\end{align}
Finally, the initial conditions for the ROM algebraic system \eqref{eq:reduced_1}-\eqref{eq:reduced_2}
are obtained performing a Galerkin projection of the initial full order condition onto the POD basis spaces:
\begin{align}
{\beta^0}_i = (\omega_0, {\varphi}_i)_{L_2(\Omega)}. %\cl
%{\overline{\beta}^0}_i = (\u(\bm{x},\bm{\pi},t_0), \bm{\overline{\varphi}}_i)_{L_2(\Omega)}. \el
\end{align}
%\anna{Queste non sarebbero le component, invece che l'intero vettore?}
%a POD strategy is exploited and is chosen to apply the POD onto the full snapshots matrices that include both
%the time and parameter dependency.%The last statement leads to the following approximation of the fields:

\begin{rem}\label{rem2}
We consider homogeneous boundary conditions \eqref{eq:bc-pz-d1}-\eqref{eq:bc-pz-d22}. So, the approximated
vorticity $\omega_r$ and stream function $\psi_r$ automatically satisfy the boundary conditions and no special treatment (such as lifting function and penalty methods \cite{Girfoglio_JCP, Lorenzi2016,Stabile2018, CiCP-30-34}) is necessary.
\end{rem}
%The Reduced Order Model (ROM) we propose is an extension of the model
%introduced in [32, 33]. In Sec 3.1 we introduce the procedure we use to construct
%a POD-Galerkin ROM and in Sec. 3.2 we present the strategy we choose for pres-
%sure stabilization at reduced order level. Finally, Sec. 3.3 describes the lifting
%function method we apply to enforce non-homogeneous Dirichlet boundary con-
%ditions for the velocity field at the reduced order level. The ROM computations
%are carried out using ITHACA-FV [47], an in-house open source C++ library.\subsection{Numerical results}

%\subsection{A POD-Galerkin projection method}

%\begin{rem}
%%REMARK PER LE BC IN ROM DA LORENZI:
% Therefore, the Boundary Conditions (BCs) are “embedded” in the Bji term and not explicit present
%in the ROM formulation. This may be of concern when dealing with parametrized ROM, e.g., in control-oriented
%applications. In particular, in the fluid dynamics field, the classic control variable is the velocity at the boundary [58]
%since it could be used to control the velocity field in the domain or an output variable of interest. In this view, if the
%reduced order model is direc
%This is out of scope 
%\end{rem}
%\div \u^{n+1} & = 0, \label{eq:disc_filter_ns-2}\\
%-\alpha^2\div \left(a(\u^{n+1}) \nabla\ubar^{n+1}\right) +\ubar^{n+1} +\nabla \lambda^{n+1} & = \u^{n+1},\label{eq:disc_filter_mom}\\
%\div \ubar^{n+1} & = 0, \label{eq:disc_filter_mass}

%w two co-rotating vortices
%\subsubsection{Channel  flow  over  a  step}
%\subsubsection{Lid driven cavity flow}
\section{Numerical results}\label{sec:num_res}

\begin{comment}
STRUTTURA DEI RISULTATI:

validazione FOM
Figura 1: condizioni iniziali per omega e u
Figura 2: confronto tra u NSE e u psi zeta per 4 tempi (t = 4, t = 8, t = 16, t = 20)
Figura 3: come sopra pero' per la psi
Figura 4: come sopra pero' per la zeta
validazione ROM ricostruzione in tempo
Figura 5: eigenvalues decay??? per psi e zeta per 250 snap
Figura 6: facciamo vedere anche i primi modi magari??
FIgura 7: errori relativi per 3 soglie energetiche?!? per psi e zeta NO SOLO UNA SOGLIA! facciamo vedere l'enstrofia!
Figura 8: ricostruzione dei campi per 4 tempi (come figura 2, 3 e 4?) direi di si!
Figura 9: coefficienti psi e zeta vs t (per far vedere che ha senso utilizzare coefficienti diversi per le due variabili visto che appunto gli andamenti in tempo sono completamente diversi) - probabilmente splittato in due figure, uno per psi e uno per zeta
validazione ROM ricostruzione parametrica Re
Figura 10: eigenvaleus decay parametrico in Re per quanti snap? 1000 per ogni Re? Quanti Re? tra 200 e 400, 5 numeri di Reynolds? Ovvero 5*500 = 2500 snap?
Figura 11: errori relativi per 3 soglie?!?
Figura 12: ricostruzione dei campi per 4 tempi (come figura 2, 3 e 4?) direi di si!
Costo computazionale molto dettagliato magari...
\end{comment}

In order to validate our FOM and ROM approaches, we consider a widely used benchmark test known as the vortex merger. It consists in fluid motion induced by a pair of co-rotating vortices separated from each other by a certain distance. One of the reasons why this test has been extensively investigated in two-dimensions is that it explains the average inverse energy and direct enstrophy cascades observed in two-dimensional turbulence \cite{Reinaud2005}.

\begin{figure}
\centering
 \begin{overpic}[width=0.34\textwidth]{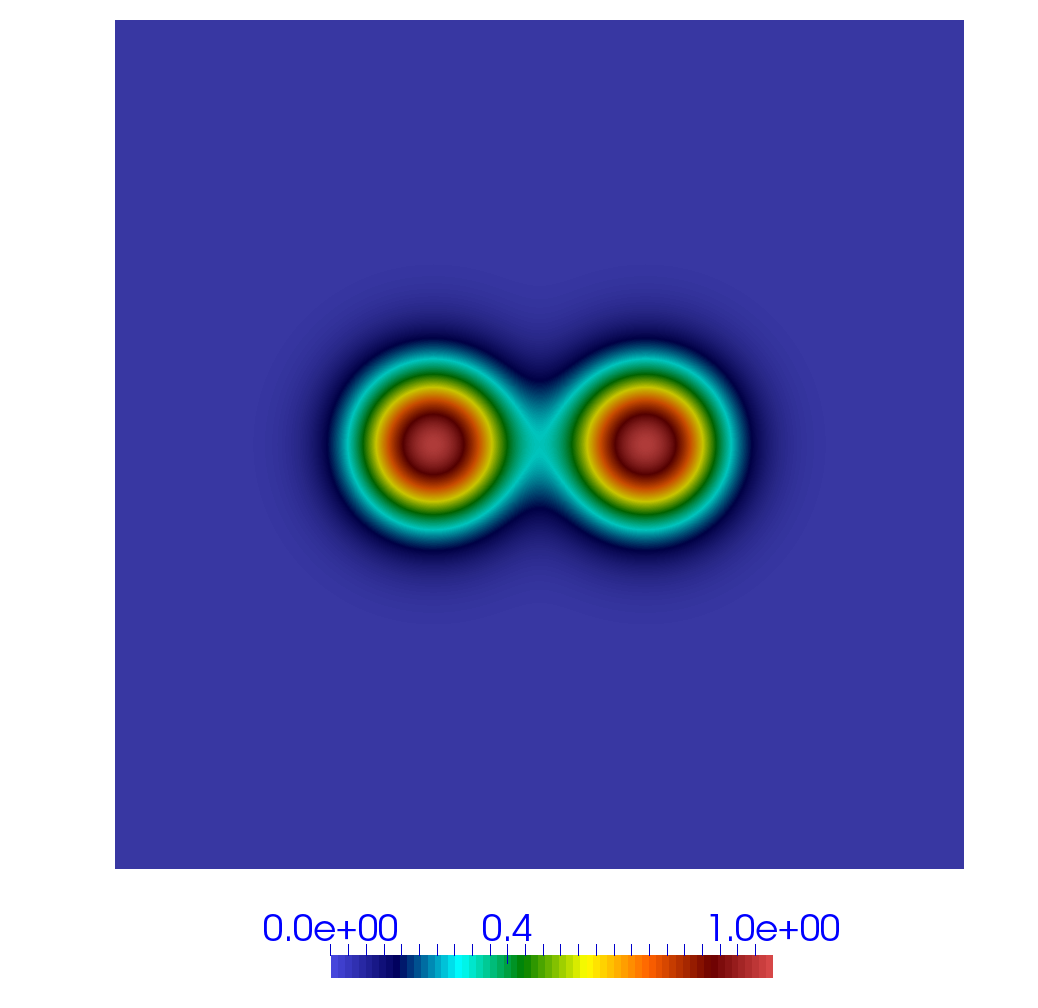}
        %\put(35,18){FOM}
        %\put(-8,7){$\u$}
      \end{overpic}
 \begin{overpic}[width=0.34\textwidth]{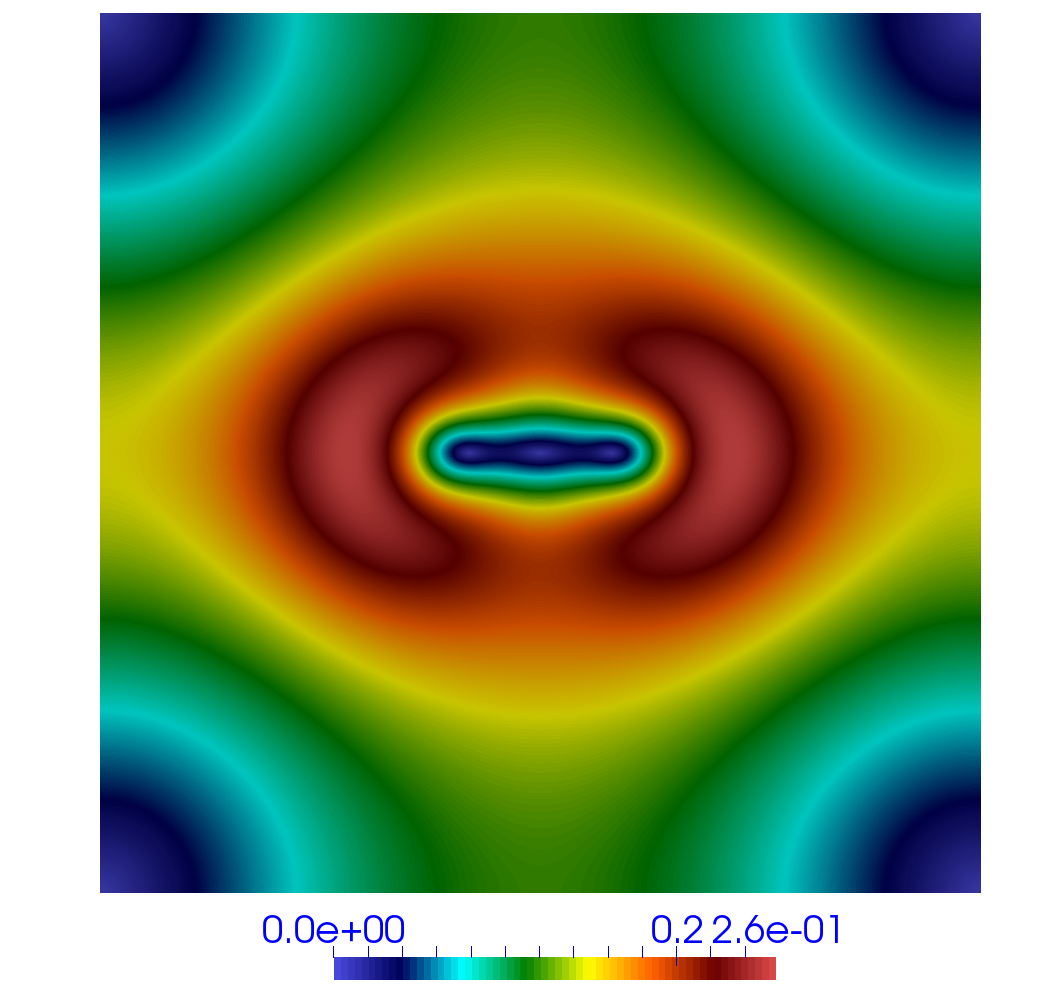}
        %\put(35,18){ROM}
      \end{overpic}\\
\caption{Initial conditions for $\omega_0$ (left) and $\u_0$ (right).}
\label{fig:IC}
\end{figure}

The computational domain is a $2\pi \times 2\pi$ rectangle. %We generated data
%snapshots for Re = [200, 400, 600, 800] with a 2562
%spatial grid and a time-step of 0.01 from time t = 0
%to t = 20. We tested the ETC framework for the out-of-sample condition at Re = 1000. The linear and
%nonlinear operators in GP equations for two-dimensional vorticity transport equation are:
%Concerning the $\psi-\omega$ formulation, 
%We impose homogenous Dirichlet and Neumann boundary conditions for $\psi$ and $\omega$, respectively. 
The initial condition for the vortex merger test case is given by:
\begin{align}
\omega(x,y,0) = \omega_0(x,y) = e^{\left(-\pi \left[\left(x-x_1\right)^2 + \left(y-y_1\right)^2\right]\right)} + e^{\left(-\pi \left[\left(x-x_2\right)^2 + \left(y-y_2\right)^2\right]\right)},
\end{align}
where ($x_1$, $y_1$) and ($x_2$, $y_2$) are the initial locations of the centers of the vortices. We set them as ($3\pi/4$, $\pi$) and ($5\pi/4$, $\pi$), respectively. Fig.~\ref{fig:IC} shows
the initial vorticity $\omega_0$ (left) and corresponding initial velocity $\u_0$ (right). We note that the initial condition $\u_0 (x,y) = \u(x,y,0) = \nabla \times \bm{\psi}_0$ is computed by solving $\Delta \psi_0 = \omega_0$. 
We let the system evolve until time $T = 20$.
%\anna{Michele, che intervallo di tempo stiamo considerando? [0, 20], giusto?} \michele{Si, poi nel caso ROM parametrico usiamo [0, 10] per evitare di collezionare troppi snapshots come riportato nel testo}
As previously mentioned, we enforce boundary conditions \eqref{eq:bc-pz-d1}-\eqref{eq:bc-pz-d22}. 
%\michele{(vogliamo evidenziare che gli altri invece usano BC periodiche per questo problema?)}
%\anna{YES! Tutti usano BC periodiche o c'e' qualcuno che usa le stesse che usiamo noi?} \michele{Direi tutti quelli che usano il vortex merger come benchmark, quindi \cite{Pawar2020, Pawar2020-2, Ahmed2020, Ahmed2019, Pawar2021} All'inizio feci una prova con le BC periodiche ma ebbi problemi in OF e quindi optai fin da subito per queste altre:) }. \anna{Io non direi niente allora, per evitare che un reviewer ci chieda di riprovare con le BC periodiche visto che tutti usano quelle. Siccome ai fini del metodo non cambia nulla, evitarei proprio di dirlo.} \michele{OK!}

Following \cite{Ahmed2019,Ahmed2020,Pawar2020-2,Pawar2020,Pawar2021},
we consider a computational grid with $256^2 = 65536$ cells for all the simulations whose results are reported next.
%\anna{Michele, a cosa si riferiscono le 2 referenze? Anche loro usano mesh con $256^2$ celle?}
%\michele{Si, anche loro, e per Re anche piu' grandi, fino a 1000/1500 mi sembra.}

\subsection{Validation of the FOM}\label{sec:FOM_validation}
Let us start with the validation of our implementation of the stream function-vorticity formulation at the FOM level. We compare the results obtained with such formulation against the results
produced by the standard Navier--Stokes solver in OpenFOAM \textit{icoFoam}, which is based on a partitioned algorithm called PISO \cite{PISO}. We set $\Delta t = 0.01$, $\f = \F = 0$ and $Re = 800$ \cite{Pawar2020-2,Pawar2020}. %For the convective term, we use a
%second-order accurate central difference scheme that features low dissipation \cite{}. We set $\Delta t = 0.01$ \michele{scriviamo anche CFL 0.14?} %which allows to obtain CF L max ≈ 0.14 at the time when the velocity reaches its
%maximum value. 
%For the FOM validation, we consider as benchmark the solution by NSE in %velocity-pressure formulation computed by using the PISO algorithm.
%\textcolor{red}{scrivere il numero di cicli PISO}.
For the simulations with \textit{icoFoam}, we enforce boundary condition \eqref{eq:bc-ns-d}. The partitioned algorithm requires also a boundary condition for the pressure problem: we set $\nabla p \cdot \n = 0$ on $\partial \Omega$. %See Figure \ref{fig:IC} for $\omega_0$ and $\u_0$.

%and set $\Deltad t = 0.01$ \cite{} that corresponds to $CFL_{max} \approx 0.14$ 

\begin{figure}[htb]
\centering
\hspace{.4cm} \begin{overpic}[width=0.2\textwidth, grid=false]{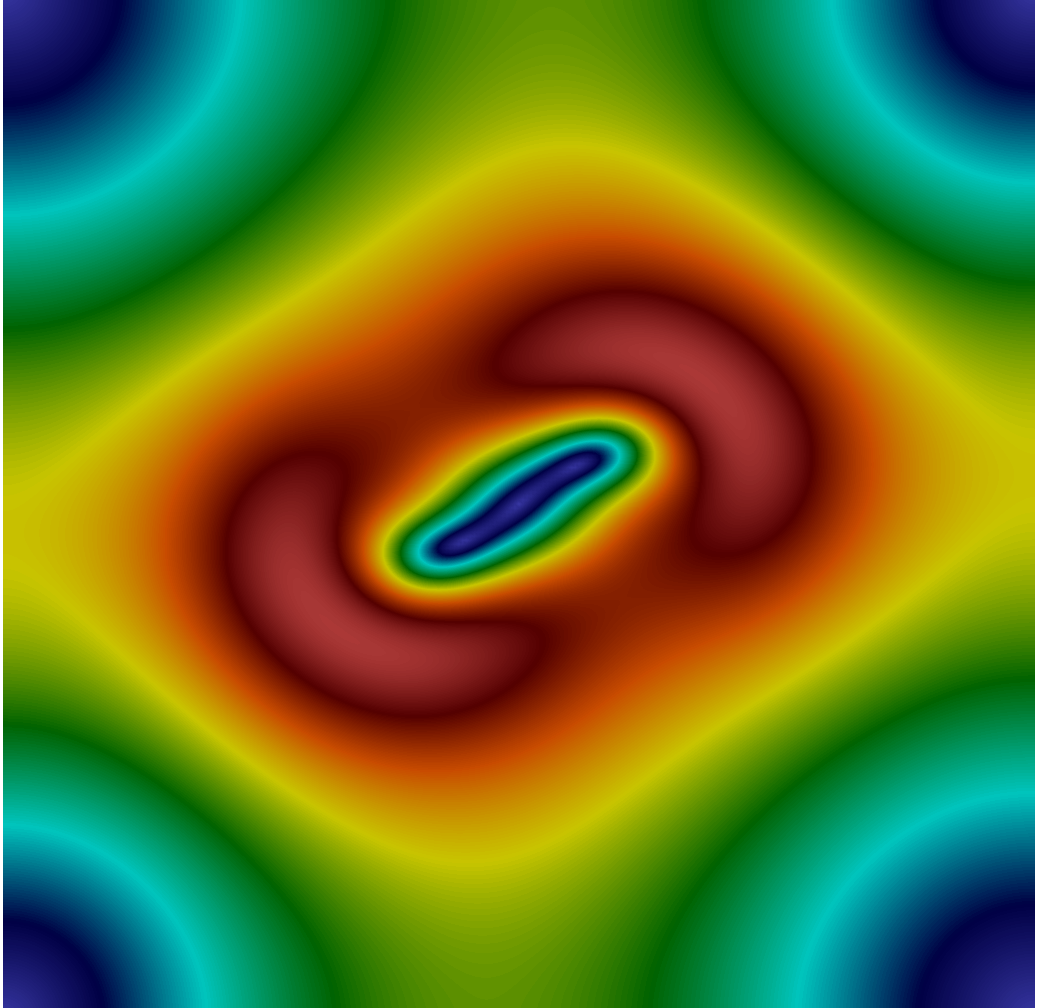}
        \put(38,101){$t = 4$}
        \put(-30,50){NSE}
        \put(-28,40){$\u, p$}
      \end{overpic}
 \begin{overpic}[width=0.2\textwidth]{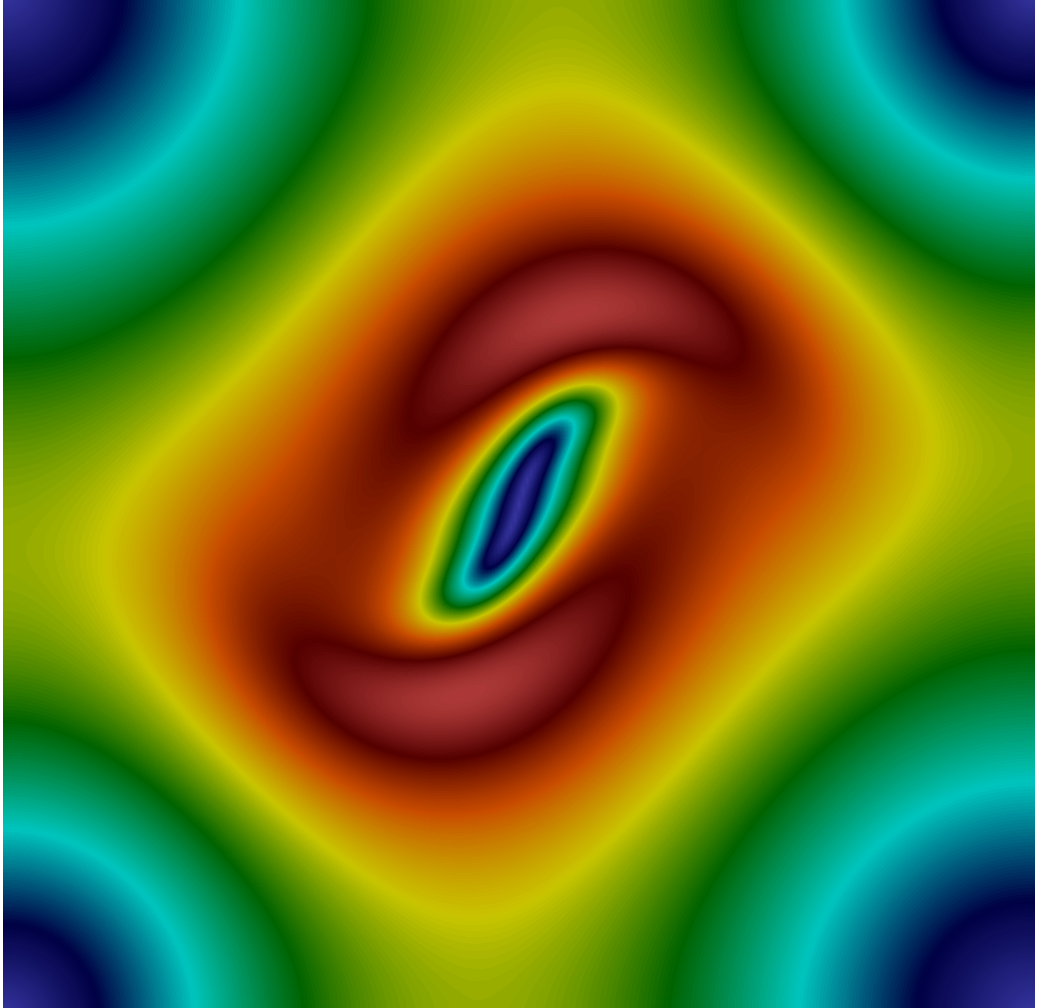}
        \put(38,101){$t = 8$}
      \end{overpic}
 \begin{overpic}[width=0.2\textwidth]{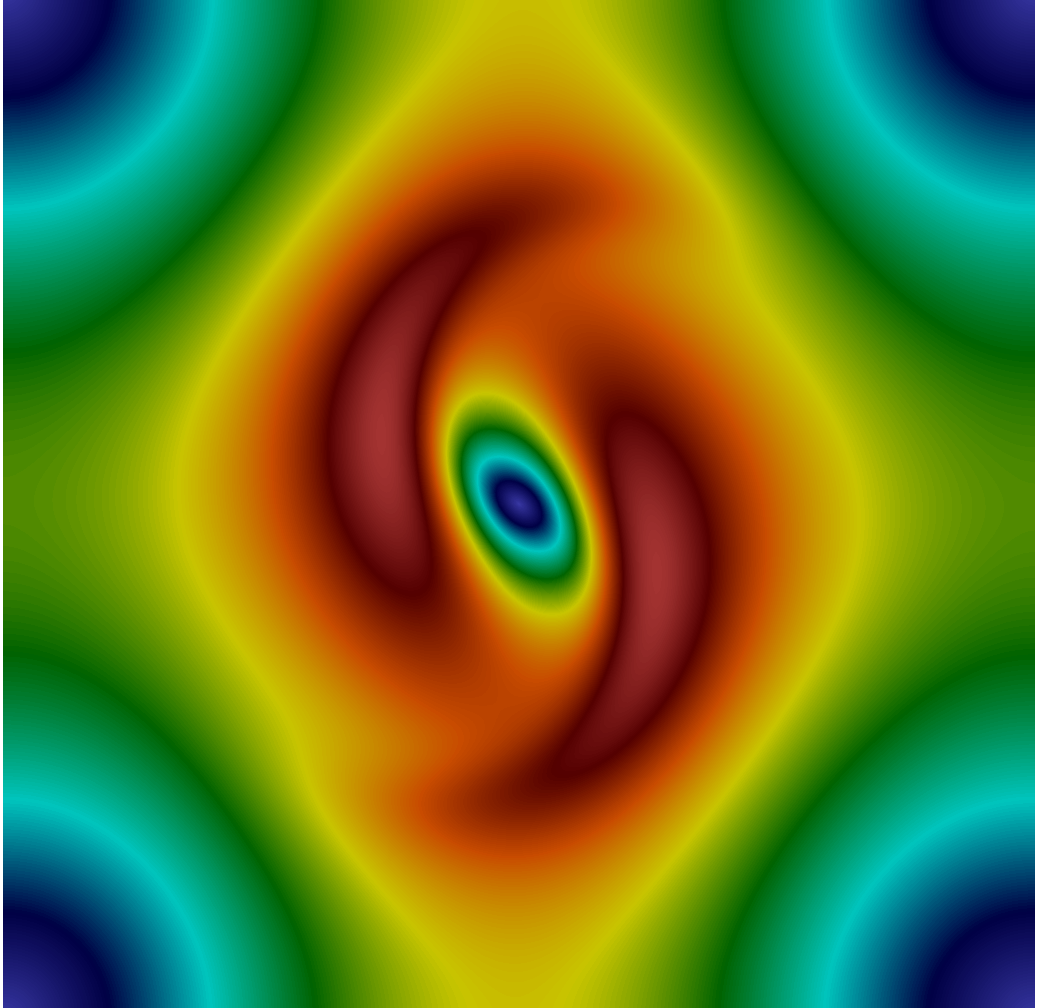}
        \put(35,101){$t = 16$}
      \end{overpic}
 \begin{overpic}[width=0.2\textwidth]{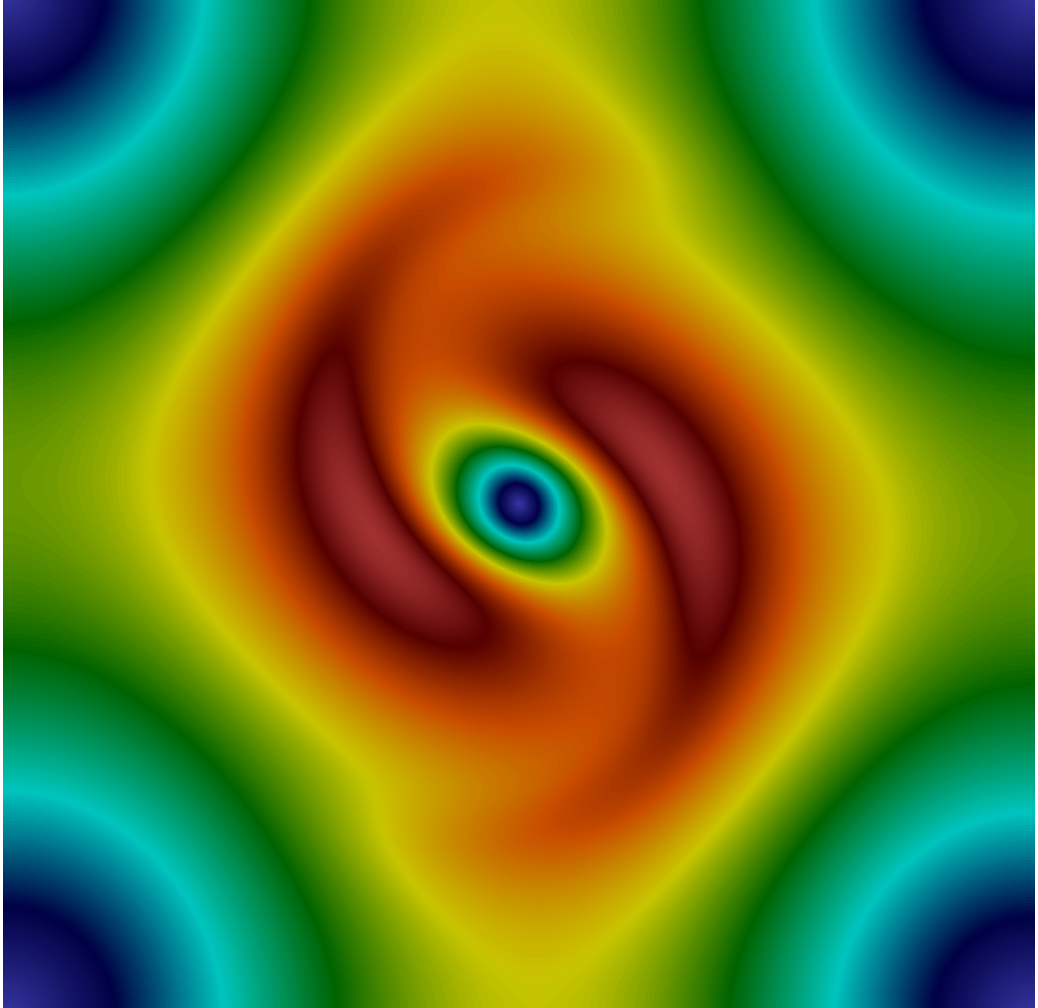}
         \put(35,101){$t = 20$}
      \end{overpic}
      \begin{overpic}[width=0.09\textwidth, grid=false]{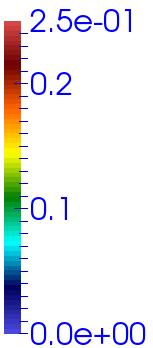}
      \end{overpic}\\
\hspace{.4cm} \begin{overpic}[width=0.2\textwidth]{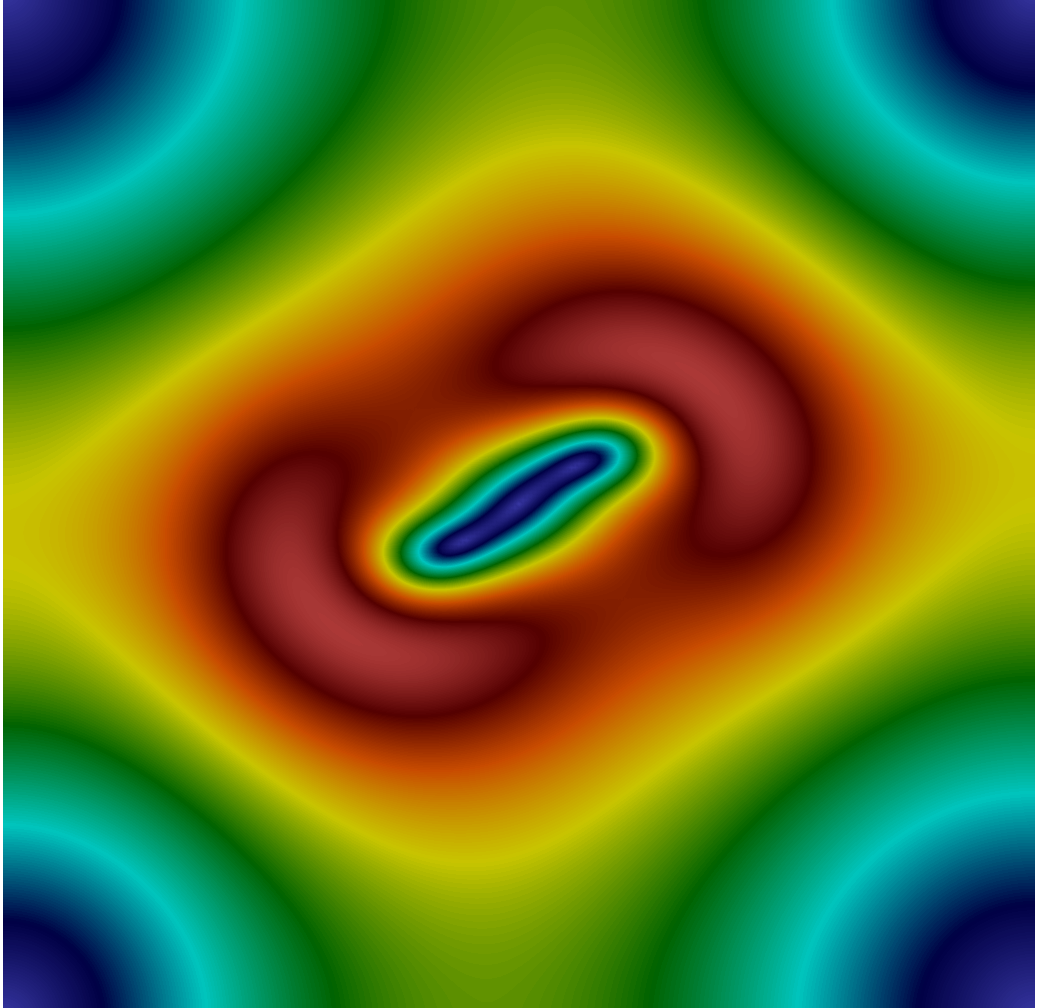}
        \put(-30,50){NSE}
        \put(-28,40){$\psi, \omega$}
      \end{overpic}
 \begin{overpic}[width=0.2\textwidth]{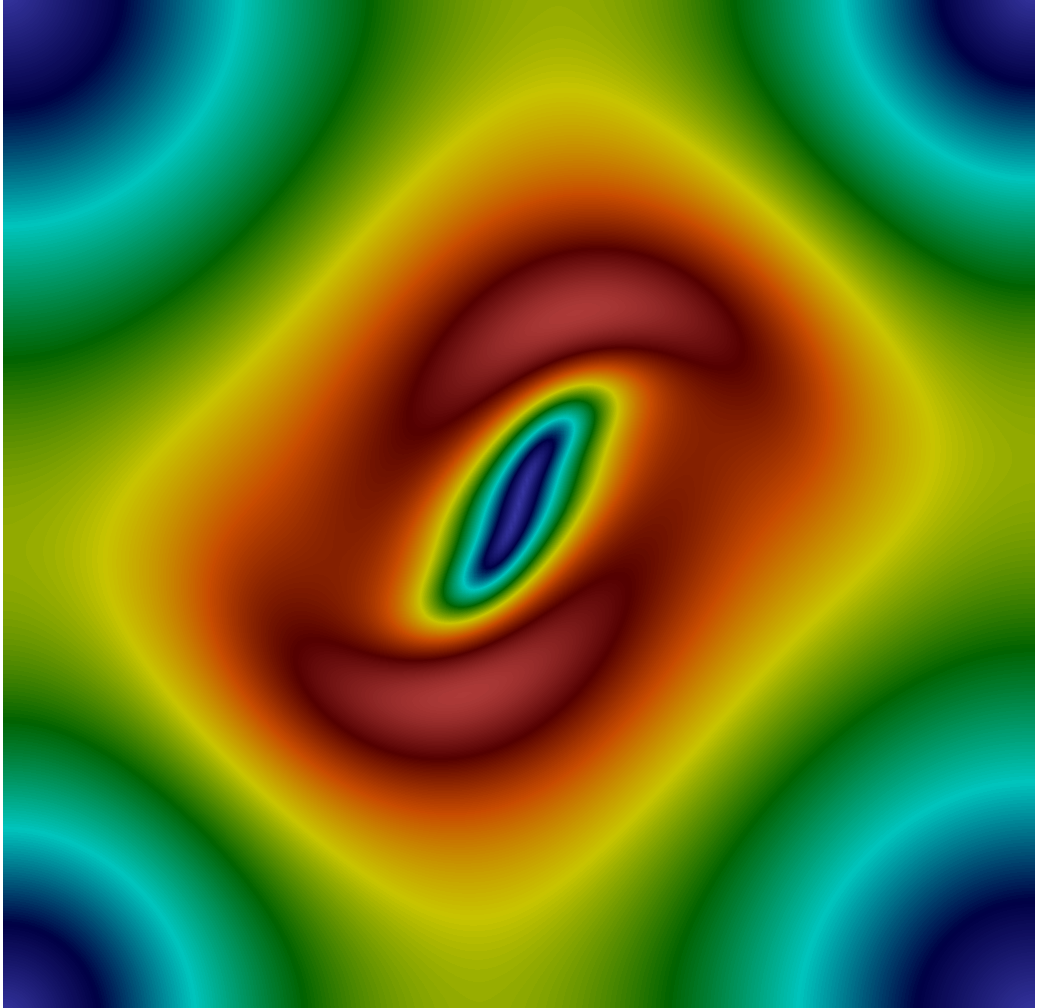}
      \end{overpic}
 \begin{overpic}[width=0.2\textwidth]{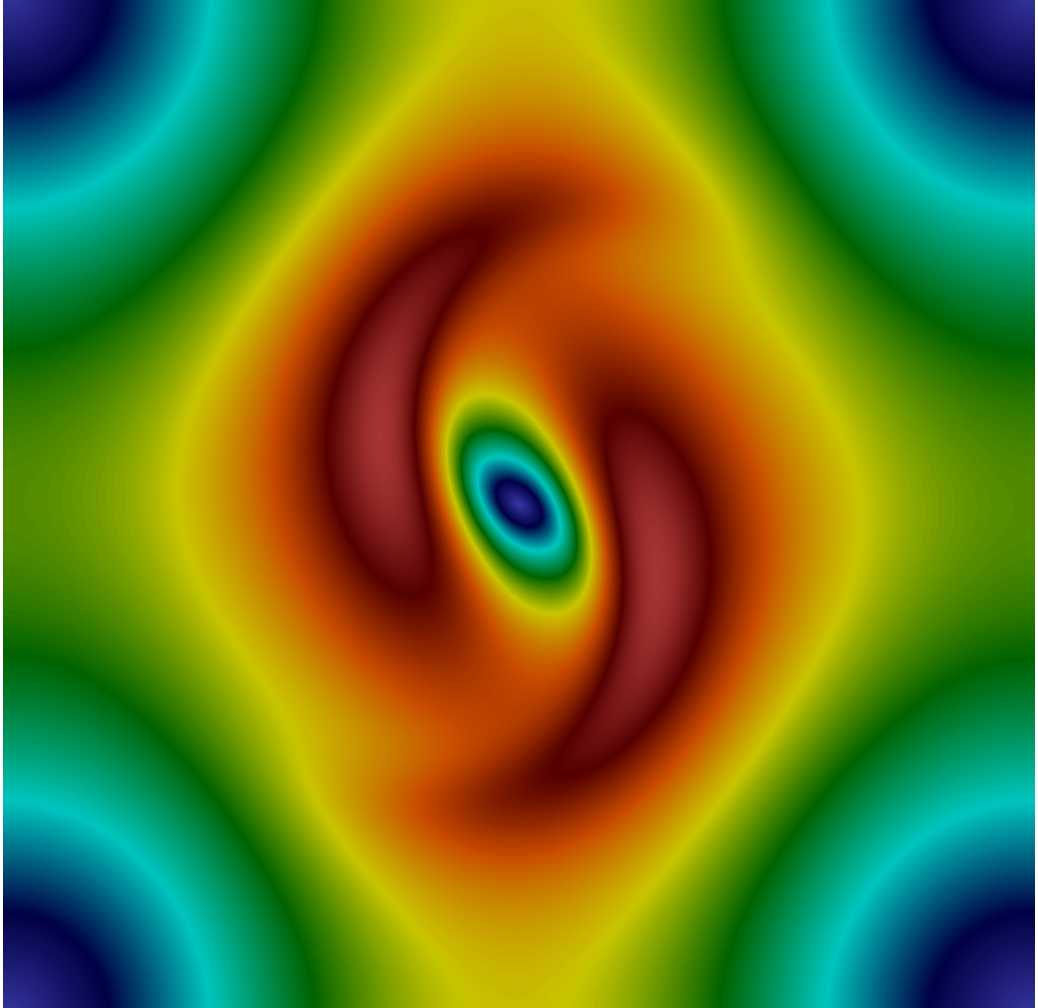}
      \end{overpic}
 \begin{overpic}[width=0.2\textwidth]{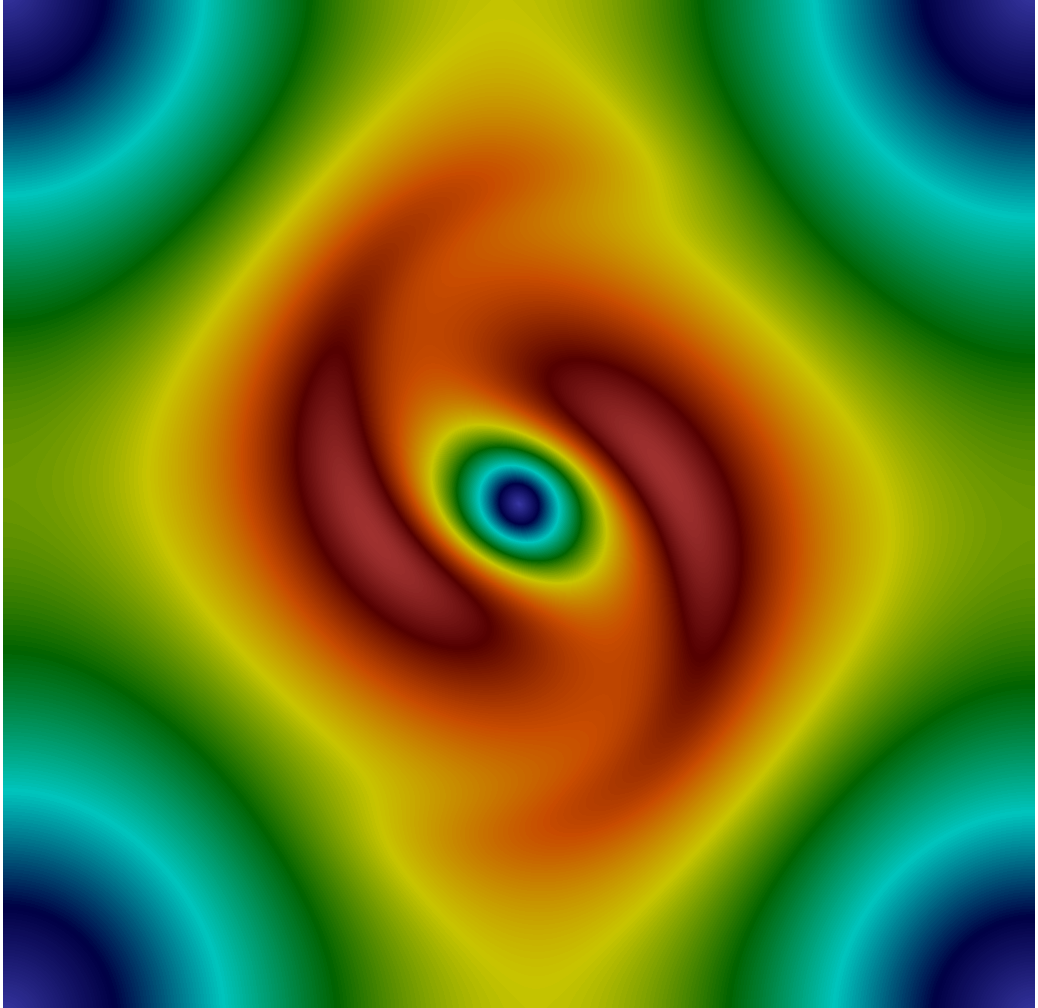}
      \end{overpic}
      \begin{overpic}[width=0.09\textwidth]{img/legendUFOM.png}
      \end{overpic}\\
      \vskip .2cm
      \hspace{.4cm}
      \begin{overpic}[width=0.2\textwidth]{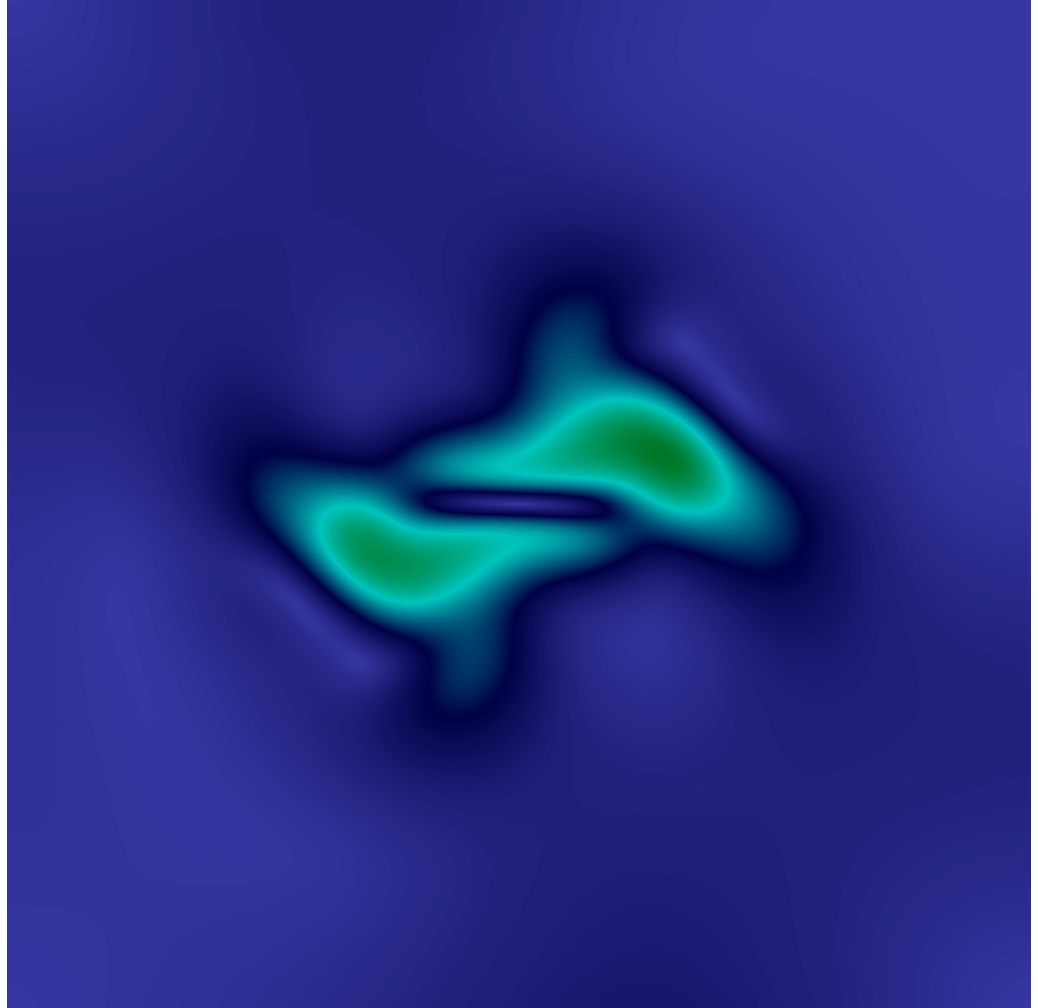}
        \put(-30,50){Diff.}
      \end{overpic}
       \begin{overpic}[width=0.2\textwidth]{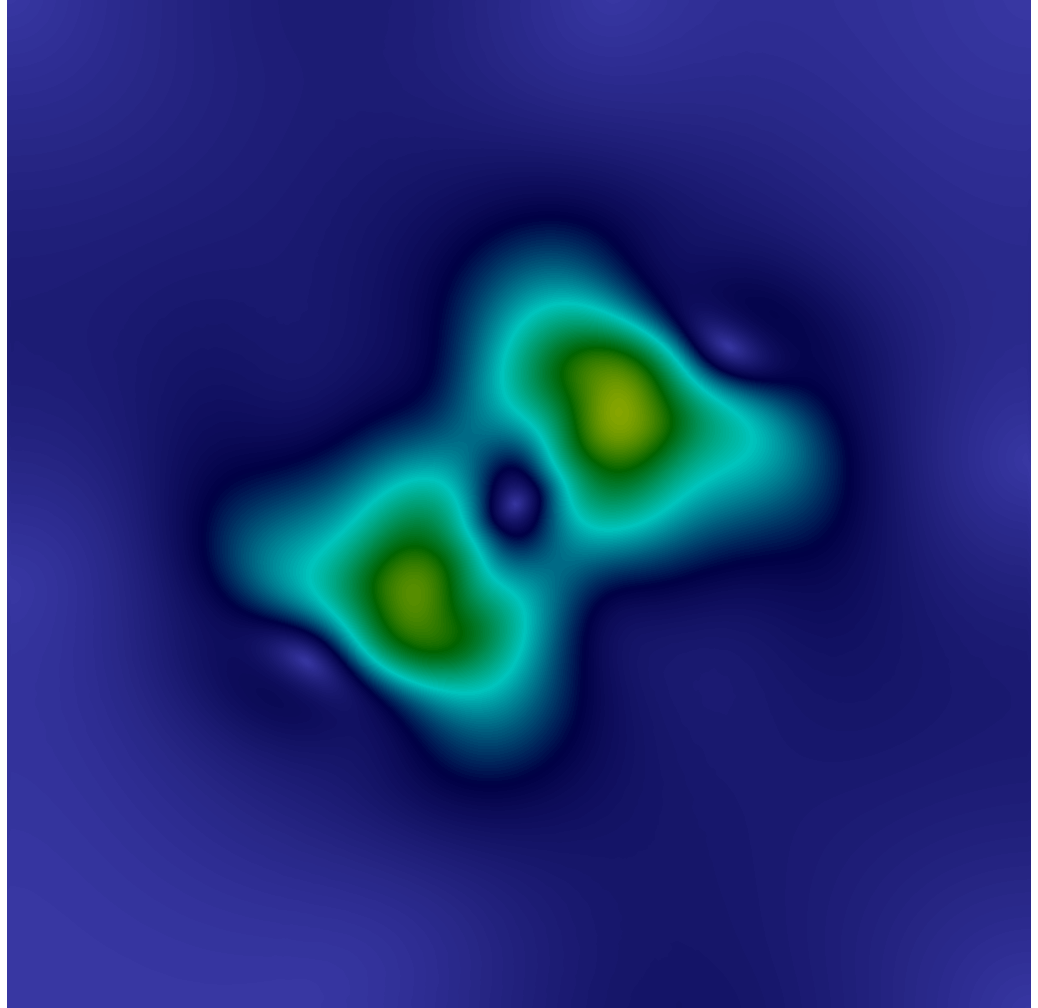}
      \end{overpic}
 \begin{overpic}[width=0.2\textwidth]{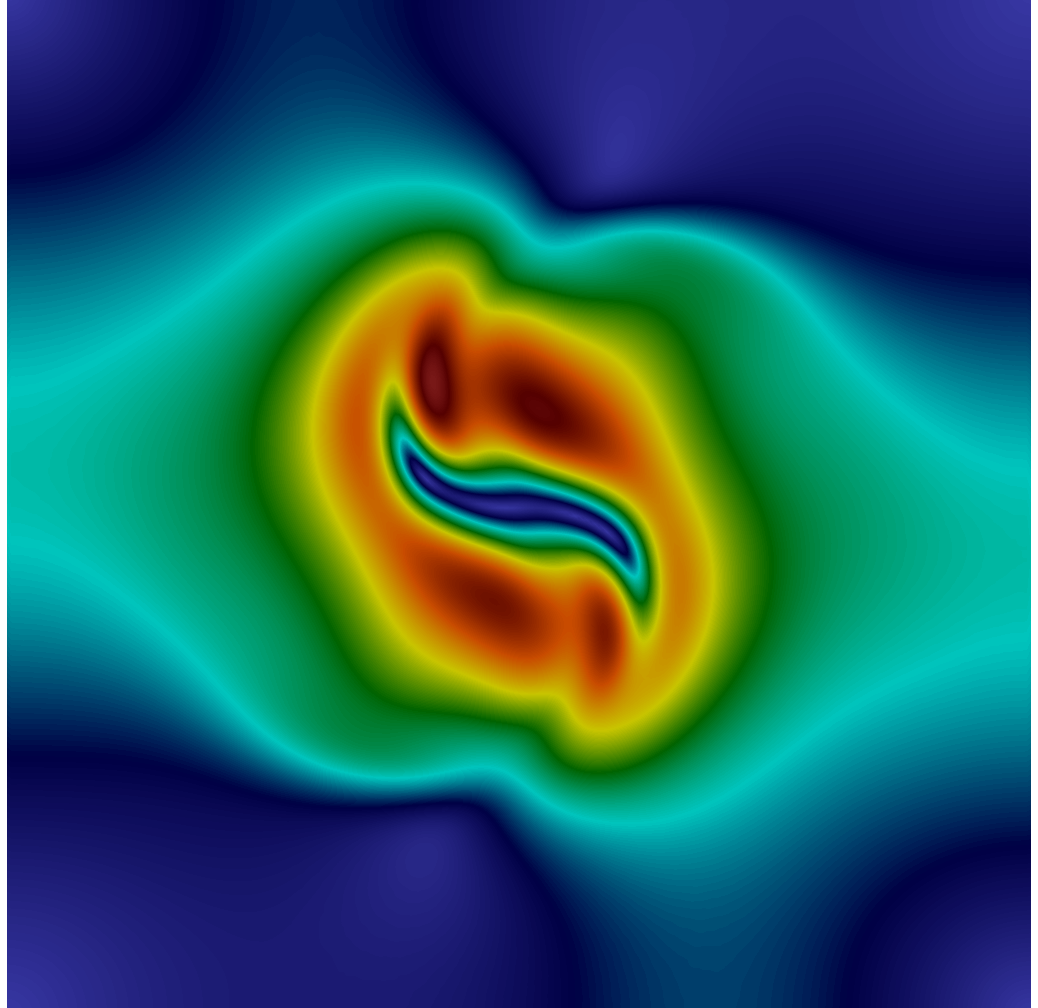}
      \end{overpic}
      \begin{overpic}[width=0.2\textwidth]{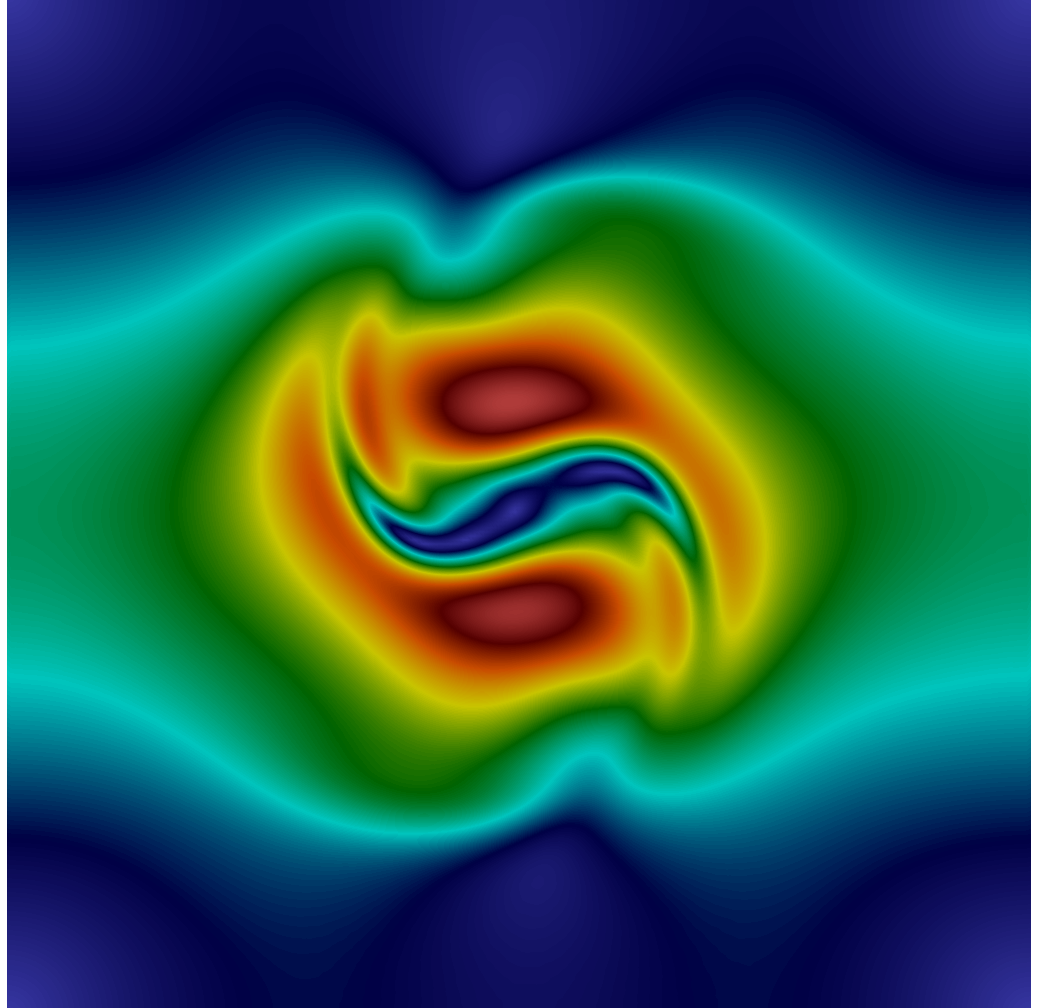}
      \end{overpic}
      \begin{overpic}[width=0.09\textwidth]{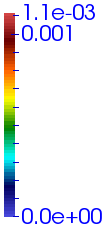}
      \end{overpic}
\caption{FOM validation: velocity $\u$ computed by the solver in velocity-pressure formulation (first row) and stream function-vorticity formulation (second row), and difference between the two fields in absolute value (third row) at $t = 4$ (first column), $t = 8$ (second column), $t = 16$ (third column) and $t = 20$ (fourth column). %\anna{Michele, io metterei una sola legenda per riga e la metterei in verticale dopo le 4 figure. Inoltre aggiugnerei una terza riga con, per es., la differenza in modulo tra i due campi cosi' il confronto risulta piu' quantitativo invece che solo ad occhio. Idem per le figure 4 e 5. Btw, queste figure mi piacciono molto :)}
} %\michele{fatto! :) si, sono molto esotiche :D}
\label{fig:comp_U_FOM}
\end{figure}

\begin{figure}[htb]
\centering
 \hspace{.4cm} \begin{overpic}[width=0.2\textwidth]{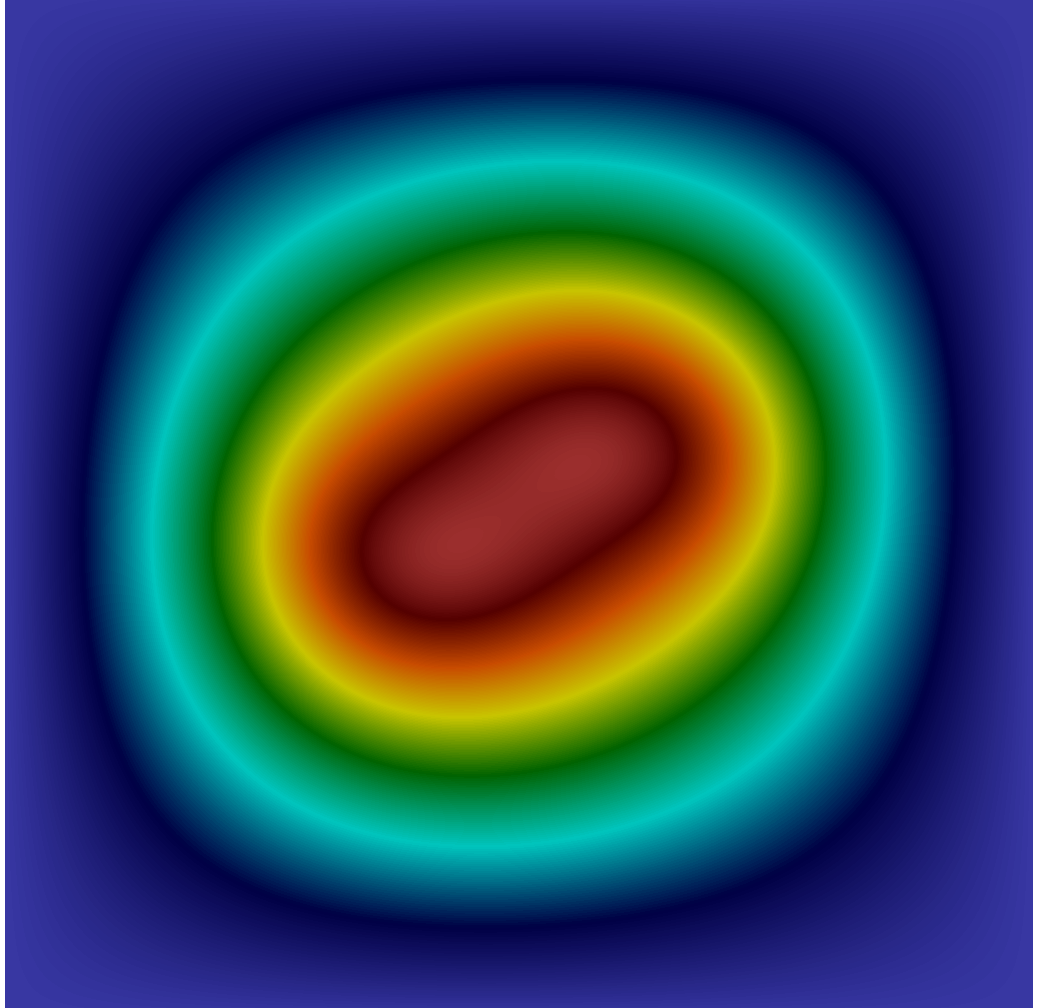}
        \put(38,101){$t = 4$}
        \put(-30,50){NSE}
        \put(-28,40){$\u, p$}
      \end{overpic}
 \begin{overpic}[width=0.2\textwidth]{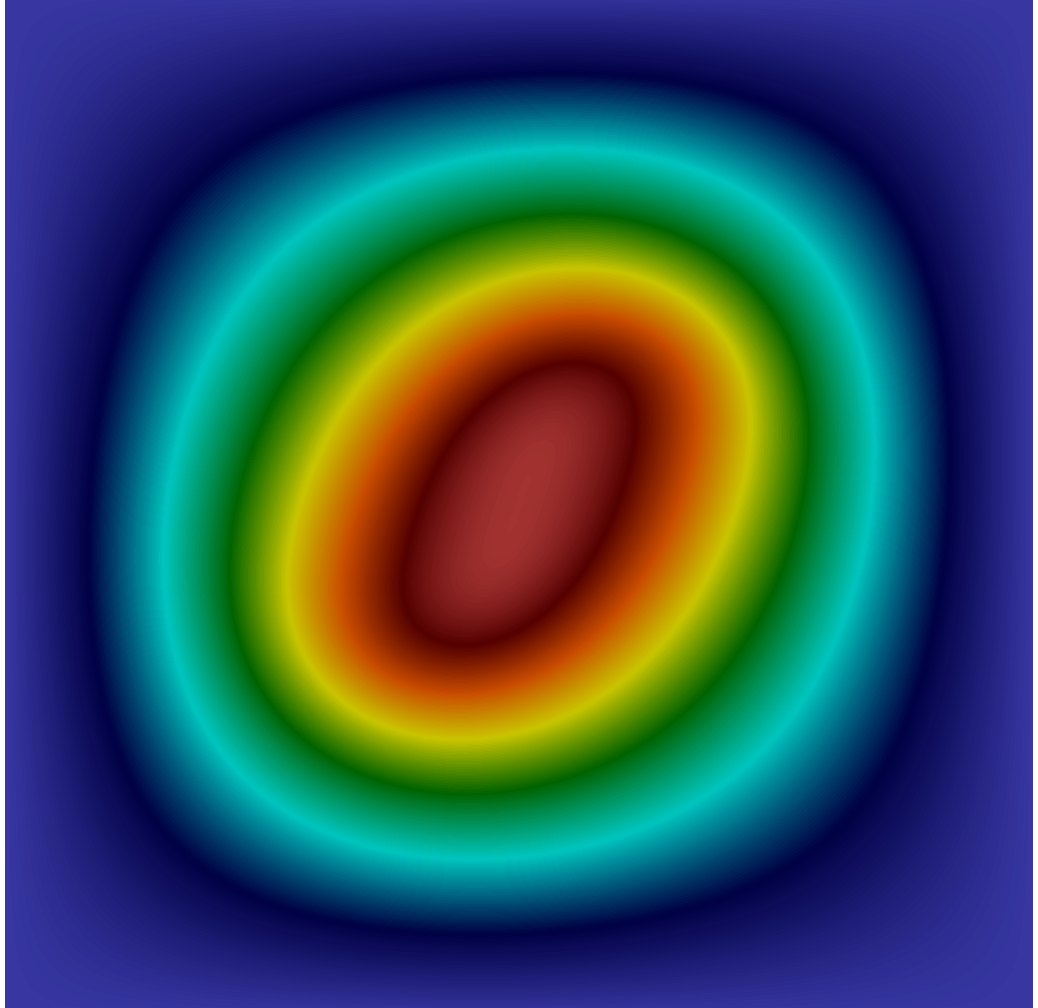}
        \put(38,101){$t = 8$}
      \end{overpic}
 \begin{overpic}[width=0.2\textwidth]{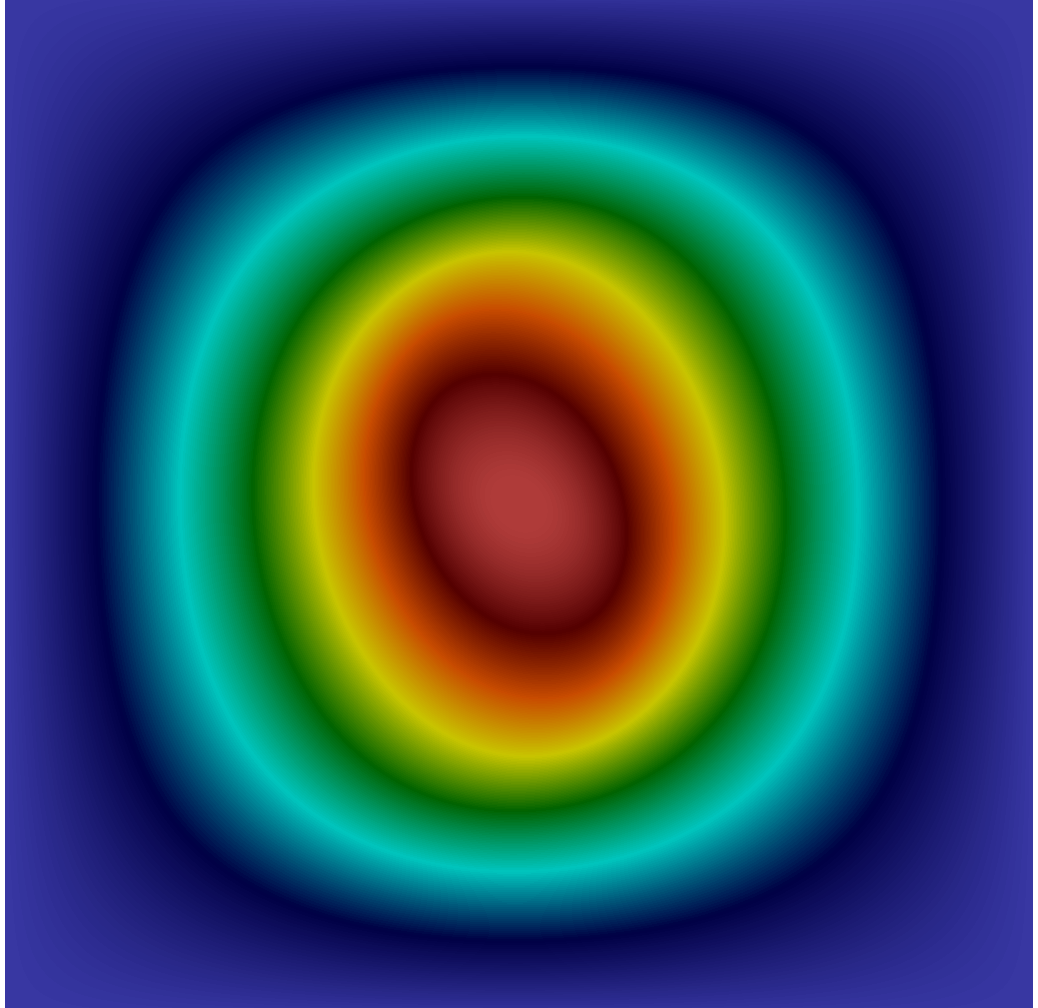}
        \put(35,101){$t = 16$}
      \end{overpic}
 \begin{overpic}[width=0.2\textwidth]{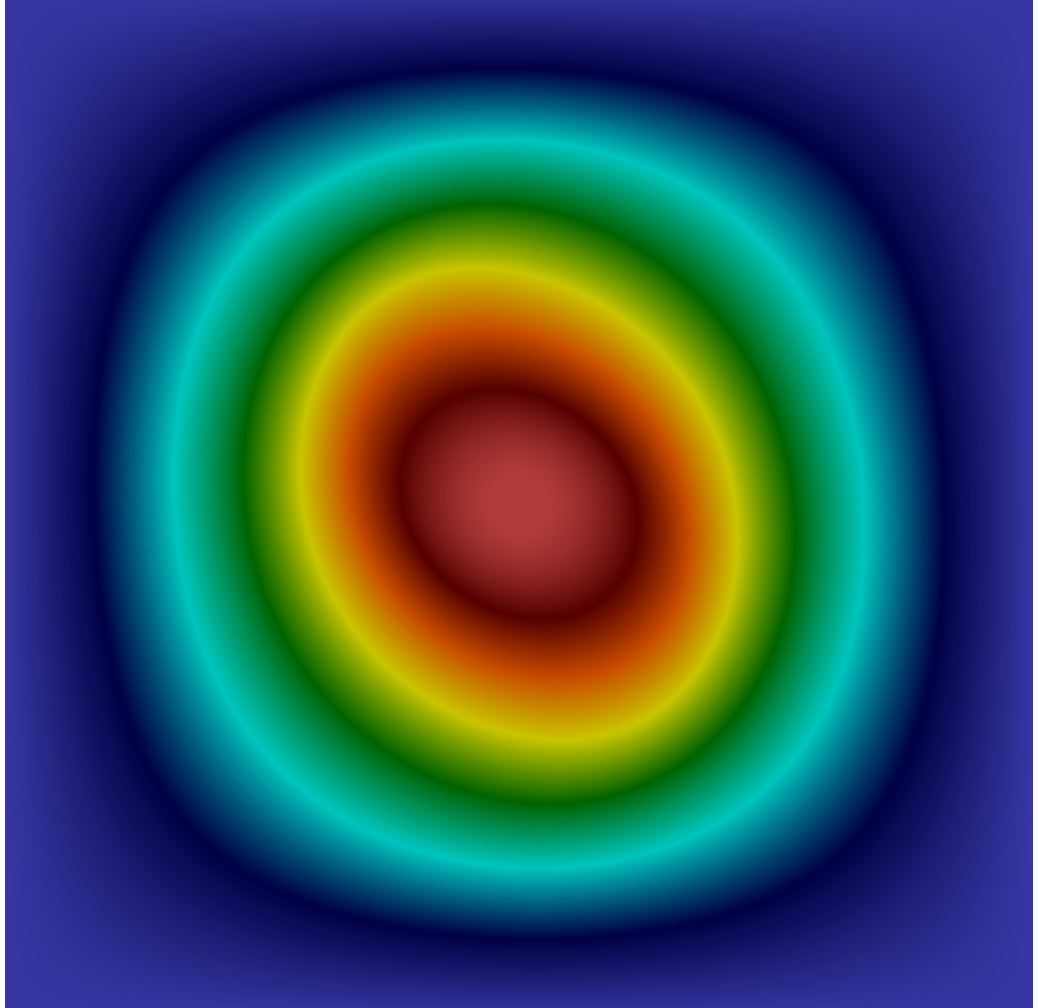}
         \put(35,101){$t = 20$}
      \end{overpic}
       \begin{overpic}[width=0.08\textwidth]{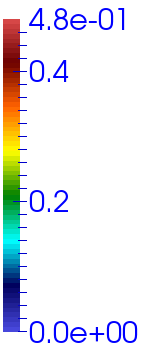}
      \end{overpic}\\
       \vskip .2cm
       \hspace{.4cm}
  \begin{overpic}[width=0.2\textwidth]{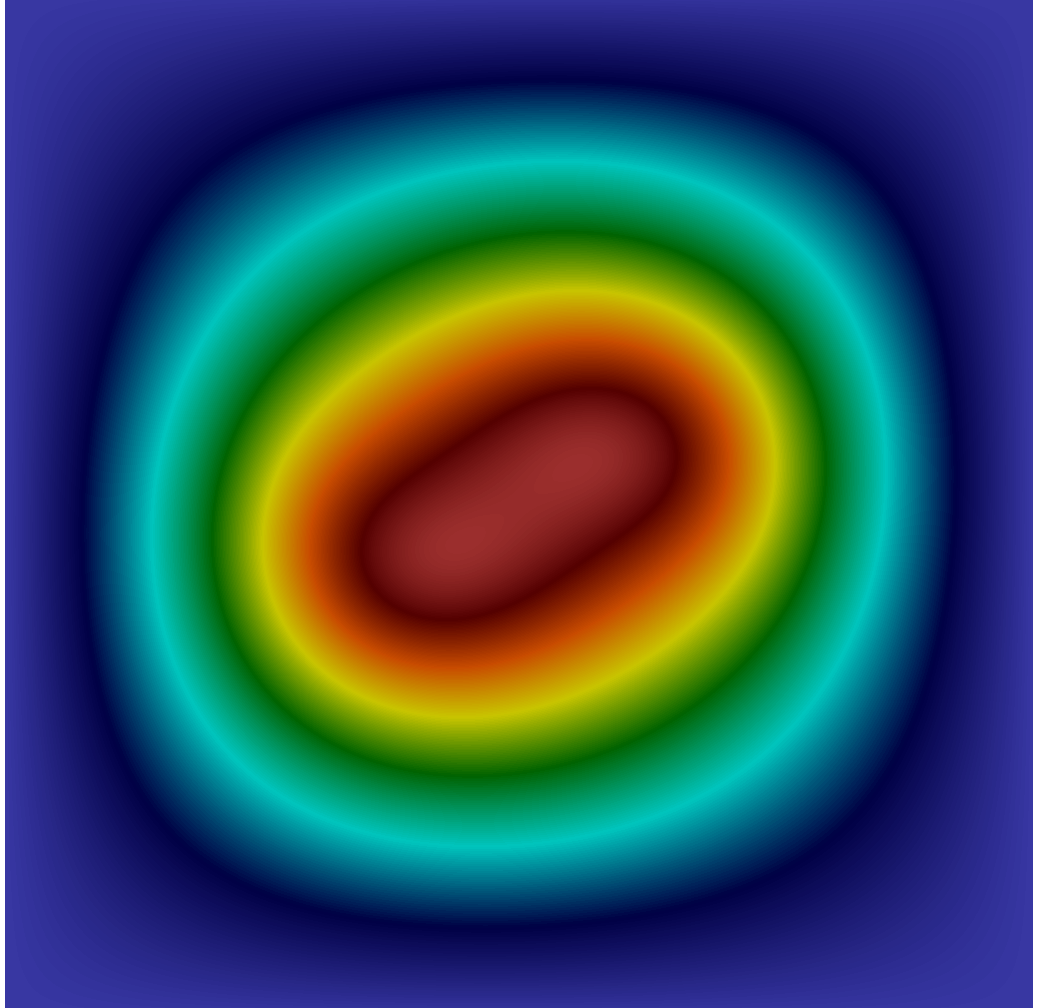}
        \put(-30,50){NSE}
        \put(-28,40){$\psi, \omega$}
      \end{overpic}
 \begin{overpic}[width=0.2\textwidth]{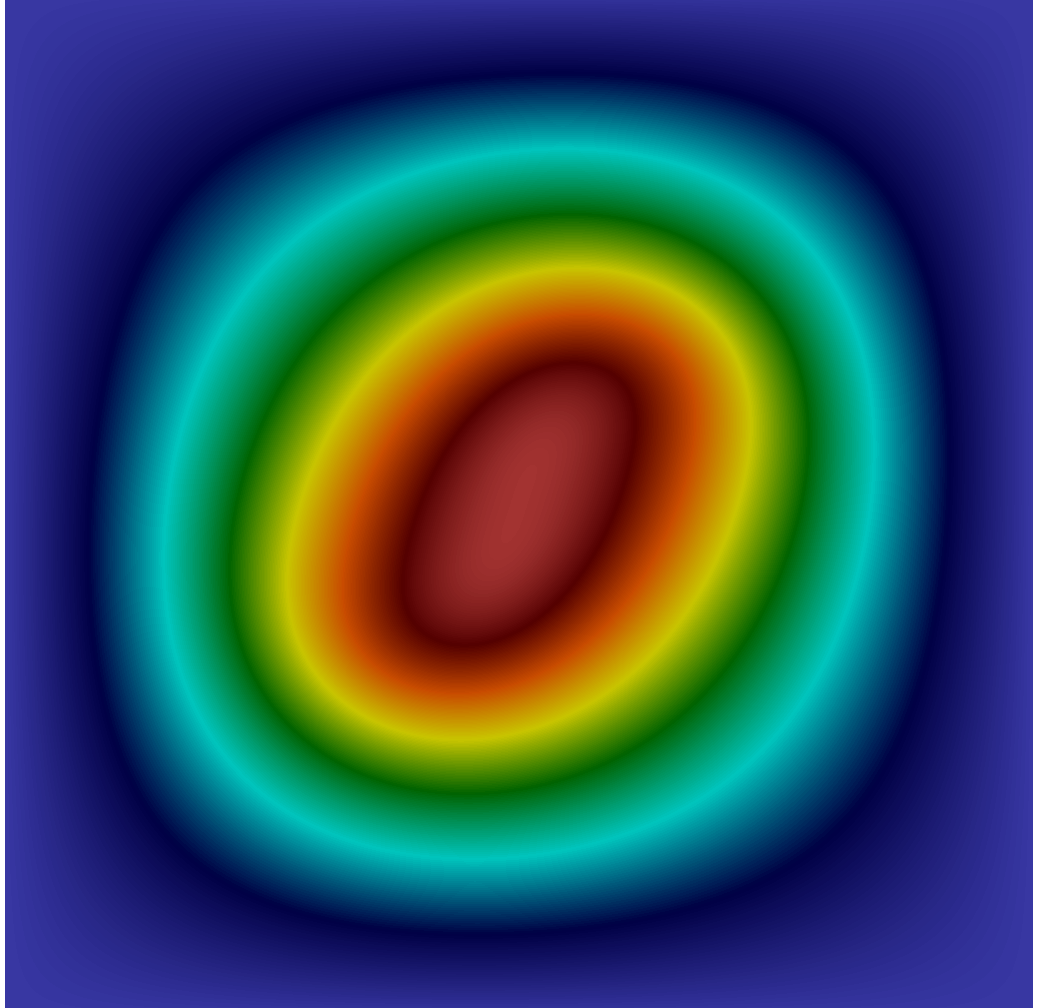}
      \end{overpic}
 \begin{overpic}[width=0.2\textwidth]{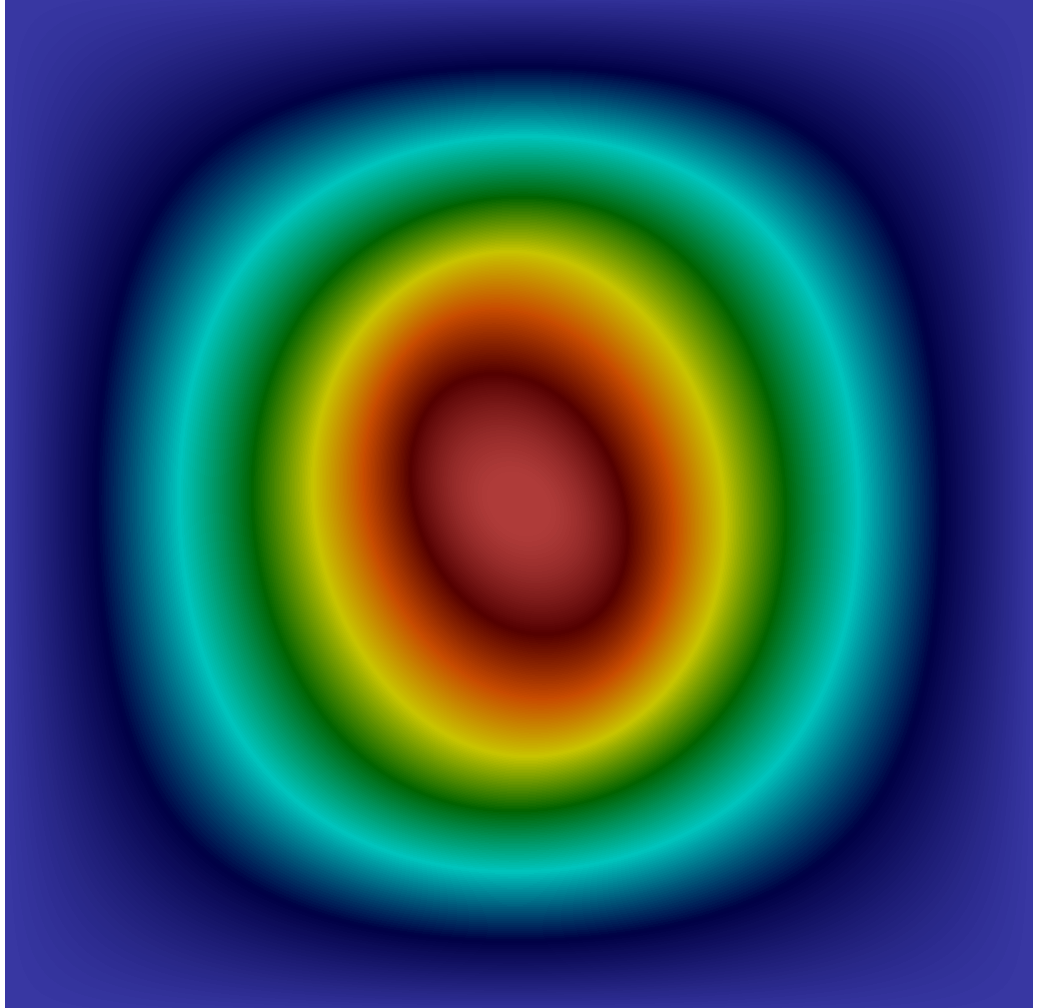}
      \end{overpic}
 \begin{overpic}[width=0.2\textwidth]{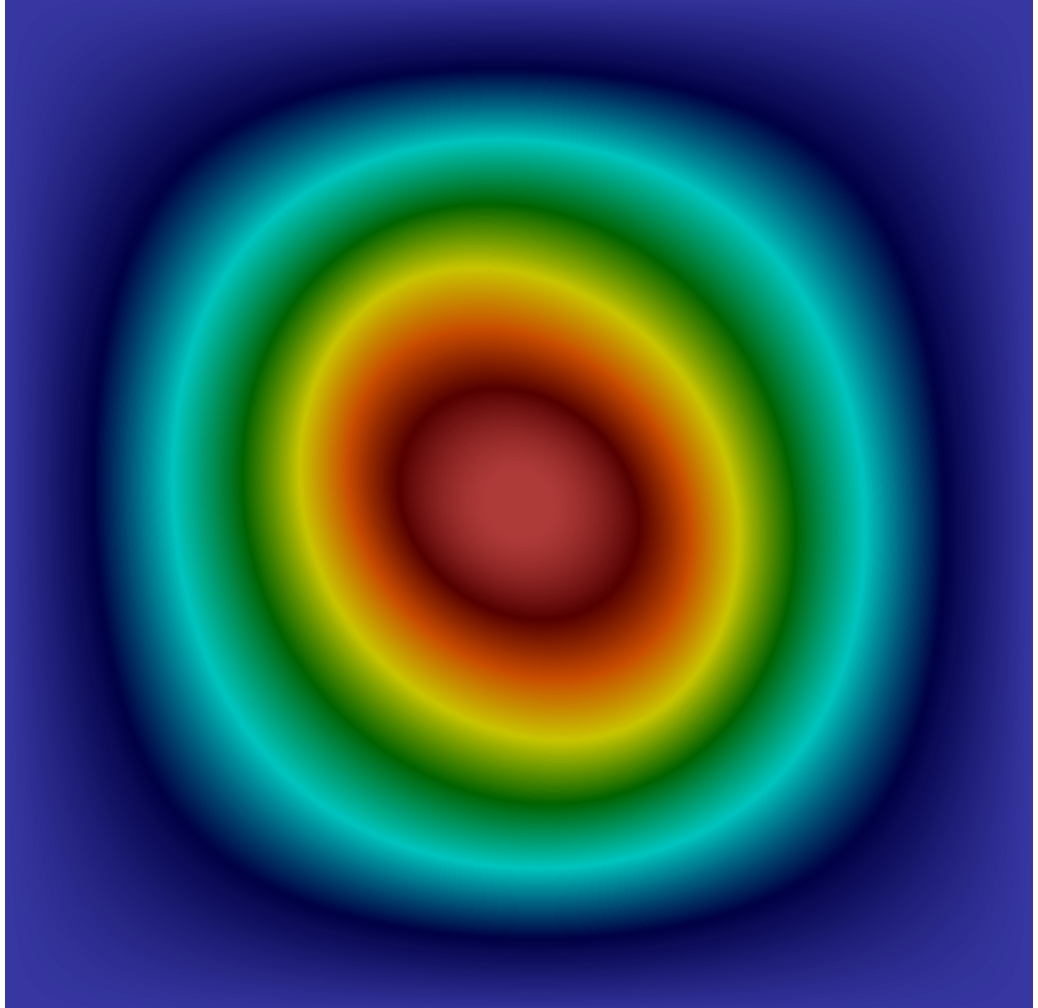}
      \end{overpic}
      \begin{overpic}[width=0.08\textwidth]{img/legendPsiFOM.png}
      \end{overpic}\\
      \vskip .2cm
      \hspace{.4cm}
 \begin{overpic}[width=0.2\textwidth]{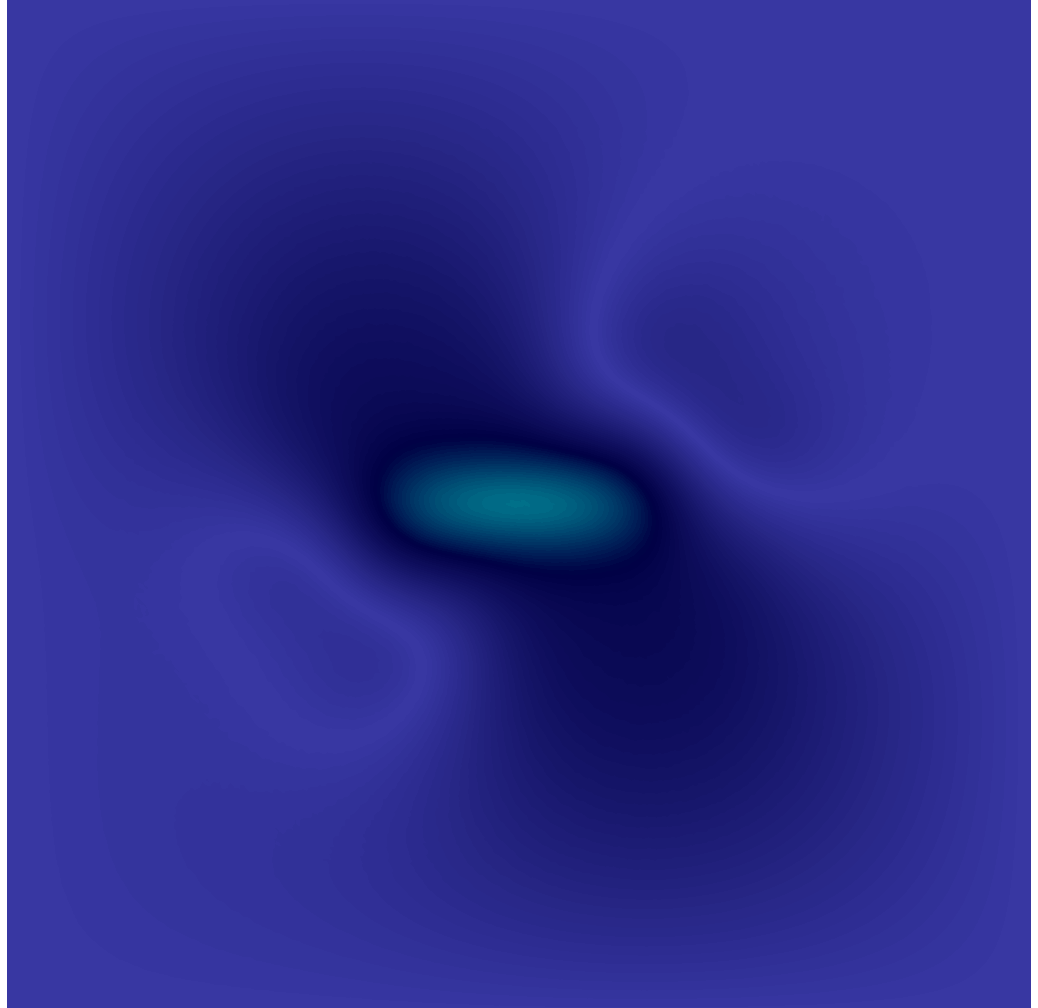}
        \put(-30,50){Diff.}
      \end{overpic}
 \begin{overpic}[width=0.2\textwidth]{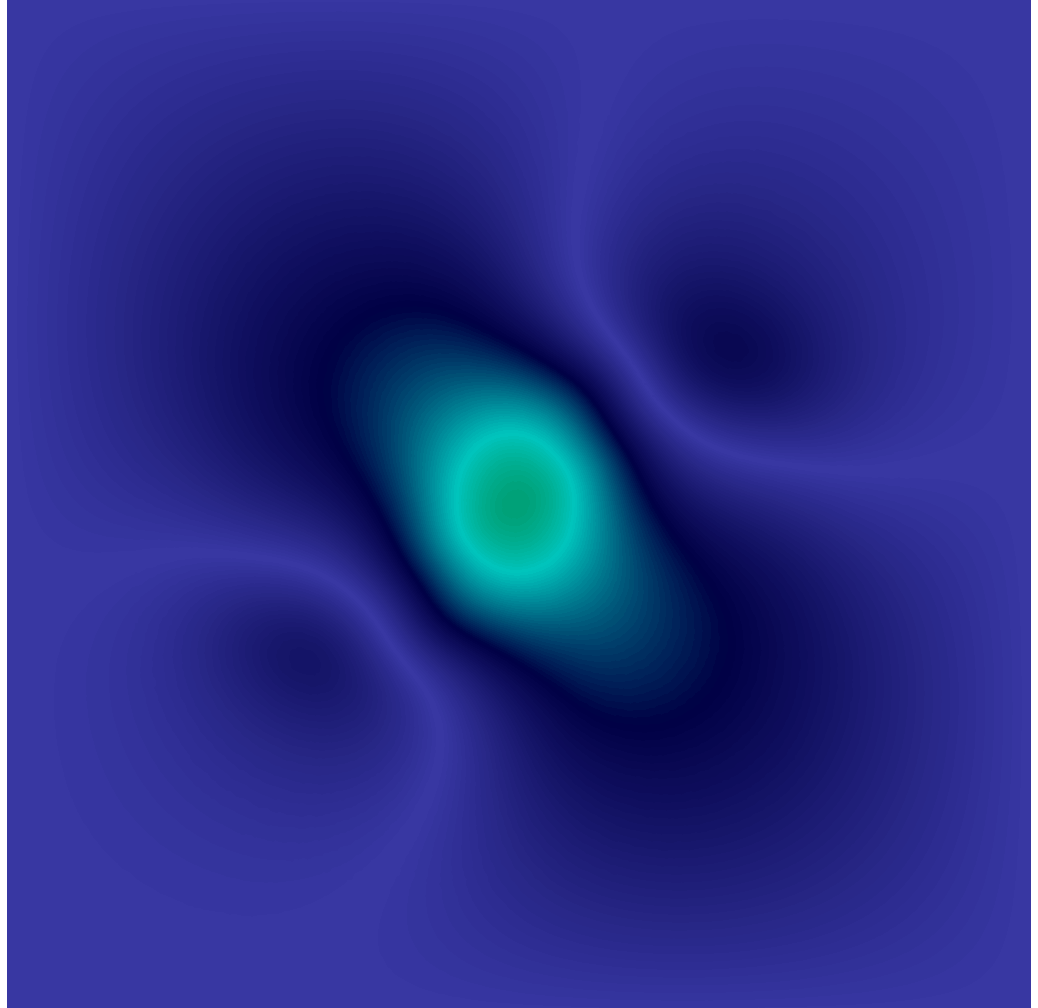}
      \end{overpic}
 \begin{overpic}[width=0.2\textwidth]{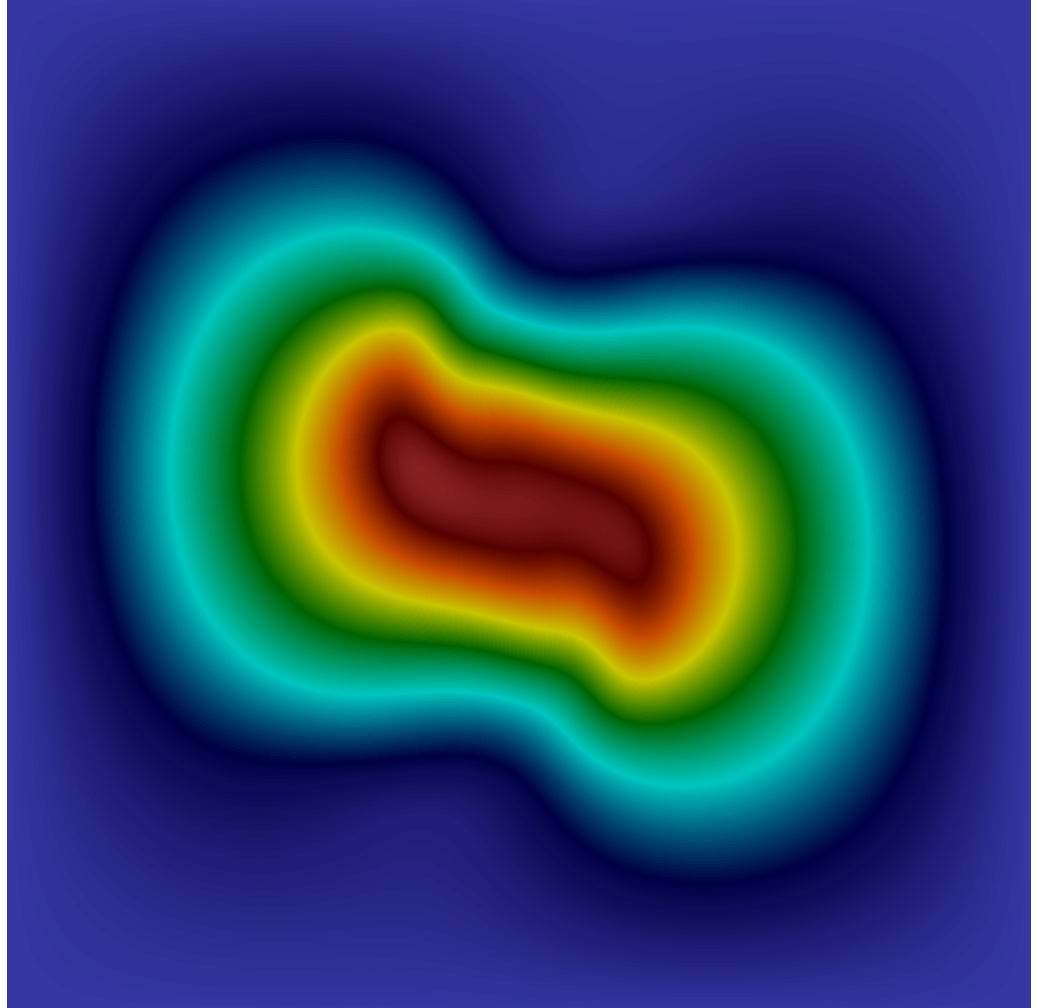}
      \end{overpic}
      \begin{overpic}[width=0.2\textwidth]{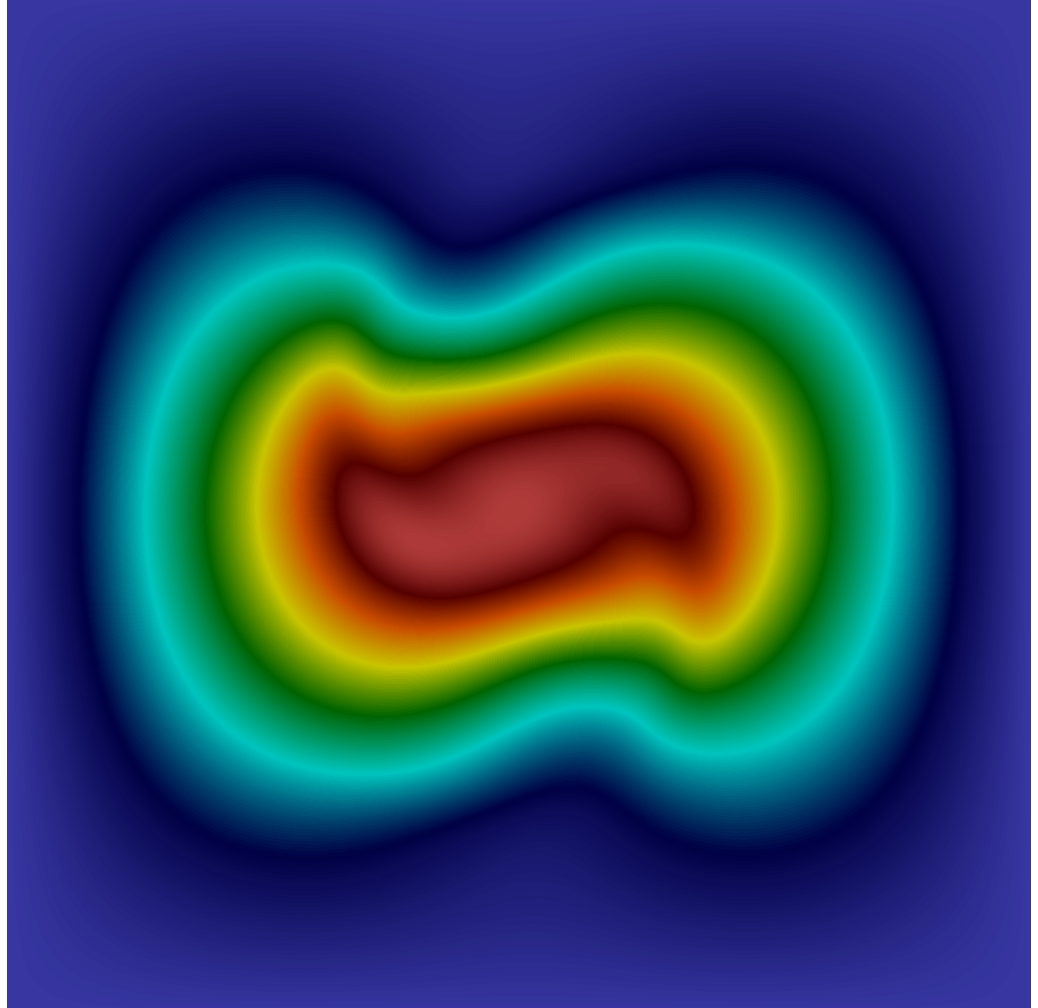}
      \end{overpic}
      \begin{overpic}[width=0.09\textwidth]{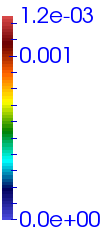}
      \end{overpic}
\caption{FOM validation: stream function $\psi$ computed by the solver in velocity-pressure formulation (first row) and stream function-vorticity formulation (second row), and difference between the two fields in absolute value (third row)  at $t = 4$ (first column), $t = 8$ (second column), $t = 12$ (third column) and $t = 20$ (fourth column).}
\label{fig:comp_psi_FOM}
\end{figure}

\begin{figure}[htb]
\centering
\hspace{.4cm}
 \begin{overpic}[width=0.2\textwidth]{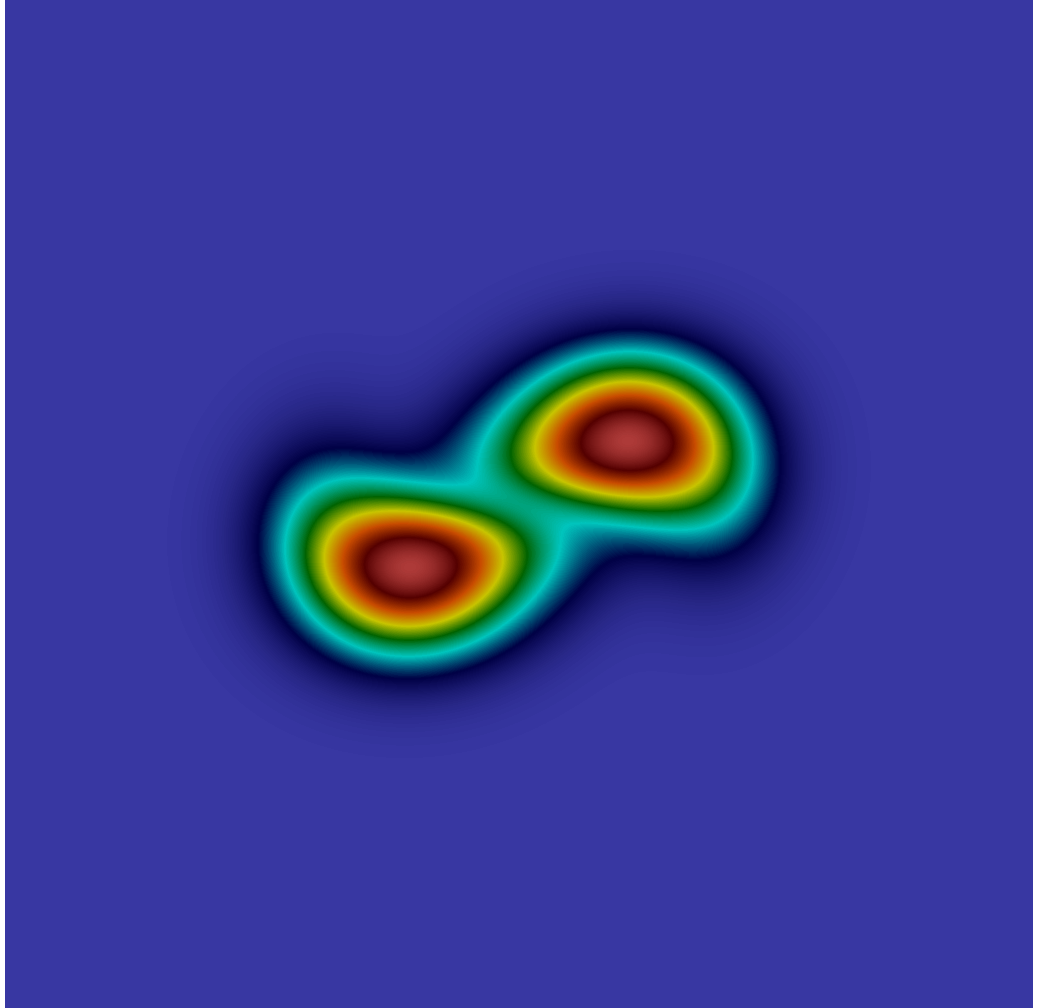}
        \put(38,101){$t = 4$}
        \put(-30,50){NSE}
        \put(-28,40){$\u, p$}
      \end{overpic}
 \begin{overpic}[width=0.2\textwidth]{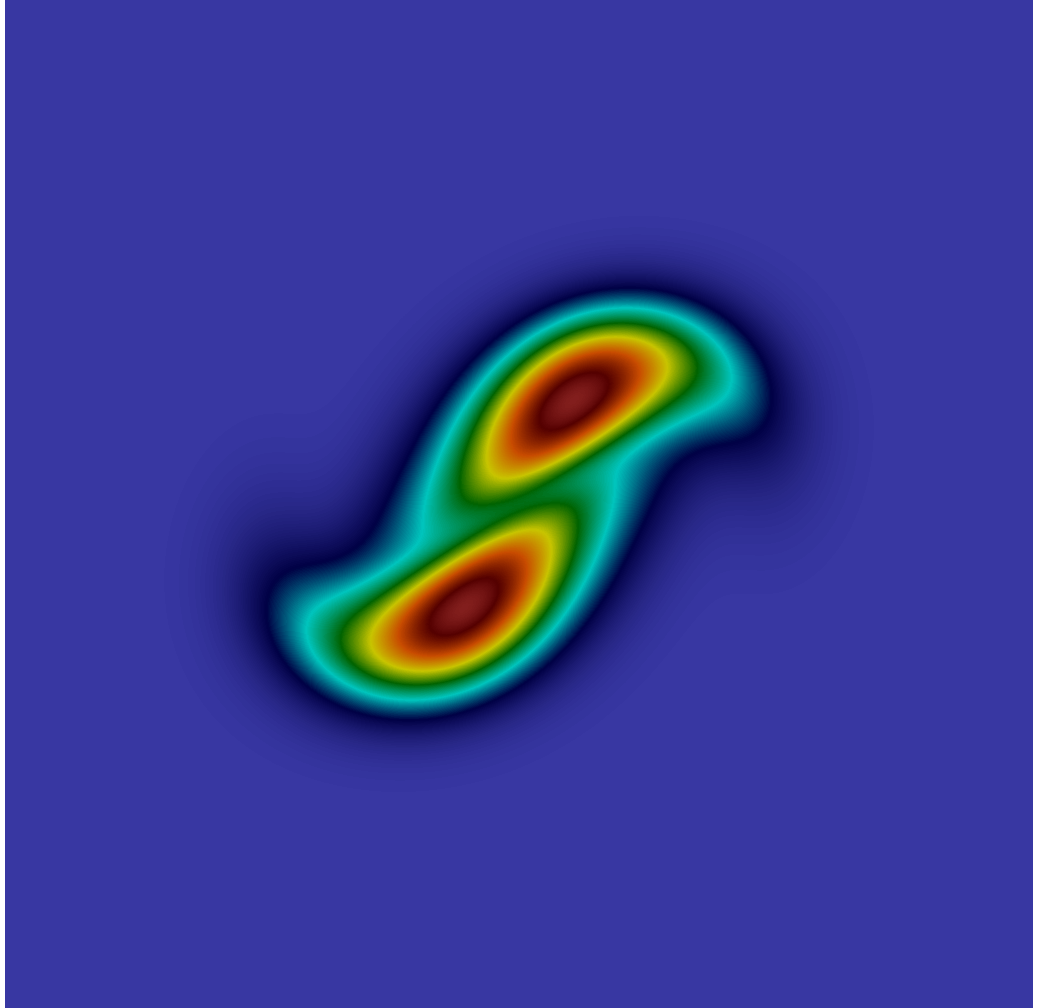}
        \put(38,101){$t = 8$}
      \end{overpic}
 \begin{overpic}[width=0.2\textwidth]{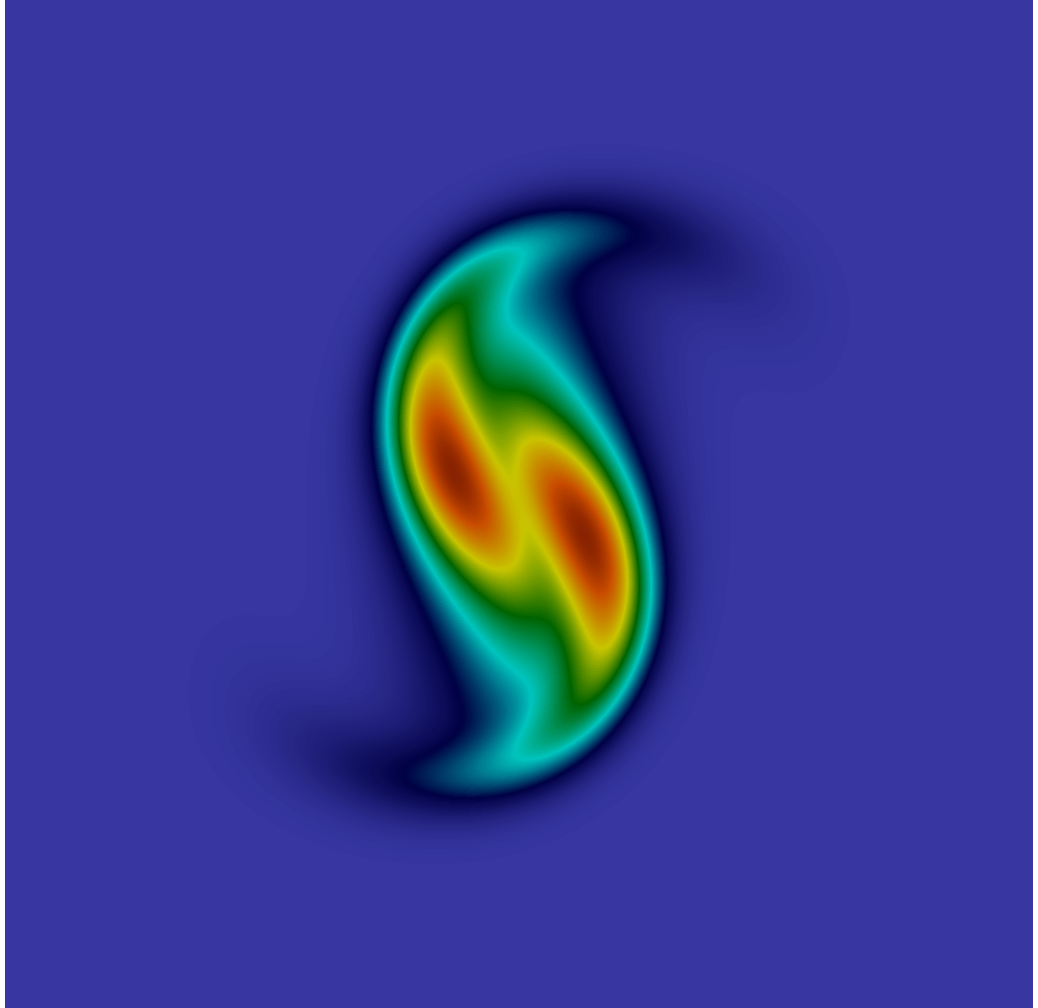}
        \put(35,101){$t = 16$}
      \end{overpic}
 \begin{overpic}[width=0.2\textwidth]{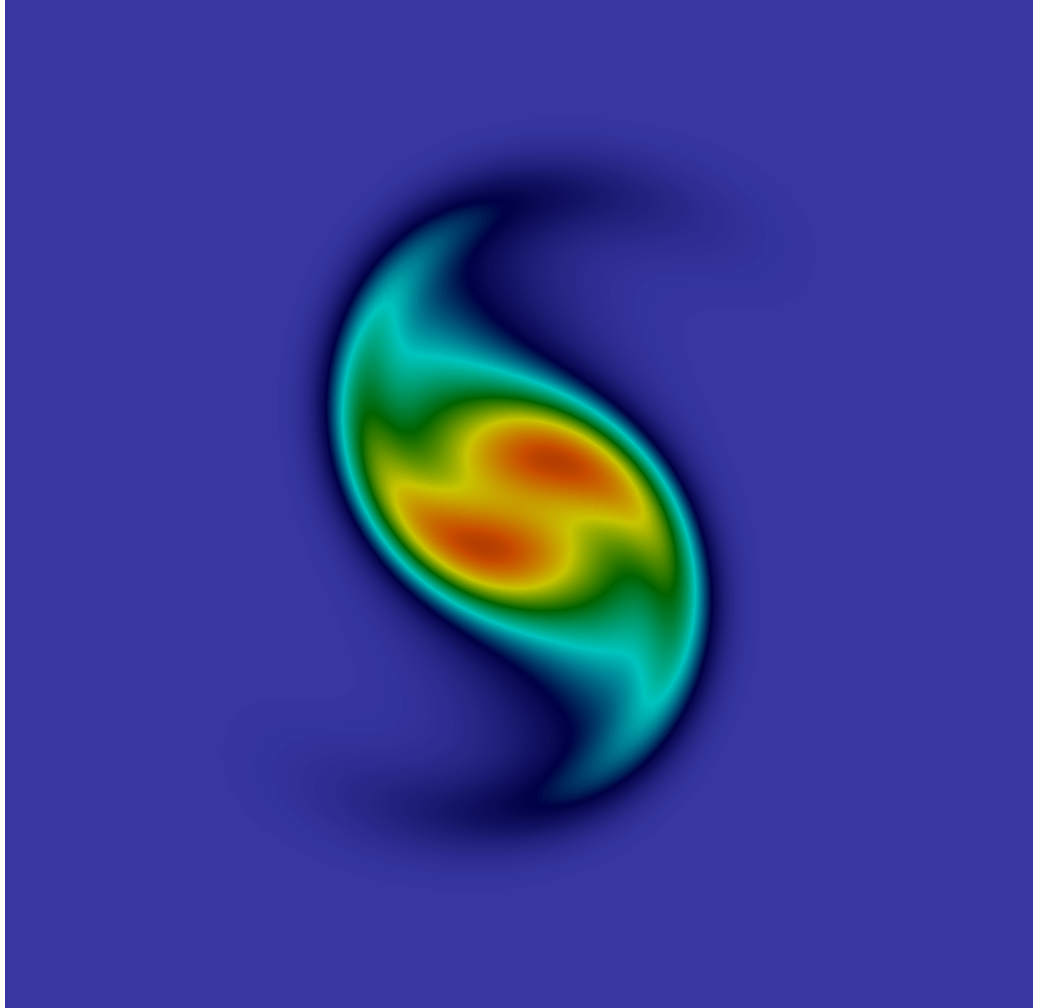}
         \put(35,101){$t = 20$}
      \end{overpic}
      \begin{overpic}[width=0.08\textwidth]{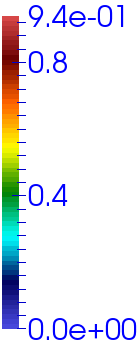}
      \end{overpic}\\
      \vskip .2cm
      \hspace{.4cm}
 \begin{overpic}[width=0.2\textwidth]{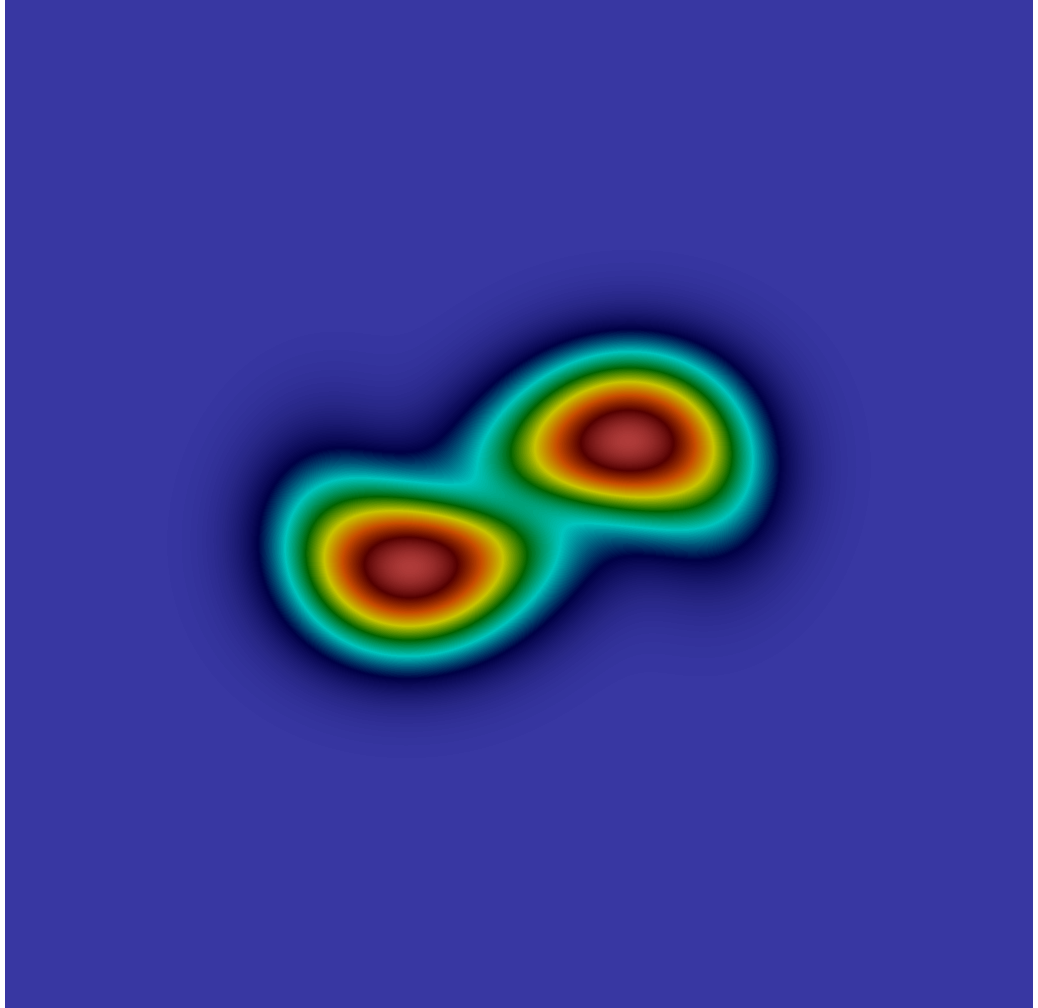}
        \put(-30,50){NSE}
        \put(-28,40){$\psi, \omega$}
      \end{overpic}
 \begin{overpic}[width=0.2\textwidth]{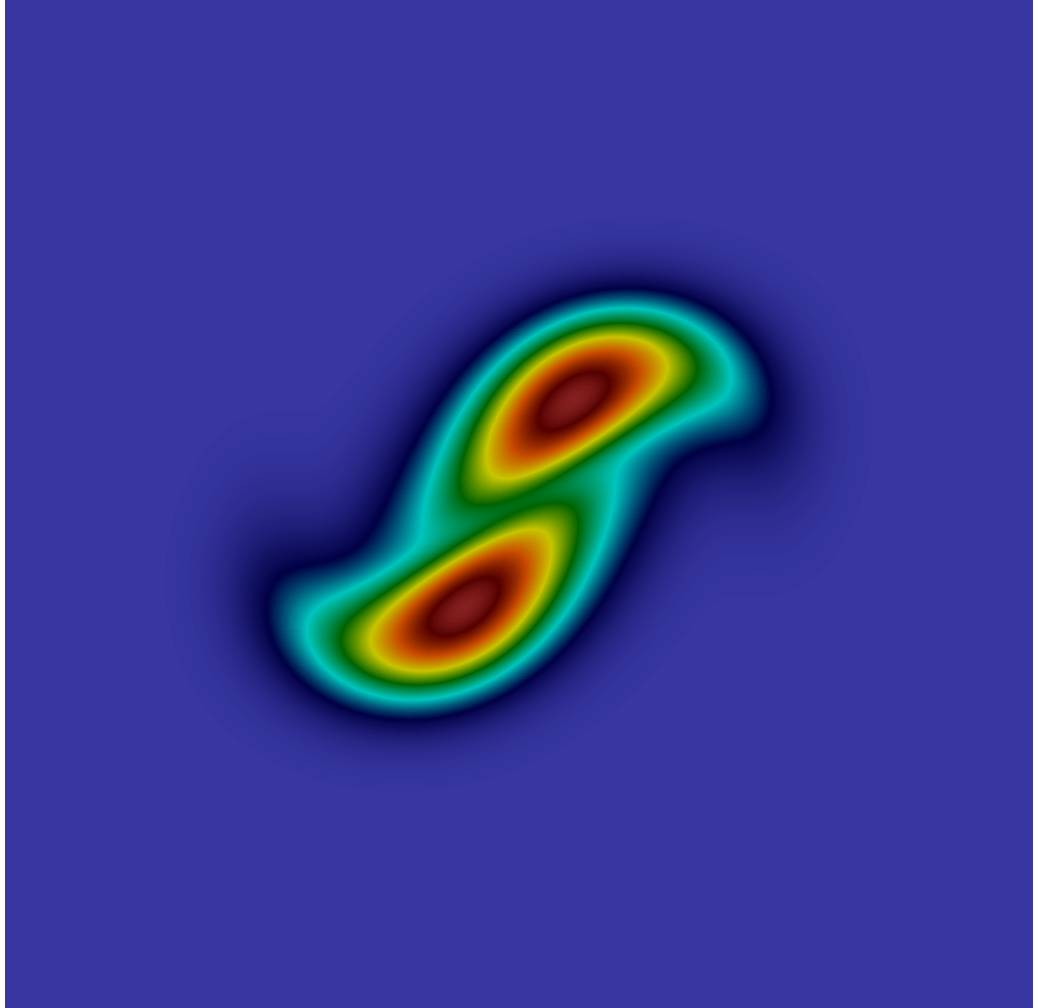}
      \end{overpic}
 \begin{overpic}[width=0.2\textwidth]{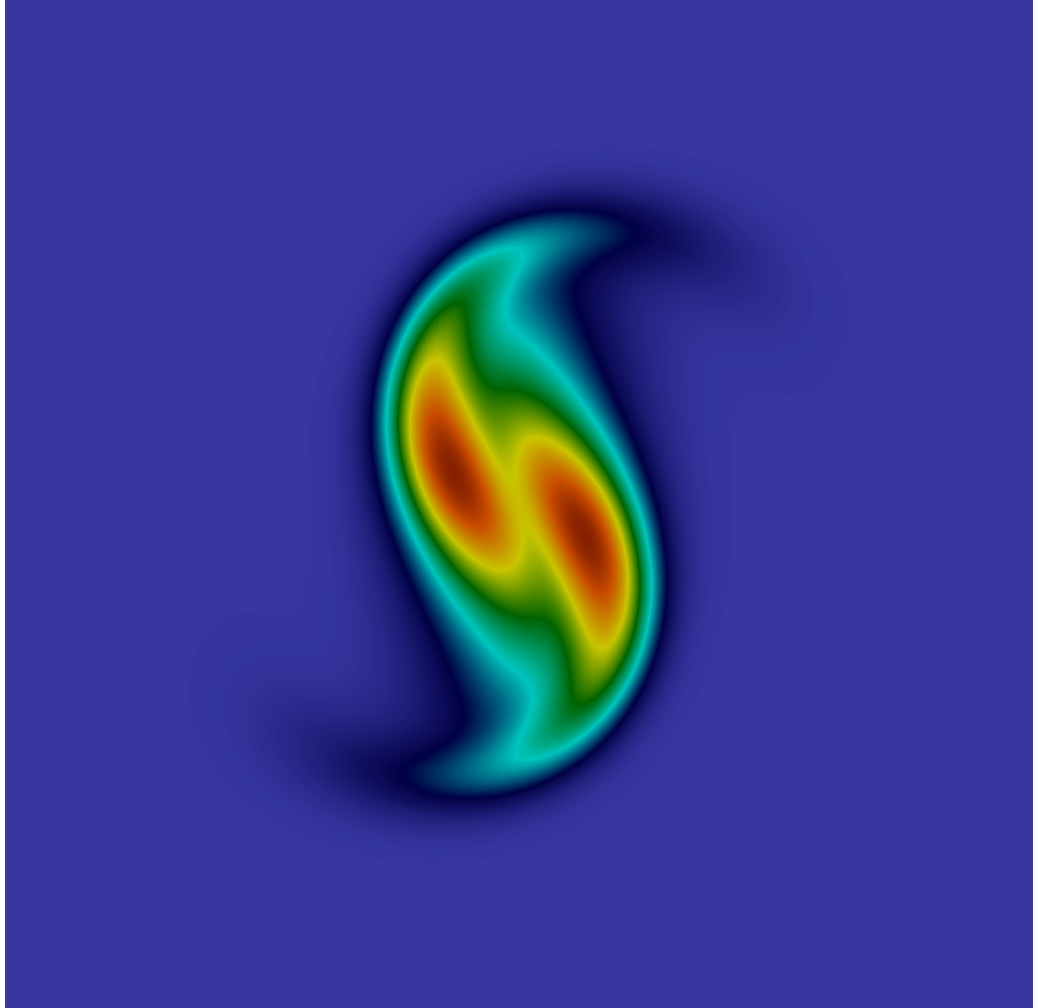}
      \end{overpic}
 \begin{overpic}[width=0.2\textwidth]{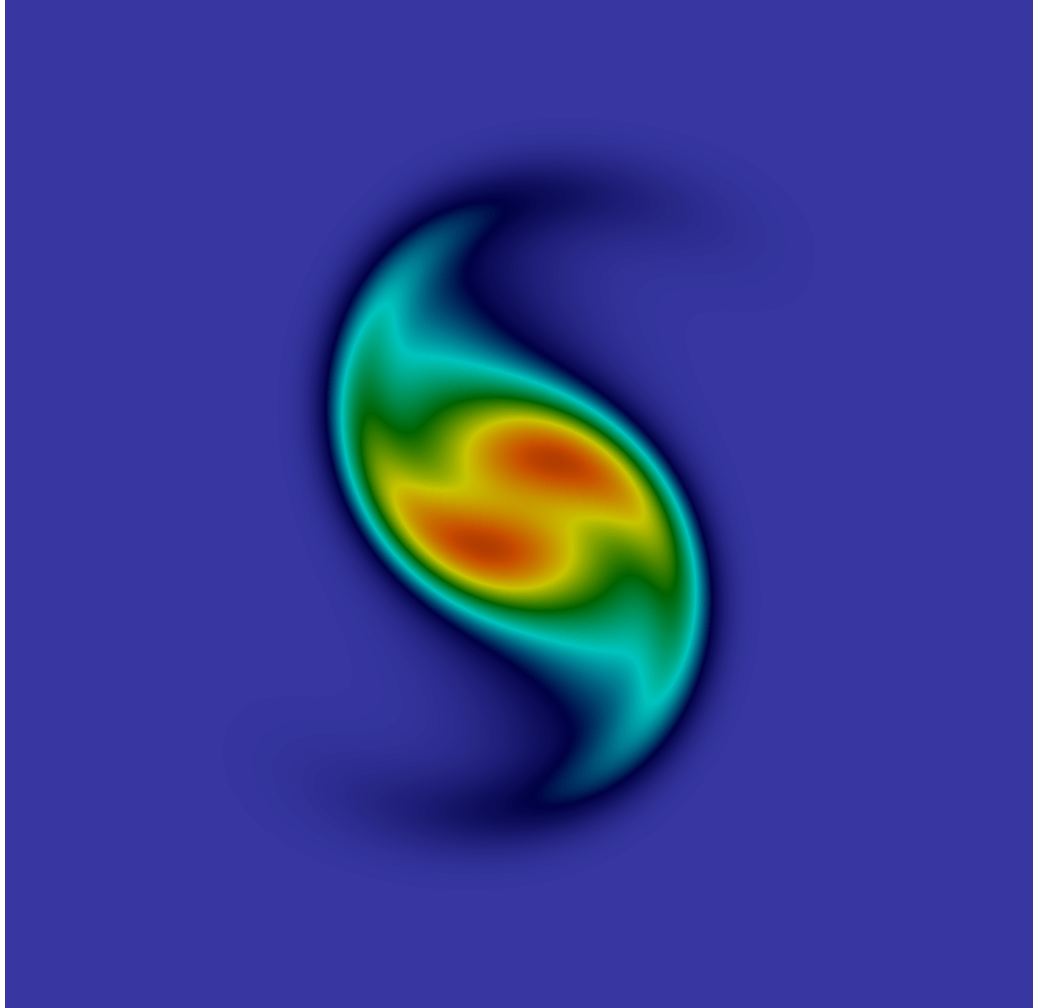}
      \end{overpic}
      \begin{overpic}[width=0.08\textwidth]{img/legendOmegaFOM.png}
      \end{overpic}\\
      \hspace{.4cm}
       ~\begin{overpic}[width=0.2\textwidth]{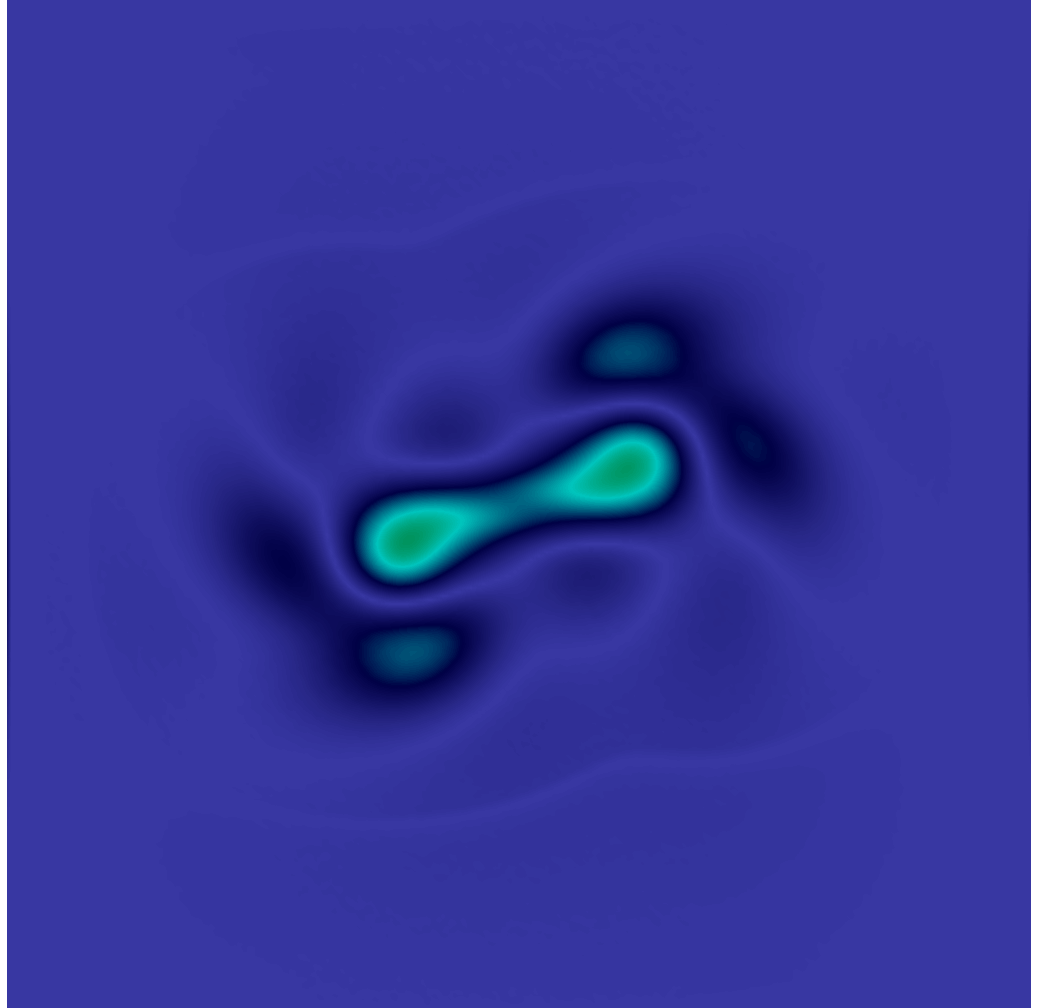}             \put(-30,50){Diff.}
      \end{overpic}
 \begin{overpic}[width=0.2\textwidth]{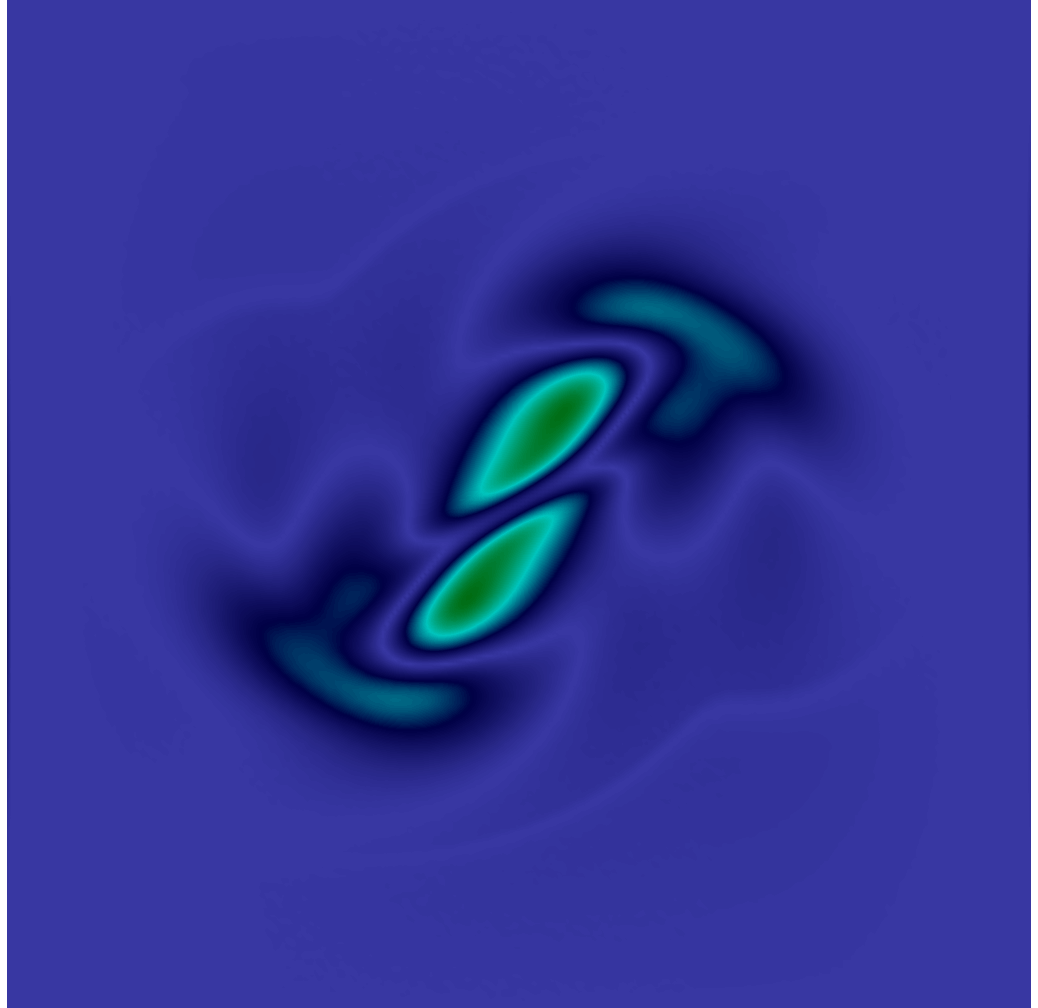}
      \end{overpic}
 \begin{overpic}[width=0.2\textwidth]{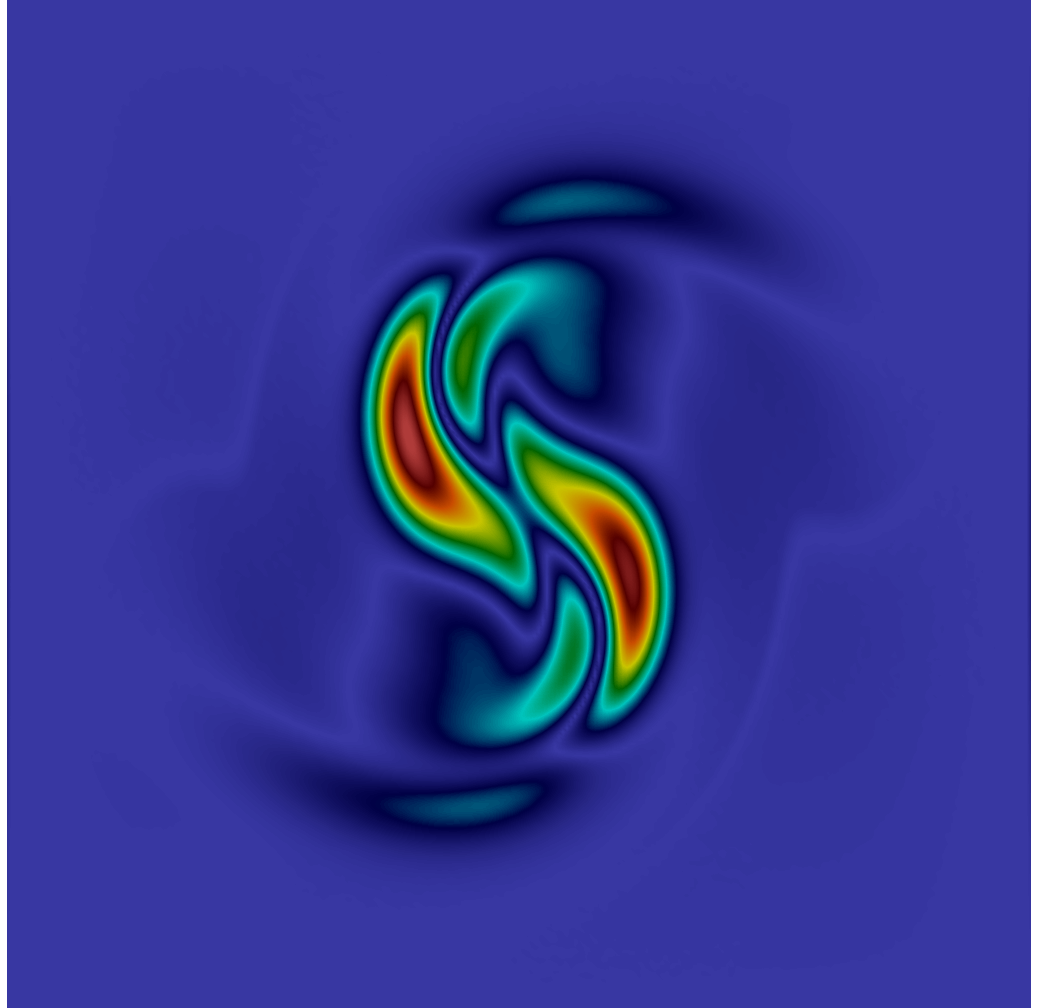}
      \end{overpic}
      \begin{overpic}[width=0.2\textwidth]{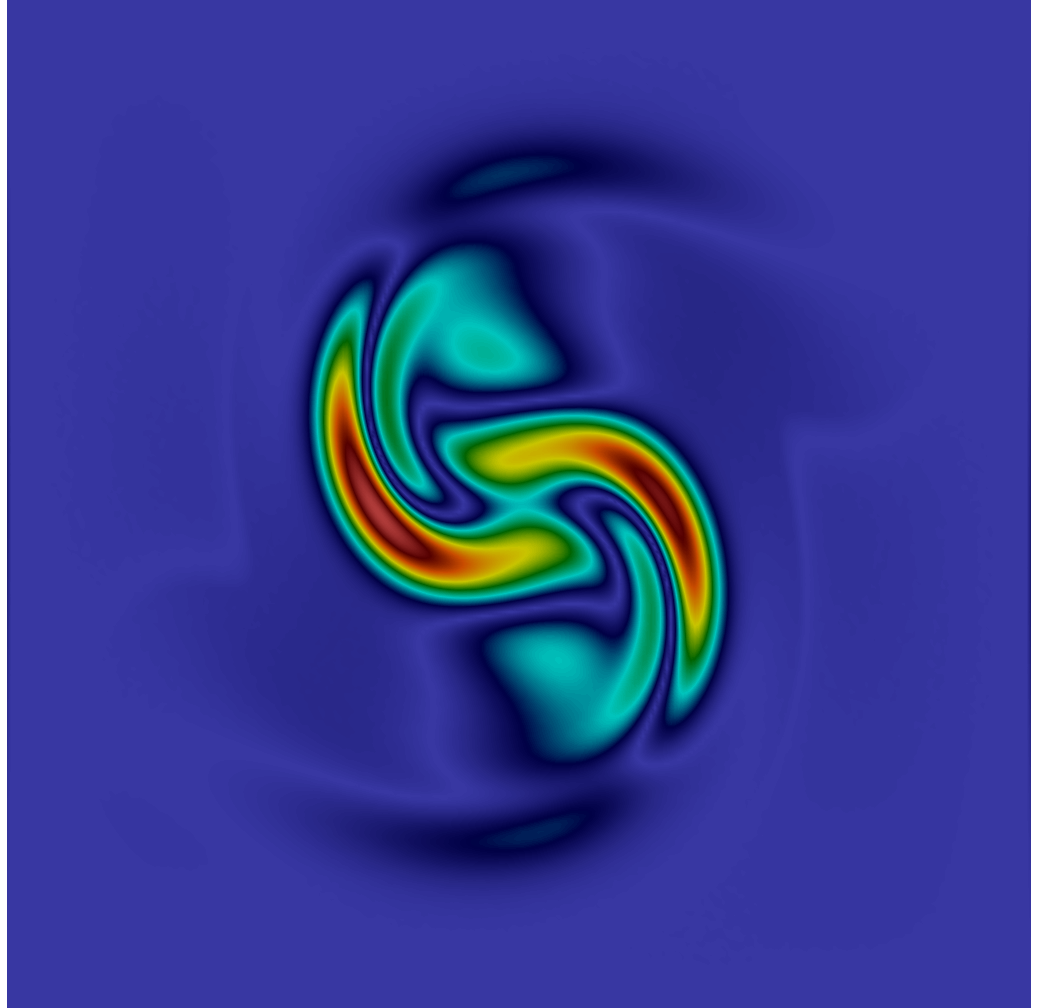}
      \end{overpic}
      \begin{overpic}[width=0.09\textwidth]{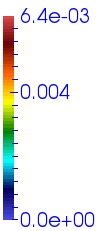}
      \end{overpic}\\
\caption{FOM validation: vorticity $\omega$ computed by the solver in velocity-pressure formulation (first row) and stream function-vorticity formulation (second row), and difference between the two fields in absolute value (third row) at $t = 4$ (first column), $t = 8$ (second column), $t = 16$ (third column) and $t = 20$ (fourth column).}
\label{fig:comp_zeta_FOM}
\end{figure}

Figures \ref{fig:comp_U_FOM}, \ref{fig:comp_psi_FOM}, and \ref{fig:comp_zeta_FOM} display a qualitative comparison
in terms of $\u$, $\psi$, and $\omega$ computed by the solvers in primitive variables and in stream function-vorticity formulation at four different times. %: $t = 4$, $t = 8$, $t = 16$, and $t = 20$. 
As we can see from these figures, the solutions are very close to each other with the maximum relative difference in absolute value not exceeding $4.4e-3$ for $\u$, $2.5e-3$ for $\psi$, and $6.8e-3$ for $\omega$.
%\michele{(vogliamo inserire qualche altro commento ora che abbiamo inserito anche gli errori?)}. %The CPU time for ...In orde FACCIAMO UN CONFRONTO SUL COSTO COMPUTAZIONALE.

%163 secondi
%56 secondi

With these results, we consider the FOM validated. Next, we are going to validate our ROM approach. Our goal is a thorough assessment of our ROM model on two fronts: (i) the reconstruction
of the time evolution of the flow field and (ii) a physical parametric setting. Let us start from the former. 

\subsection{Validation of the ROM: time reconstruction}\label{sec:time}

\begin{figure}[h]
\centering
\includegraphics[height=0.35\textwidth]{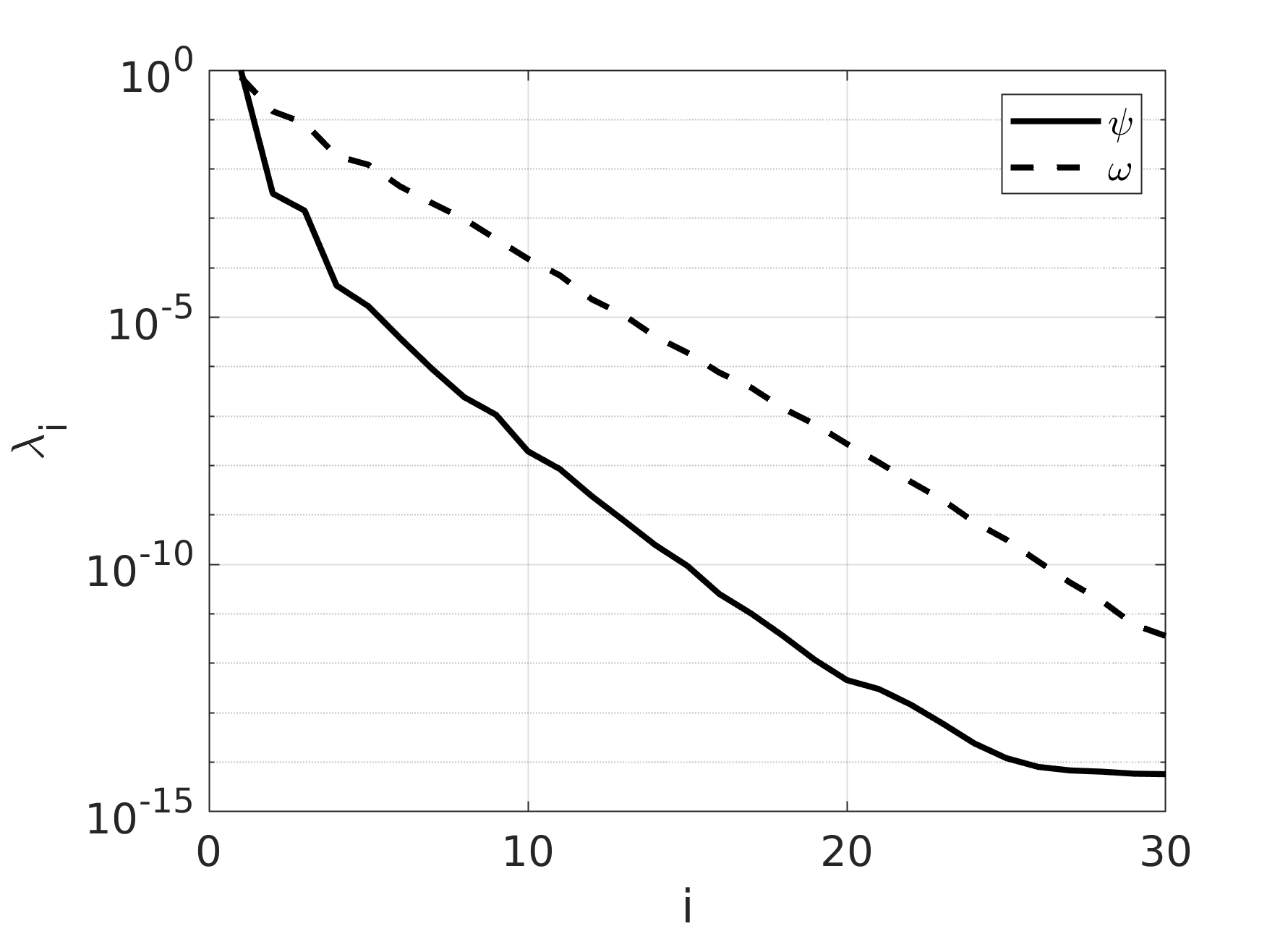}
\caption{ROM validation - time reconstruction: eigenvalue decay for the stream function and the vorticity.}
\label{fig:eig_t}
\end{figure}

We collect 250 FOM snapshots, one every 0.08 s, i.e.~we use an equispaced grid in time. Fig.~\ref{fig:eig_t} shows the eigenvalues decay for the stream function and the vorticity.
We observe that the eigenvalues decay for $\omega$ is much slower than the eigenvalues decay for $\psi$. We set the threshold for the selection of the eigenvalues to $1e-5$, resulting in 6 modes for $\psi$ and 14 modes for $\omega$. We suspect that this larger number of modes for $\omega$ is due to the richer structure of $\omega$.
%\michele{(forse quest'ultima osservazione potrebbe essere oggetto di critica, nel senso che forse è difficile trovare una connessione diretta tra il rateo di decrescita degli autovalori e l'errore tra le due formulazioni a livello FOM, che ne pensi?)}.

%\end{comment}

We calculate the relative $L^2$ error in percentage: 
%\anna{In \eqref{eq:space} abbiamo usato $\Phi$ mentre adesso usiamo $X$, sono la stessa cosa? In caso, scegliamo come uniformare la notazione.}
\begin{equation}\label{eq:error1}
E_{\Phi}(t) = 100 \cdot \dfrac{||\Phi_h(t) - \Phi_r(t)||_{L^2(\Omega)}}{||{\Phi_h}(t)||_{L^2(\Omega)}},
\end{equation}
where $\Phi_h$ is a field computed with the FOM ($\psi_h$ or $\omega_h$) and $\Phi_r$ is the corresponding field computed with the ROM ($\psi_r$ or $\omega_r$). Moreover, we evaluate the relative error in percentage for the enstrophy $e = \int_\Omega \omega^2 d\Omega$, i.e.:
\begin{equation}\label{eq:error2}
E_{e}(t) = 100 \cdot \dfrac{e_h(t) - e_r(t)}{e_h(t)},
\end{equation}
where $e_h$ and $e_r$ are the values of the enstrophy 
%\anna{(ho visto che hai commentato la definizione di enstofia, secondo me vale la pena di lasciarla spt perche' non mi pare che al momento abbiamo problemi di spazio)} \michele{in realta' dovrebbe esserci, sebbene non cm equazione ma nel testo, un rigo prima eq. (35)}
%\begin{equation}\label{eq:enstrophy}
%e = \int_\Omega \omega^2 d\Omega,
%\end{equation}
computed by the FOM and the ROM, respectively.

\begin{figure}
\centering
 \begin{overpic}[width=0.45\textwidth]{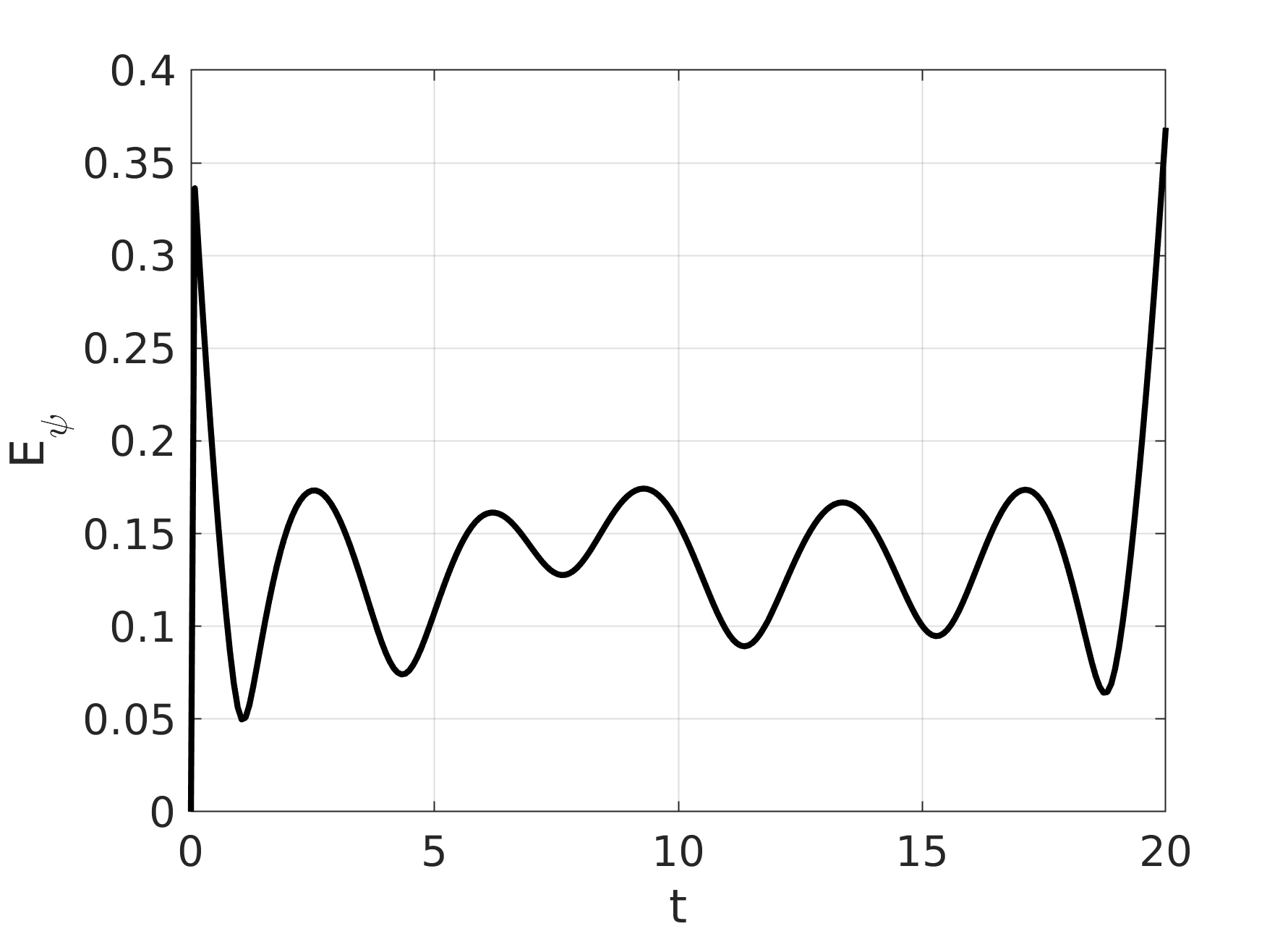}
        %\put(35,18){FOM}
        %\put(-8,7){$\u$}
      \end{overpic}
 \begin{overpic}[width=0.45\textwidth]{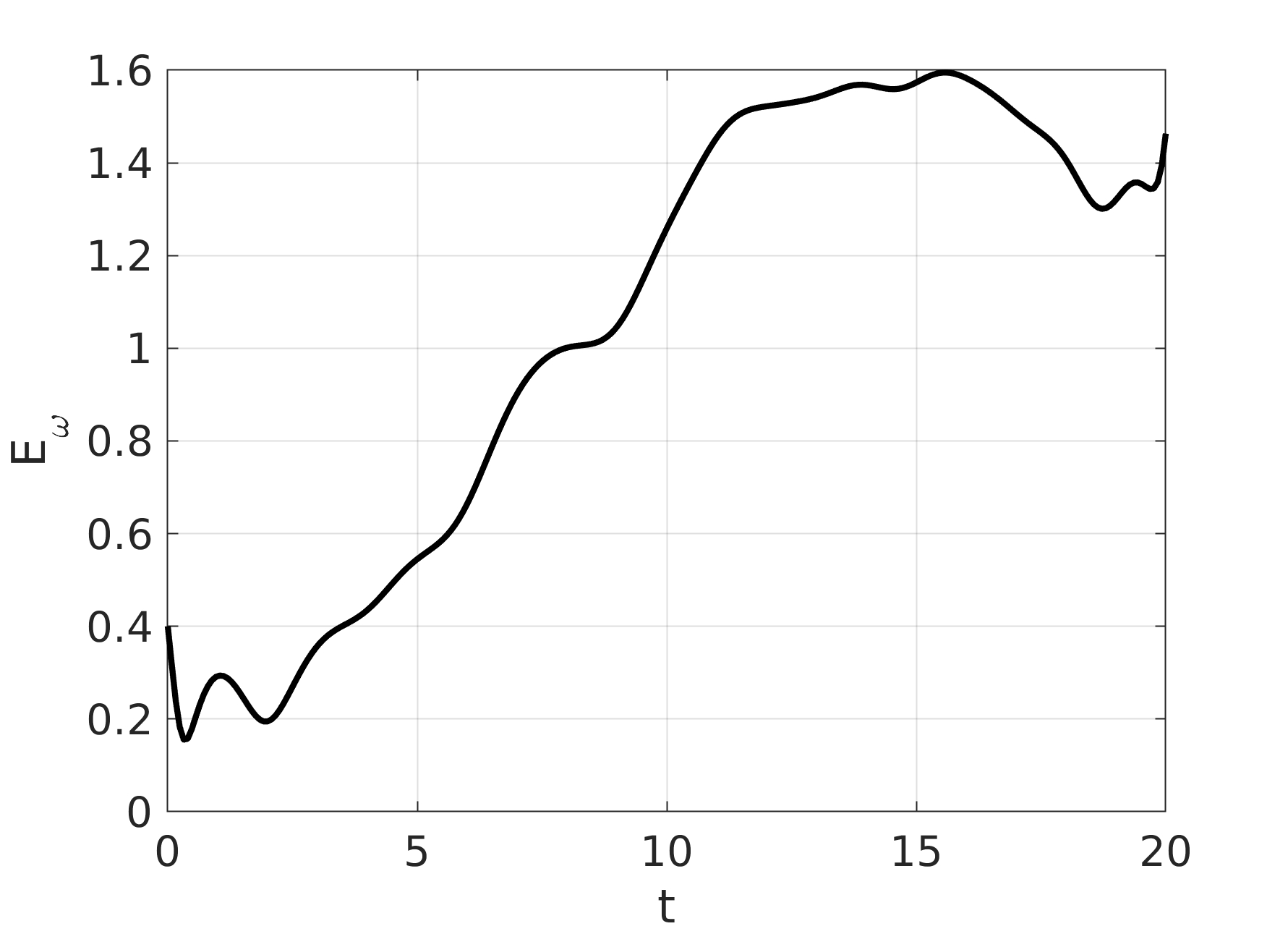}
        %\put(35,18){ROM}
      \end{overpic}
      %\vskip .2cm
       %\begin{overpic}[width=0.3\textwidth]{img/U_diff_1s.png}
        %\put(40,40){ROM}
      %\end{overpic}\\
 \begin{overpic}[width=0.45\textwidth]{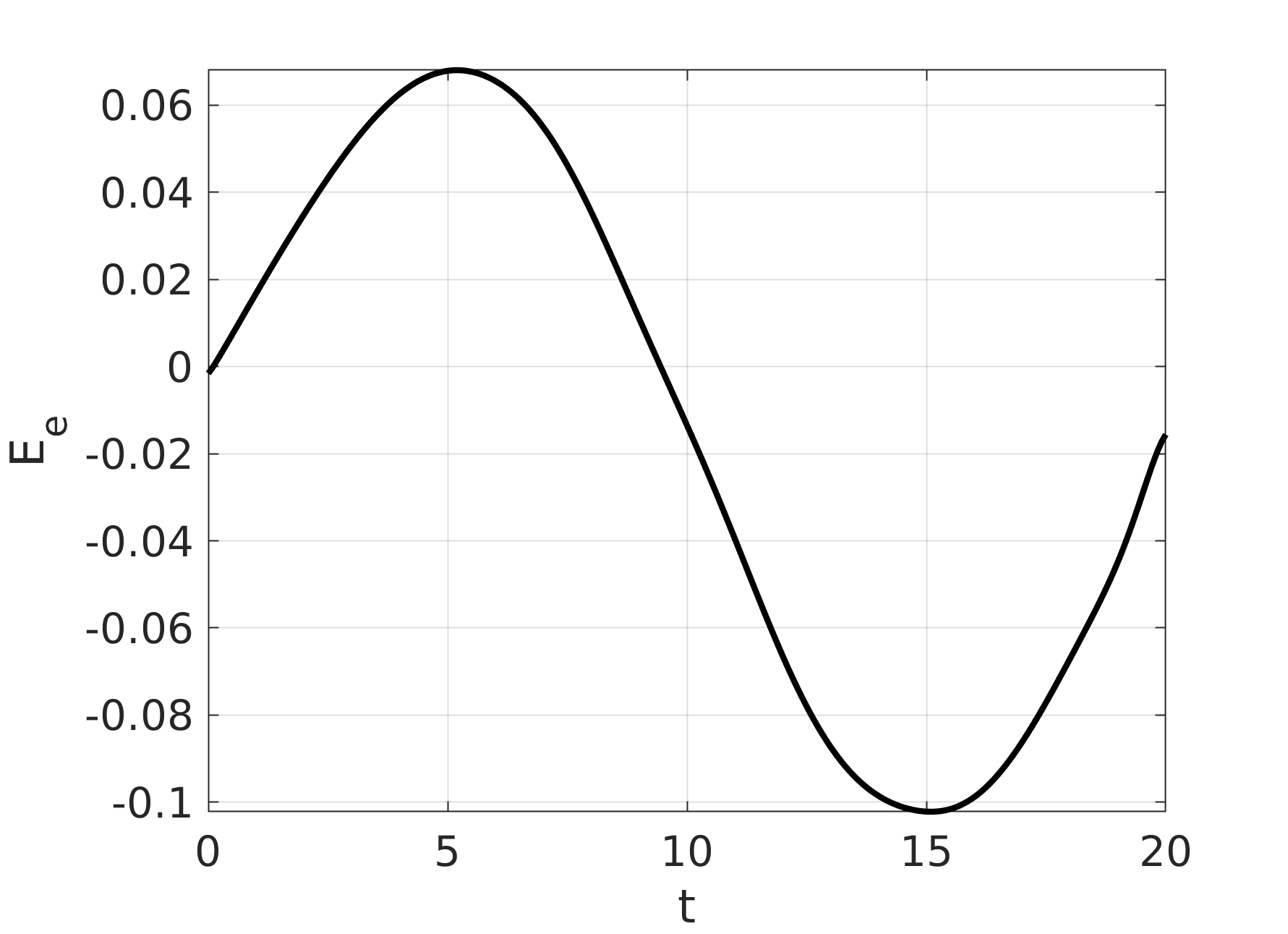}
        %\put(-8,7){$p$}
      \end{overpic}\\
\caption{ROM validation - time reconstruction: evolution of error \eqref{eq:error1} for stream function $\psi$ (top left) and vorticity $\omega$ (top right), and error for enstrophy $e$ \eqref{eq:error2} (bottom).}
\label{fig:errors}
\end{figure}

Fig.~\ref{fig:errors} shows error \eqref{eq:error1} for the stream function and vorticity and error \eqref{eq:error2} for the enstrophy over time. We see that all relative errors in percentage achieve very low values. In particular, over the entire time interval the error for $\psi$ is lower than 0.4\%, the error for $\omega$ is lower than 1.6\%, and the error for the enstrophy is lower than 0.1\% in absolute value. One observation is in order: the relative errors are significantly lower than the values obtained by ROM for the Navier--Stokes equations in primitive variables (up to 3 \% for the velocity and 15 \% for the pressure in the 2D flow past a cylinder benchmark \cite{Stabile2018}) 
%\michele{poi verifico per sicurezza i numeri} \anna{Si', non vogliamo fare gaffes :)} \michele{verificato, tutto ok! :)}.  \michele{forse qui ci sta una incomprensione. Nel riferimento biblio citato non viene trattato il problema del vortice con NSE u p ma il problema del flow past a cylinder che, se uno volesse trattare con NSE psi omega, anche a livello full non sarebbe banalissimo visto il problema di dover assegnare delle BC alla vorticità. In ogni caso, al di la' del benchmark, secondo me il commento cosi' come la spiegazione sotto rimangono cmq valide, perchè piu' che il benchmark in sè riguarda il tipo di formulazione, pero' forse è da rifrasare un attimo}. 
%\anna{OK, rifrasato.} \michele{OK! :)}
%\anna{(Questa e' una cosa che si sapeva gia'?? Intendo, nel paper \cite{Stabile2018} vengono confrontate le 2 formulazioni o no? Altrimenti e' una cosa carina da sottolineare)} \michele{No, no, nessun confronto. Dovremmo essere i primi, almeno to the best of my knowledge ahah :D}. 
We speculate that the reason for this difference could be the fact that the development of a ROM for the Navier--Stokes equations in primitive variables requires strategies for the stabilization of the pressure (e.g. supremizer enrichment and Poisson pressure equation) that makes the ROM framework more complex and error-prone \cite{Ballarin2014,Girfoglio2021b,Girfoglio_JCP,Rozza2007,Stabile2018}.
%\anna{Michele, cita anche il nostro lavoro sul confronto tra metodi di stabilizzazione.} \michele{ok!}
%\michele{inserire un po di citazioni}. %(Ladyzhenskaya-Brezzi-Babuska) condition \michele{inserire un po di citazioni}.
%Then, we observe that the smallest relative error are achieved
%for the enstrophy of the system.

\begin{figure}
\centering
 \begin{overpic}[width=0.45\textwidth]{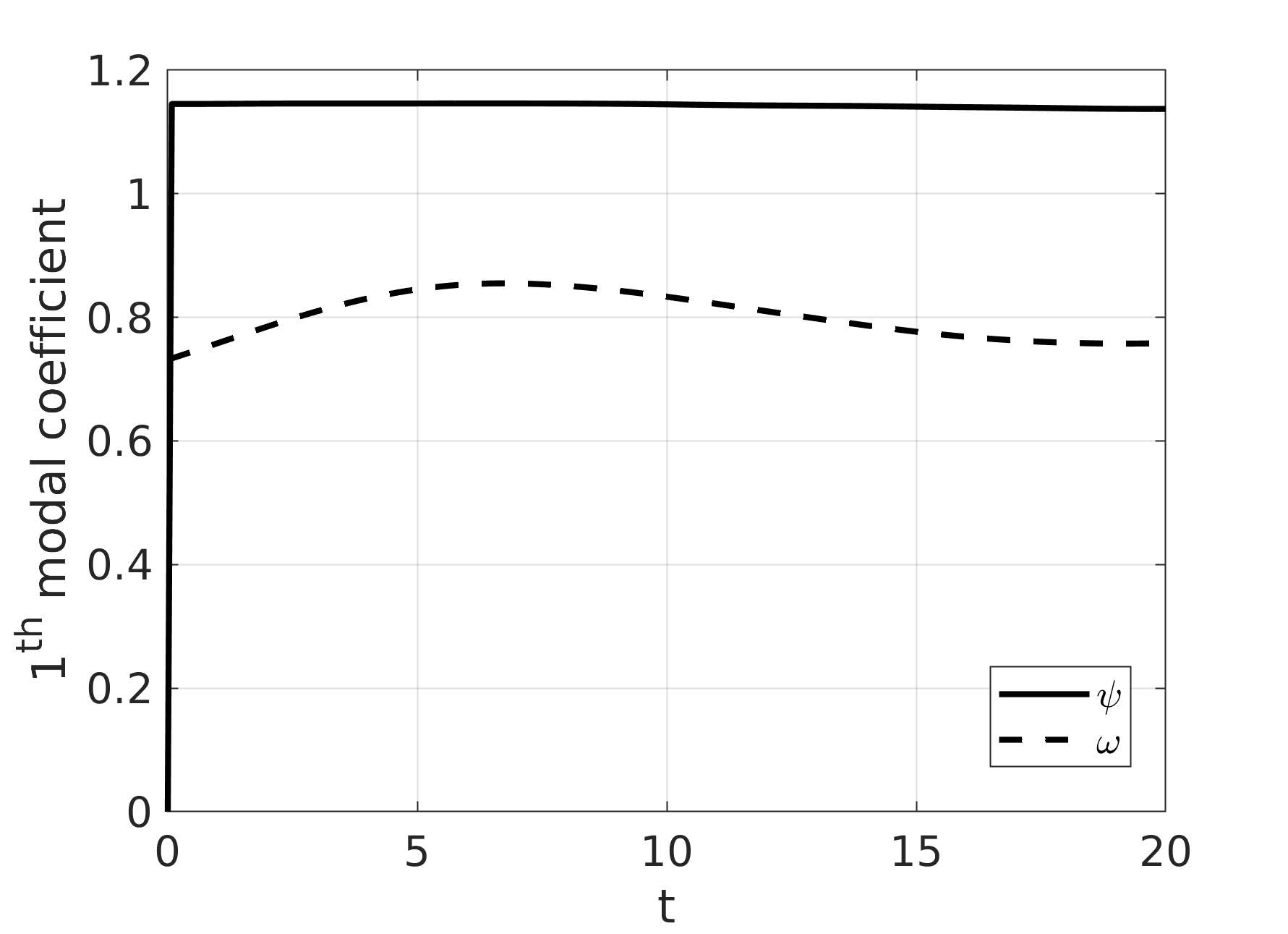}
        %\put(35,18){FOM}
        %\put(-8,7){$\u$}
      \end{overpic}
 \begin{overpic}[width=0.45\textwidth]{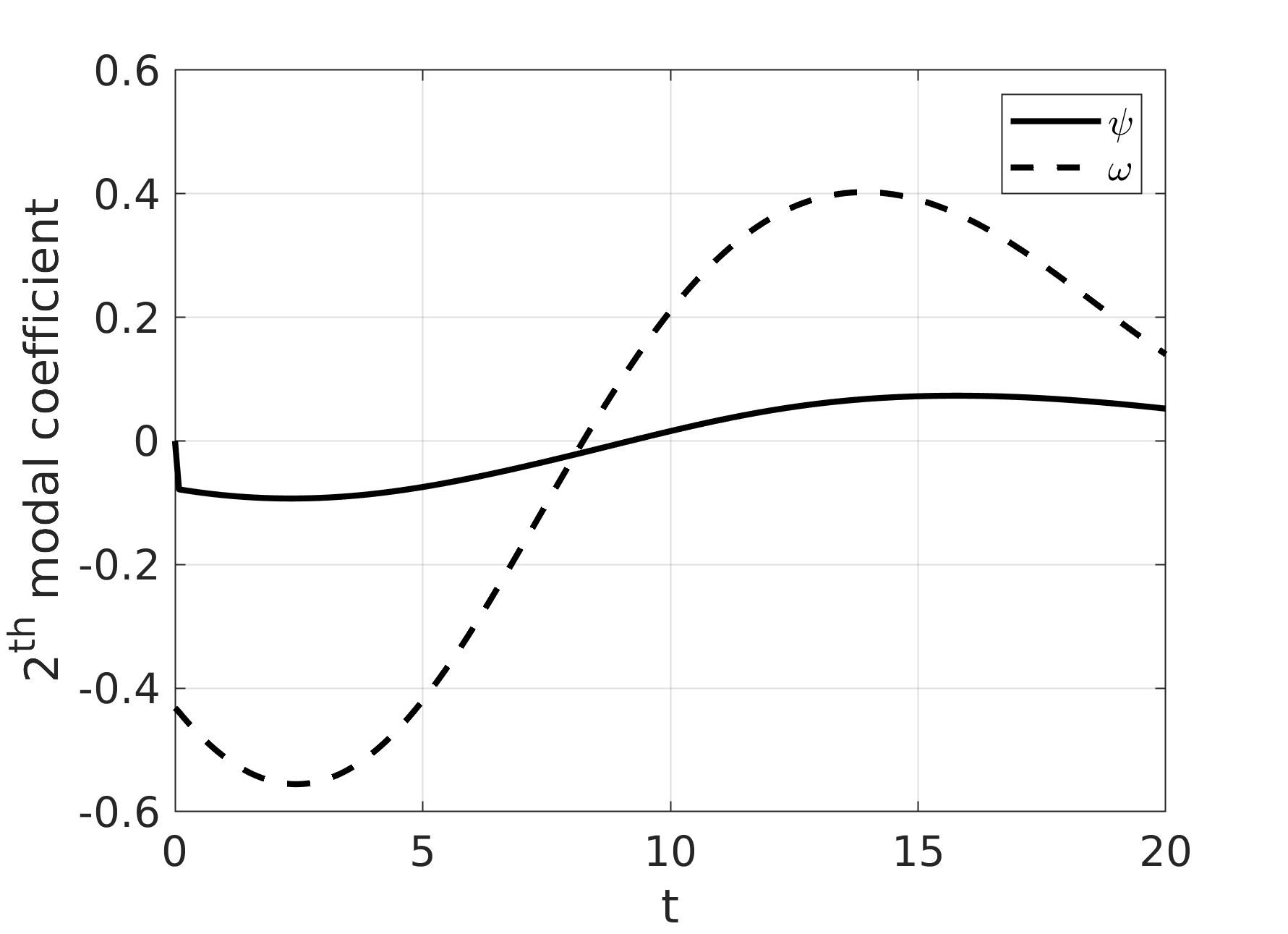}
        %\put(35,18){ROM}
      \end{overpic}
      %\vskip .2cm
       %\begin{overpic}[width=0.3\textwidth]{img/U_diff_1s.png}
        %\put(40,40){ROM}
      %\end{overpic}\\
 \begin{overpic}[width=0.45\textwidth]{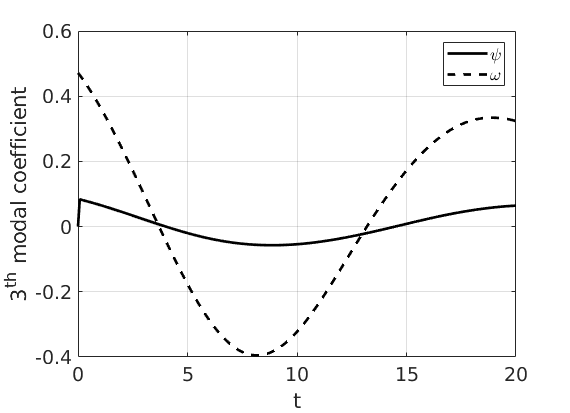}
        %\put(-8,7){$p$}
      \end{overpic}\\
\caption{ROM validation - time reconstruction: evolution of the first three reduced coefficients for $\psi$ and $\omega$.}
\label{fig:modal_coeff_t}
\end{figure}

In order to justify our choice to use different reduced coefficients to approximate stream-function and vorticity in \eqref{eq:ROM_1}, we show in 
Fig.~\ref{fig:modal_coeff_t} the time evolution of the first three reduced coefficients for $\psi$ and $\omega$: the differences
are significant. Thus, using the same reduced coefficients would lead to a less accurate reconstruction of $\psi$ and $\omega$.
%We note that wemade two important choices in our approach:  (i) We enforced the coupling betweenthe POD vorticity and streamfunction basis functions in (15); and (ii) We used thesame ROM coefficients in the ROM vorticity approximation (16) and in the ROMstreamfunction approximation (17). 

\begin{figure}
\centering
\hspace{.4cm}
 \begin{overpic}[width=0.19\textwidth]{img/psi_psiomega_4.png}
        \put(38,101){$t = 4$}
        \put(-38,45){FOM}  %-28
      \end{overpic}
 \begin{overpic}[width=0.19\textwidth]{img/psi_psiomega_8.png}
        \put(38,101){$t = 8$} %0.2
      \end{overpic}
 \begin{overpic}[width=0.19\textwidth]{img/psi_psiomega_16.png}
        \put(35,101){$t = 16$}
      \end{overpic}
 \begin{overpic}[width=0.19\textwidth]{img/psi_psiomega_20.png}
         \put(35,101){$t = 20$}
      \end{overpic}
       \begin{overpic}[width=0.08\textwidth]{img/legendPsiFOM.png} %0.09
      \end{overpic}\\
       \vskip .2cm
       \hspace{.4cm}
 \begin{overpic}[width=0.19\textwidth]{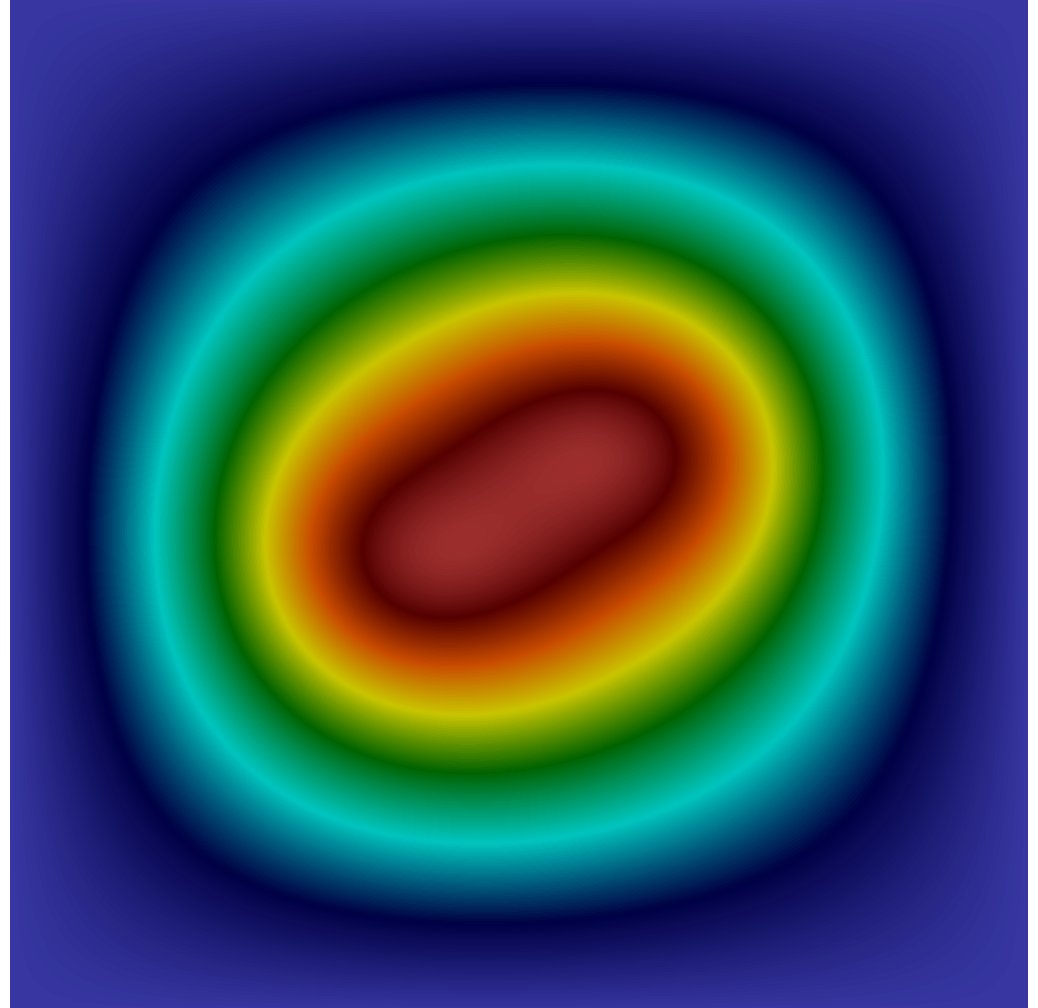}
        \put(-38,45){ROM}
      \end{overpic}
 \begin{overpic}[width=0.19\textwidth]{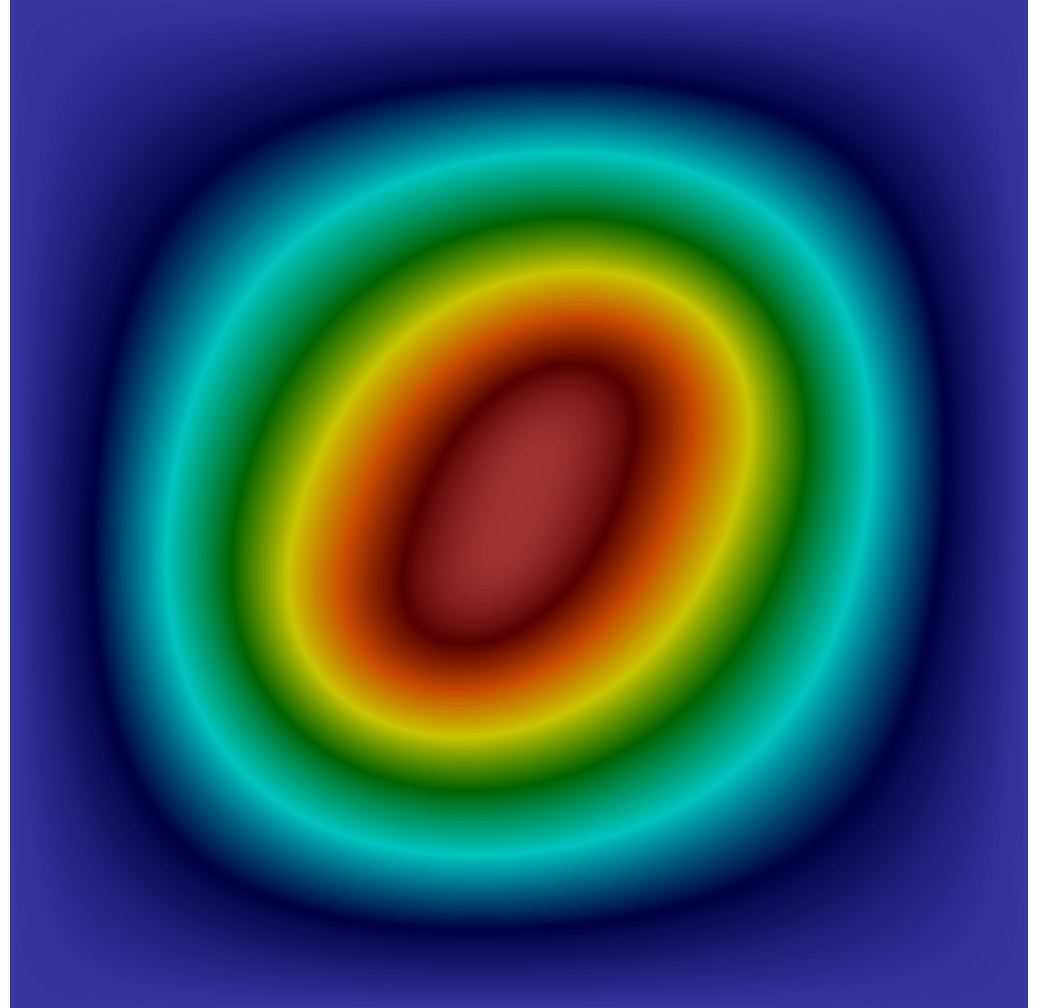}
      \end{overpic}
 \begin{overpic}[width=0.19\textwidth]{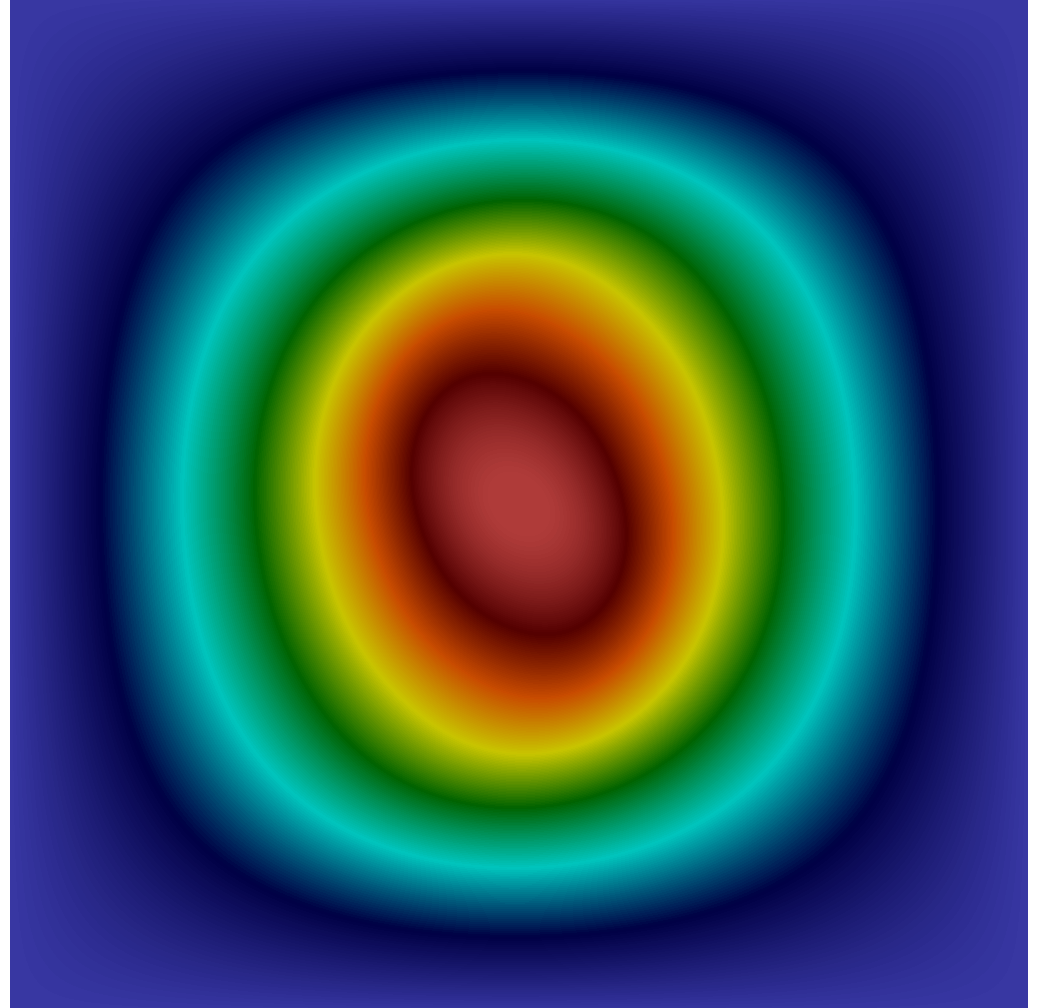}
      \end{overpic}
 \begin{overpic}[width=0.19\textwidth]{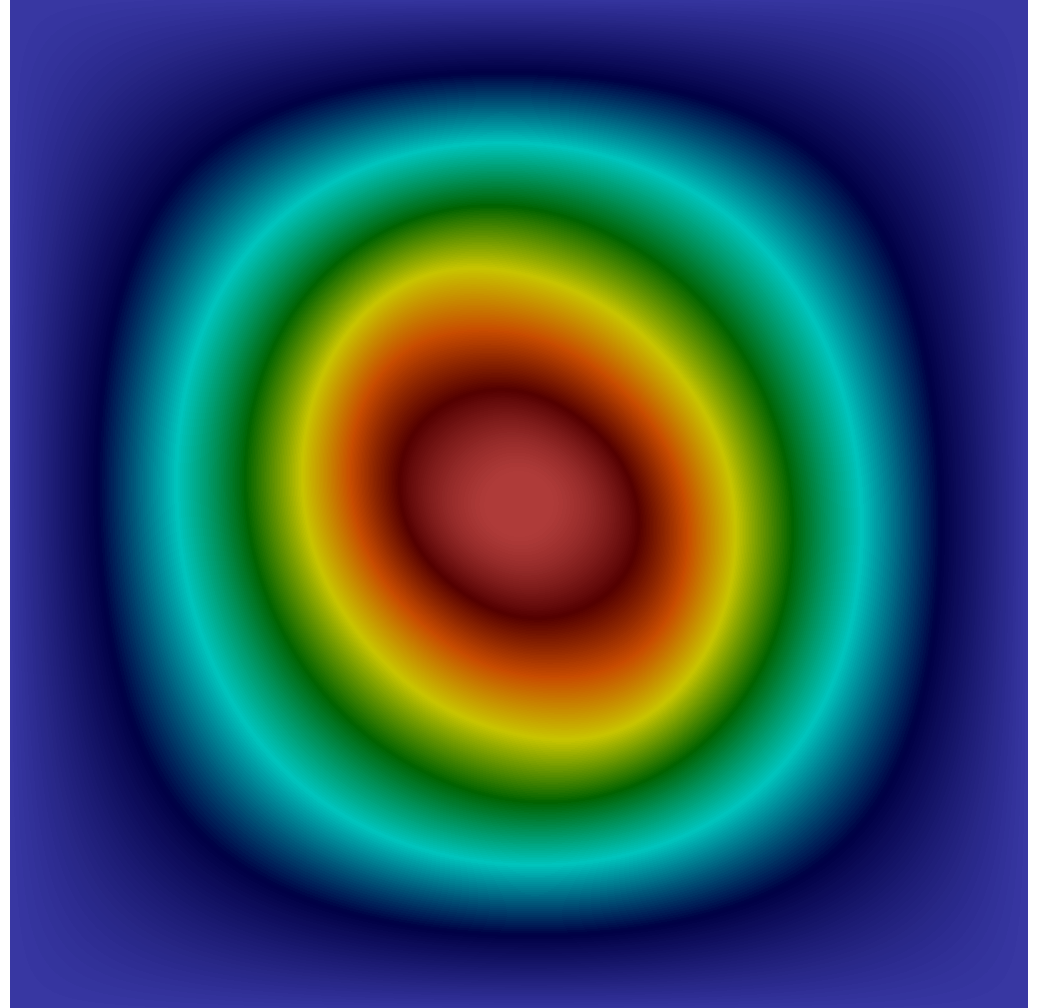}
      \end{overpic}
      \begin{overpic}[width=0.08\textwidth]{img/legendPsiFOM.png}
      \end{overpic}
 \vskip .2cm
 \hspace{.4cm}
       \begin{overpic}[width=0.19\textwidth]{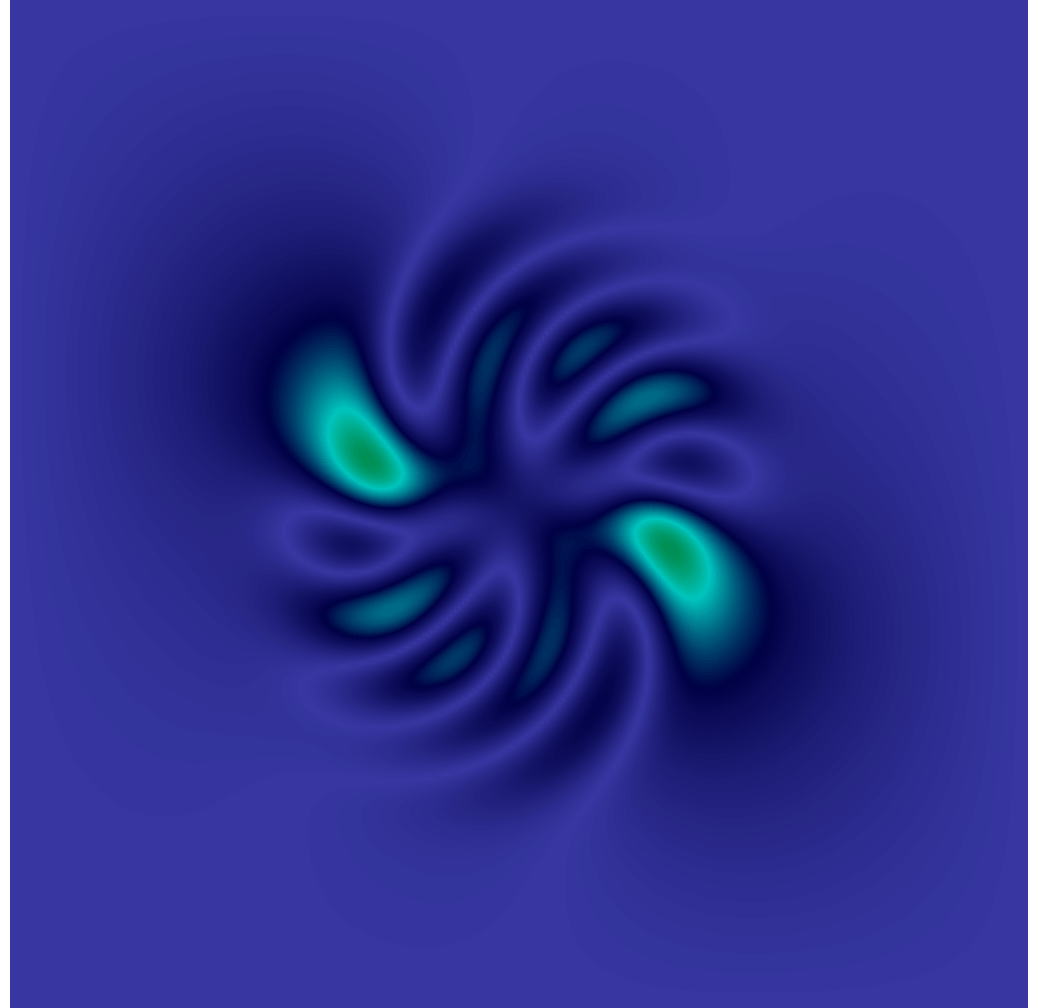}
       \put(-35,45){Diff.}
      \end{overpic}
 \begin{overpic}[width=0.19\textwidth]{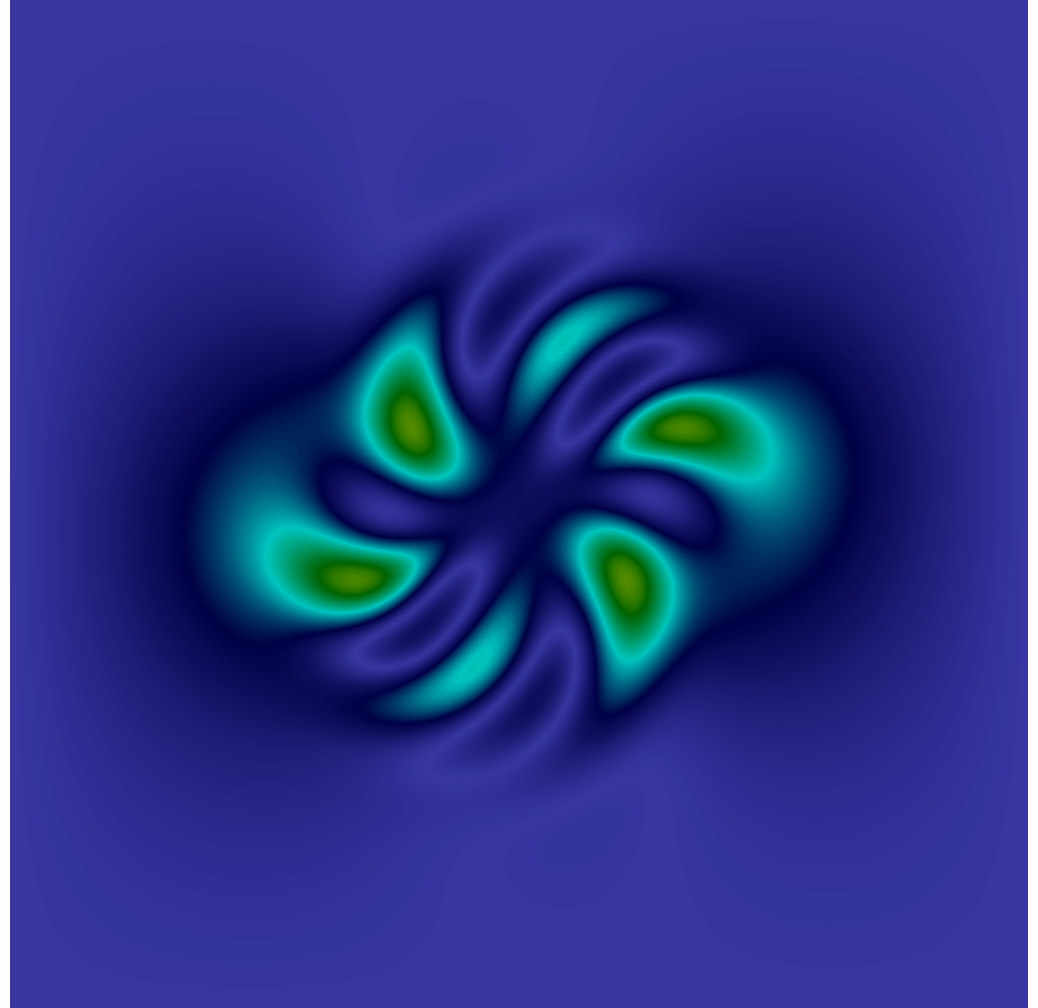}
      \end{overpic}
 \begin{overpic}[width=0.19\textwidth]{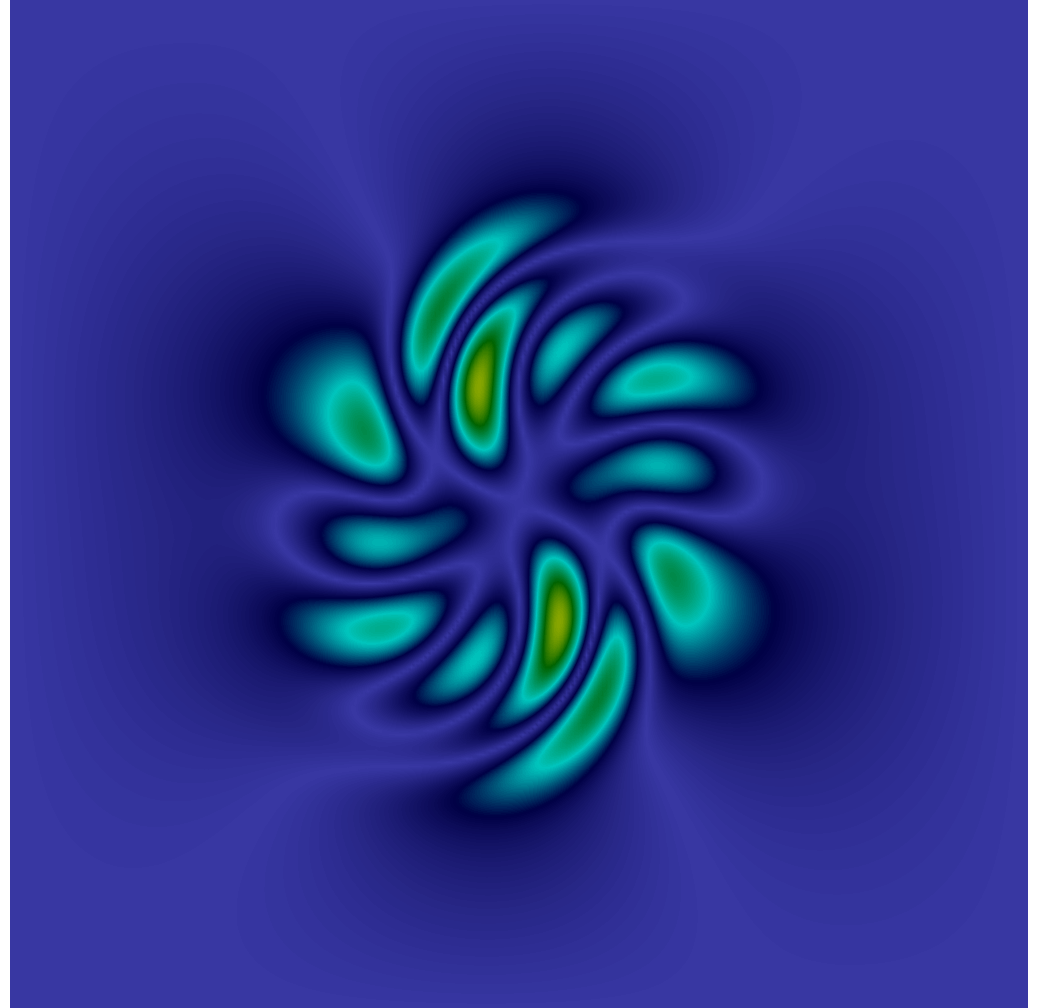}
      \end{overpic}
 \begin{overpic}[width=0.19\textwidth]{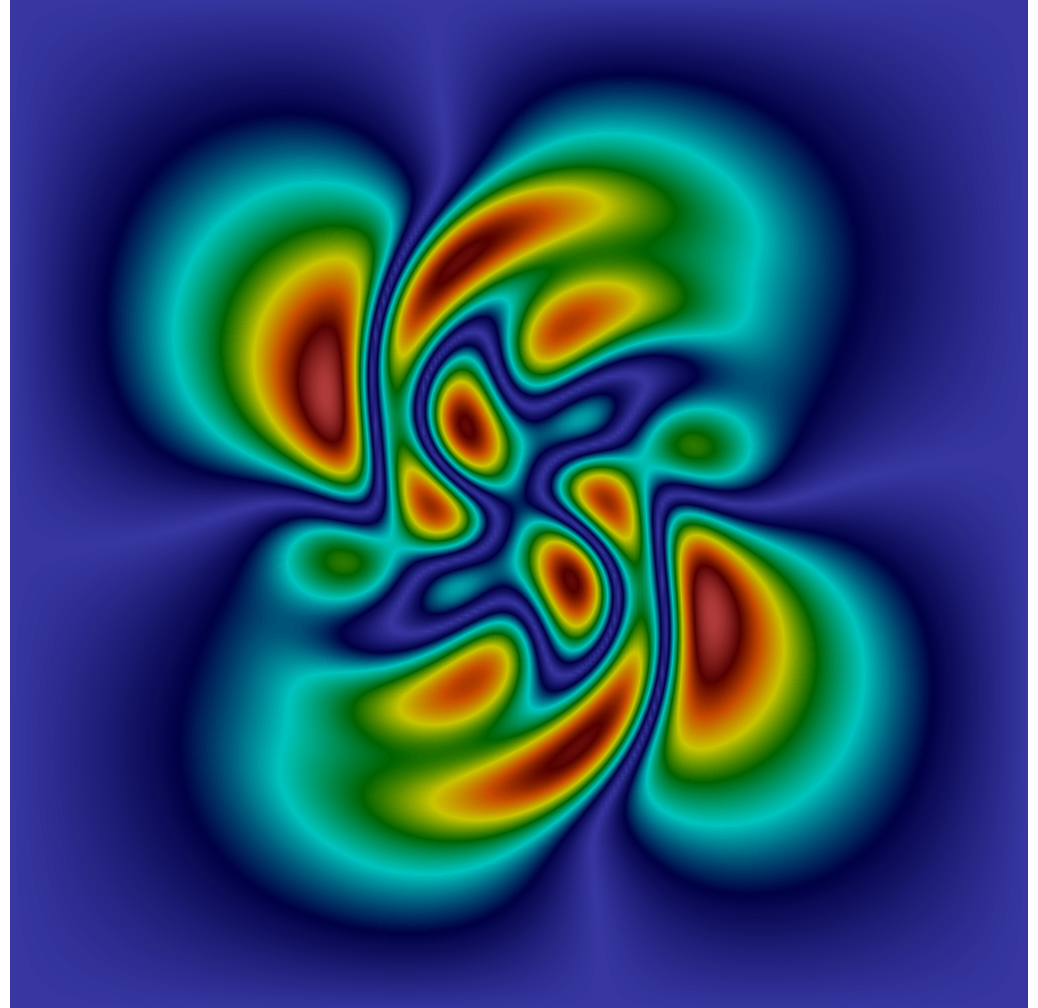}
      \end{overpic}
      \begin{overpic}[width=0.08\textwidth]{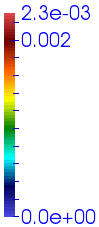}
      \end{overpic}
\caption{ROM validation - time reconstruction: stream function $\psi$ computed by the FOM (first row) and the ROM (second row), and difference between the two fields in absolute value (third row) at times $t = 4$ (first column), $t = 8$ (second column), $t = 16$ (third column) and $t = 20$ (fourth column). We consider 6 modes. %\anna{Il commento della Fig.~3 e' valido anche per questa figura e la prossima.} \michele{fatto! :)}
}
\label{fig:comp_psi_ROM}
\end{figure}

\begin{figure}
\centering
\hspace{.4cm}
 \begin{overpic}[width=0.19\textwidth]{img/omega_psiomega_4.png}
        \put(38,101){$t = 4$}
        \put(-38,45){FOM}
      \end{overpic}
 \begin{overpic}[width=0.19\textwidth]{img/omega_psiomega_8.png}
        \put(38,101){$t = 8$}
      \end{overpic}
 \begin{overpic}[width=0.19\textwidth]{img/omega_psiomega_16.png}
        \put(35,101){$t = 16$}
      \end{overpic}
 \begin{overpic}[width=0.19\textwidth]{img/omega_psiomega_20.png}
         \put(35,101){$t = 20$}
      \end{overpic}
      \begin{overpic}[width=0.08\textwidth]{img/legendOmegaFOM.png}
      \end{overpic}\\
       \vskip .1cm
       \hspace{.4cm}
 \begin{overpic}[width=0.19\textwidth]{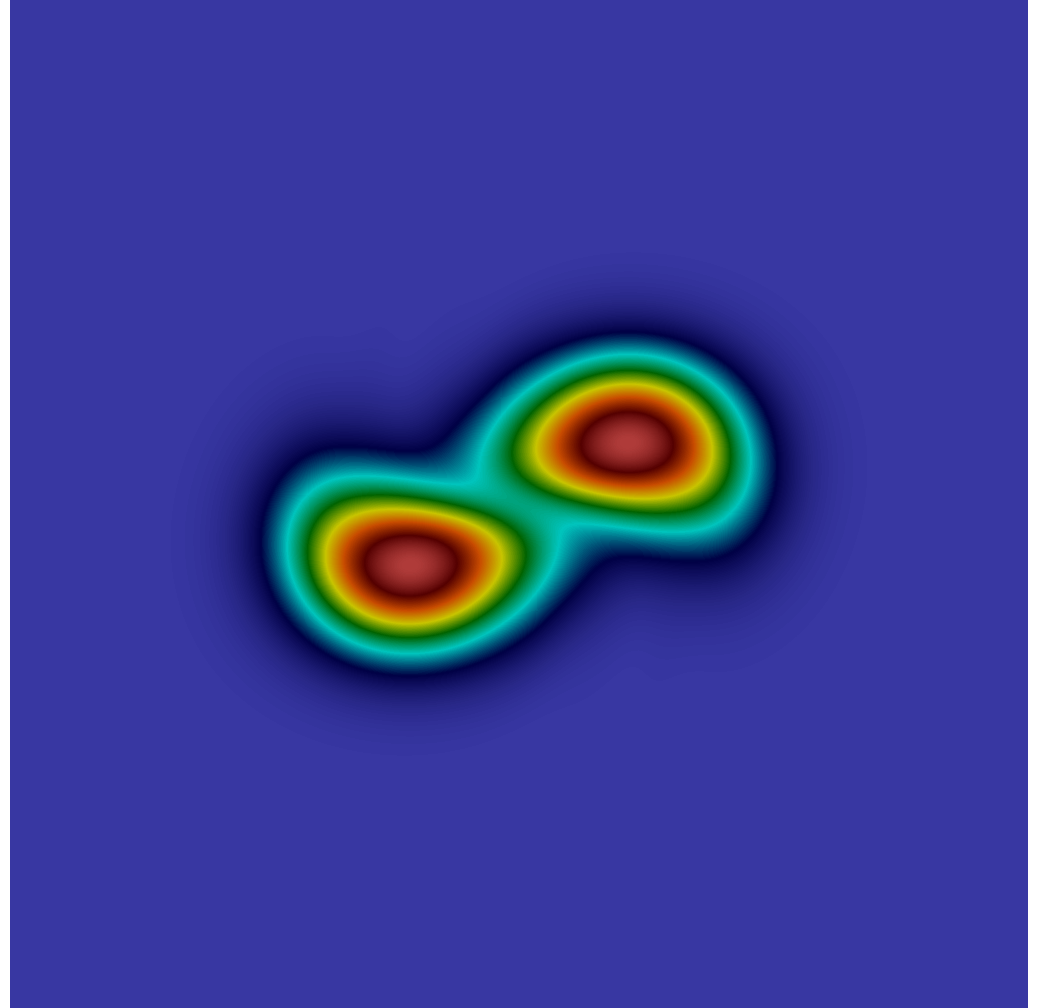}
        \put(-38,45){ROM}
      \end{overpic}
 \begin{overpic}[width=0.19\textwidth]{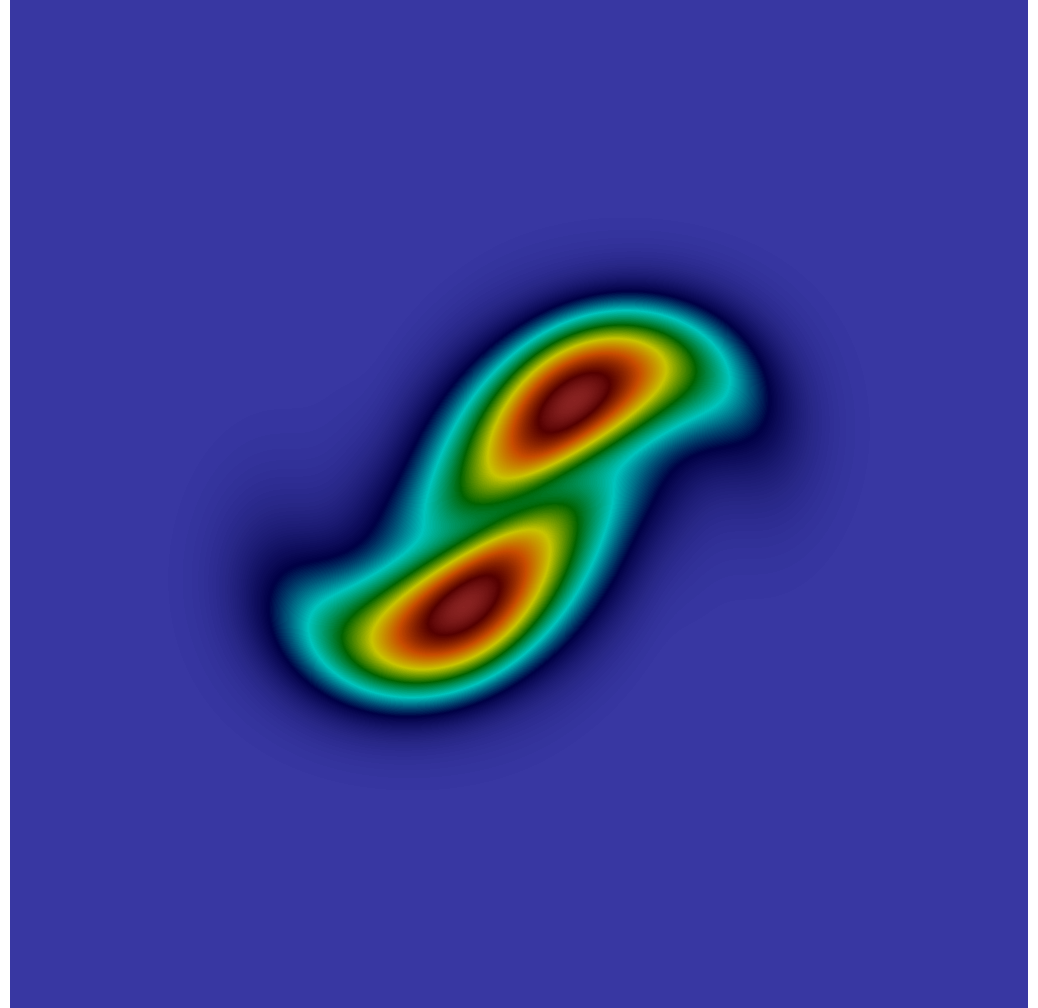}
      \end{overpic}
 \begin{overpic}[width=0.19\textwidth]{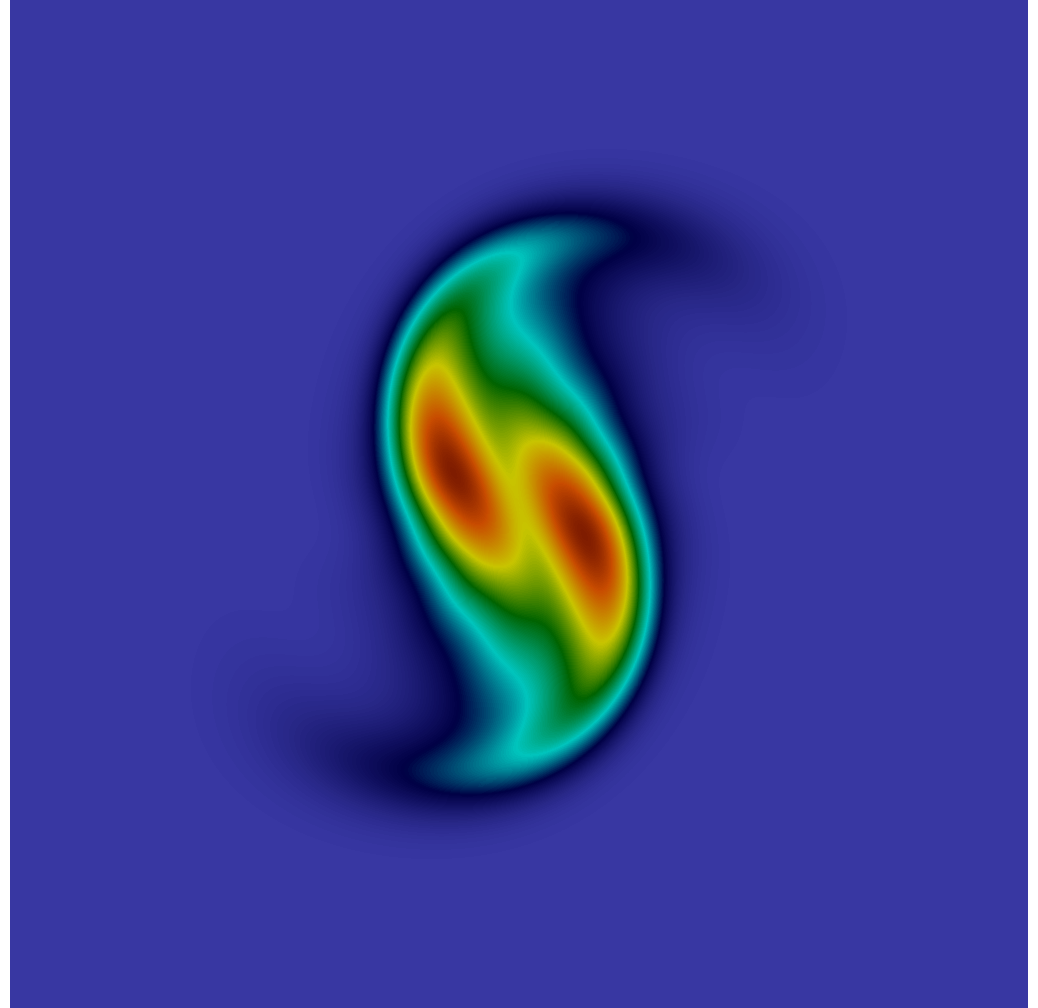}
      \end{overpic}
 \begin{overpic}[width=0.19\textwidth]{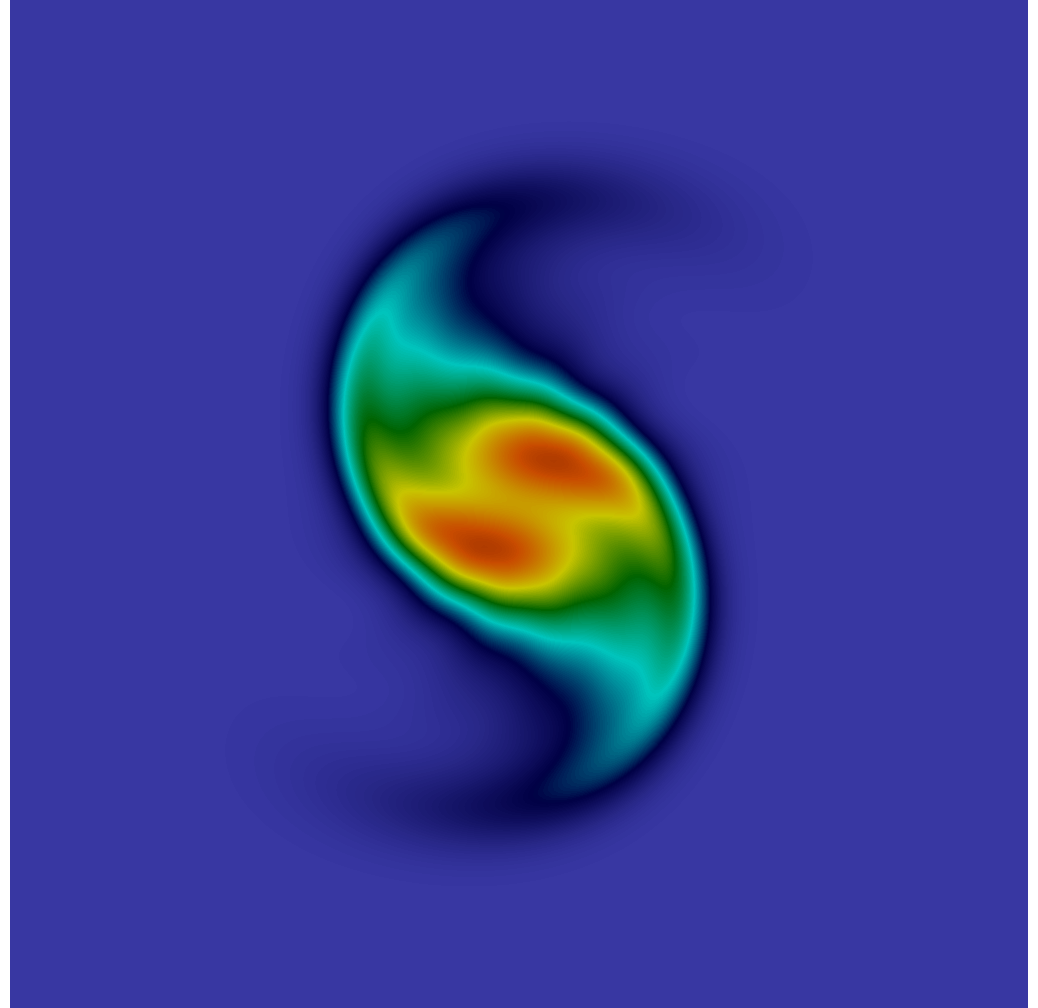}
      \end{overpic}
      \begin{overpic}[width=0.08\textwidth]{img/legendOmegaFOM.png}
      \end{overpic}
       \vskip .2cm 
       \hspace{.4cm}
       \begin{overpic}[width=0.19\textwidth]{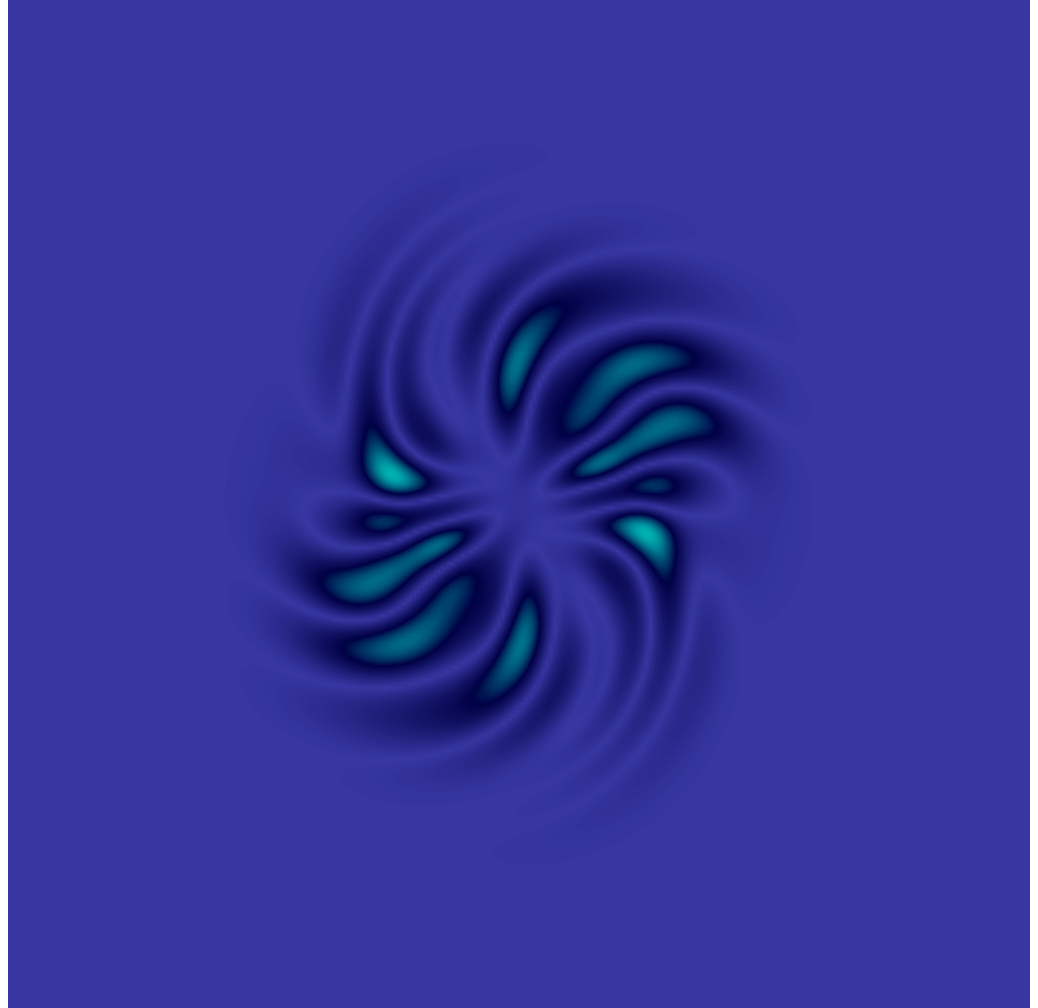}
       \put(-35,45){Diff.}
      \end{overpic}
 \begin{overpic}[width=0.19\textwidth]{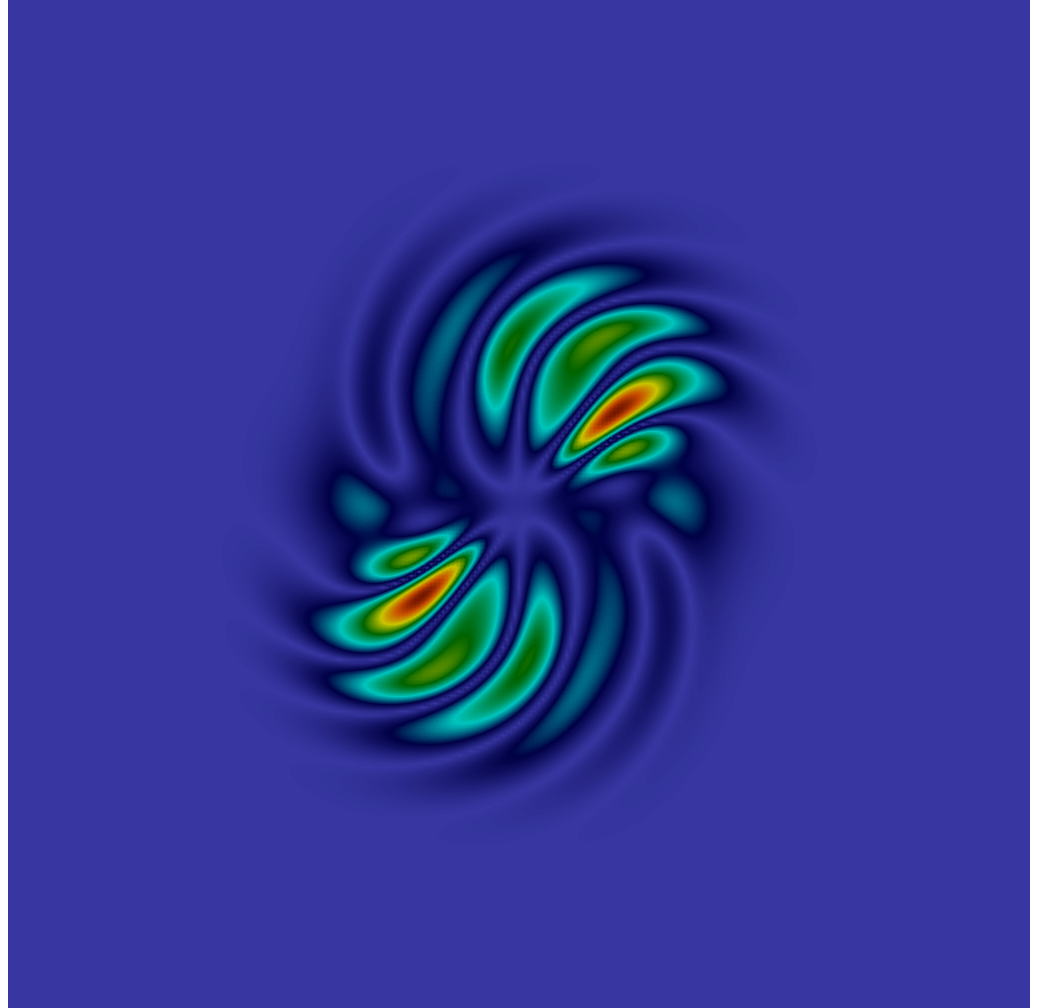}
      \end{overpic}
 \begin{overpic}[width=0.19\textwidth]{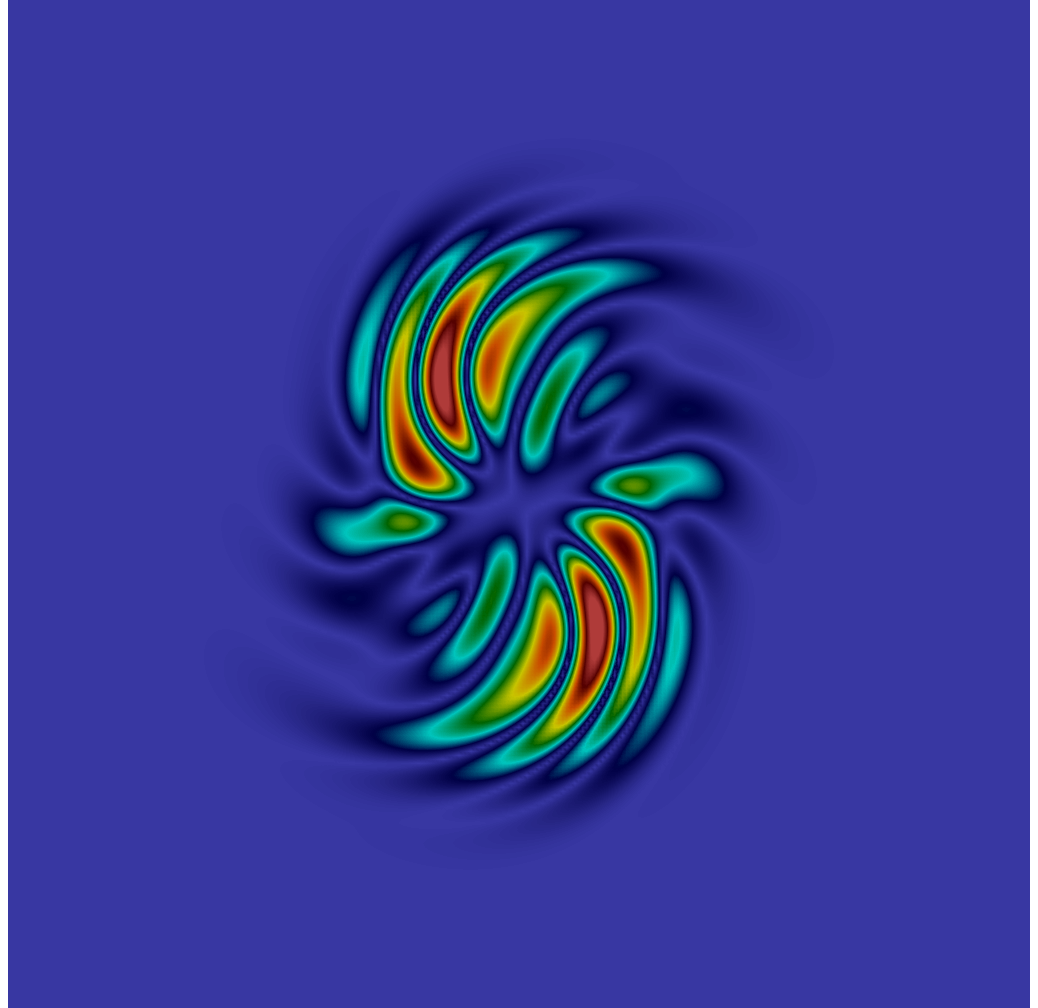}
      \end{overpic}
 \begin{overpic}[width=0.19\textwidth]{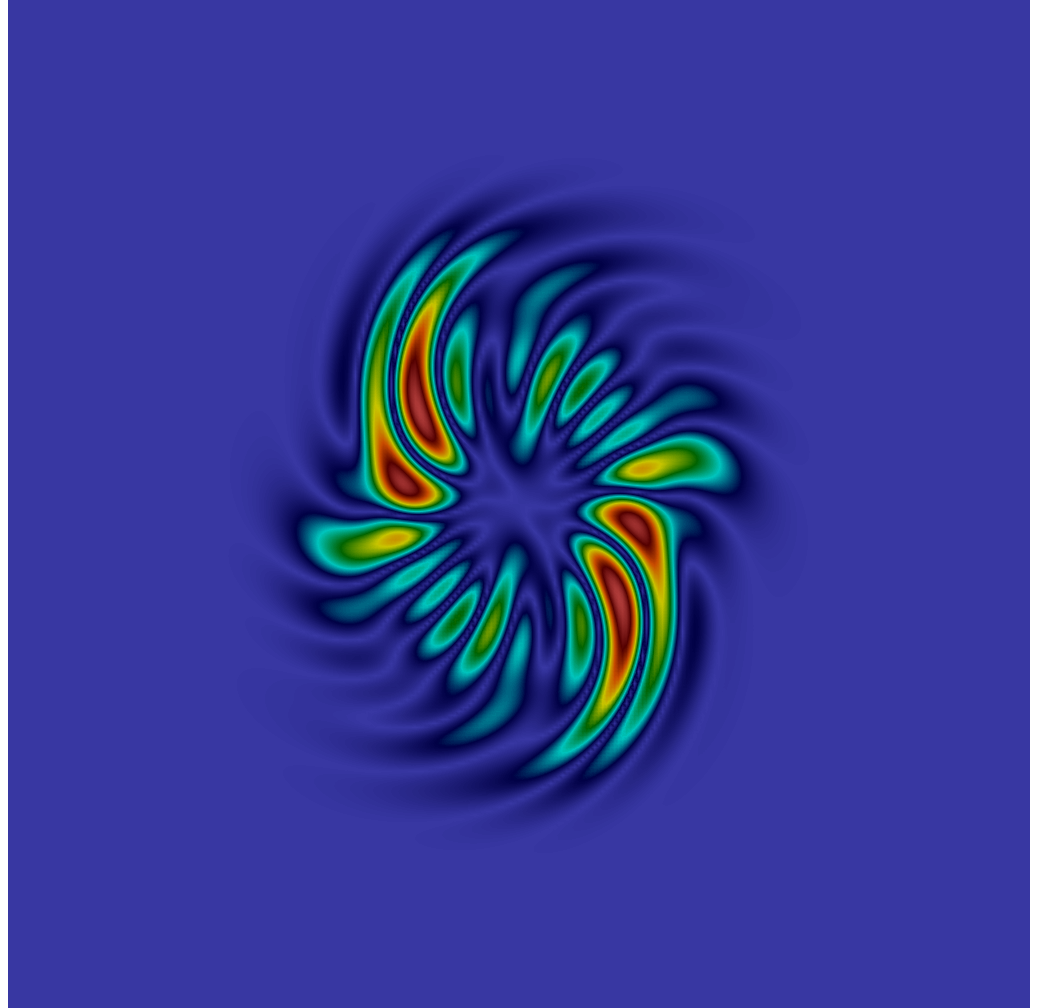}
      \end{overpic}
      \begin{overpic}[width=0.08\textwidth]{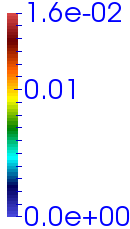}
      \end{overpic}
\caption{ROM validation - time reconstruction: vorticity $\omega$ computed by
the FOM (first row) and the ROM (second row), and difference between the two fields in absolute value (third row) at times $t = 4$ (first column), $t = 8$ (second column), $t = 16$ (third column) and $t = 20$ (fourth column). We consider 14 modes.}
\label{fig:comp_omega_ROM}
\end{figure}

Finally, we compare the solutions computed by FOM and ROM. 
Fig.~\ref{fig:comp_psi_ROM} and \ref{fig:comp_omega_ROM} display such comparison for $\psi$ and $\omega$ at four different times, respectively. 
From these figures, we see that our ROM approach provides a good global reconstruction of both stream function and vorticity. 
In fact, the maximum relative difference in absolute value does not exceed $4.7e-3$ for $\psi$ and $1.7e-2$ for $\omega$. 
%\anna{Mi chiedo come mai l'errore in $\omega$ e' sempre piu' grande di quello in $\psi$.} \michele{bella domanda! Non ho elementi certi al riguardo, è pur vero che la struttura di omega è ben piu' complessa di quella di psi per questo benchmark, (testimoniato pure dal fatto che abbiamo bisogno di piu' modi per omega), l'evoluzione in tempo è molto marcata, e dipende in maniera autonoma dal tempo, a differenza di psi che è f(omega(t))}
%\michele{anche qui, vogliamo inserire qualche altro commento visto che ci sono gli errori? Credo cmq che sia esaustivo cosi'}%Next, we make the comparison more quantitative.

We conclude this subsection by proving some information about the efficiency of our ROM approach.
The total CPU time required by a FOM simulation is about 64 s, while the computation of the
modal coefficients over the entire time window of interest takes 0.47s. 
The resulting speed-up is about 136.

%\michele{DA QUI}

\subsection{Validation of the ROM: physical parametrization}
In this section, we are going to consider a physical parametric setting. We set an arbitrary array of decaying Taylor–Green vortices as source term in the vorticity equation \eqref{eq:psizeta-mom} given by
\begin{equation}\label{eq:source}
F = - \gamma e^{-t/Re} \cos(3x) \cos(3y),
\end{equation}
where $\gamma$ is the strength of the source term. 
We consider $Re$ and $\gamma$ as the control parameters. We remark that a $Re$ parameterization has been considered in \cite{Ahmed2020} and parameterization with respect to both $Re$ and $\gamma$ has been studied in  \cite{Pawar2020-2,Pawar2020,Pawar2021} %\michele{(anche in [2] è stato fatto uno studio parametrico con una tecnica ibrida ML/Galerkin projection ma solo rispetto a Re, il termine di forcing non vien considerato. Citiamo e specifichiamo la cosa per completezza? Se no lo possiamo semplicitare citare insieme a tutti gli altri come fatto in intro, senza differenziare cosi' in dettaglio, che forse è la cosa più semplice :))}, %although within a hybrid ROM framework. 
where the aim was to develop a ROM framework to account for hidden physics through data-driven strategies based on machine learning techniques. In order to reduce the offline cost, 
we will focus on the first half of the time interval of interest considered in Sec.~\ref{sec:time}, i.e. $(0, 10]$.

%, so in the Galerkin projection the source term in the vorticity equation is treated as an unkown quantity to be recovered through data-driven strategies based on machine learning techniques. On the other hand, here we are going to develop an efficient and accurate fully intrusive ROM, so the source term in the vorticity equation is also considered in the Galerkin projection.% and one at a time, i.e. two different training phases have been performed, the one with respect to a sample distribution of $Re$ (at a given $\gamma$) and the other one with respect to a sample distribution of $\gamma$ (at a given $Re$). 
%Moreover, we consider, for each parameter, a global POD basis space obtained by a database re .... %  training phase whose samples are the elements of the set $\bm{Re} \times \bm{\gamma}$ where $\bm{Re} = [200, 400, 600, 800]$ and $\bm{\gamma} = [0.06, 0.07, 0.08, 0.09]$.

\begin{figure}
\centering
 \begin{overpic}[width=0.2\textwidth, grid = false]{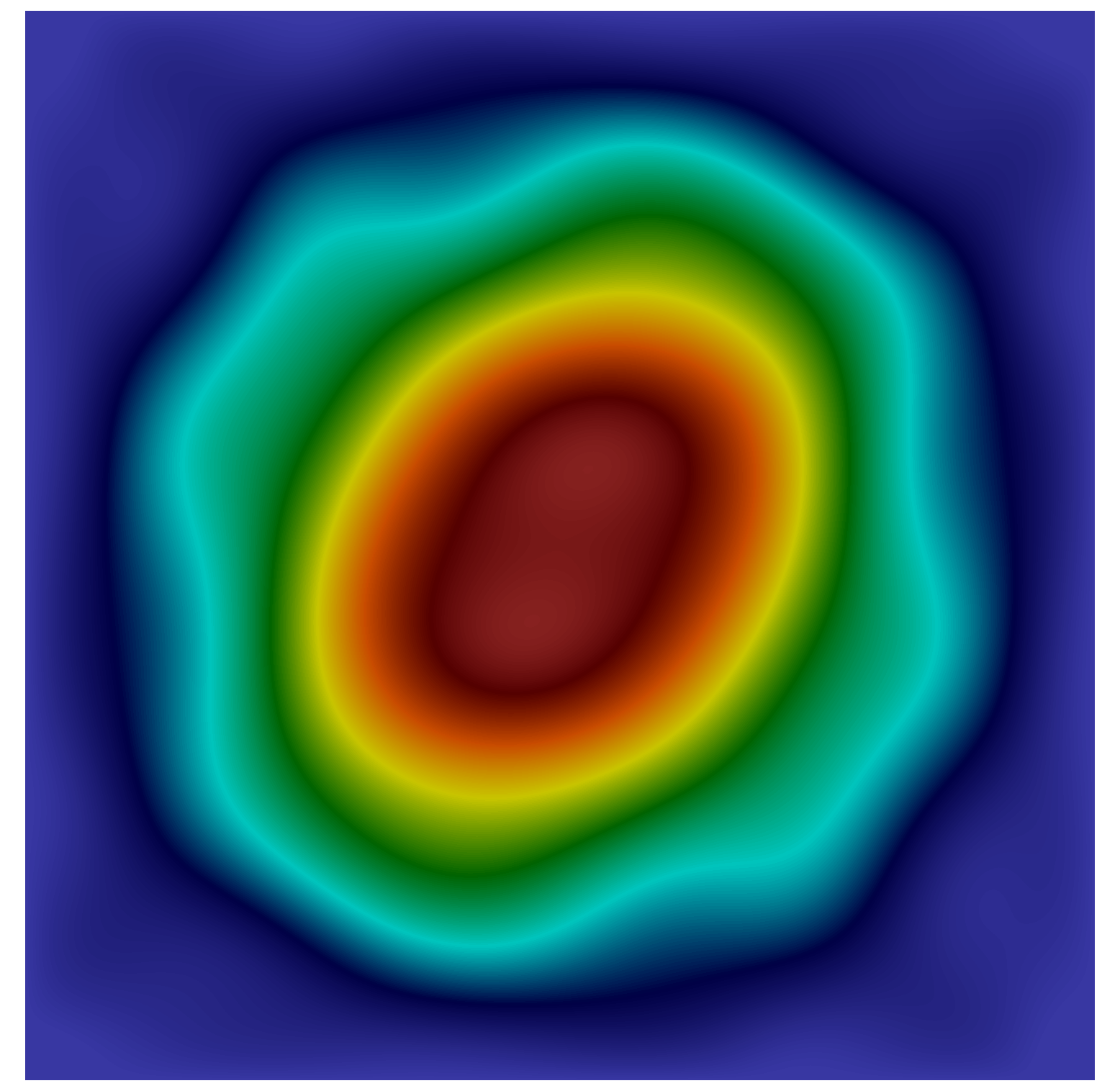}
        \put(-10,50){$\psi$}
        \put(30,100){$Re = 200$}
      \end{overpic} ~~
 \begin{overpic}[width=0.2\textwidth]{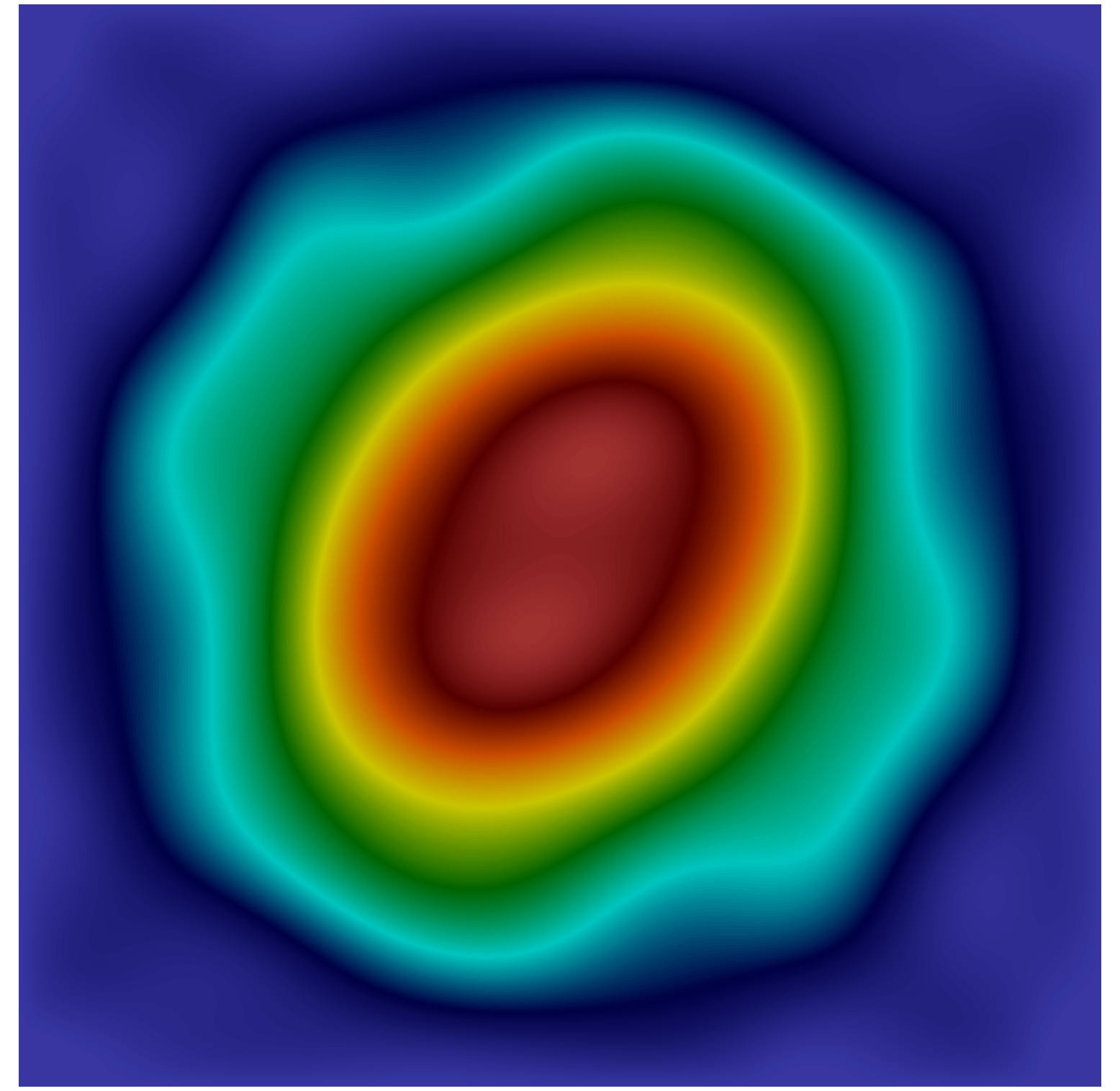}
        \put(30,100){$Re = 400$}
      \end{overpic} ~~
 \begin{overpic}[width=0.2\textwidth]{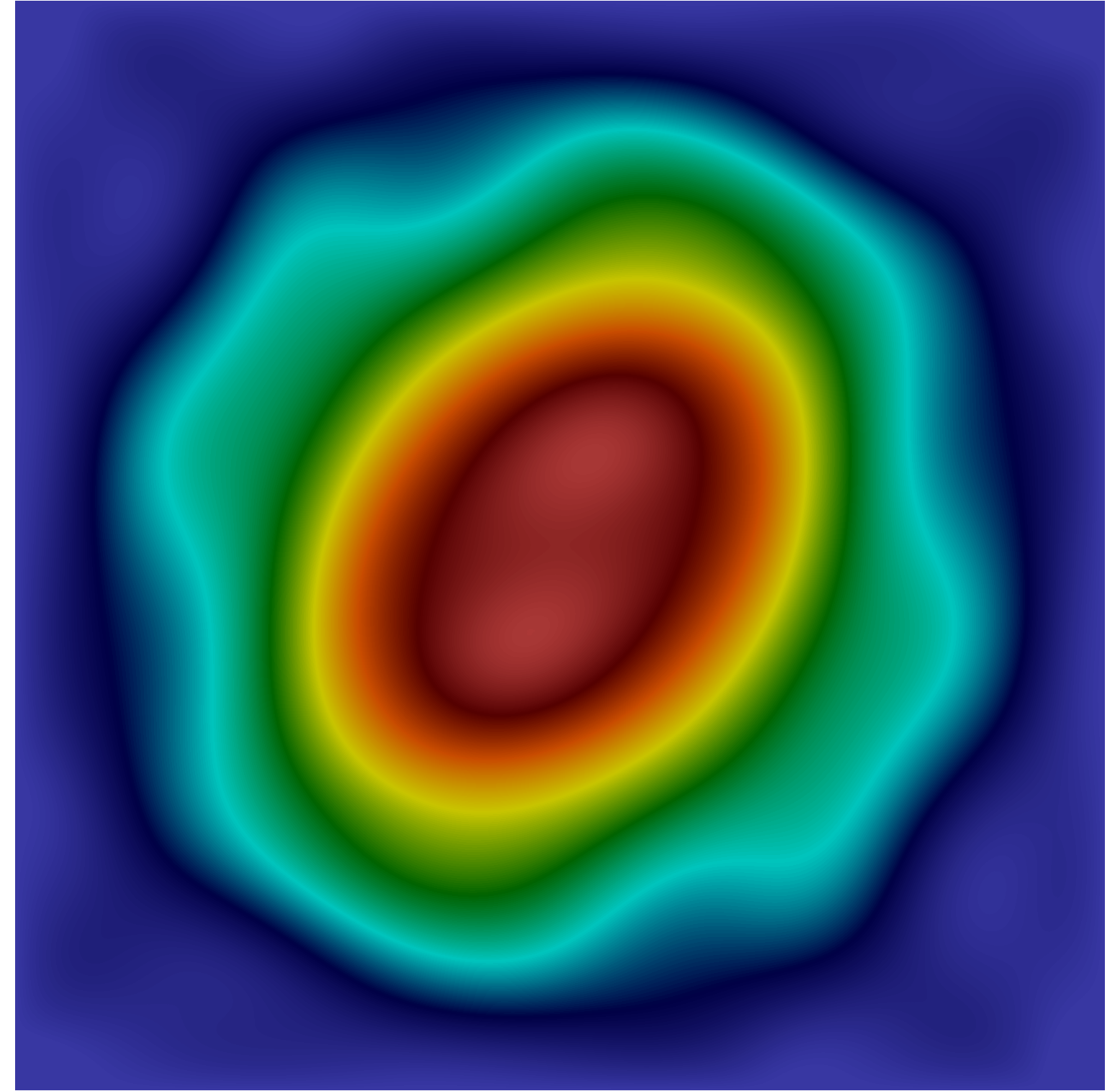}
        \put(30,100){$Re = 600$}
      \end{overpic} ~~
 \begin{overpic}[width=0.2\textwidth]{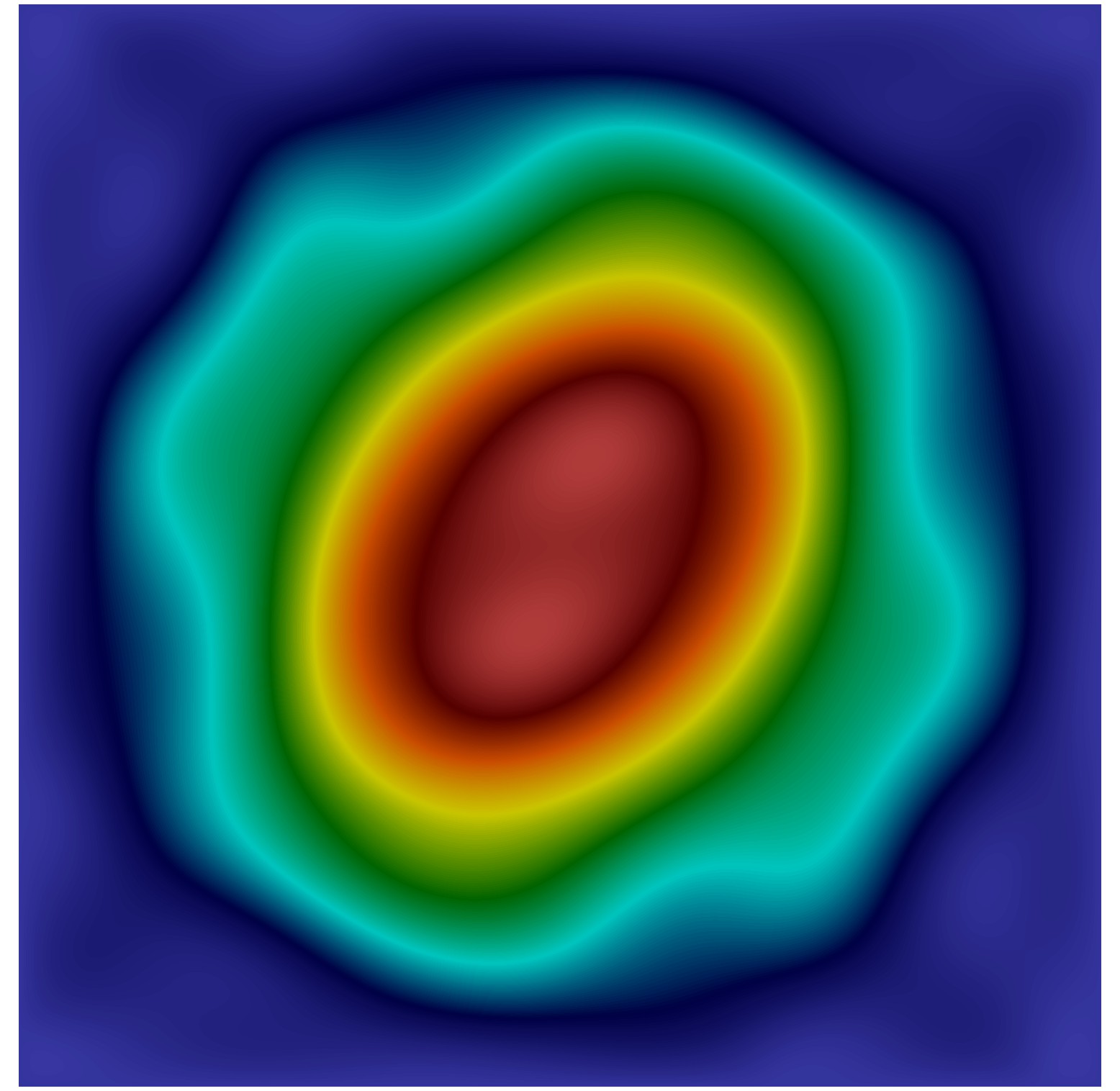}
       \put(30,100){$Re = 800$}
      \end{overpic}
            \begin{overpic}[width=0.073\textwidth]{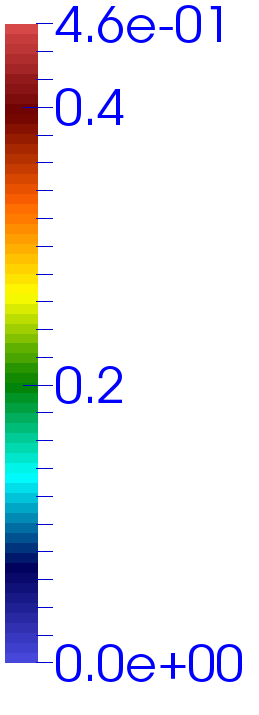}
      \end{overpic}\\
 \begin{overpic}[width=0.22\textwidth]{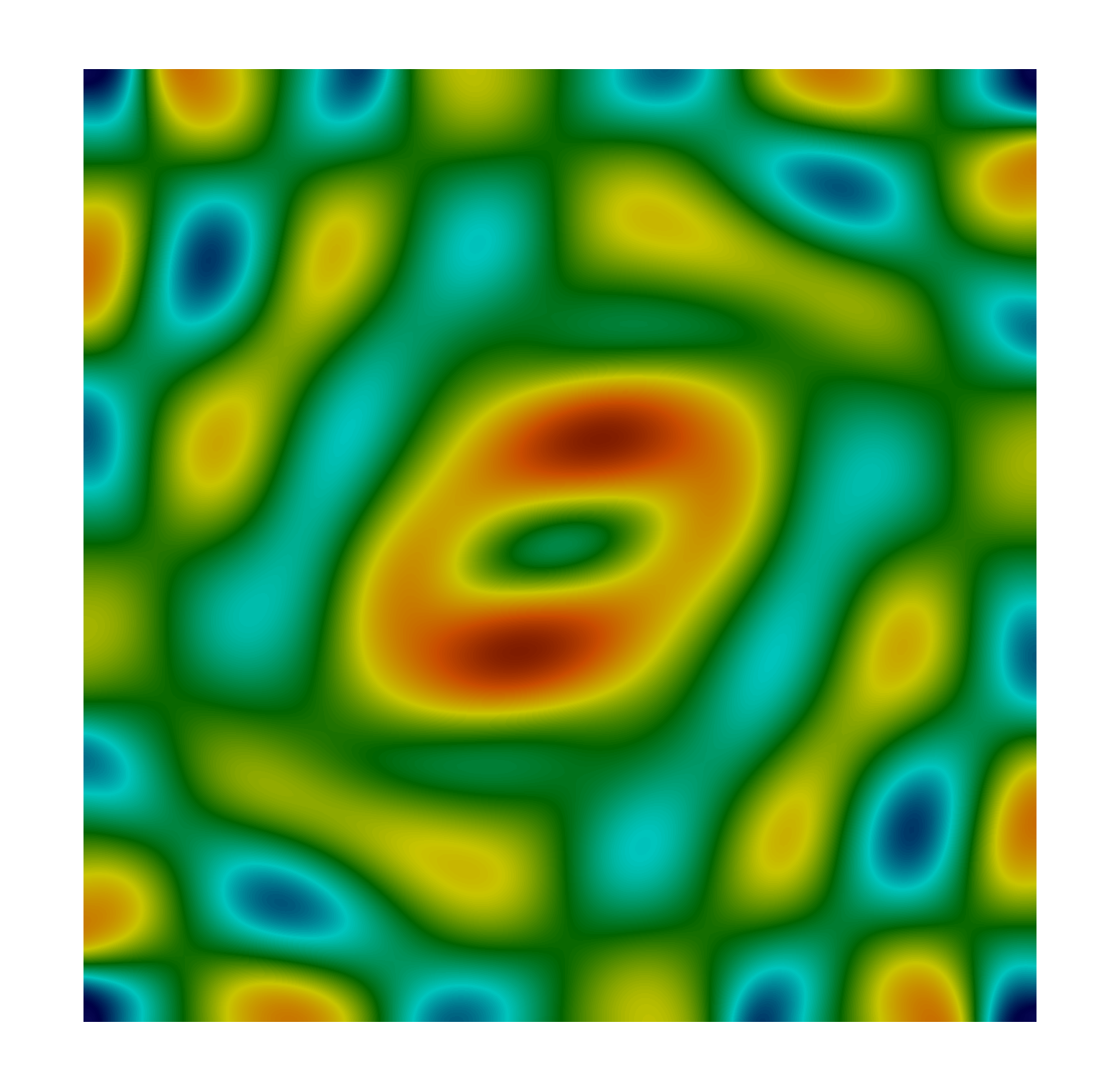}
 \put(-5,50){$\omega$}
      \end{overpic}
 \begin{overpic}[width=0.22\textwidth]{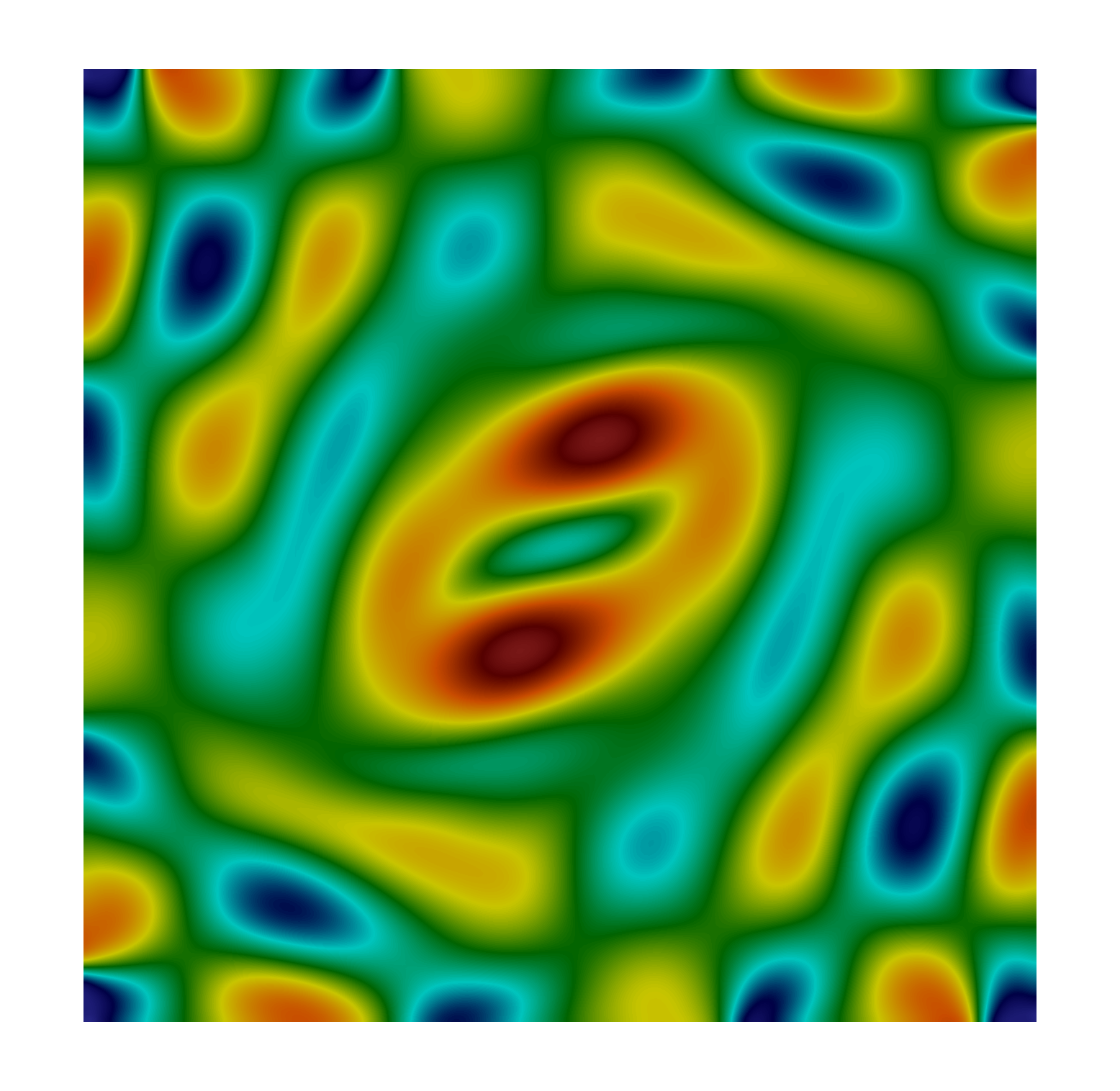}
      \end{overpic}
 \begin{overpic}[width=0.22\textwidth]{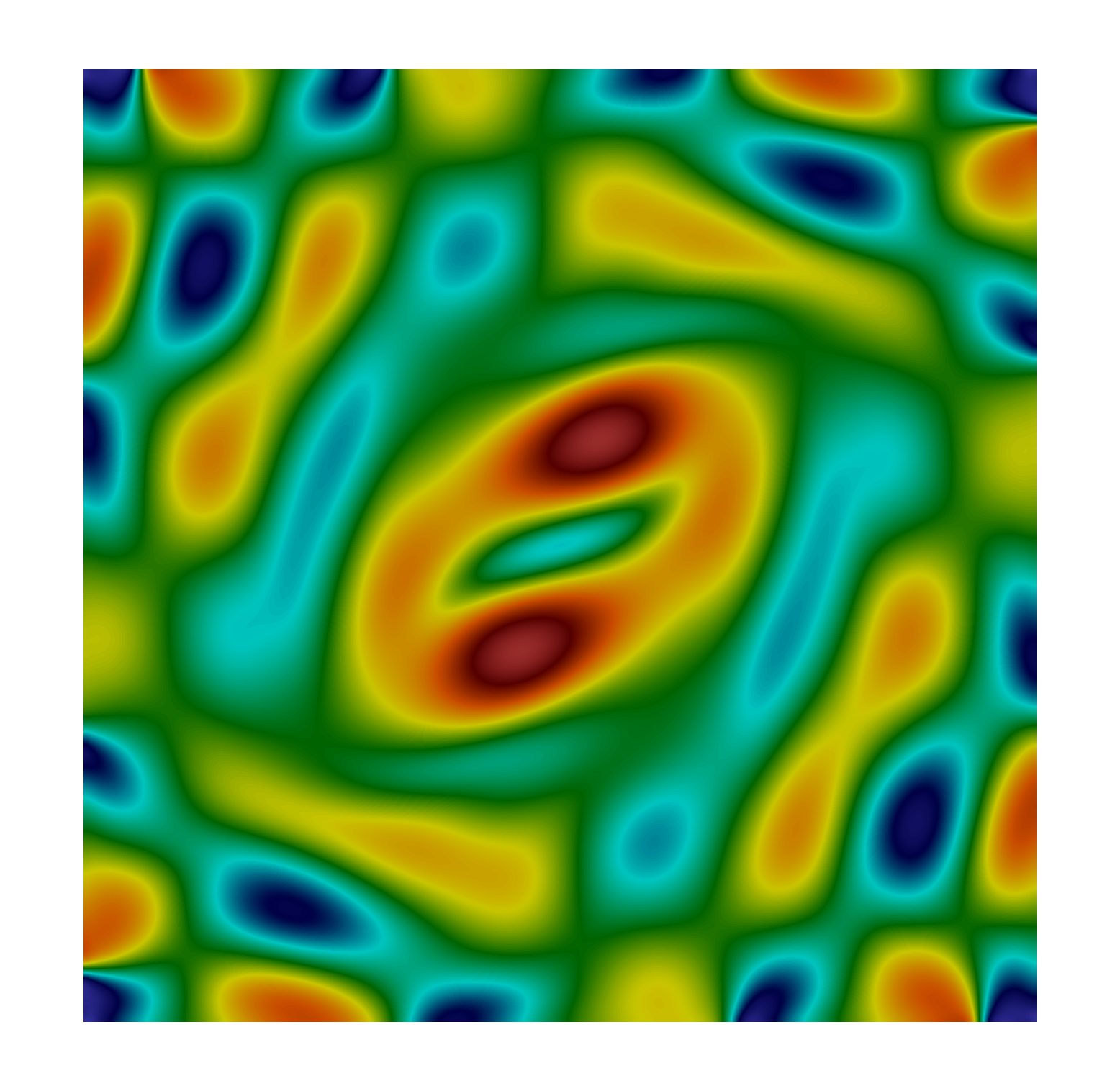}
      \end{overpic}
 \begin{overpic}[width=0.22\textwidth]{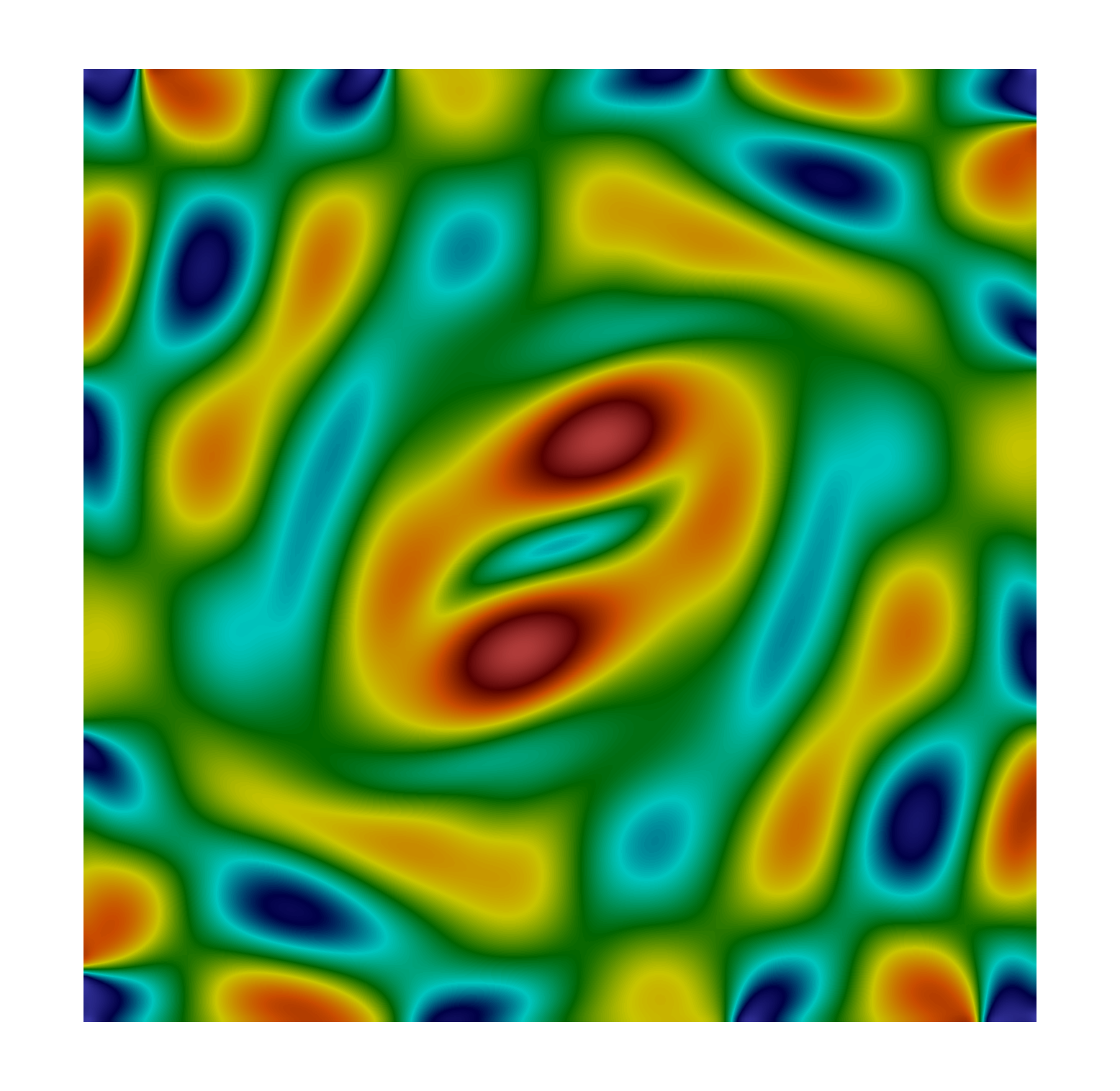}
      \end{overpic}
            \begin{overpic}[width=0.073\textwidth]{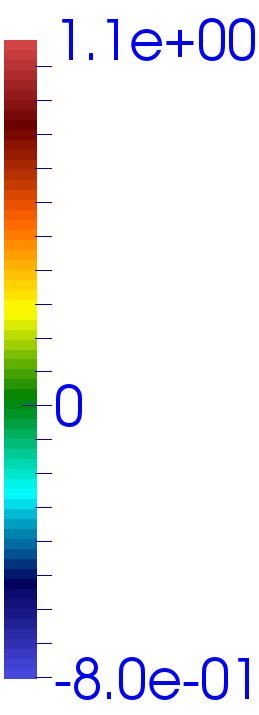}
      \end{overpic}\\
\caption{ROM validation - $Re$ parameterization: stream function $\psi$ (first row) and vorticity $\omega$ (second row) computed by the FOM at $Re = 200$ (first column), $Re = 400$ (second column), $Re = 600$ (third column), and $Re = 800$ (fourth column), at time $t = 10$ for $\gamma = 0.09$.}
\label{fig:param_FOM_Re}
\end{figure}

Let us start with the parametrization with respect to $Re$ and set $\gamma = 0.09$. To train the ROM, we choose a uniform sample distribution in the range $Re \in [200, 800]$ with 4 sampling points: 200, 400, 600, and 800. For each value of the Reynolds number in the training set, a simulation is run over time interval $(0, 10]$. Fig.~\ref{fig:param_FOM_Re} displays the stream function and vorticity fields computed by the FOM throughout the training set under consideration at $t = 10$. We observe that the stream function does not significantly vary as $Re$ changes, 
while the differences in the vorticity are more evident. 

\begin{figure}
\centering
 \begin{overpic}[width=0.45\textwidth]{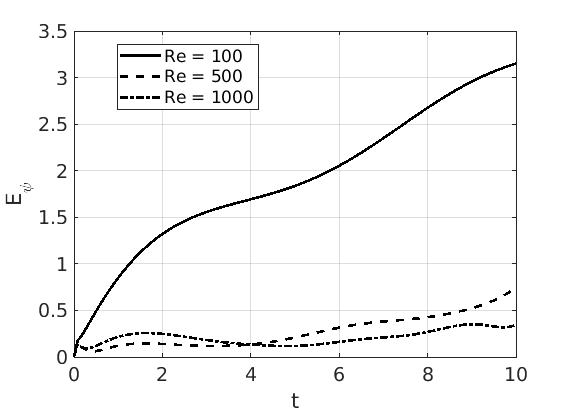}
        %\put(35,18){FOM}
        %\put(-8,7){$\u$}
      \end{overpic}
 \begin{overpic}[width=0.45\textwidth]{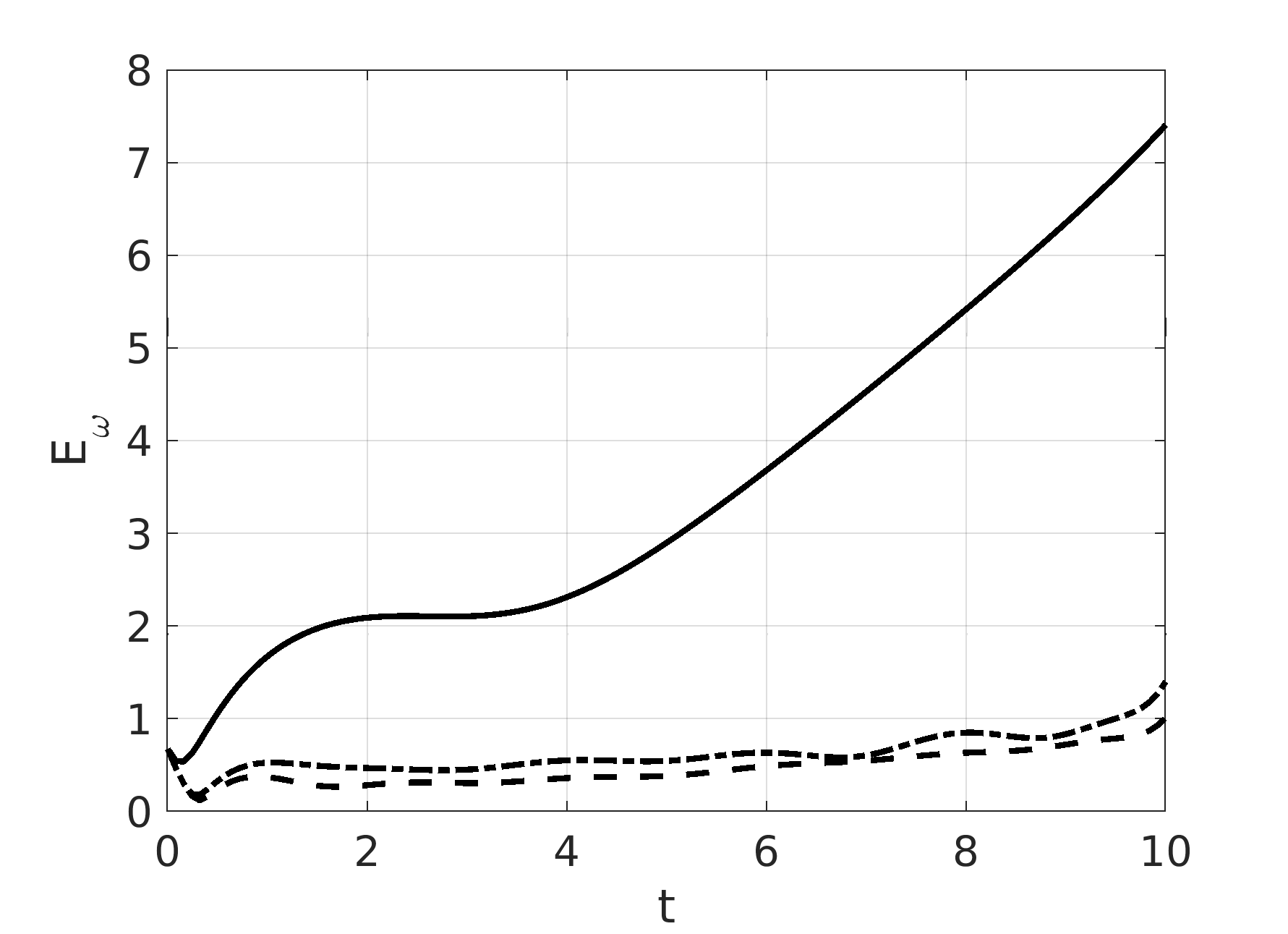}
        %\put(35,18){ROM}
      \end{overpic}
      %\vskip .2cm
       %\begin{overpic}[width=0.3\textwidth]{img/U_diff_1s.png}
        %\put(40,40){ROM}
      %\end{overpic}\\
\caption{ROM validation - $Re$ parameterization: time history of error \eqref{eq:error1} for stream function $\psi$ (left) and vorticity $\omega$ (right) for three different test values: $Re = 100$, $Re = 500$ and $Re = 1000$.}
\label{fig:errors_Re}
\end{figure}

Based on the results presented for the time reconstruction (Sec.~\ref{sec:time}), the snapshots are collected every 0.08 s. So, we collect a total of $4 \times 125 = 500$ snapshots. We set the threshold for the selection of the eigenvalues to $1e-5$, resulting in 6 modes for $\psi$ and 11 modes for $\omega$. %\michele{(Anna come noterai abbiamo 11 modi anziche' 14 stavolta per la omega, probabilmente perche' non consideriamo la seconda parte dell'intervallo di tempo. E' pur vero che consideriamo pero' la variabilita' in Reynolds che, pero', evidentemente, non è tanto "pronunciata" come quella in tempo. Non so se mi sono spiegato bene :), in ogni caso eviterei di commentare al riguardo essendo una nota un po' troppo qualitativa, che ne pensi? Idem per il caso parametrico in gamma.)} \anna{Sono d'accordo! Mentre leggevo non mi era nemmeno sorto il dubbio, quindi magari non sorge nemmeno in un altro lettore.}\michele{OK! :) }
We take three different test values to evaluate the performance of the parametrized ROM: one value ($Re = 500$) in the range under consideration but not in the training set and two values
($Re = 100, 1000$) outside the range under consideration. The latter cases are more challenging. 
Fig.~\ref{fig:errors_Re} shows error \eqref{eq:error1} for the stream function and vorticity over time for the three values of $Re$. We see that for the interpolatory test value $Re = 500$ both the errors are below 1\% over the entire time interval. %\michele{(come fatto in tempo riporto giusto gli upper bound degli errori relativi "a occhio",  direi che visto il contesto e' la cosa migliore, informativa e smooth allo stesso tempo :) Altrimenti possiamo fare come di consueto delle tabelle! Dimmi tu :) )}. \anna{Mi sembra che vada bene cosi'!} \michele{Va bene!} 
As for the extrapolatory test values, the errors
for $Re = 1000$ are comparable to the errors for $Re = 500$, while the errors for $Re = 100$ are much larger for $Re = 100$ (up to about 3\% for $\psi$ and up to about 7\% for $\omega$). %\textcolor{red}{(Addirittura per 1000 andiamo leggermente meglio come puoi ben vedere ma non sottolinerei la cosa, anche perche' siamo davvero li' e giustificare una differenza cosi' piccola è difficile, che ne dici?)} \anna{Esatto, non lo farei notare.} \michele{D'accordo!}
The poorer performance of our ROM at $Re = 100$ could be due to the fact that the vorticity computed for $Re = 200$ looks pretty different from the vorticity at the higher $Re$ included in the offline database (see Fig.~\ref{fig:param_FOM_Re}). Thus, we suspect that more solutions for lower values of $Re$ would have to be included into the training set in order to obtain a
more accurate reconstruction of the flow field at $Re = 100$.

\begin{figure}
\centering
 \begin{overpic}[width=0.175\textwidth]{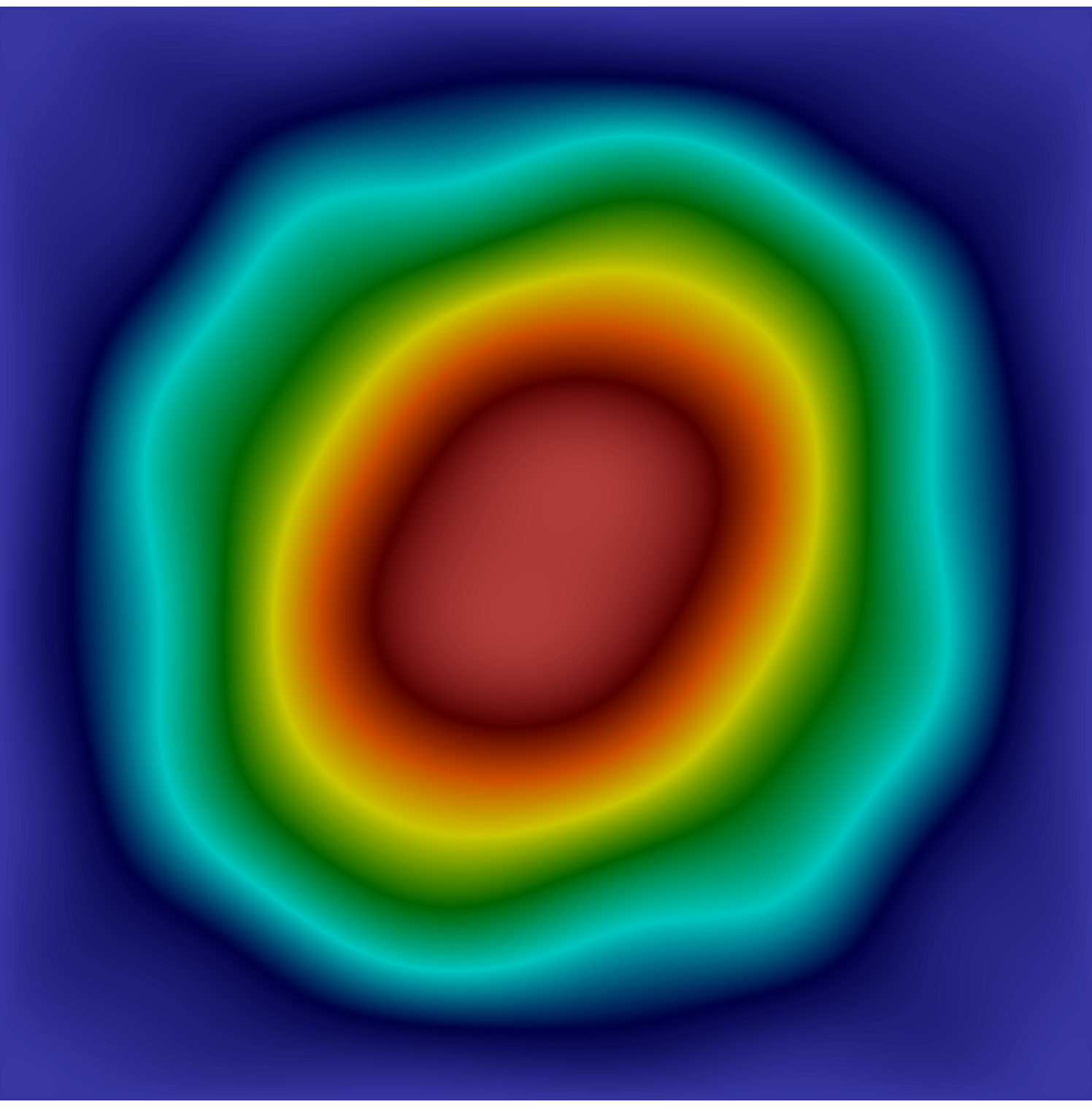}
        \put(20,103){$Re = 100$}
        \put(-42,45){FOM}   %-32
      \end{overpic}
 \begin{overpic}[width=0.075\textwidth]{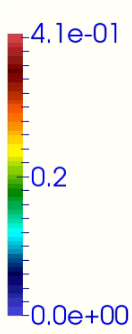}
      \end{overpic}
       \begin{overpic}[width=0.175\textwidth]{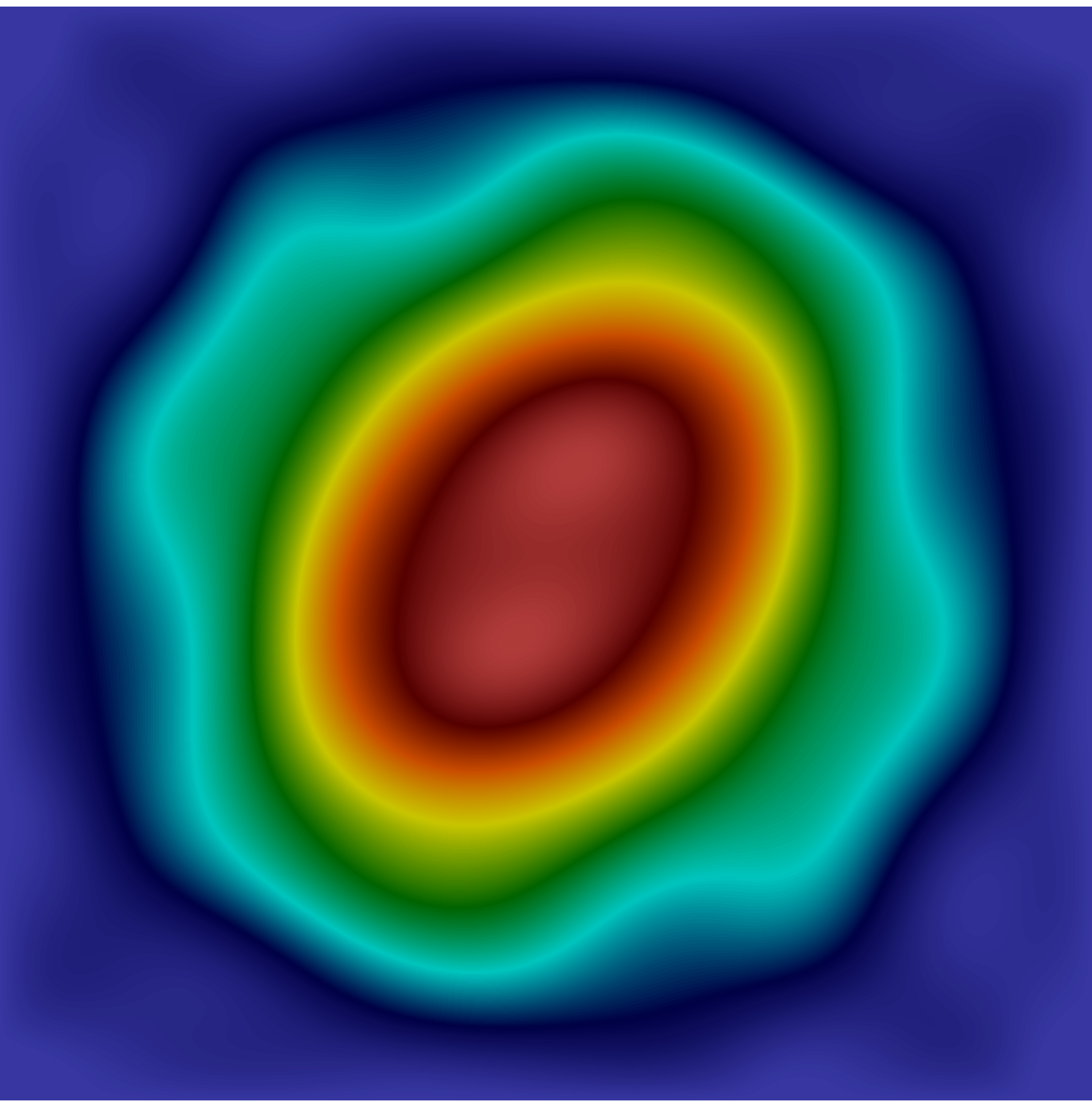}
              \put(20,103){$Re = 500$}
      \end{overpic}
 \begin{overpic}[width=0.075\textwidth]{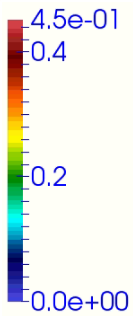}
      \end{overpic}
   \begin{overpic}[width=0.175\textwidth]{img/psiFOM_500_0_09.png}
     \put(20,103){$Re = 1000$}
      \end{overpic}
 \begin{overpic}[width=0.075\textwidth]{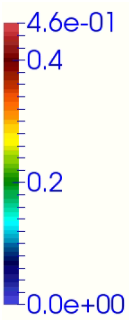}
      \end{overpic}\\
      \vspace{0.5cm}
 \begin{overpic}[width=0.175\textwidth]{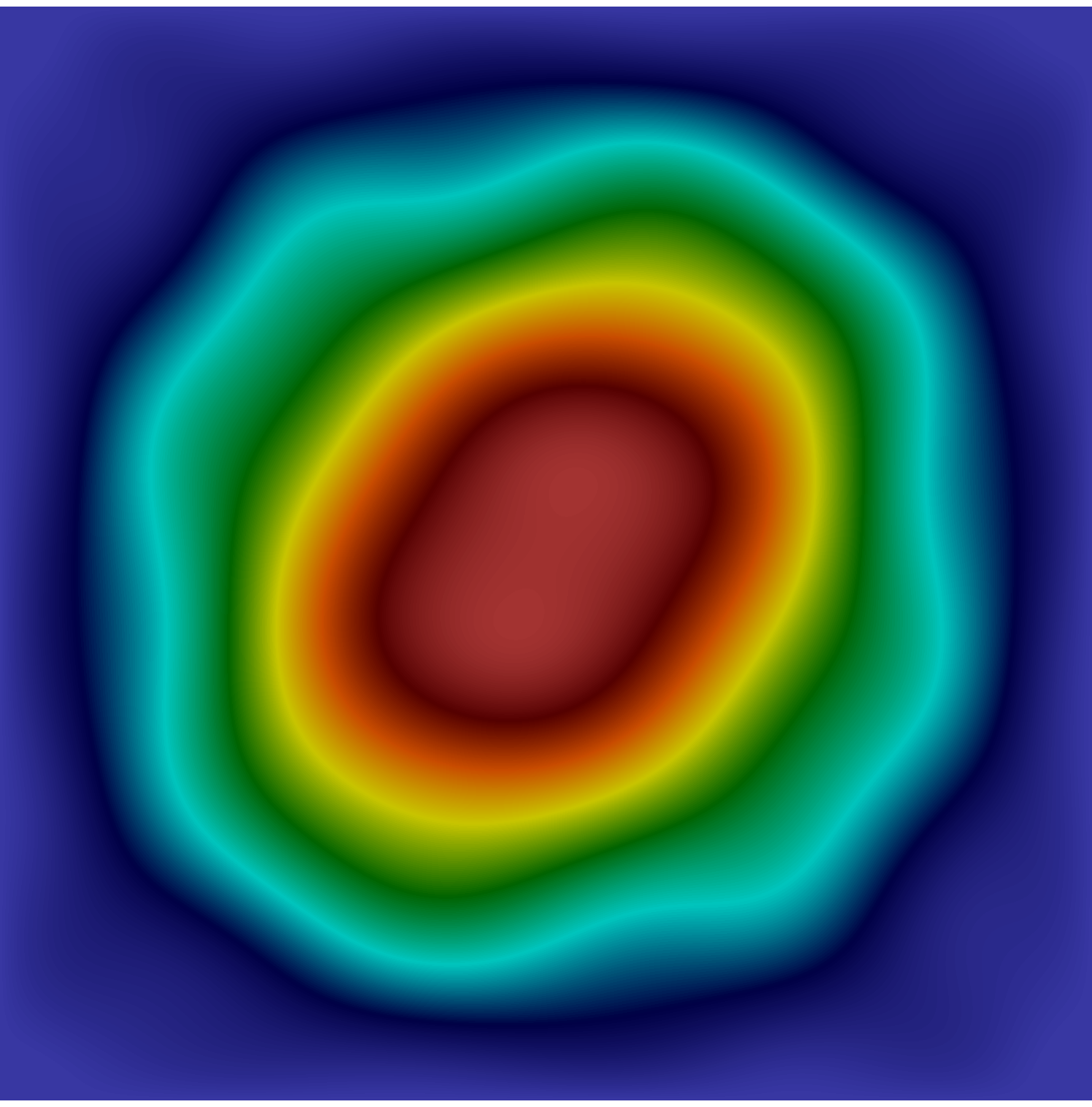}
\put(-42,45){ROM}
      \end{overpic}
 \begin{overpic}[width=0.075\textwidth]{img/legenda_100_0_09_psi.png}
      \end{overpic}
   \begin{overpic}[width=0.175\textwidth]{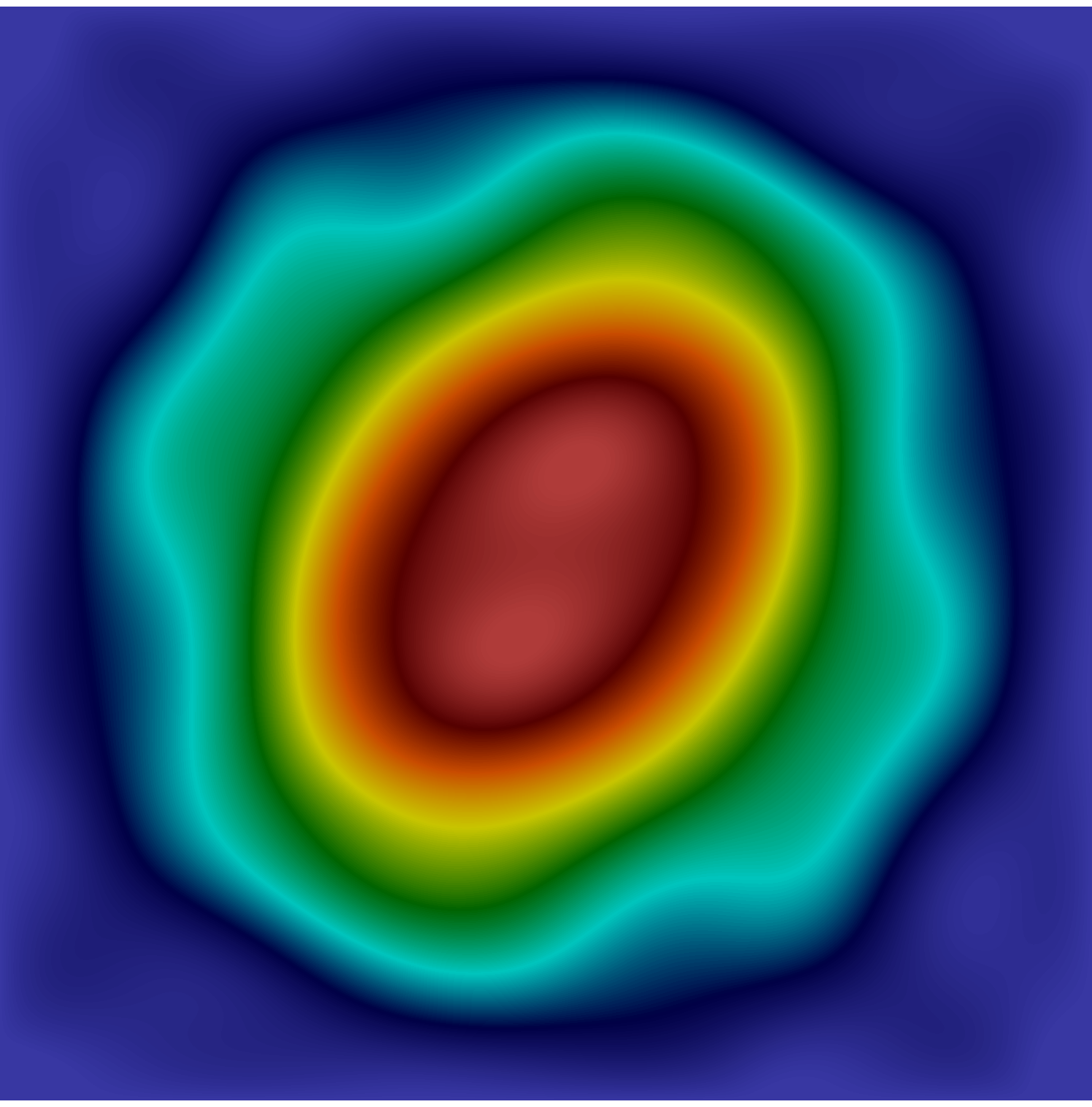}
      \end{overpic}
 \begin{overpic}[width=0.075\textwidth]{img/legenda_500_0_09_psi.png}
      \end{overpic}
      \begin{overpic}[width=0.175\textwidth]{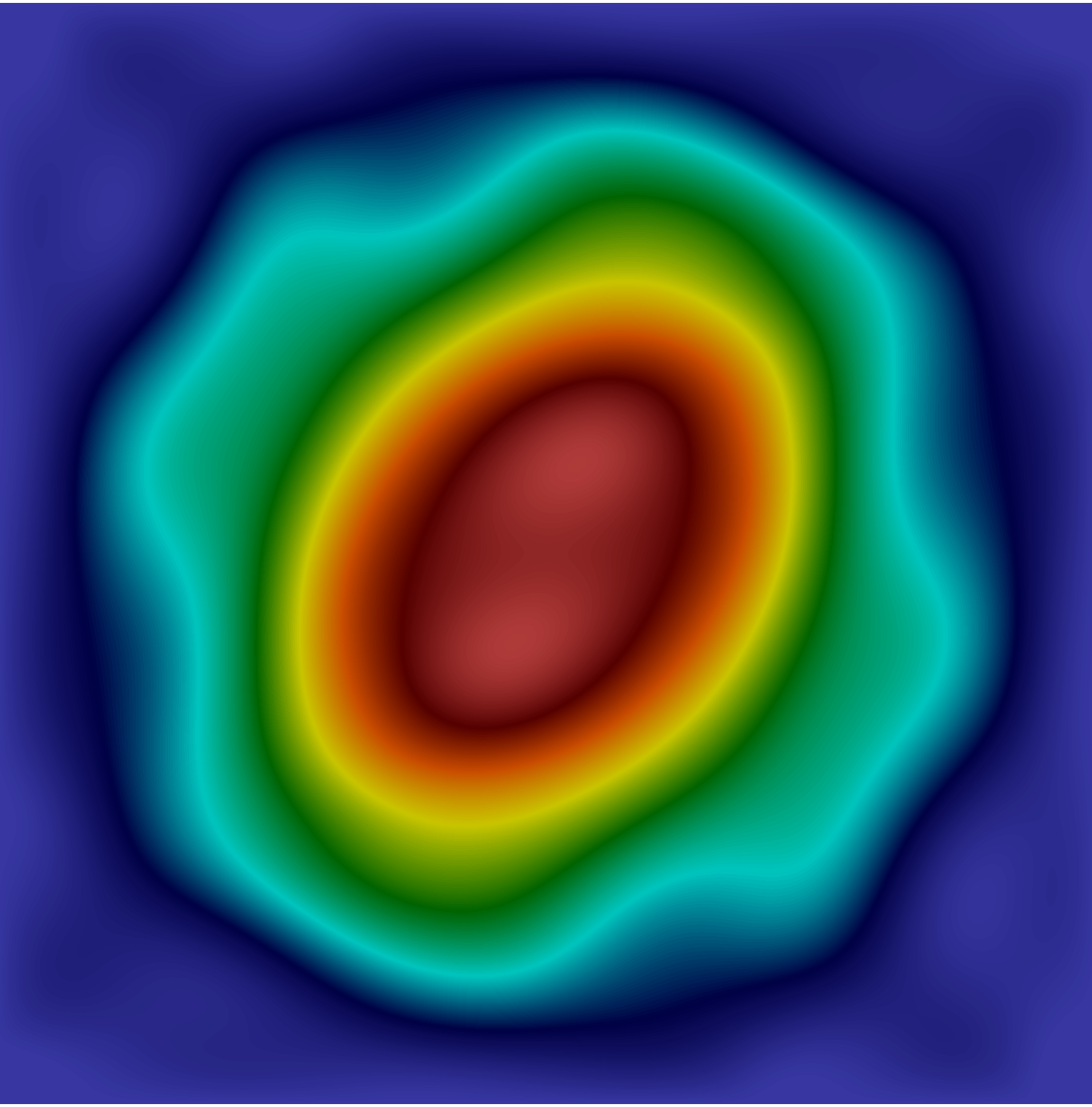}
      \end{overpic}
 \begin{overpic}[width=0.075\textwidth]{img/legenda_1000_0_09_psi.png}
      \end{overpic}\\
            \vspace{0.5cm}
\begin{overpic}[width=0.18\textwidth]{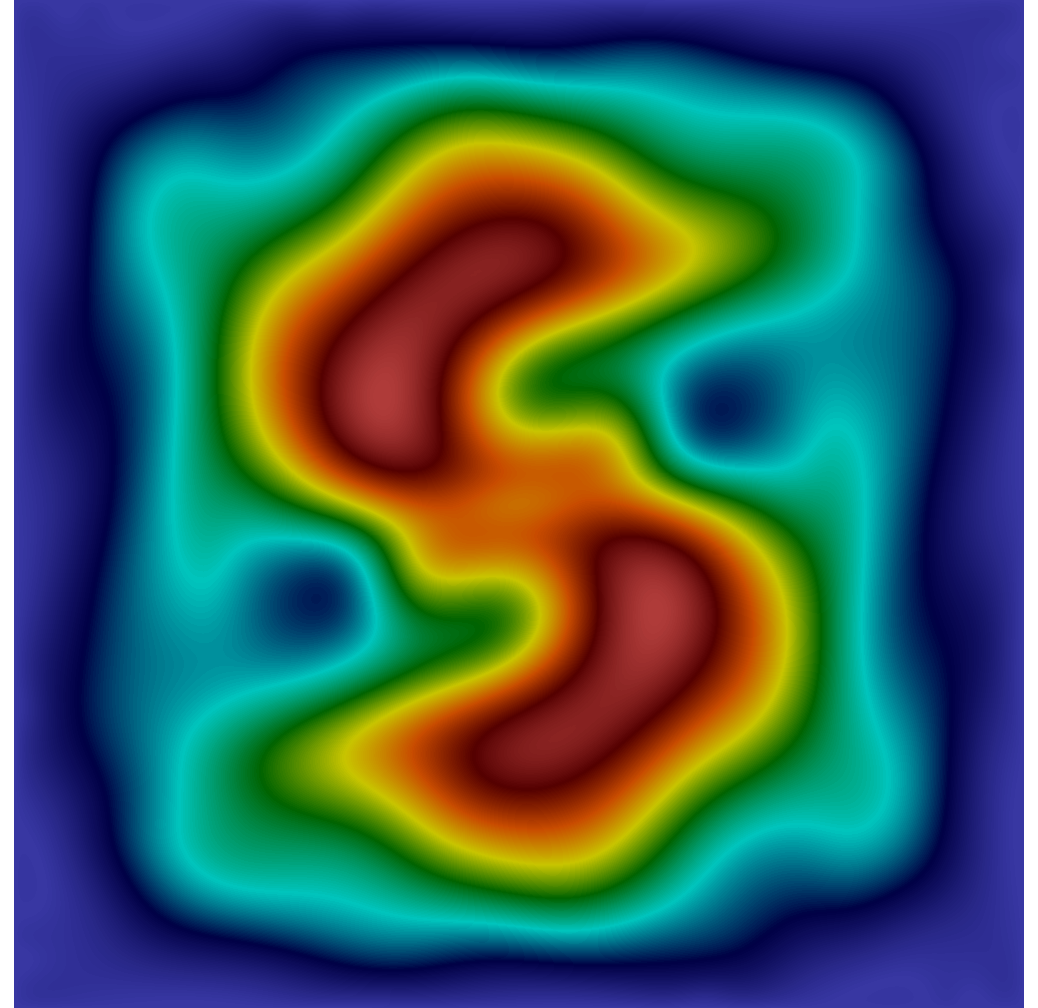}
\put(-32,45){Diff.}
      \end{overpic}
 \begin{overpic}[width=0.077\textwidth]{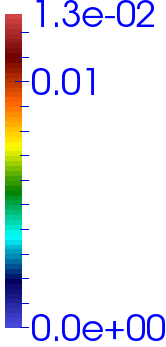}
        %\put(35,18){ROM}
      \end{overpic}
   \begin{overpic}[width=0.18\textwidth]{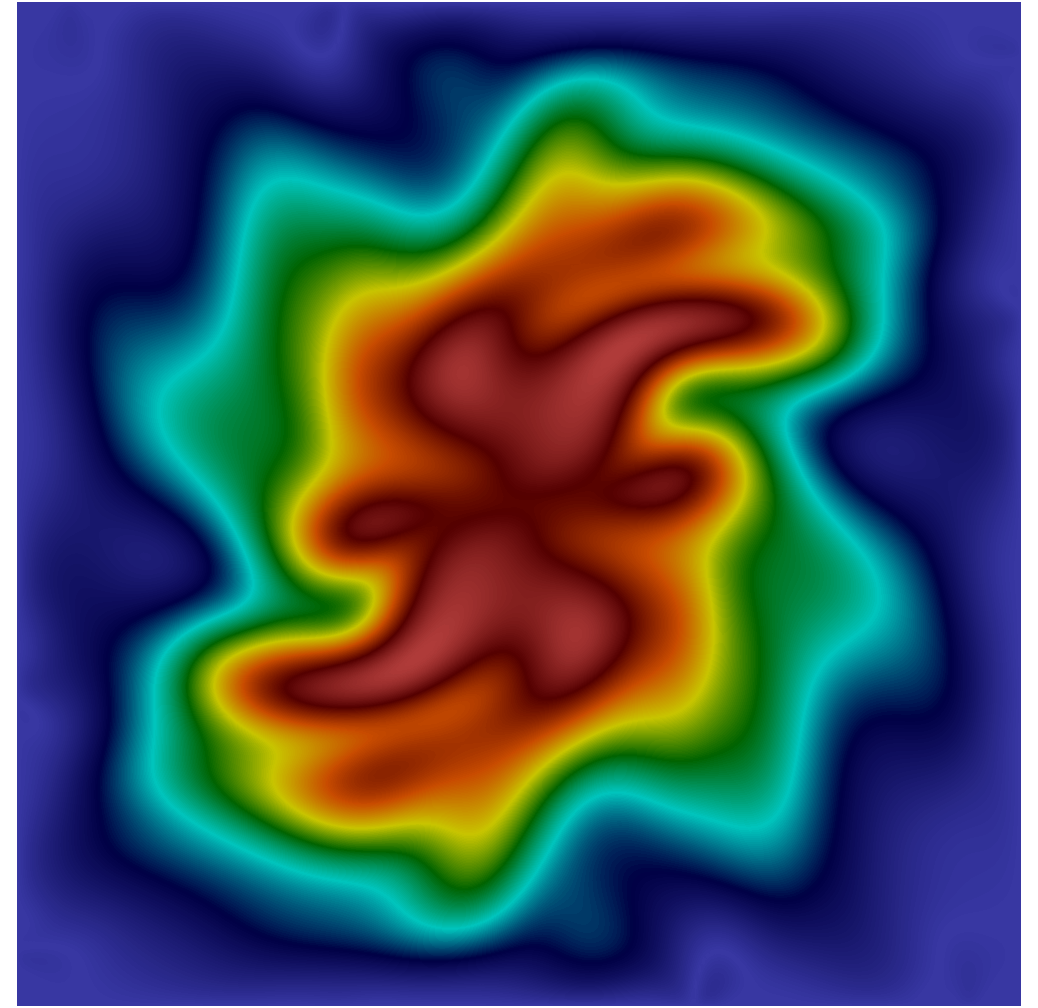}
        %\put(35,18){FOM}
        %\put(-8,7){$\u$}
      \end{overpic}
 \begin{overpic}[width=0.075\textwidth]{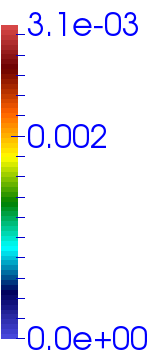}
        %\put(35,18){ROM}
      \end{overpic}
      \begin{overpic}[width=0.18\textwidth]{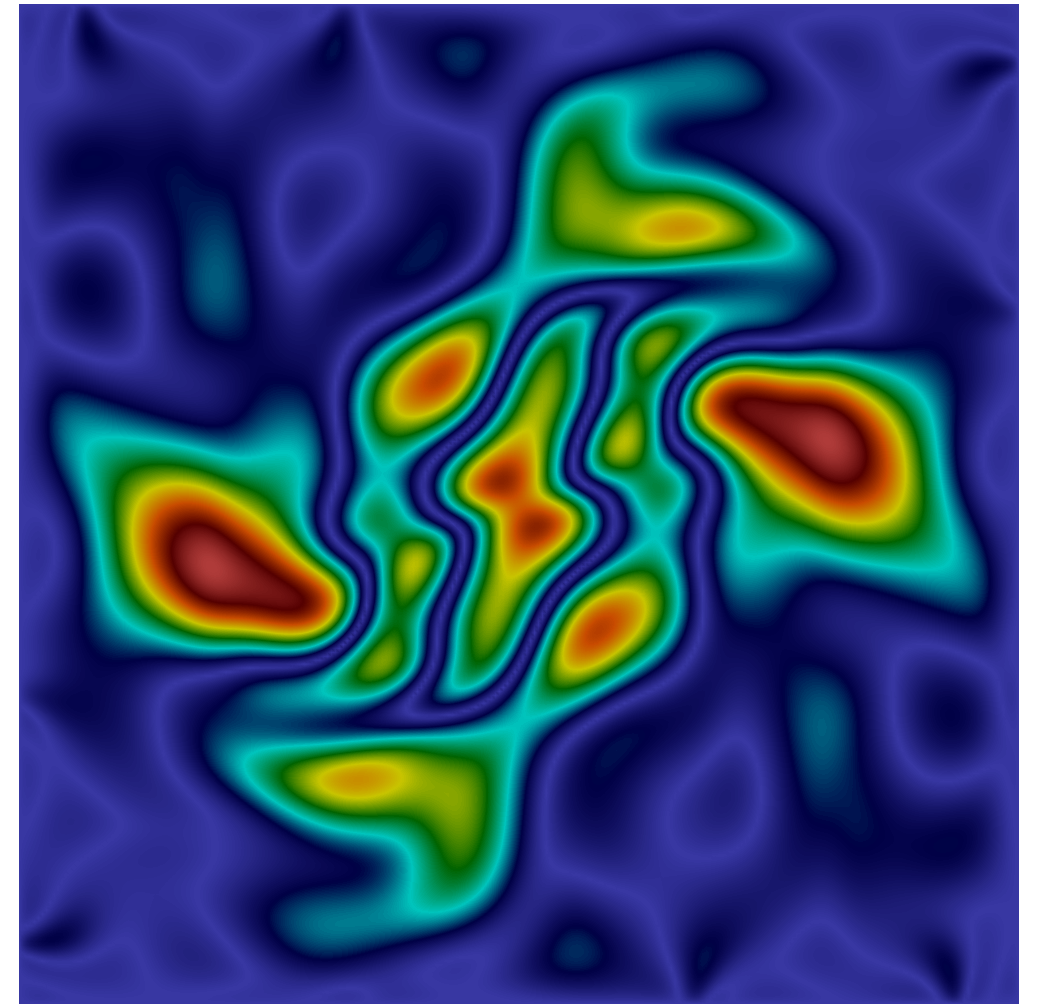}
        %\put(35,18){FOM}
        %\put(-8,7){$\u$}
      \end{overpic}
 \begin{overpic}[width=0.075\textwidth]{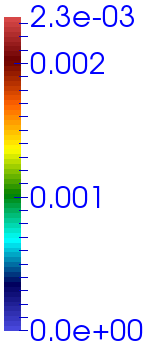}
        %\put(35,18){ROM}
      \end{overpic}\\
\caption{ROM validation - $Re$ parameterization: stream function  $\psi$ computed by the FOM (first row) and the ROM (second row), and difference between the two fields in absolute value (third row) for $Re = 100$ (first column), $Re = 500$ (second column) and $Re = 1000$ (third one) at time $t = 10$. Six modes for $\psi$ 
were considered.}
\label{fig:errors_psi_absolute_Re}
\end{figure}

\begin{figure}
\centering
 \begin{overpic}[width=0.175\textwidth]{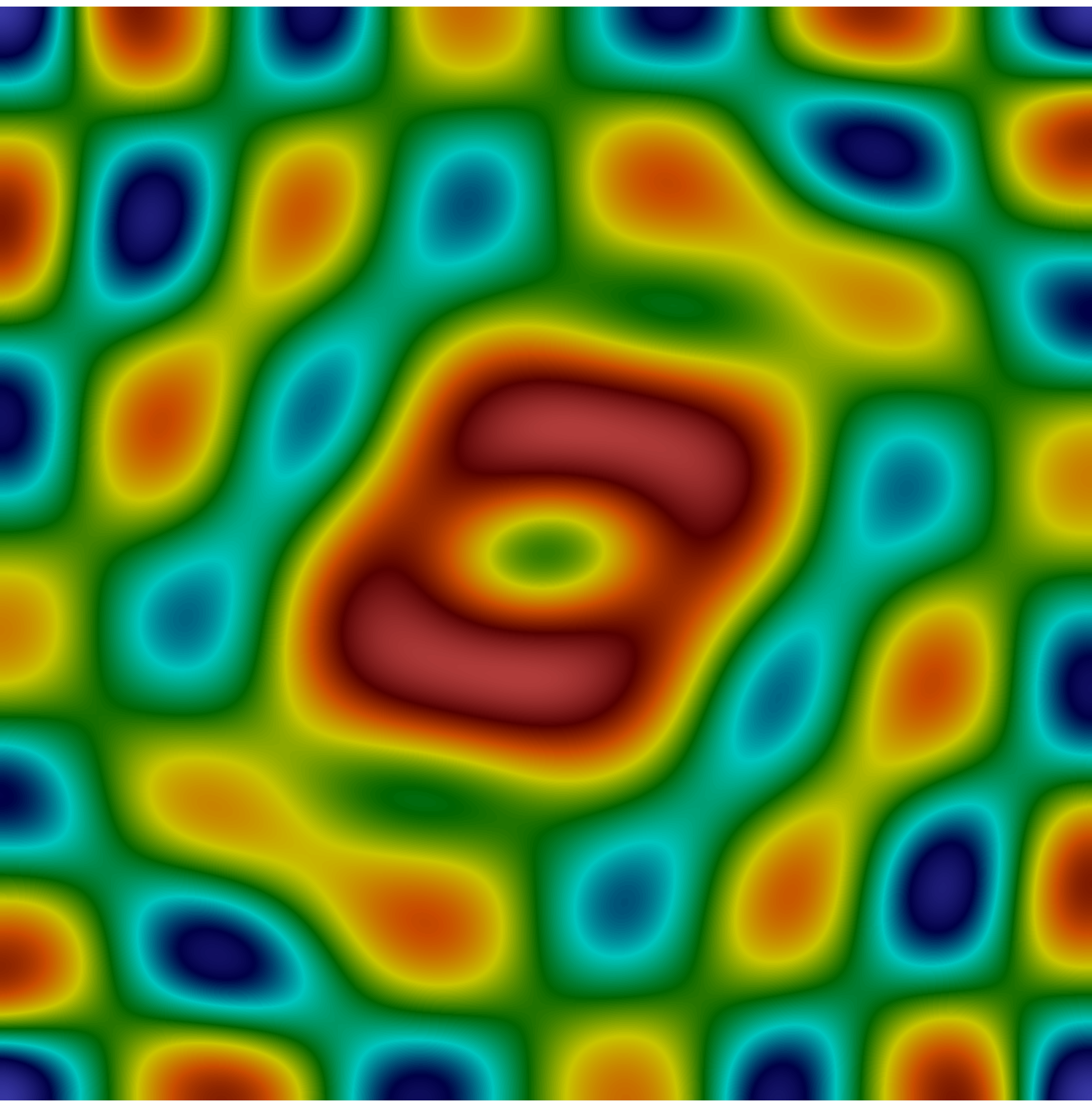}
        \put(25,103){$Re = 100$}
        \put(-42,45){FOM}
      \end{overpic}
 \begin{overpic}[width=0.065\textwidth]{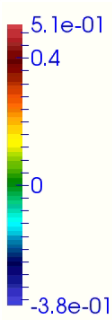}
      \end{overpic}\quad
       \begin{overpic}[width=0.175\textwidth]{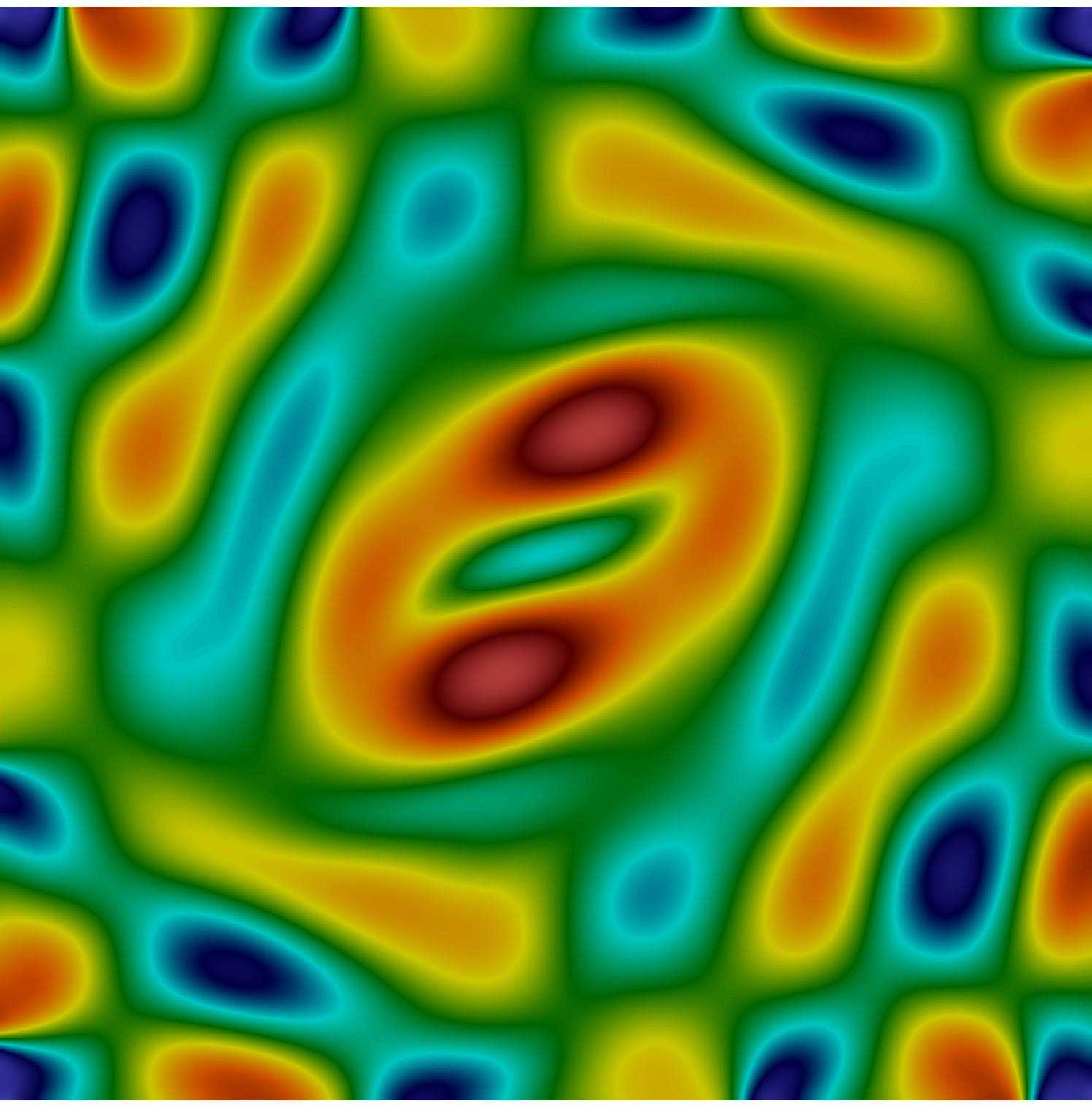}
        \put(25,103){$Re = 500$}
      \end{overpic}
 \begin{overpic}[width=0.065\textwidth]{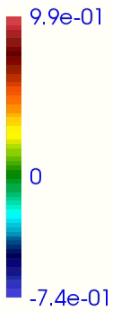}
      \end{overpic} \quad
   \begin{overpic}[width=0.175\textwidth]{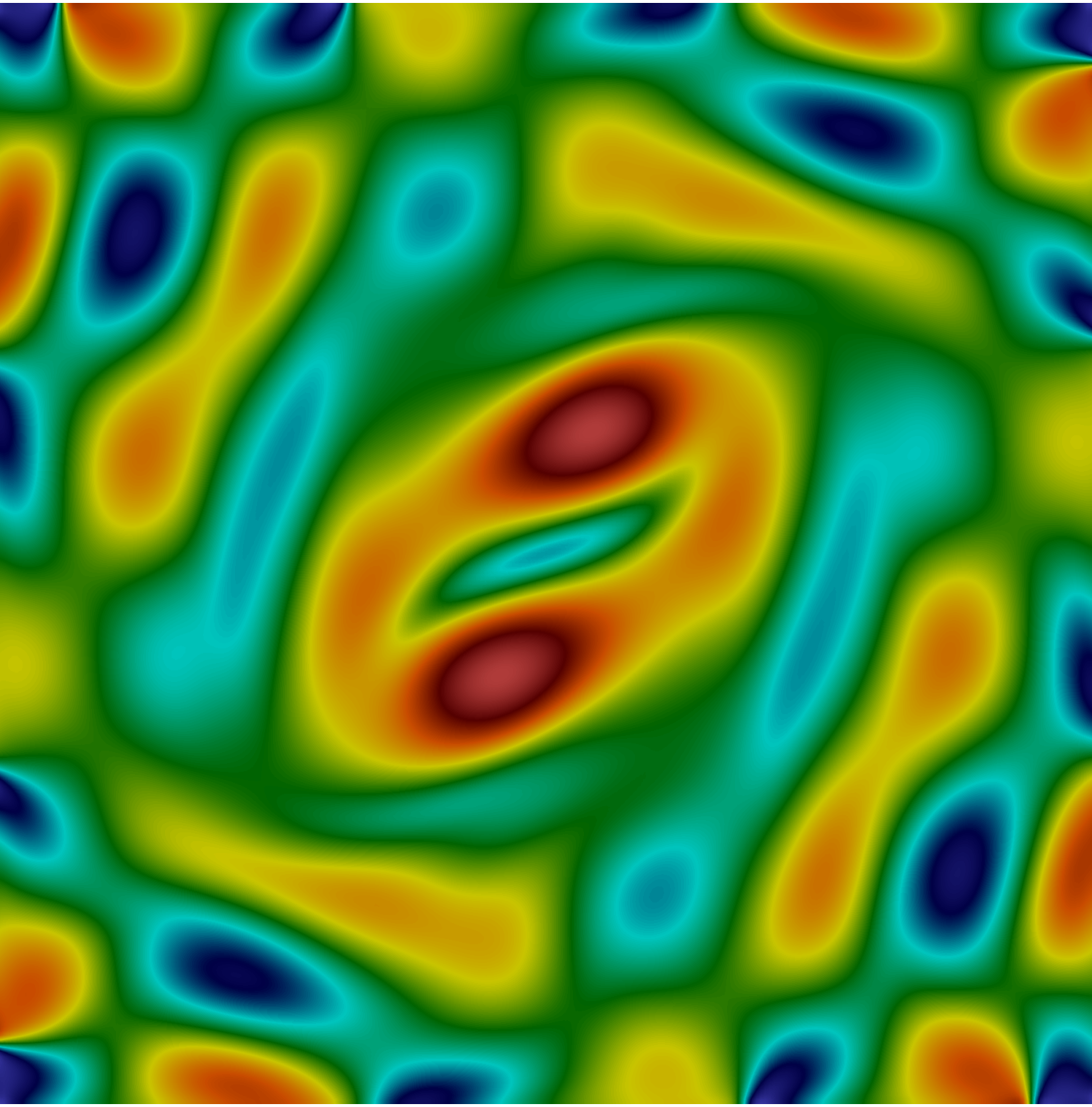}
        \put(25,103){$Re = 1000$}
      \end{overpic}
 \begin{overpic}[width=0.065\textwidth]{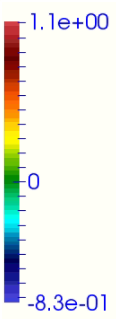}
      \end{overpic}\\
      \vspace{0.5cm}
 \begin{overpic}[width=0.175\textwidth]{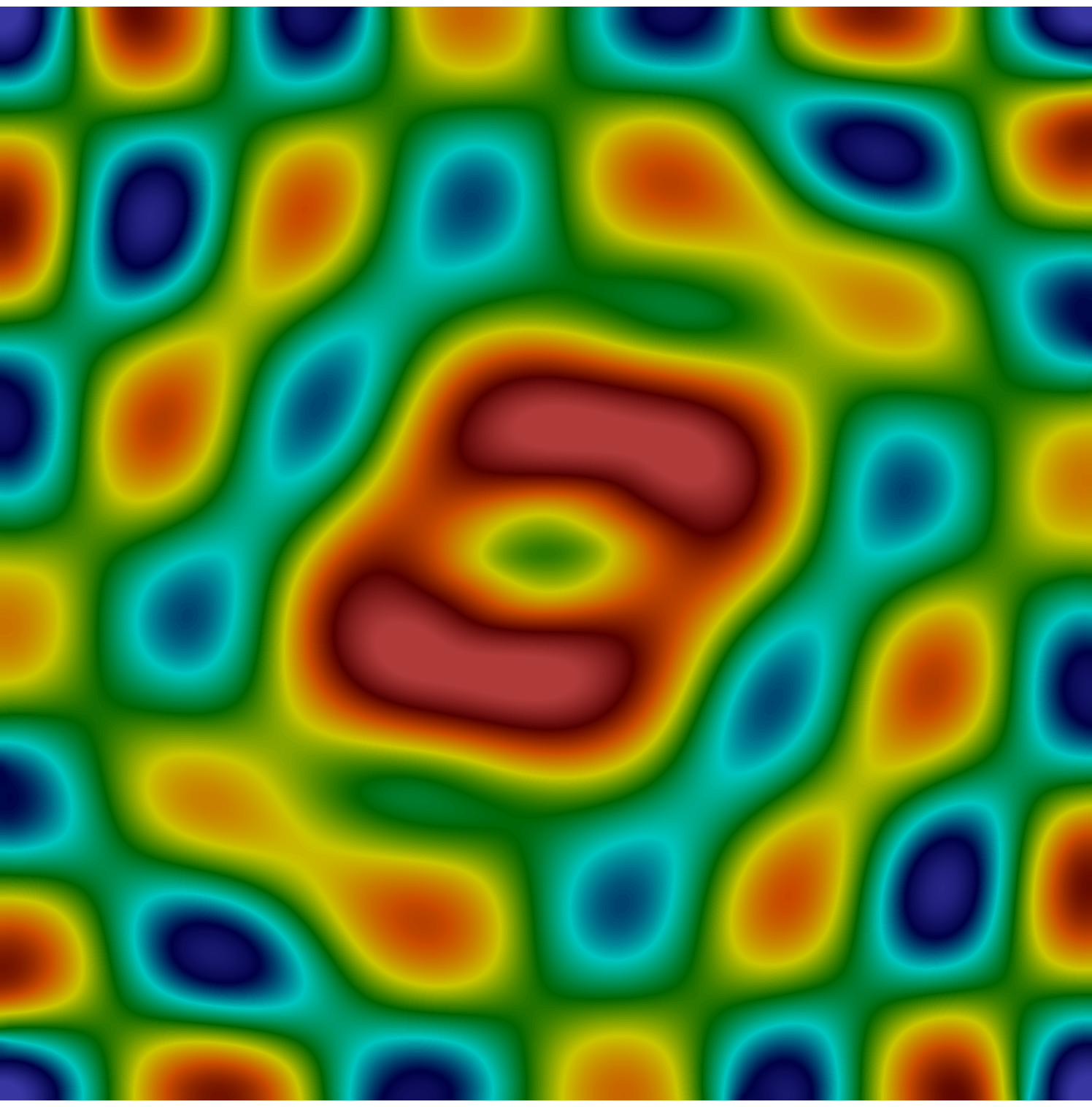}
        \put(-42,45){ROM}
      \end{overpic}
 \begin{overpic}[width=0.065\textwidth]{img/legenda_100_0_09_zeta.png}
      \end{overpic} \quad
   \begin{overpic}[width=0.175\textwidth]{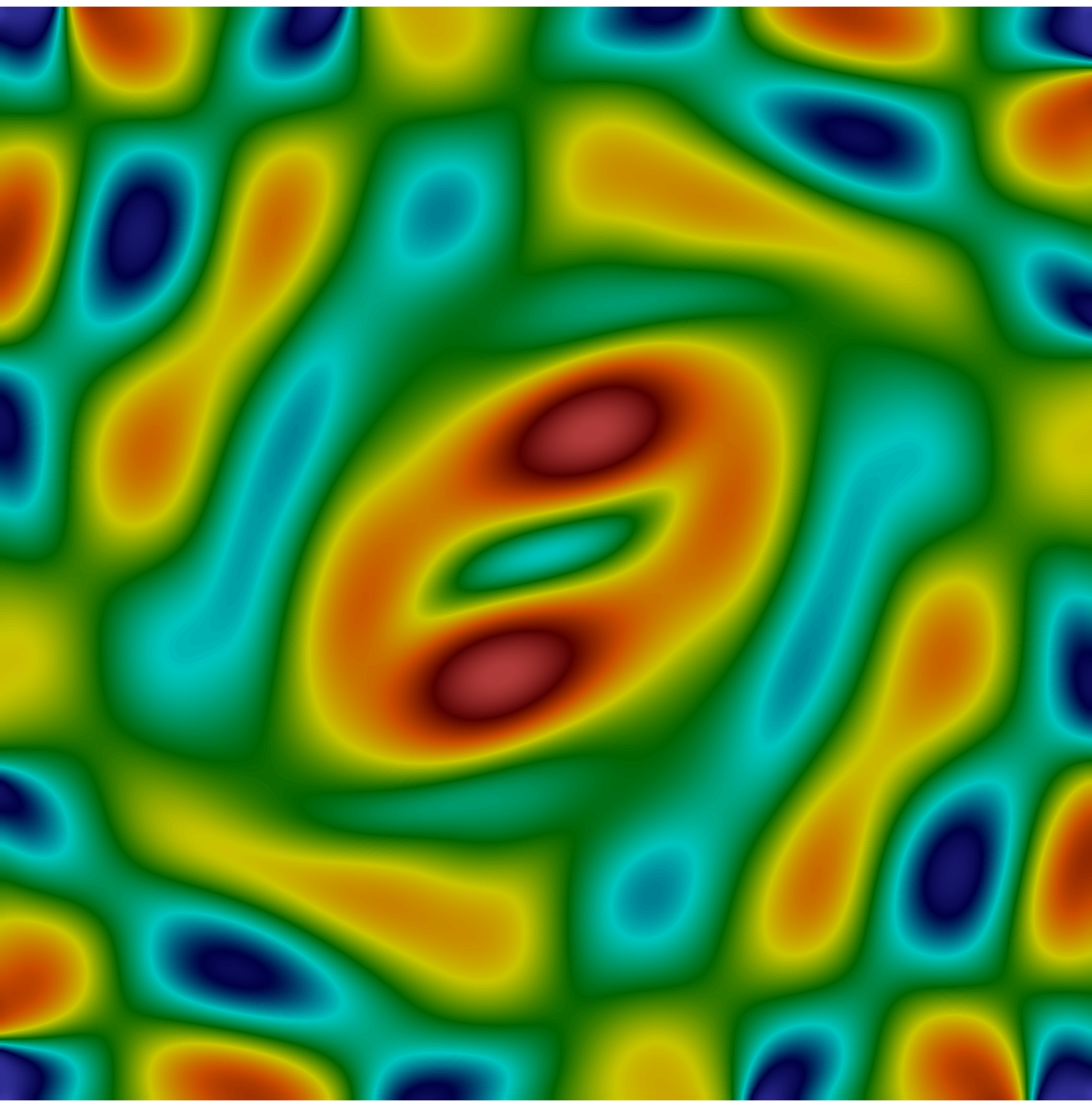}
      \end{overpic}
 \begin{overpic}[width=0.065\textwidth]{img/legenda_500_0_09_zeta.png} 
      \end{overpic} \quad
      \begin{overpic}[width=0.175\textwidth]{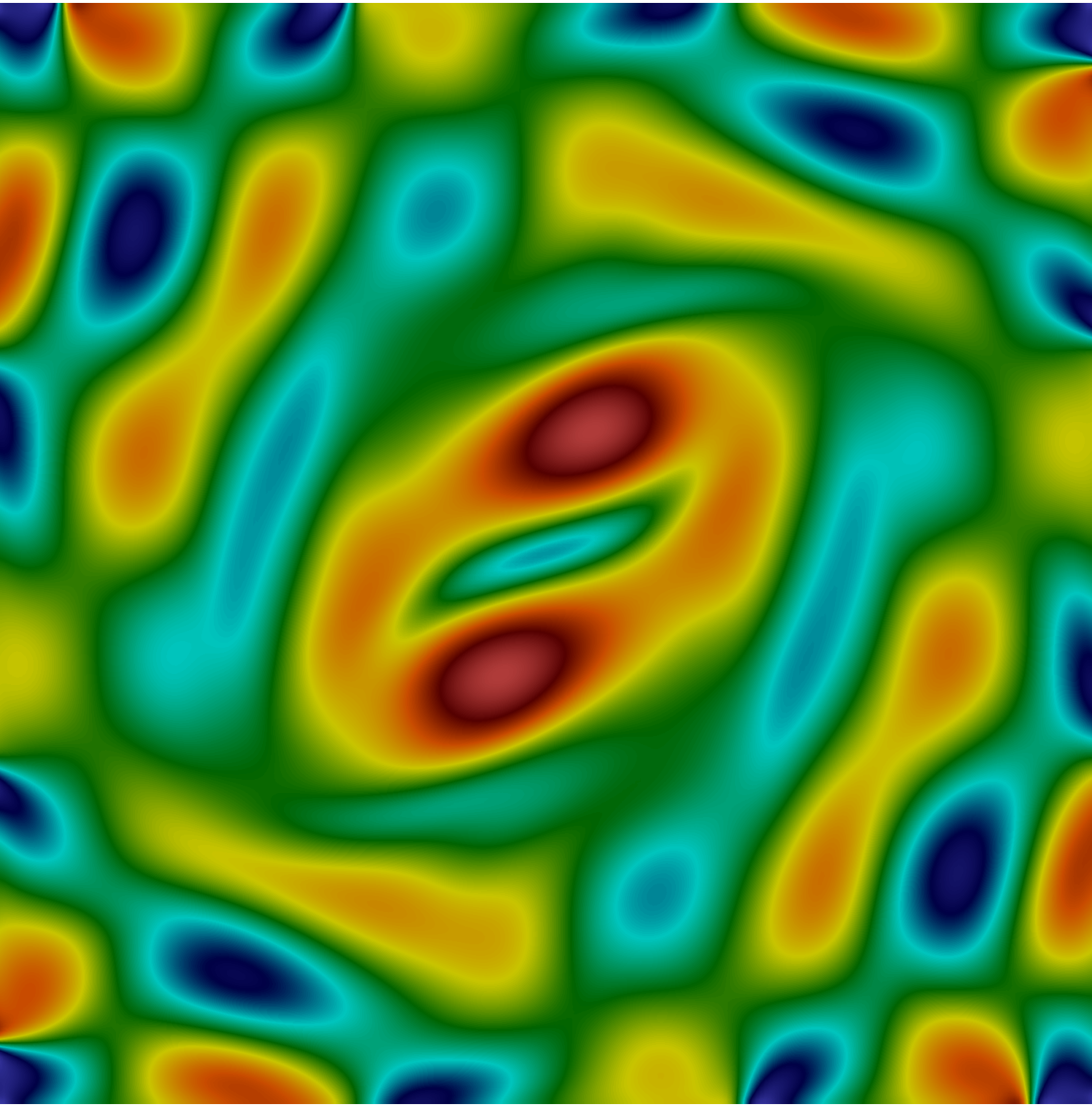}
      \end{overpic}
 \begin{overpic}[width=0.065\textwidth]{img/legenda_1000_0_09_zeta.png}
      \end{overpic}\\
            \vspace{0.5cm}

\begin{overpic}[width=0.185\textwidth]{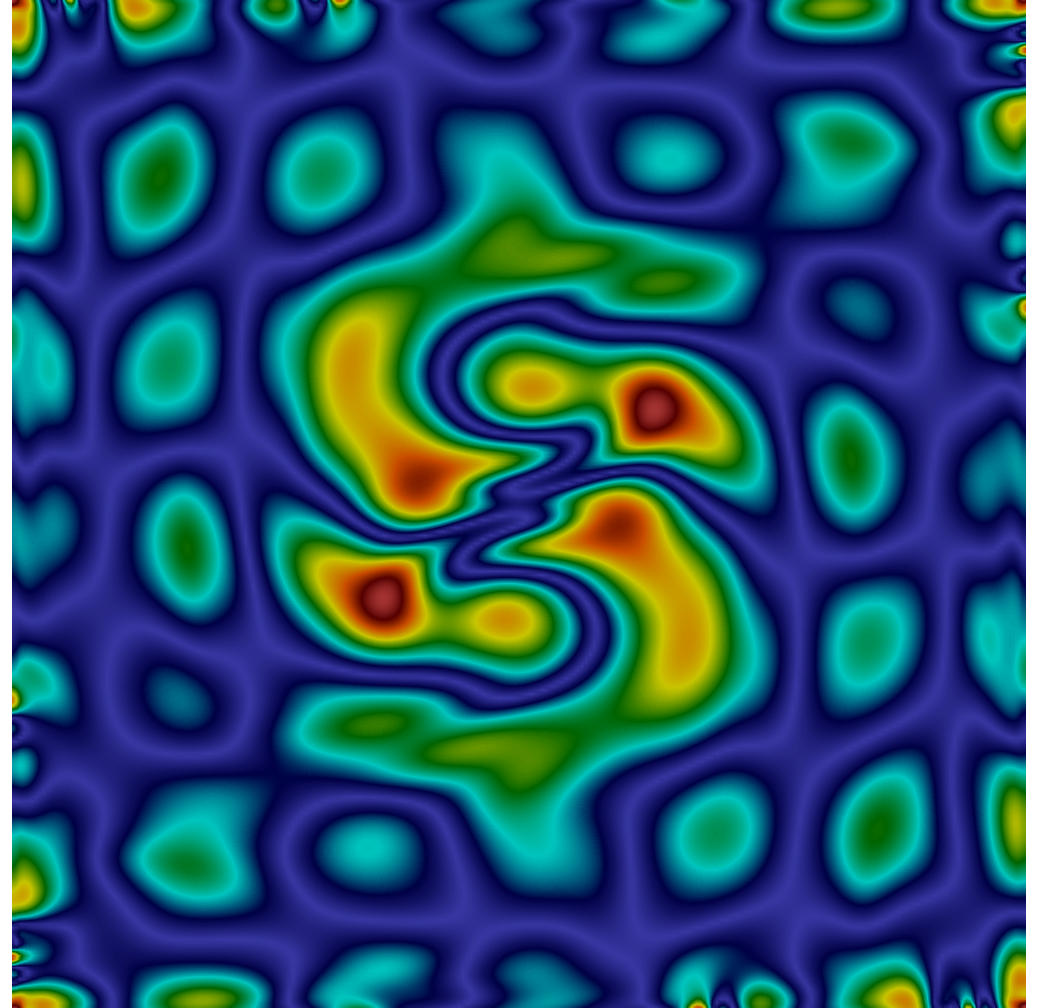}
        \put(-32,45){Diff.}
      \end{overpic}
 \begin{overpic}[width=0.077\textwidth]{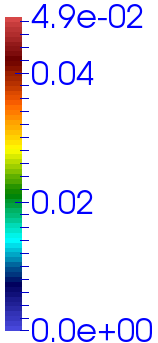}
      \end{overpic}
   \begin{overpic}[width=0.185\textwidth]{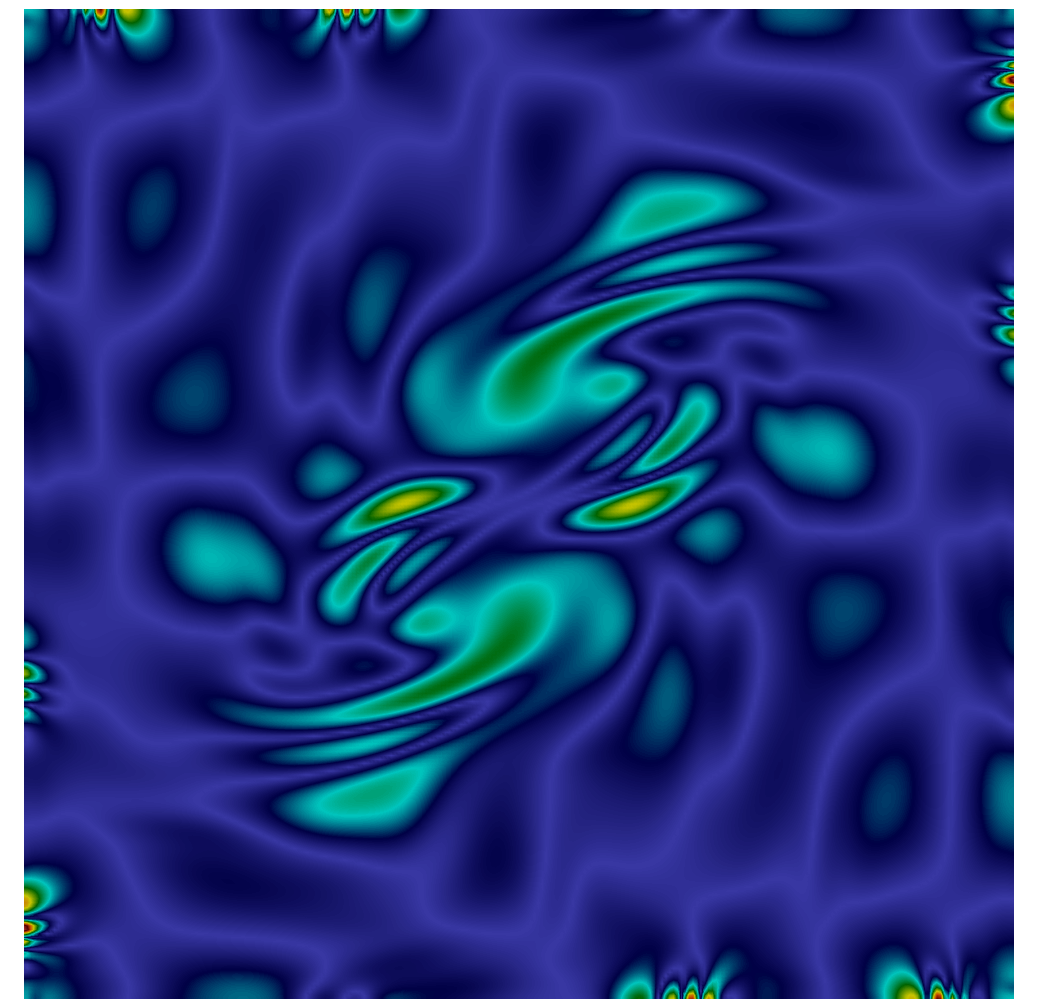}
        %\put(35,18){FOM}
        %\put(-8,7){$\u$}
      \end{overpic}
 \begin{overpic}[width=0.075\textwidth]{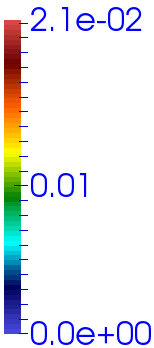}
        %\put(35,18){ROM}
      \end{overpic}
      \begin{overpic}[width=0.185\textwidth]{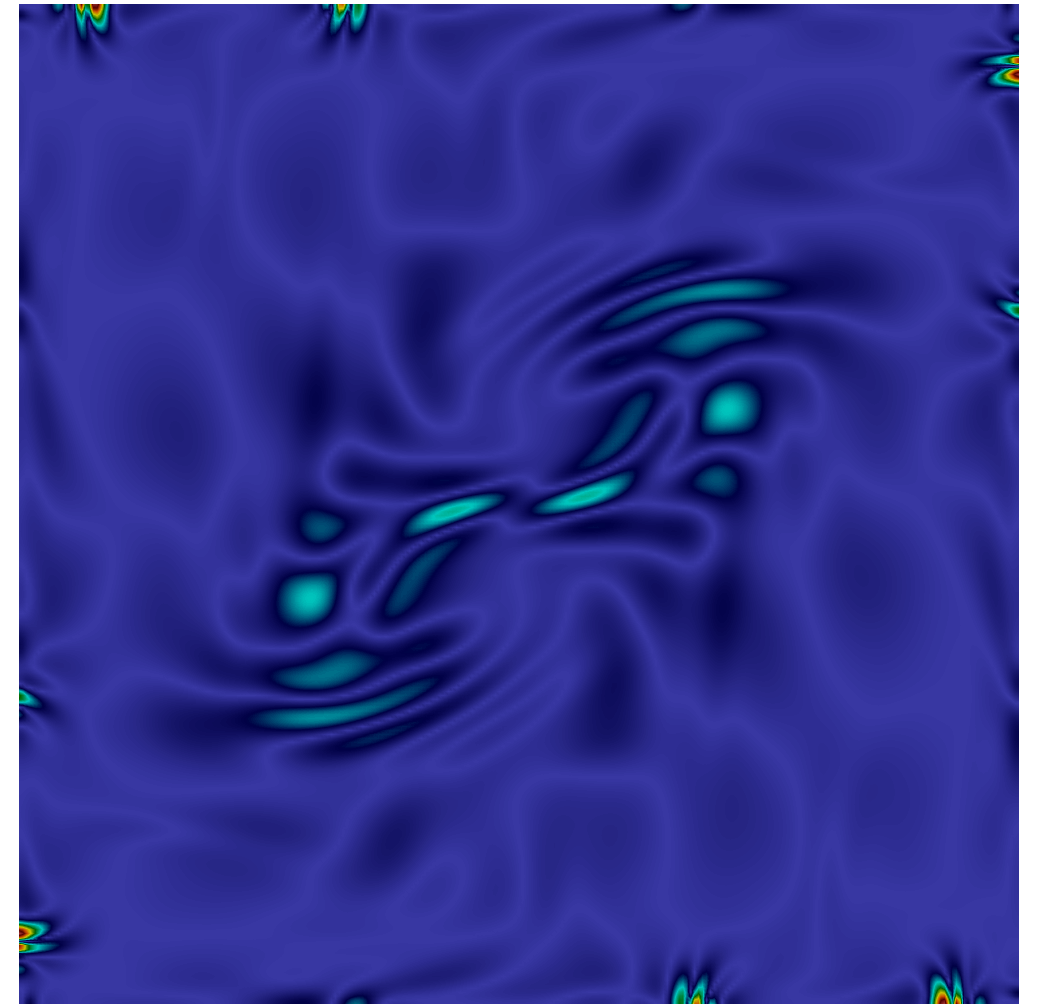}
        %\put(35,18){FOM}
        %\put(-8,7){$\u$}
      \end{overpic}
 \begin{overpic}[width=0.075\textwidth]{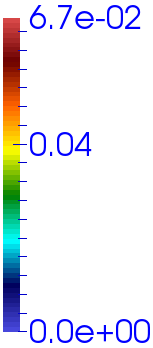}
        %\put(35,18){ROM}
      \end{overpic}\\
      
      %\vskip .2cm
       %\begin{overpic}[width=0.3\textwidth]{img/U_diff_1s.png}
        %\put(40,40){ROM}
      %\end{overpic}\\
\caption{ROM validation - $Re$ parameterization: vorticity $\omega$ computed by the FOM (first row) and the ROM (second row), and difference between the two fields in absolute value (third row) for $Re = 100$ (first column), $Re = 500$ (second column) and $Re = 1000$ (third one) at time $t = 10$. Eleven modes for $\omega$ 
were considered.}
\label{fig:errors_zeta_absolute_Re}
\end{figure}

Fig.~\ref{fig:errors_psi_absolute_Re} and \ref{fig:errors_zeta_absolute_Re} present a qualitative
comparison of the solutions computed by FOM and ROM at $t = 10$
for the the three test value of $Re$.
We observe that our ROM approach provides a good global reconstruction of both stream function and vorticity. 
In fact, the maximum relative difference in absolute value does not exceed $3.1e-2$ for $\psi$ and $6.7e-2$ for $\omega$. 

\begin{figure}
\centering
 \begin{overpic}[width=0.21\textwidth]{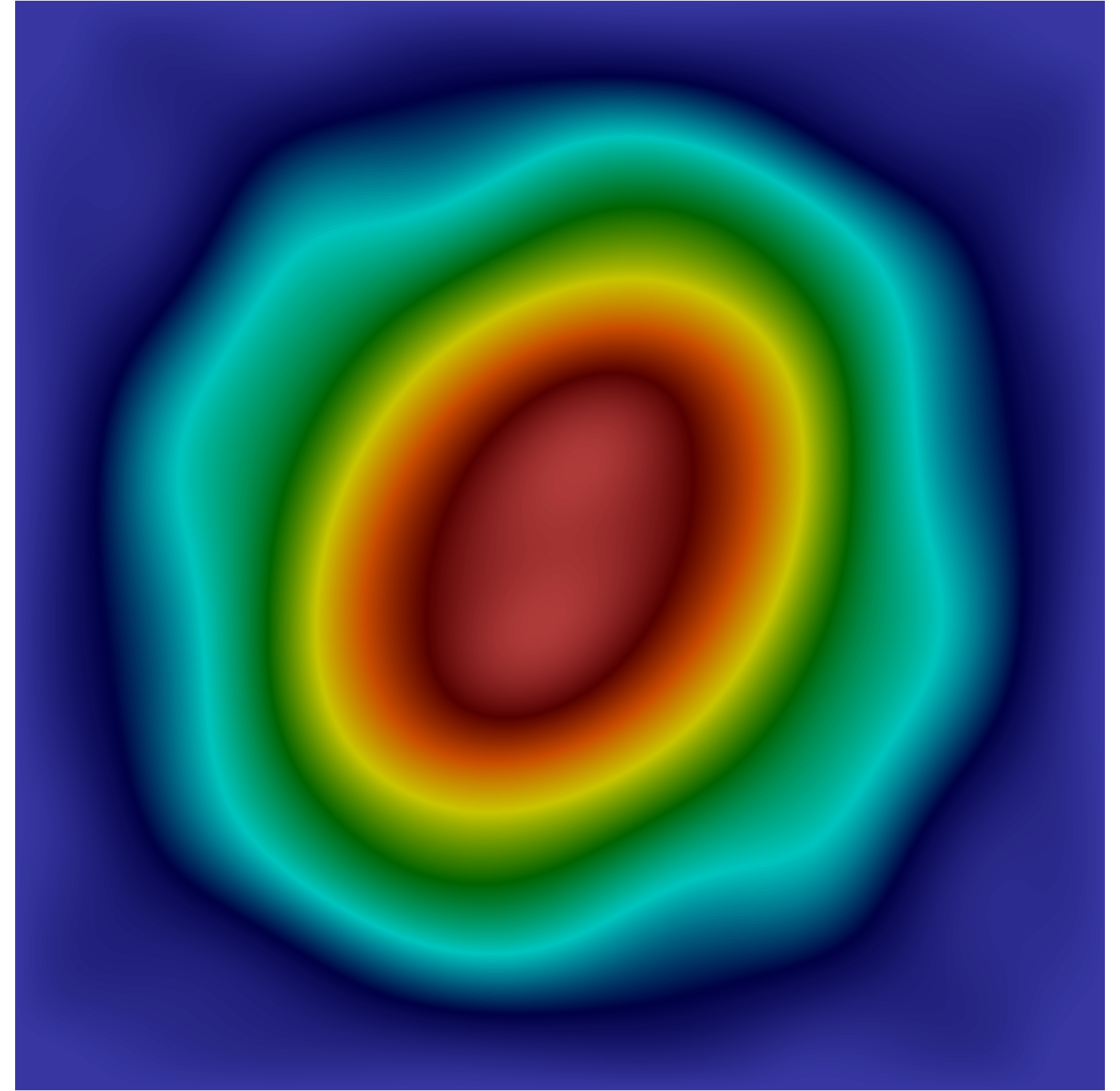}
        \put(-10,50){$\psi$}
        \put(30,100){$\gamma = 0.06$}
      \end{overpic}
 \begin{overpic}[width=0.21\textwidth]{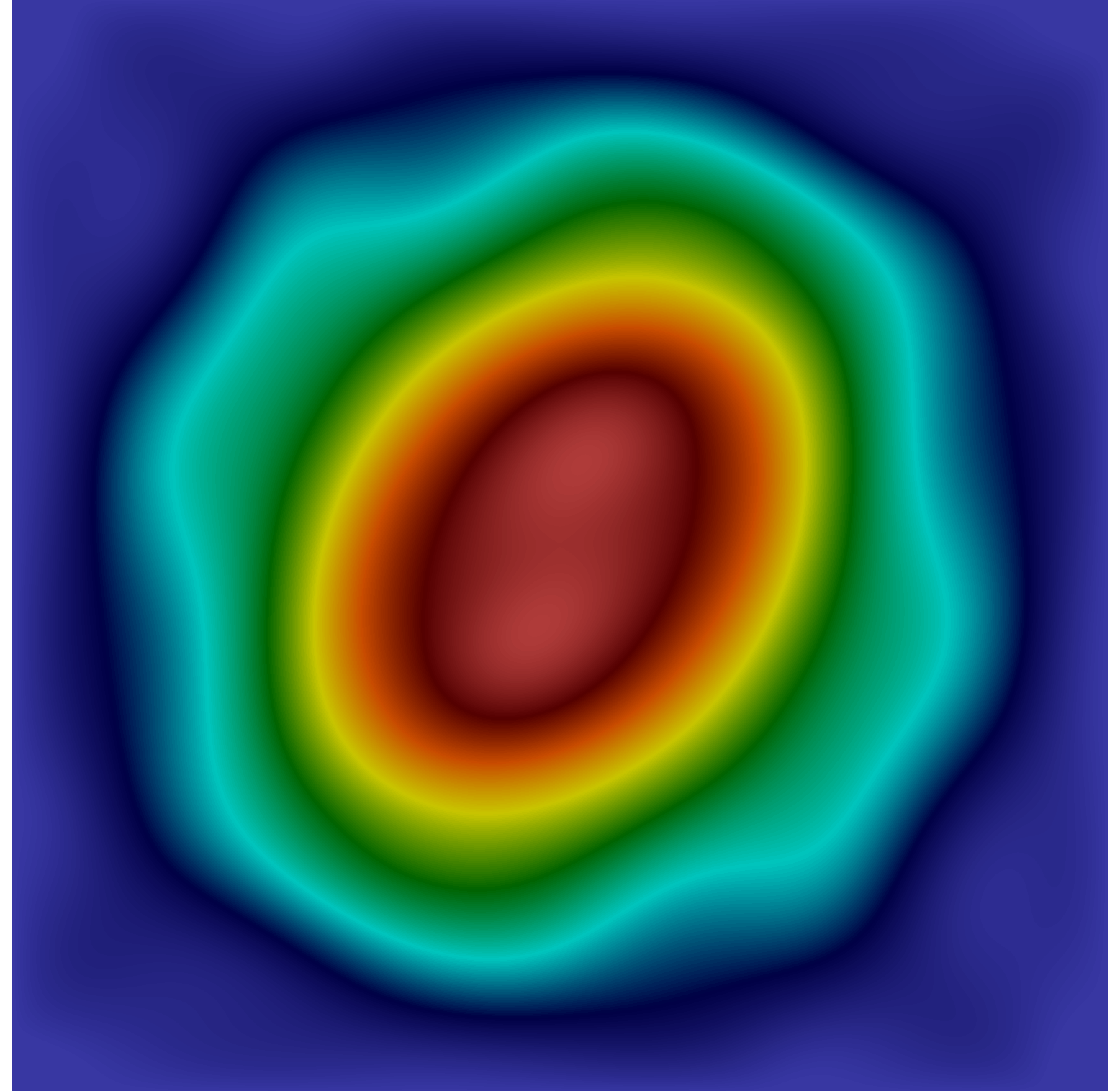}
      \put(30,100){$\gamma = 0.07$}
      \end{overpic}
 \begin{overpic}[width=0.21\textwidth]{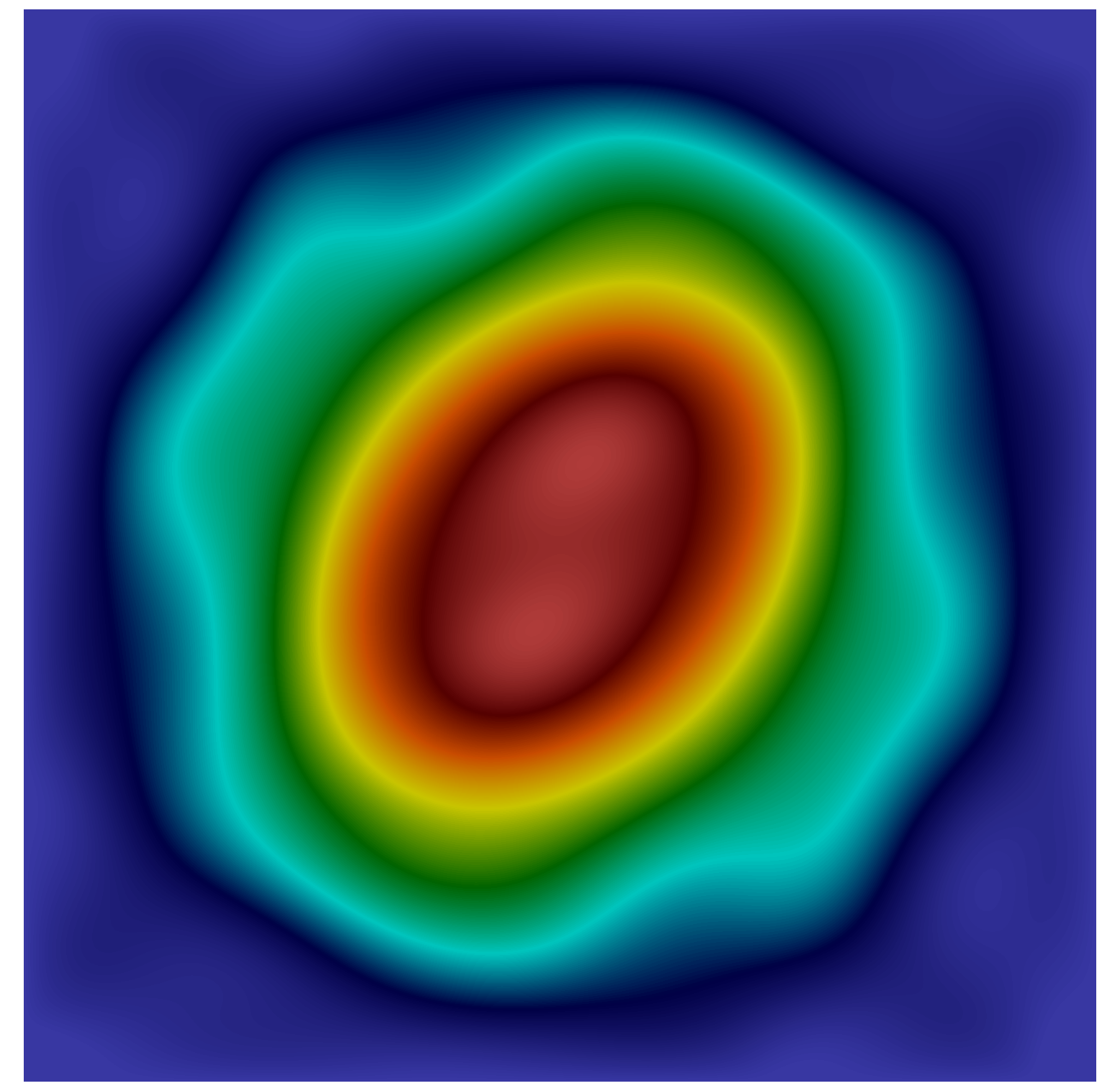}
        \put(30,100){$\gamma = 0.08$}
      \end{overpic}
      \begin{overpic}[width=0.21\textwidth]{img/psiRe_800.png}
       \put(30,100){$\gamma = 0.09$}
      \end{overpic}
       \begin{overpic}[width=0.073\textwidth]{img/scale_psi_Re.png}
      \end{overpic}\\
 \begin{overpic}[width=0.21\textwidth]{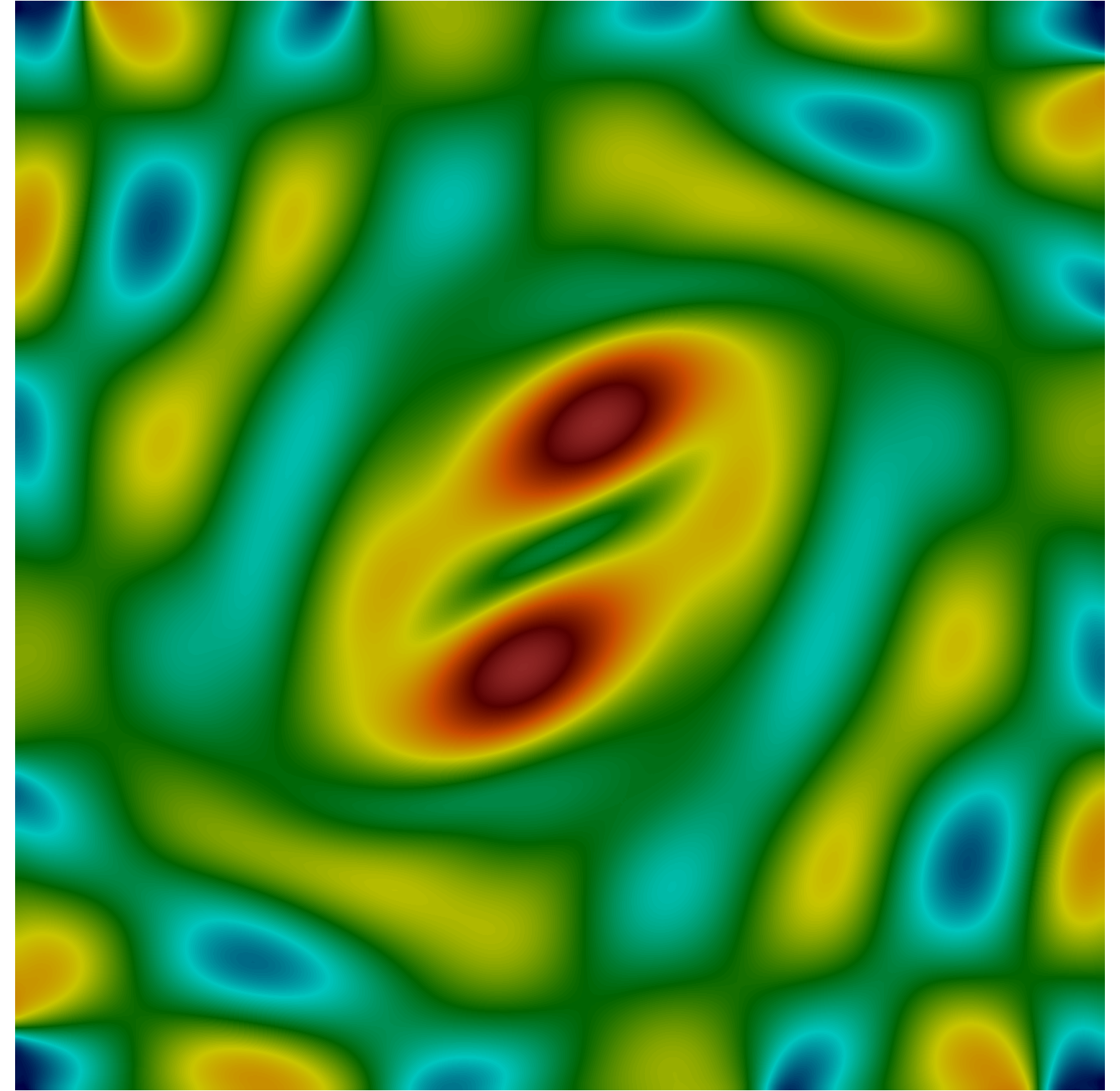}
        \put(-10,50){$\omega$}
      \end{overpic}
 \begin{overpic}[width=0.21\textwidth]{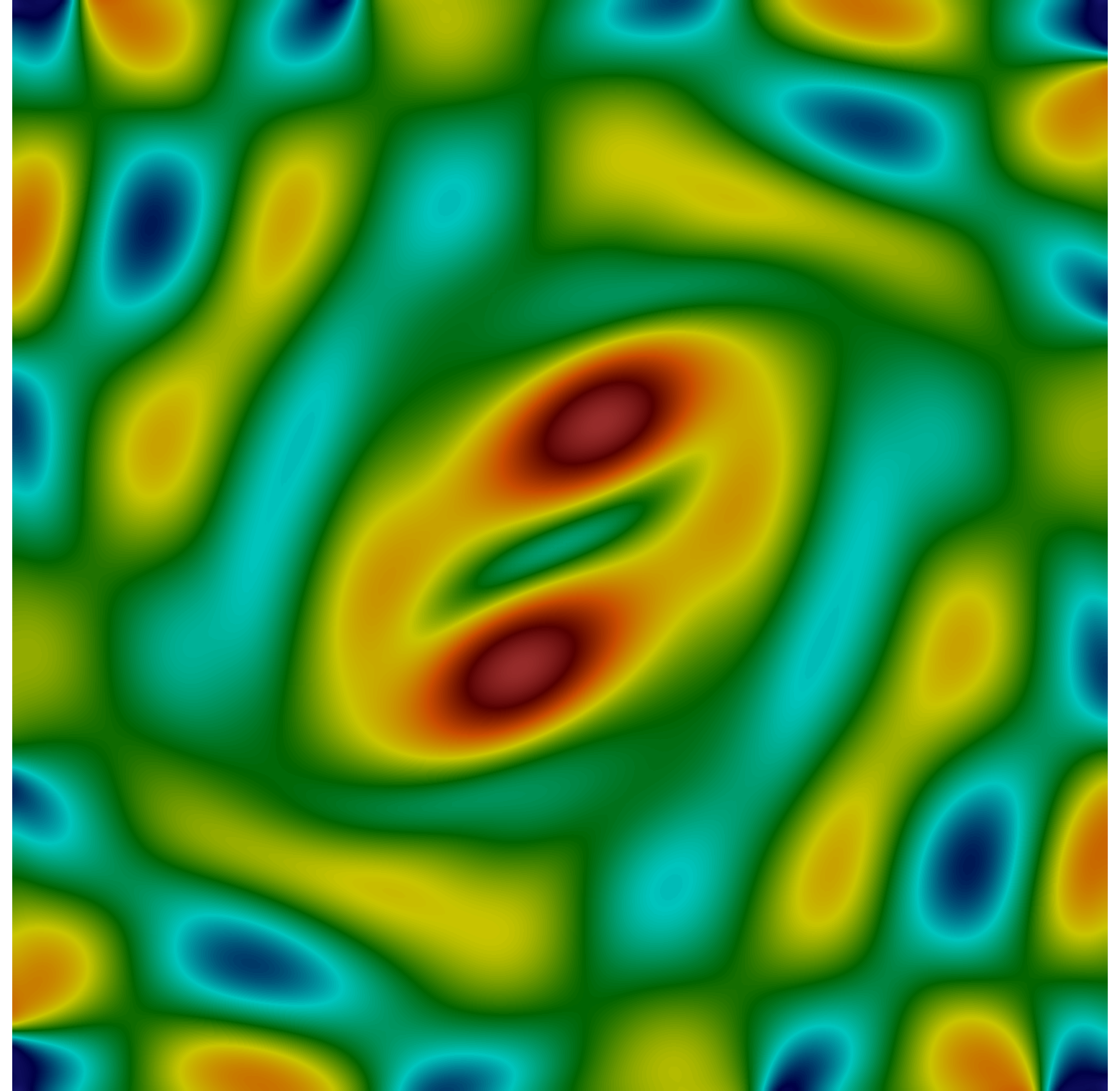}
      \end{overpic}
 \begin{overpic}[width=0.21\textwidth]{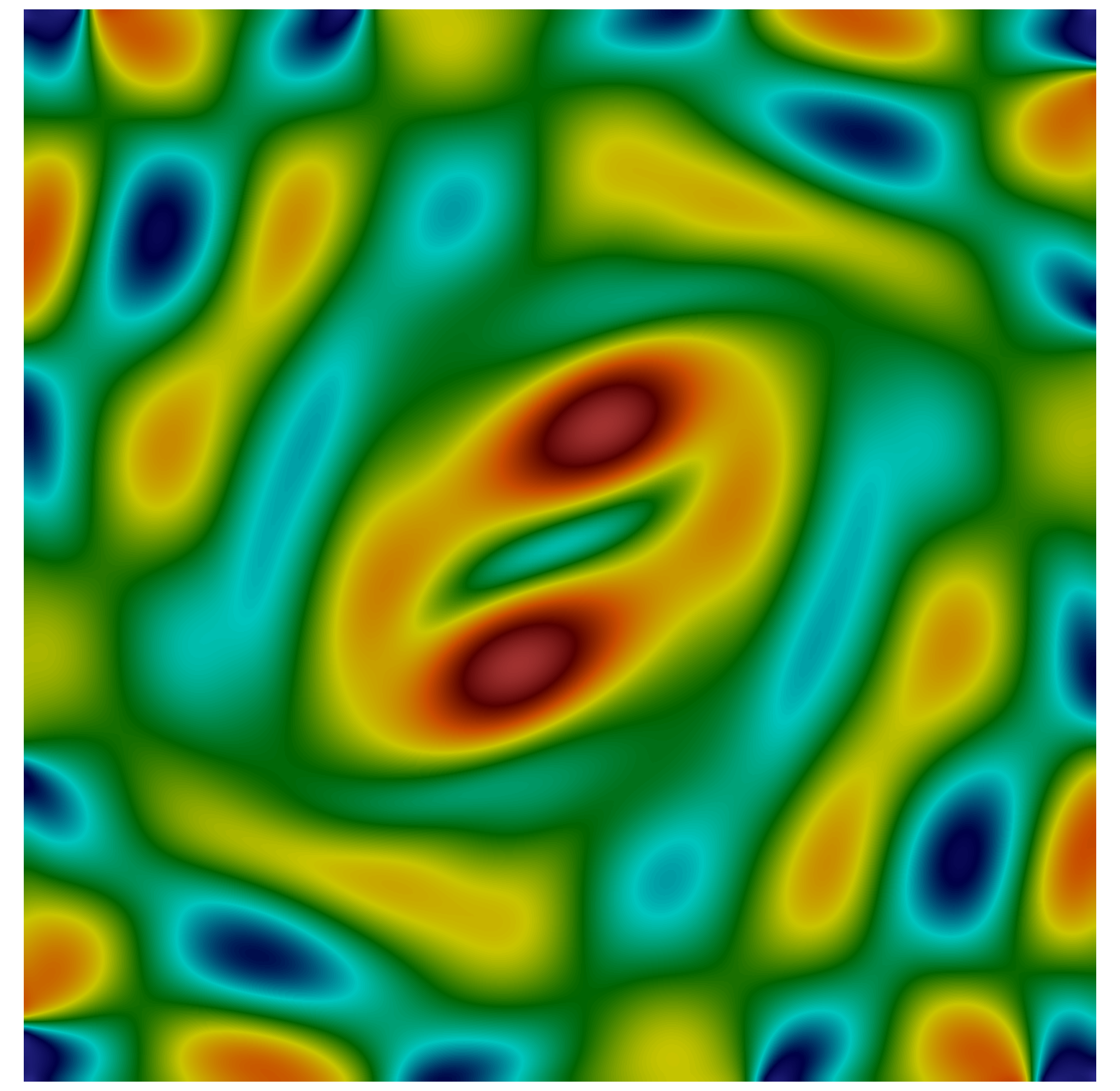}
      \end{overpic}
 \begin{overpic}[width=0.21\textwidth]{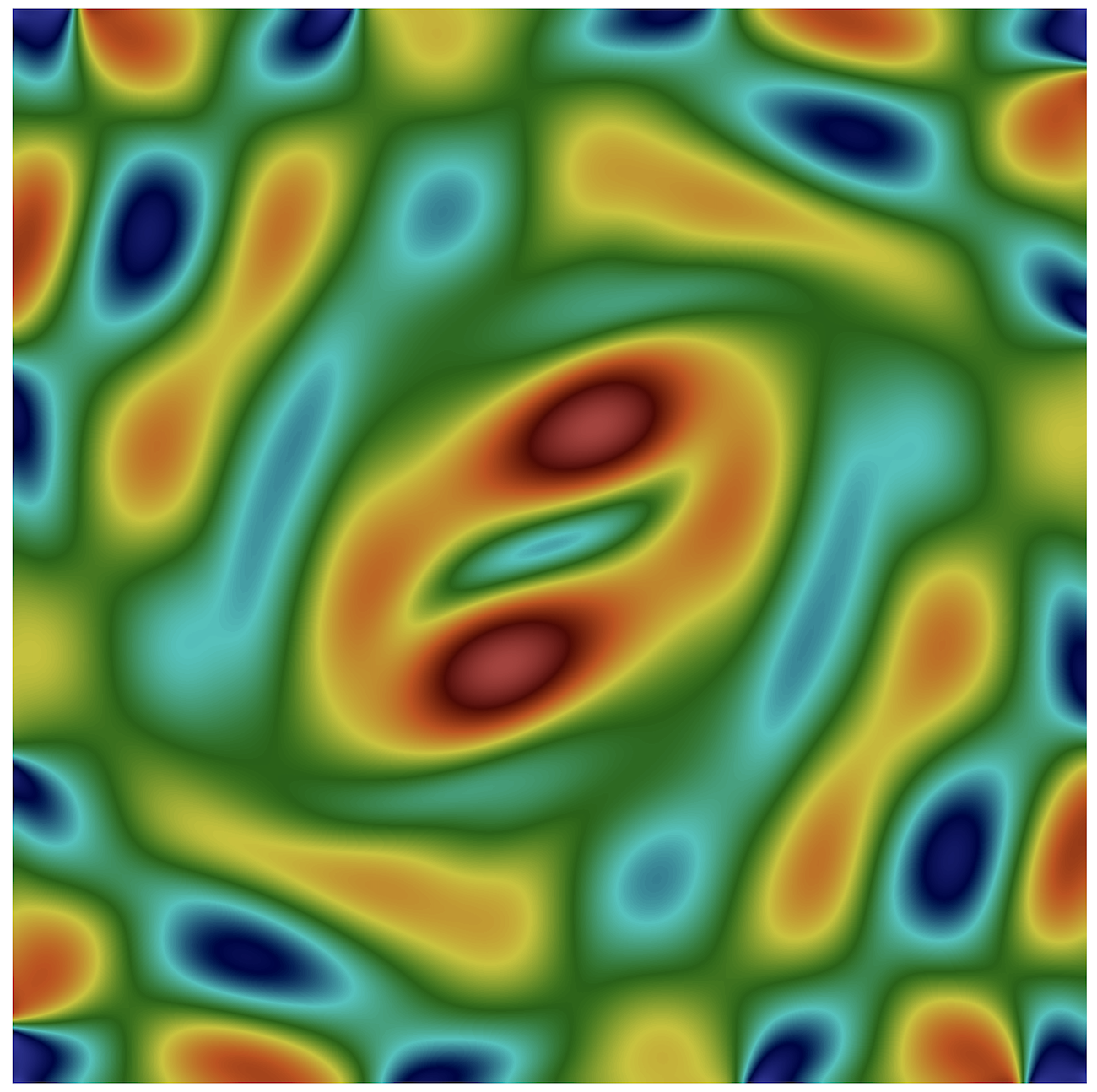}
      \end{overpic}
      \begin{overpic}[width=0.073\textwidth]{img/scale_zeta_Re.png}
      \end{overpic}\\
\caption{ROM validation - $\gamma$ parameterization: 
stream function $\psi$ (first row) and vorticity $\omega$  (second row) computed by the FOM for $Re = 800$ and $\gamma = 0.06$ (first column), $\gamma = 0.07$ (second column), $\gamma = 0.08$ (third column), and $\gamma = 0.09$ (fourth column) at time $t = 10$.}
\label{fig:param_FOM_gamma}
\end{figure}

Next, we consider $\gamma$ as a variable parameter and fix $Re$ to 800. 
Similarly to the $Re$ parameterization,
the training of the ROM is carried out by a uniform sample distribution in the range $\gamma \in [0.06, 0.09]$ consisting of 4 sampling points: 0.06, 0.07, 0.08 and 0.09. For each value of $\gamma$ inside the training set, a simulation is run over time interval $(0, 10]$. Fig. \ref{fig:param_FOM_gamma} shows the stream function and vorticity fields computed by the FOM at $t = 10$ for the four sampling values of $\gamma$. We see
that both stream function and vorticity 
do not significantly vary as $\gamma$ changes value. %On the other hand, for what concerns the vorticity, we observe that in t$\gamma = 0.06$ is not so different from the one at $\gamma = 0.07, 0.08,$ and $0.09$.

\begin{figure}
\centering
 \begin{overpic}[width=0.45\textwidth]{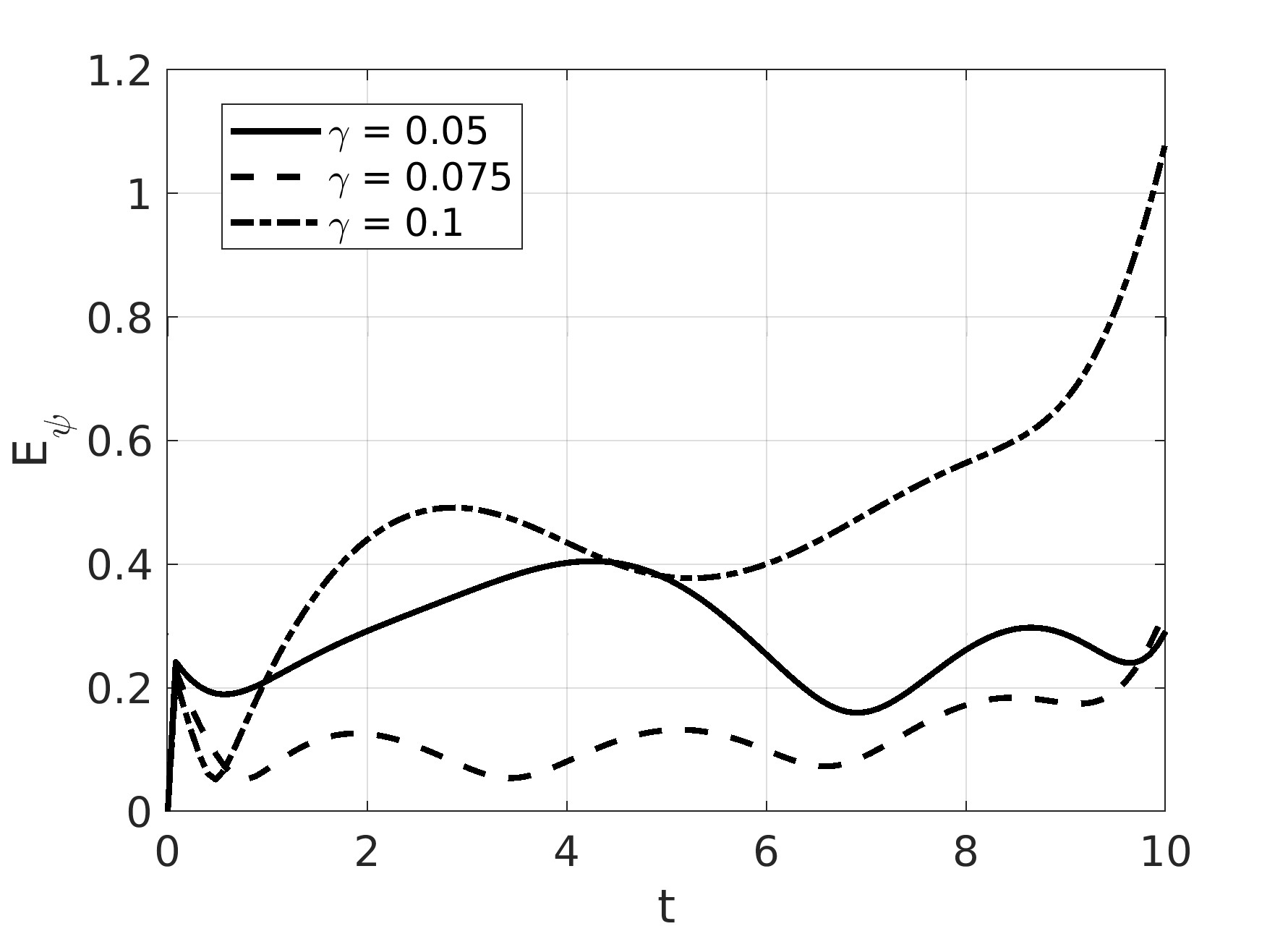}
        %\put(35,18){FOM}
        %\put(-8,7){$\u$}
      \end{overpic}
 \begin{overpic}[width=0.45\textwidth]{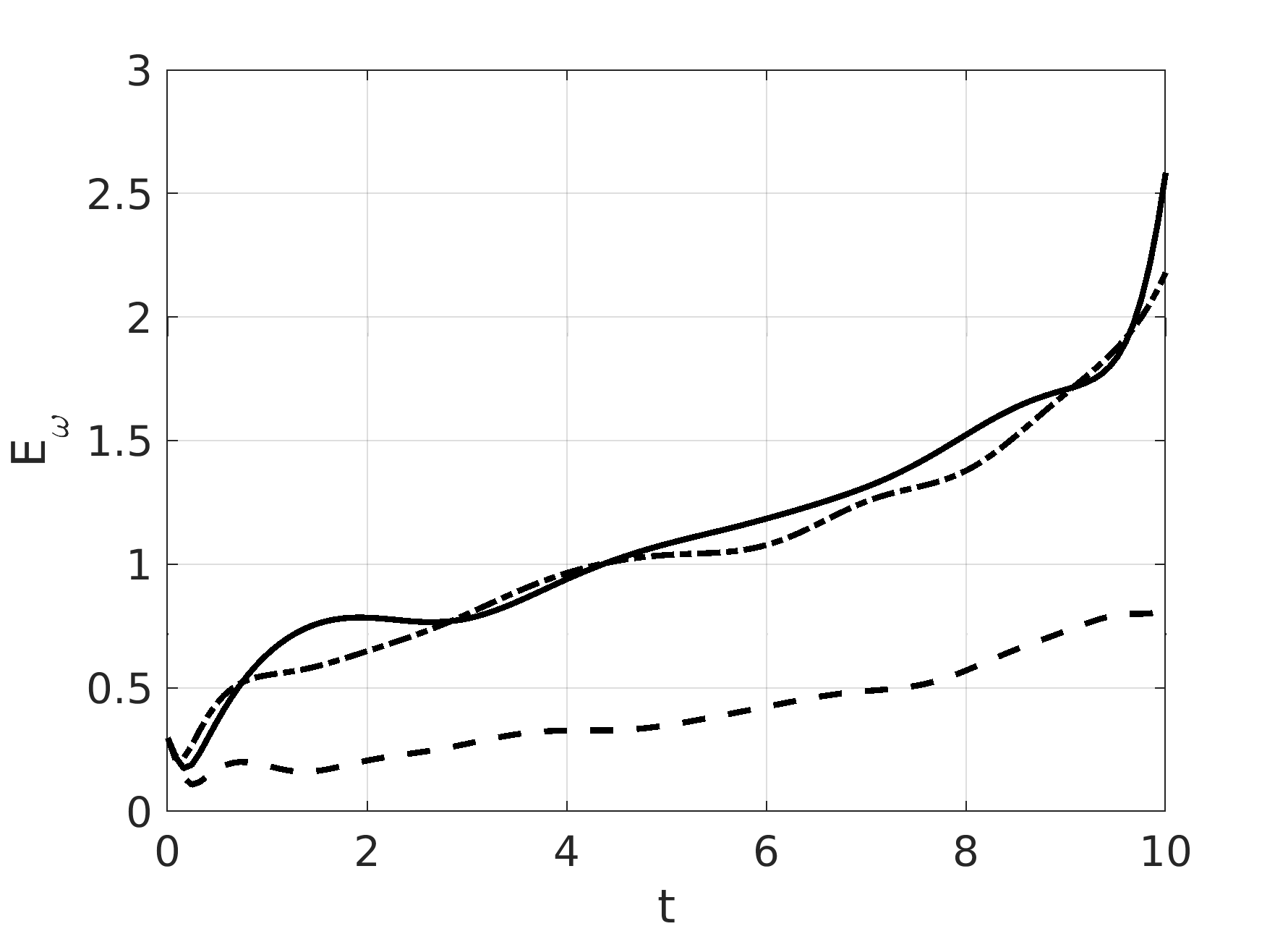}
        %\put(35,18){ROM}
      \end{overpic}
      %\vskip .2cm
       %\begin{overpic}[width=0.3\textwidth]{img/U_diff_1s.png}
        %\put(40,40){ROM}
      %\end{overpic}\\
\caption{ROM validation - $\gamma$ parameterization: time history of error \eqref{eq:error1} for stream function $\psi$ (left) and vorticity $\omega$ (right) for the three different test values $\gamma = 0.05$, $\gamma = 0.075$ and $\gamma = 0.1$ at $Re = 800$.}
\label{fig:errors_gamma}
\end{figure}

Analogously to what we have done for the previous parametric test case, the snapshots are collected every 0.08 s for a total of $4 \times 125 = 500$ snapshots. We set the threshold for the selection of the eigenvalues to $1e-5$, which results in 6 modes for $\psi$ and 12 modes for $\omega$. 
Once again, 
we take three different test values to evaluate the performance of the parametrized ROM: $\gamma = 0.075$ (in the range under consideration but not in the training set) and 
$\gamma = 0.05, 0.1$ (outside the range under consideration).
Fig.~\ref{fig:errors_gamma} shows error \eqref{eq:error1} for the stream function and vorticity over time for these three values of $\gamma$. For the interpolatory test value ($\gamma = 0.075$), the error for $\psi$ is lower than 0.4\% and the error for $\omega$ is lower than 1\% over the entire time interval of interest. As for the extrapolatory test values ($\gamma = 0.05, 0.1$), the error for $\psi$ is lower than 1.2\% and the error for $\omega$ is lower than 2.5\%. So, unlike the $Re$ parametric case, the errors obtained for all test values are comparable. This is obviously due to the fact that there is no abrupt change in the solution as $\gamma$ varies in $[0.06, 0.09]$
(see Fig.~\ref{fig:param_FOM_gamma}) and the POD space seems to include enough information for a very good reconstruction of the flow field also at values of $\gamma$ right outside the training set.
% $Re = 1000$ are comparable with those related to the interpolary test value $Re = 500$ whilst for $Re = 100$
%the ROM performs significantly worse (although not so bad in absolute) for $\omega$. 
%for a correct reconstruction of the flow field at $Re = 100$ it could be necessary to consider in the training set lower values of $Re$.

\begin{figure}
\centering
 \begin{overpic}[width=0.175\textwidth]{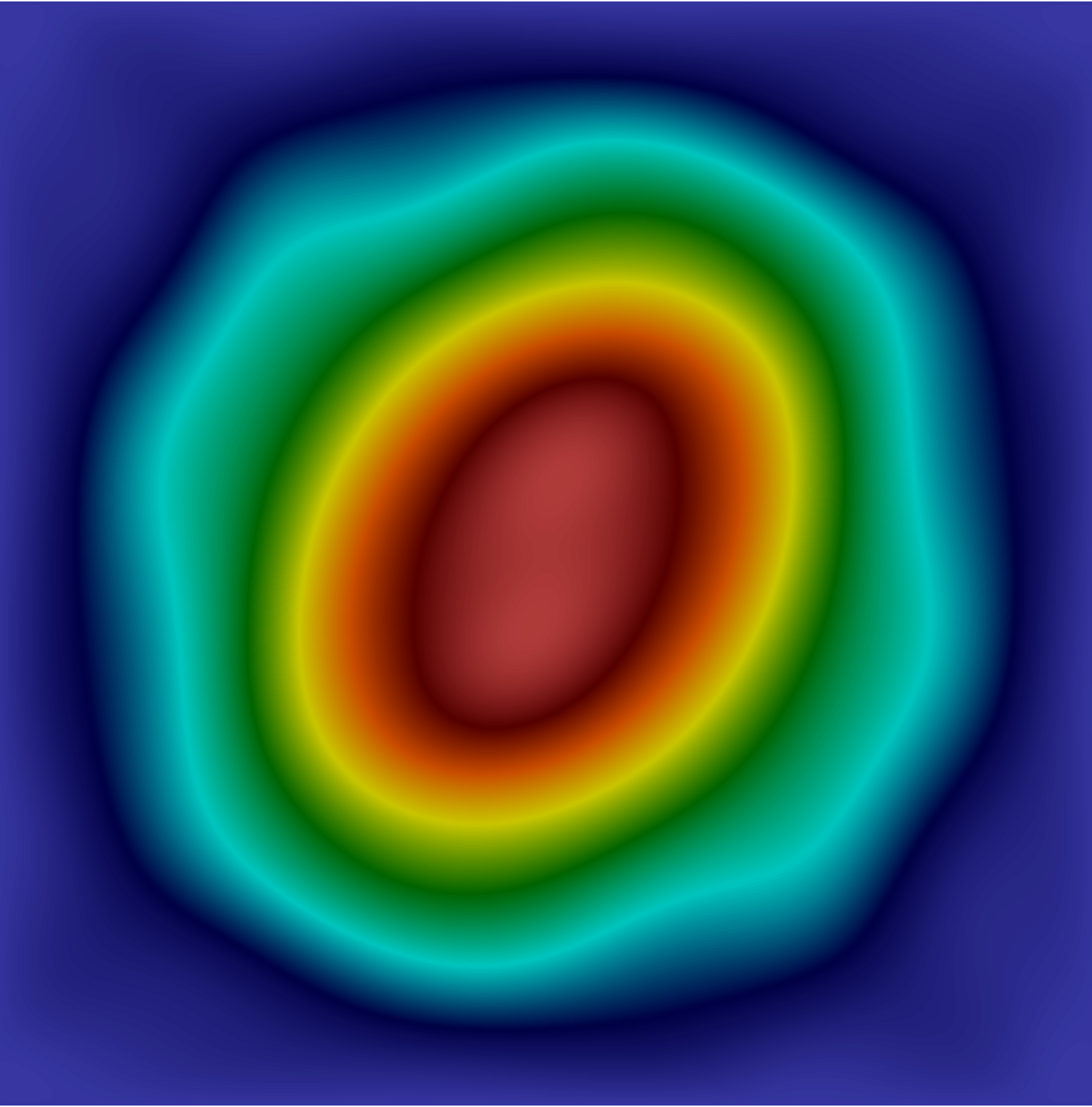}
        \put(22,103){$\gamma = 0.05$}
        \put(-42,48){FOM}
      \end{overpic}
 \begin{overpic}[width=0.075\textwidth]{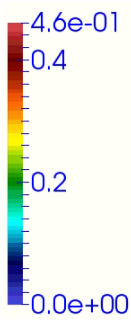}
      \end{overpic}
       \begin{overpic}[width=0.175\textwidth]{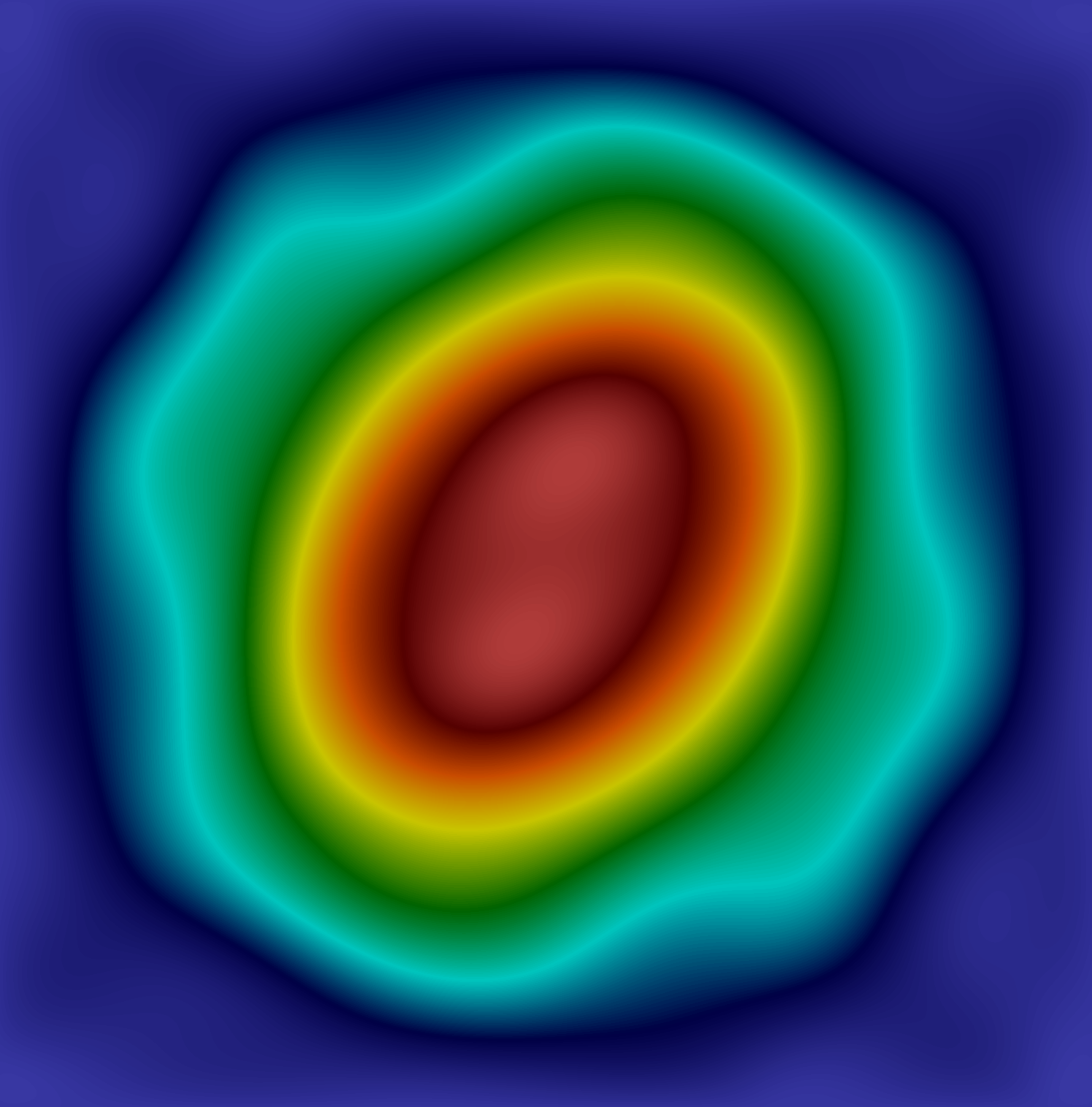}
        \put(20,103){$\gamma = 0.075$}
      \end{overpic}
 \begin{overpic}[width=0.075\textwidth]{img/legenda_800_0_05_psi.png}
      \end{overpic}
   \begin{overpic}[width=0.175\textwidth]{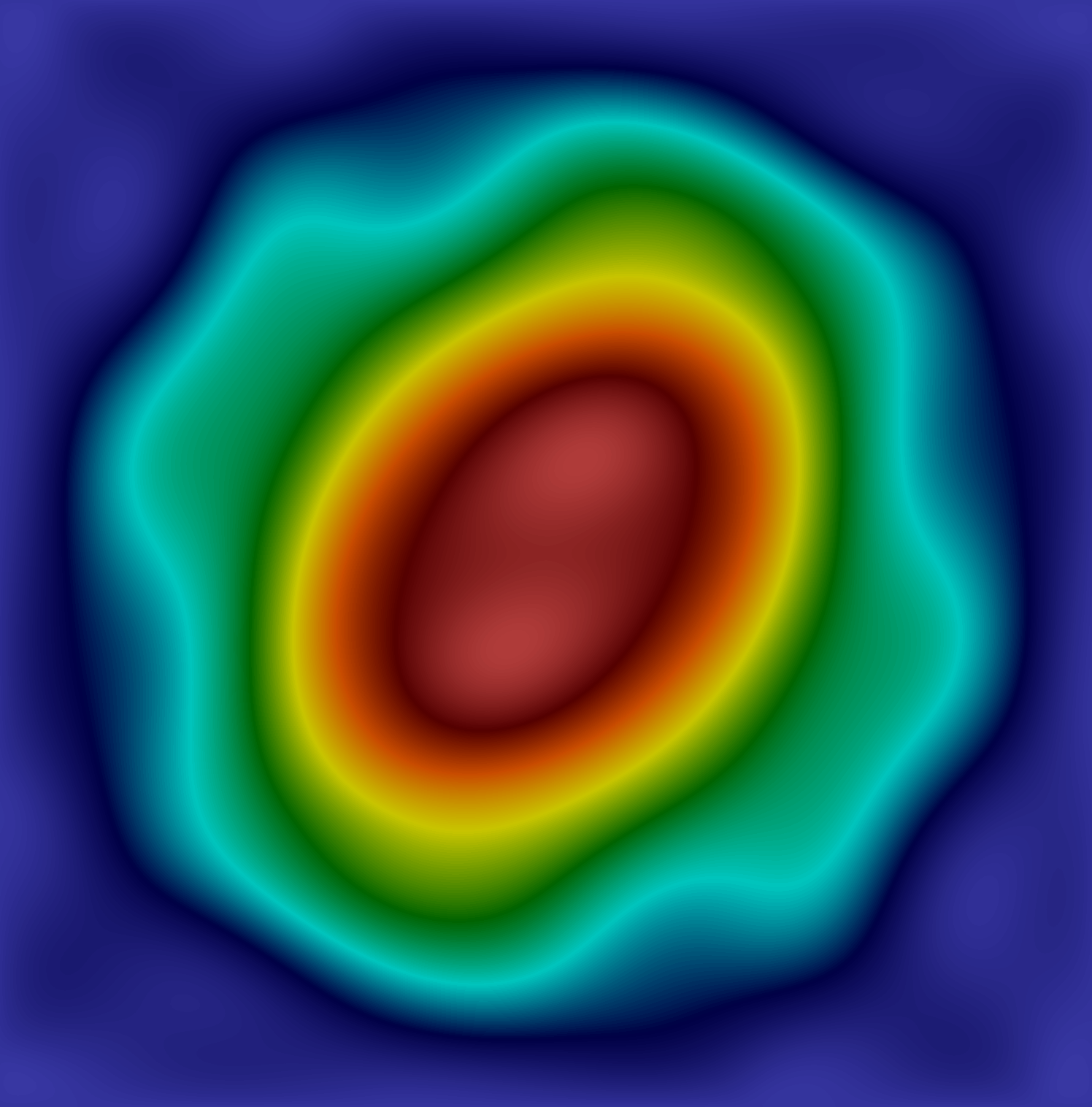}
        \put(24,103){$\gamma = 0.1$}
      \end{overpic}
 \begin{overpic}[width=0.075\textwidth]{img/legenda_800_0_05_psi.png}
      \end{overpic}\\
      \vspace{0.5cm}
 \begin{overpic}[width=0.175\textwidth]{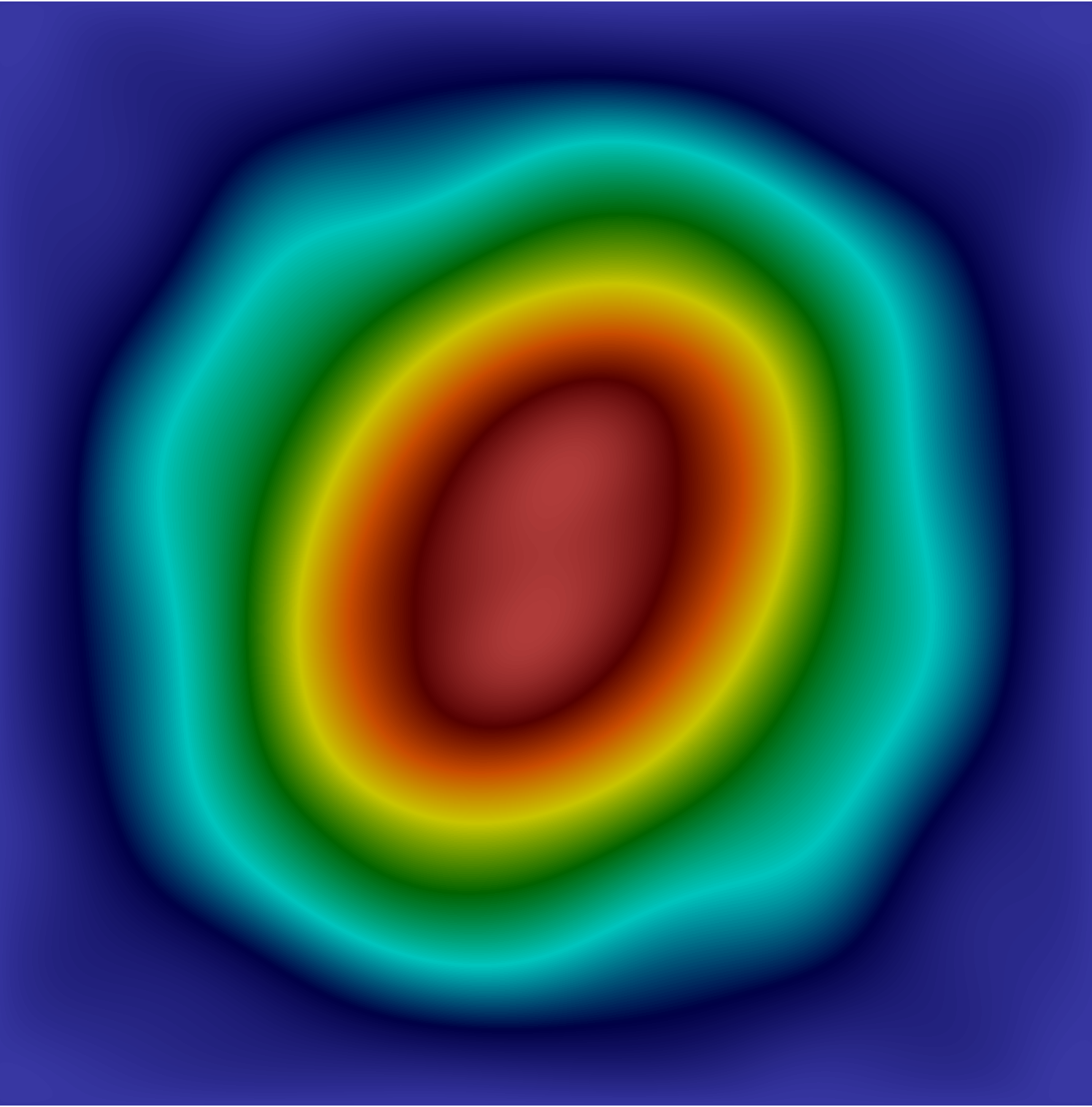}
        \put(-42,48){ROM}
      \end{overpic}
 \begin{overpic}[width=0.075\textwidth]{img/legenda_800_0_05_psi.png}
      \end{overpic}
   \begin{overpic}[width=0.175\textwidth]{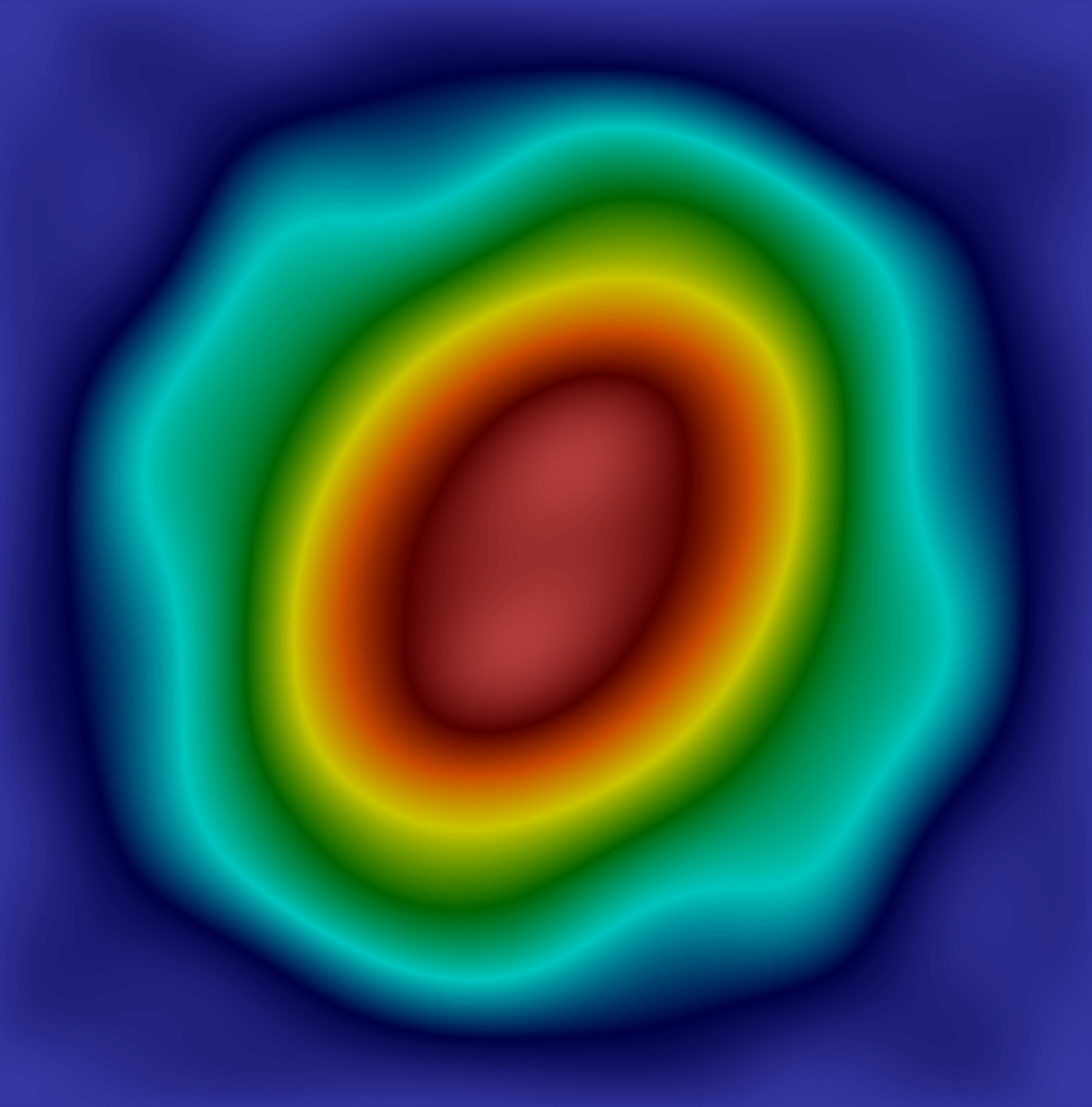}
      \end{overpic}
 \begin{overpic}[width=0.075\textwidth]{img/legenda_800_0_05_psi.png}
      \end{overpic}
      \begin{overpic}[width=0.175\textwidth]{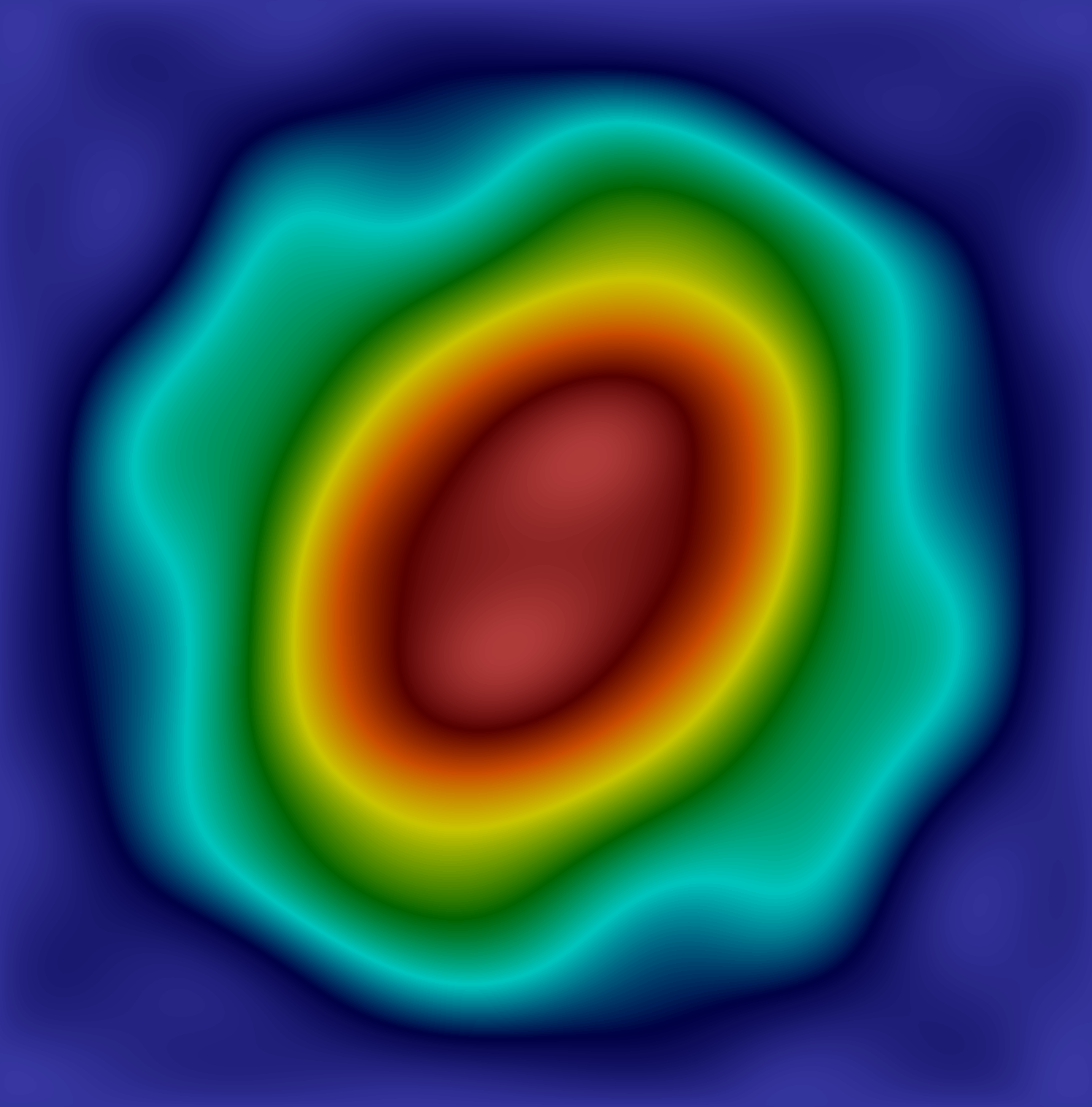}
      \end{overpic}
 \begin{overpic}[width=0.075\textwidth]{img/legenda_800_0_05_psi.png}
      \end{overpic}\\
            \vspace{0.5cm}
\begin{overpic}[width=0.185\textwidth]{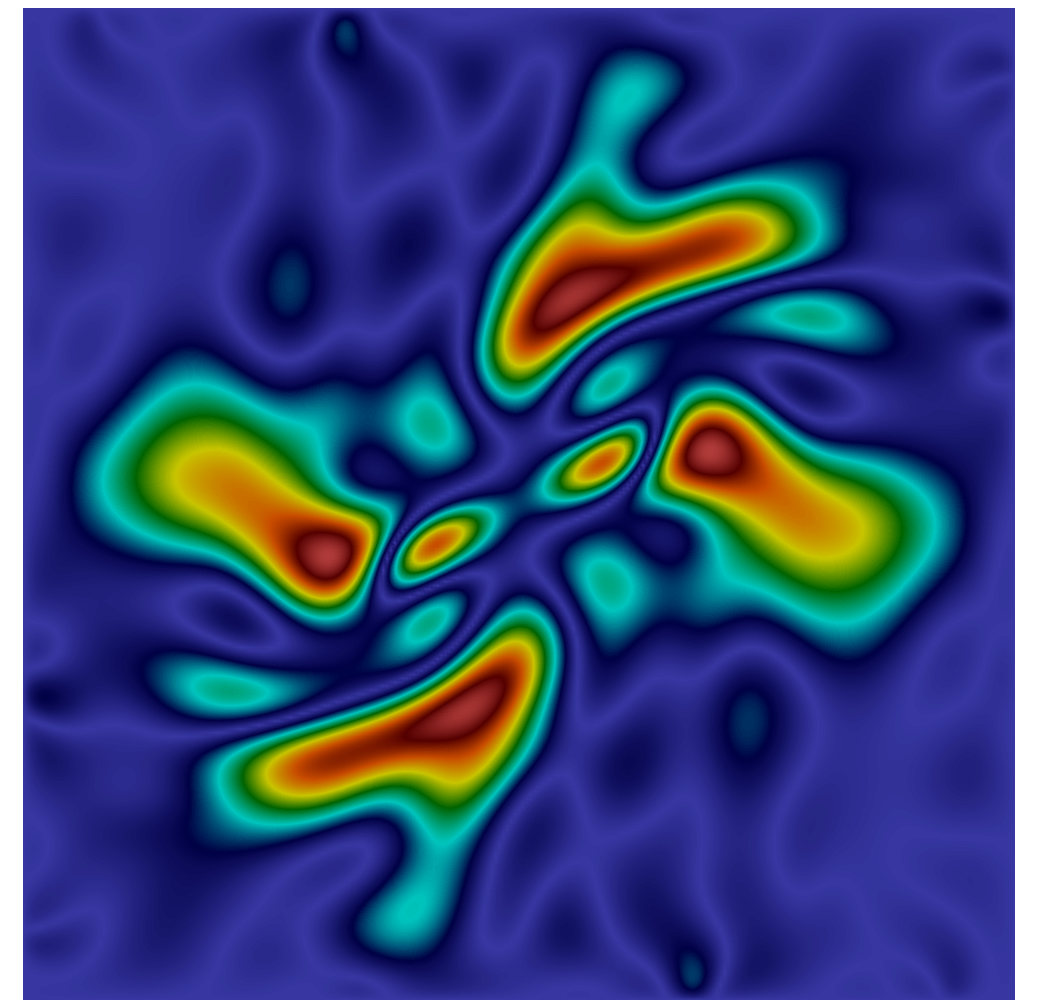}
        \put(-32,48){Diff.}
      \end{overpic}
 \begin{overpic}[width=0.065\textwidth]{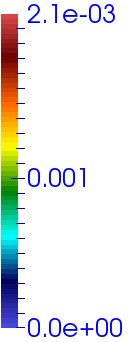}
        %\put(35,18){ROM}
      \end{overpic}
   \begin{overpic}[width=0.175\textwidth]{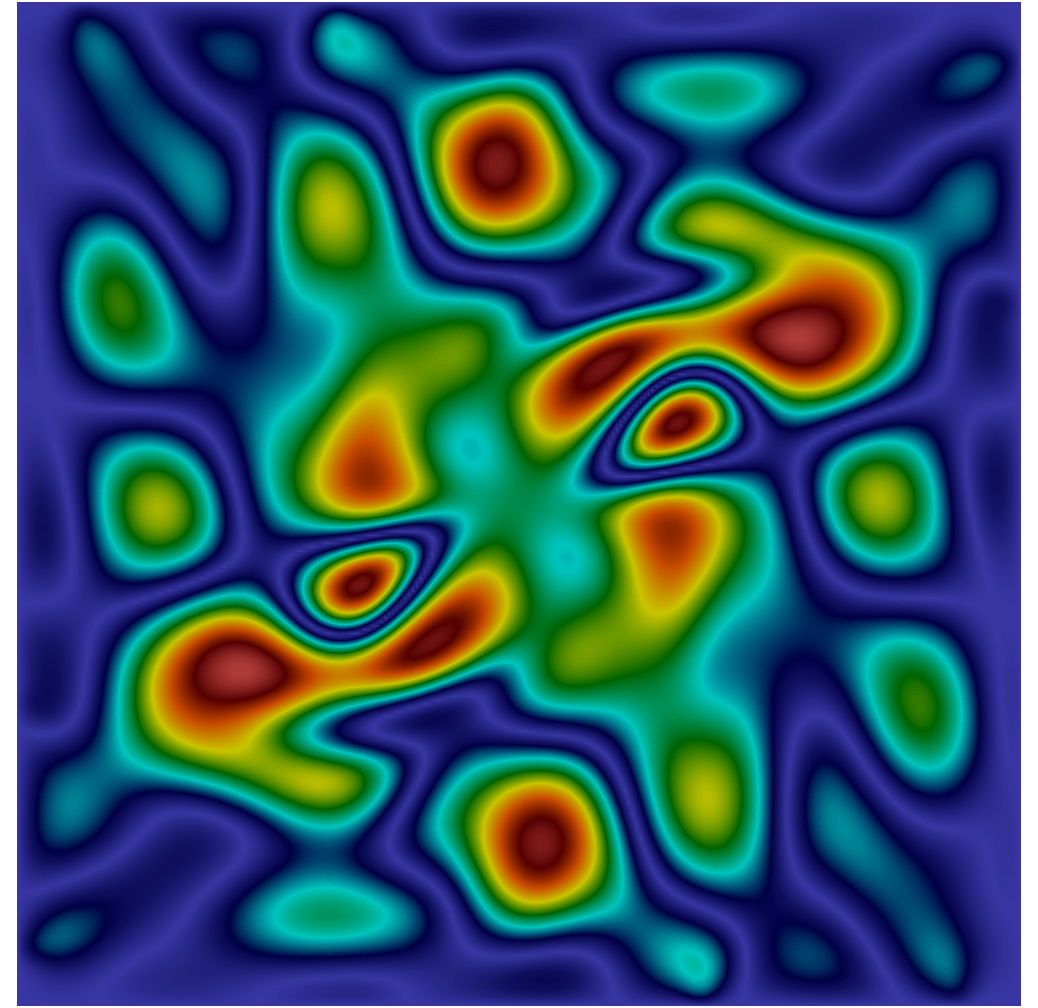}
        %\put(35,18){FOM}
        %\put(-8,7){$\u$}
      \end{overpic}
 \begin{overpic}[width=0.065\textwidth]{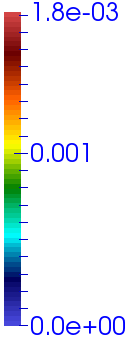}
        %\put(35,18){ROM}
      \end{overpic}
      \begin{overpic}[width=0.175\textwidth]{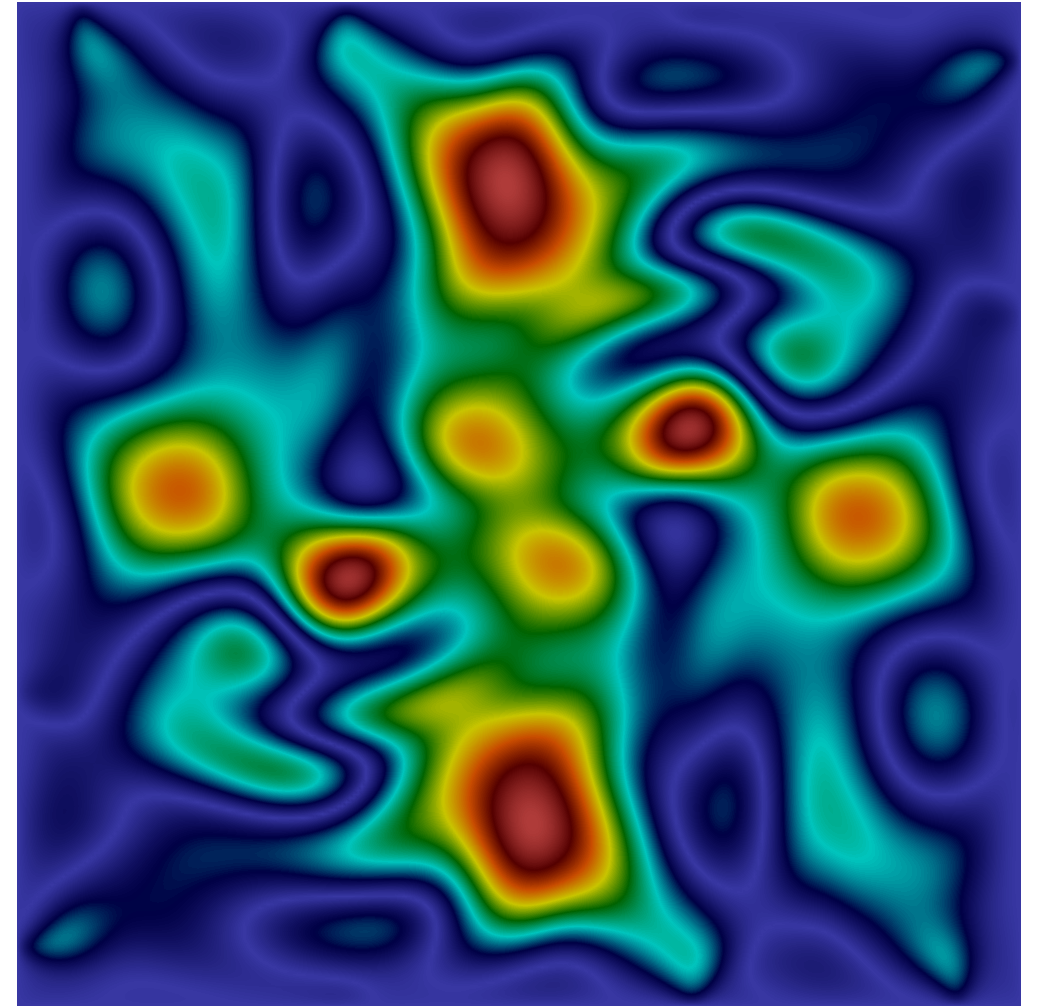}
        %\put(35,18){FOM}
        %\put(-8,7){$\u$}
      \end{overpic}
 \begin{overpic}[width=0.065\textwidth]{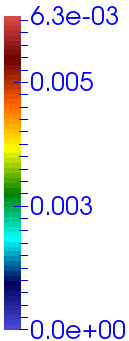}
        %\put(35,18){ROM}
      \end{overpic}\\
      
      %\vskip .2cm
       %\begin{overpic}[width=0.3\textwidth]{img/U_diff_1s.png}
        %\put(40,40){ROM}
      %\end{overpic}\\
\caption{ROM validation - $\gamma$ parameterization:
stream function  $\psi$ computed by the FOM (first row) and the ROM (second row), and difference between the two fields in absolute value (third row) for $\gamma = 0.05$ (first column), $\gamma = 0.075$ (second column) and $\gamma = 0.1$ (third one) at time $t = 10$. Six modes for $\psi$ were considered.}
\label{fig:errors_psi_absolute_gamma}
\end{figure}

\begin{figure}
\centering
 \begin{overpic}[width=0.175\textwidth]{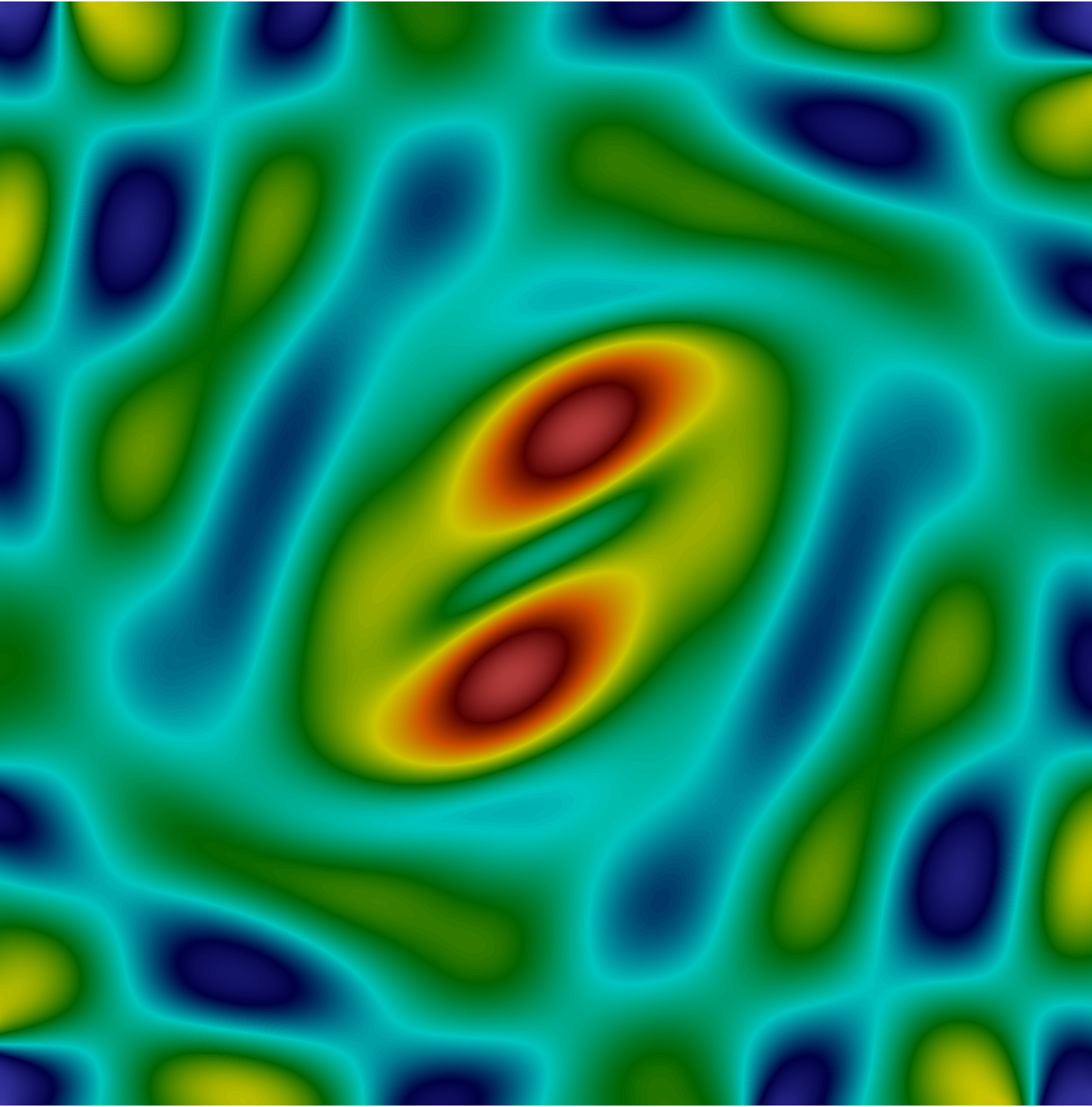}
        \put(22,103){$\gamma = 0.05$}
        \put(-42,48){FOM}
      \end{overpic}
 \begin{overpic}[width=0.065\textwidth]{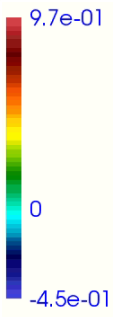}
      \end{overpic}
       \begin{overpic}[width=0.175\textwidth]{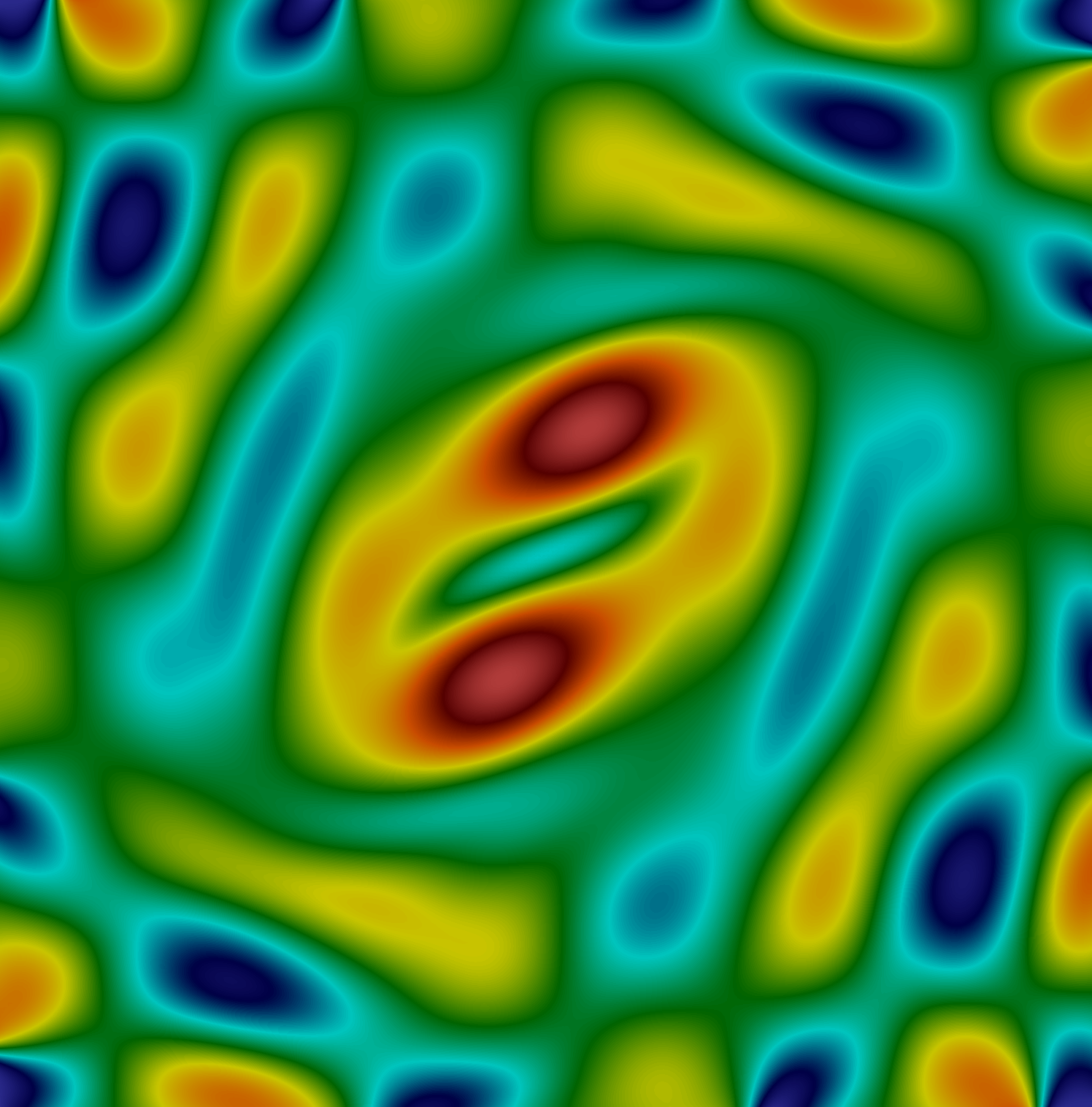}
        \put(20,103){$\gamma = 0.075$}
      \end{overpic}
 \begin{overpic}[width=0.07\textwidth]{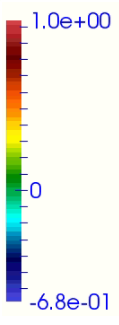}
      \end{overpic}
   \begin{overpic}[width=0.175\textwidth]{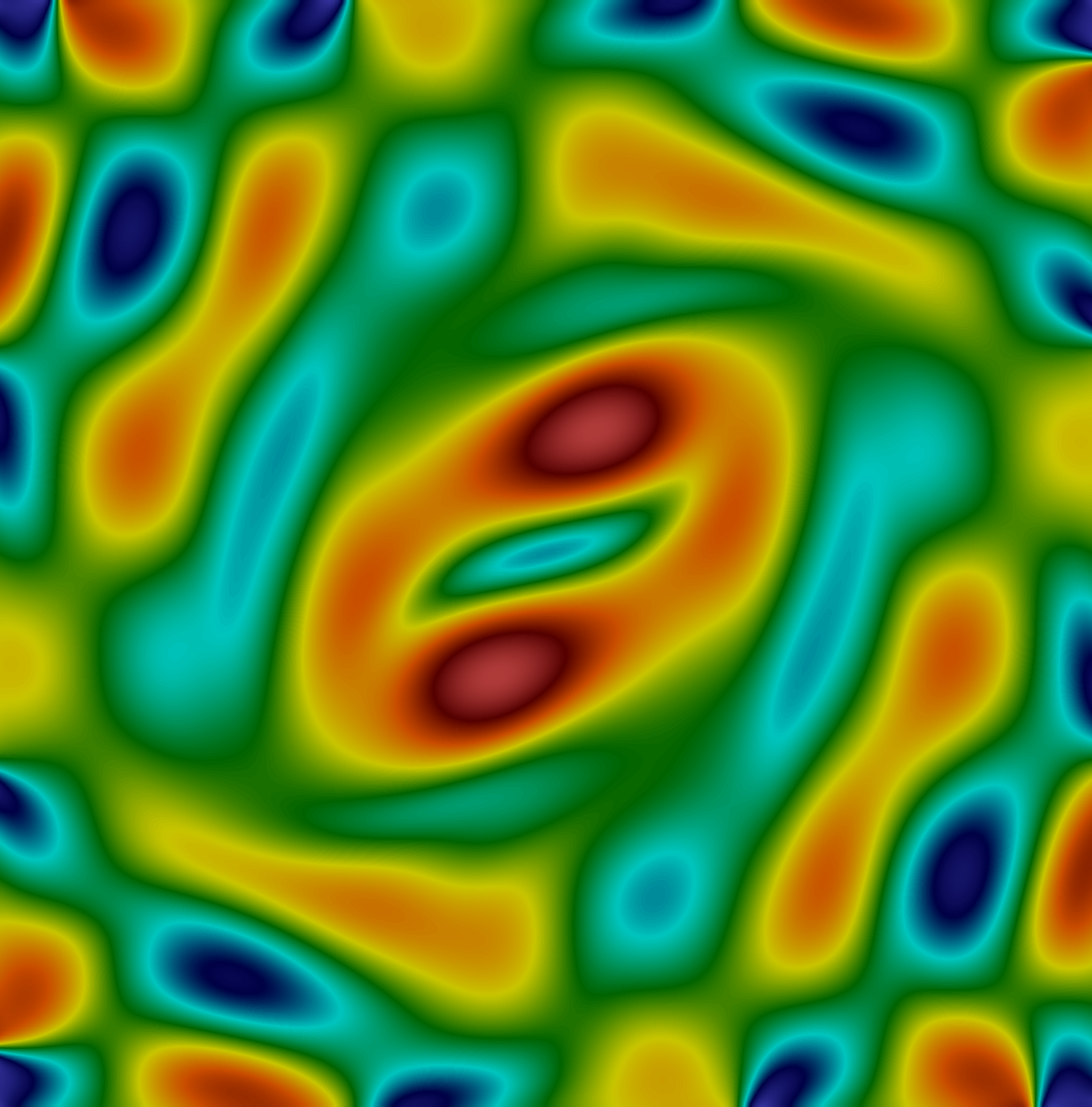}
        \put(20,103){$\gamma = 0.1$}
      \end{overpic}
 \begin{overpic}[width=0.065\textwidth]{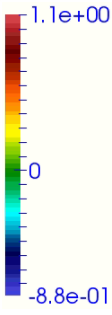}
      \end{overpic}\\
      \vspace{0.5cm}
 \begin{overpic}[width=0.175\textwidth]{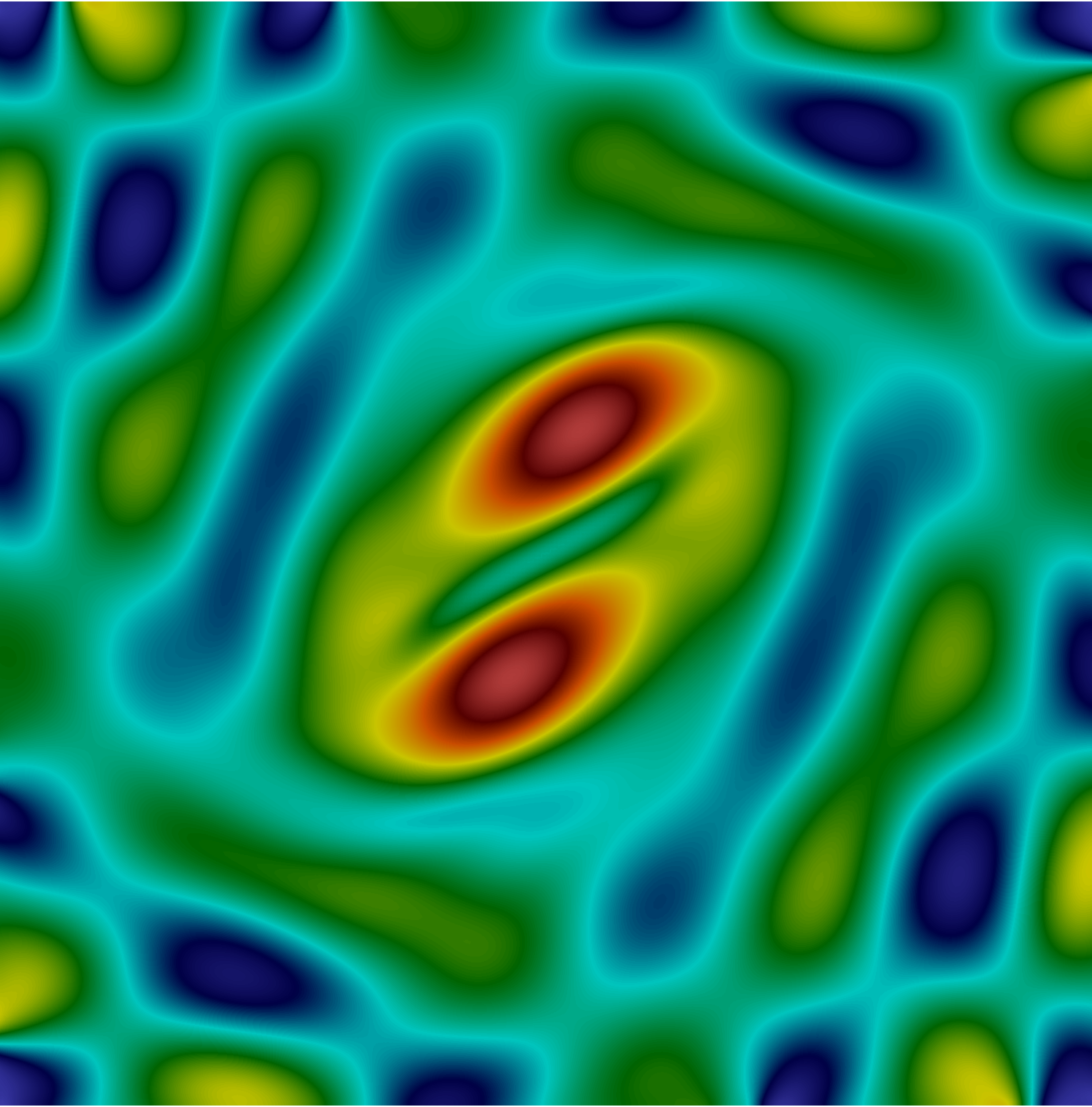}
        \put(-42,48){ROM}
      \end{overpic}
 \begin{overpic}[width=0.065\textwidth]{img/legenda_800_0_05_zeta.png}
      \end{overpic}
   \begin{overpic}[width=0.175\textwidth]{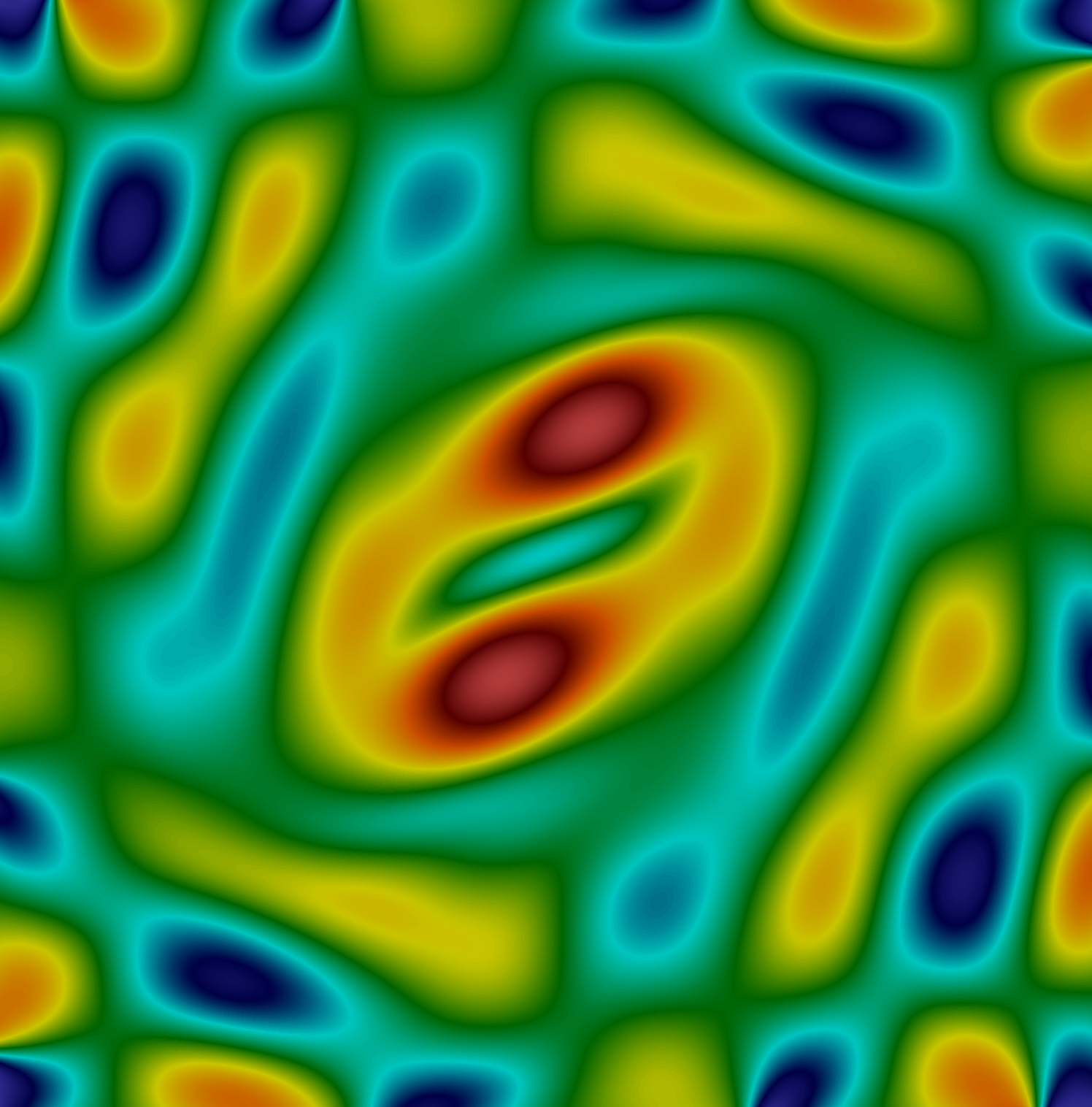}
      \end{overpic}
 \begin{overpic}[width=0.07\textwidth]{img/legenda_800_0_075_zeta.png}
      \end{overpic}
      \begin{overpic}[width=0.175\textwidth]{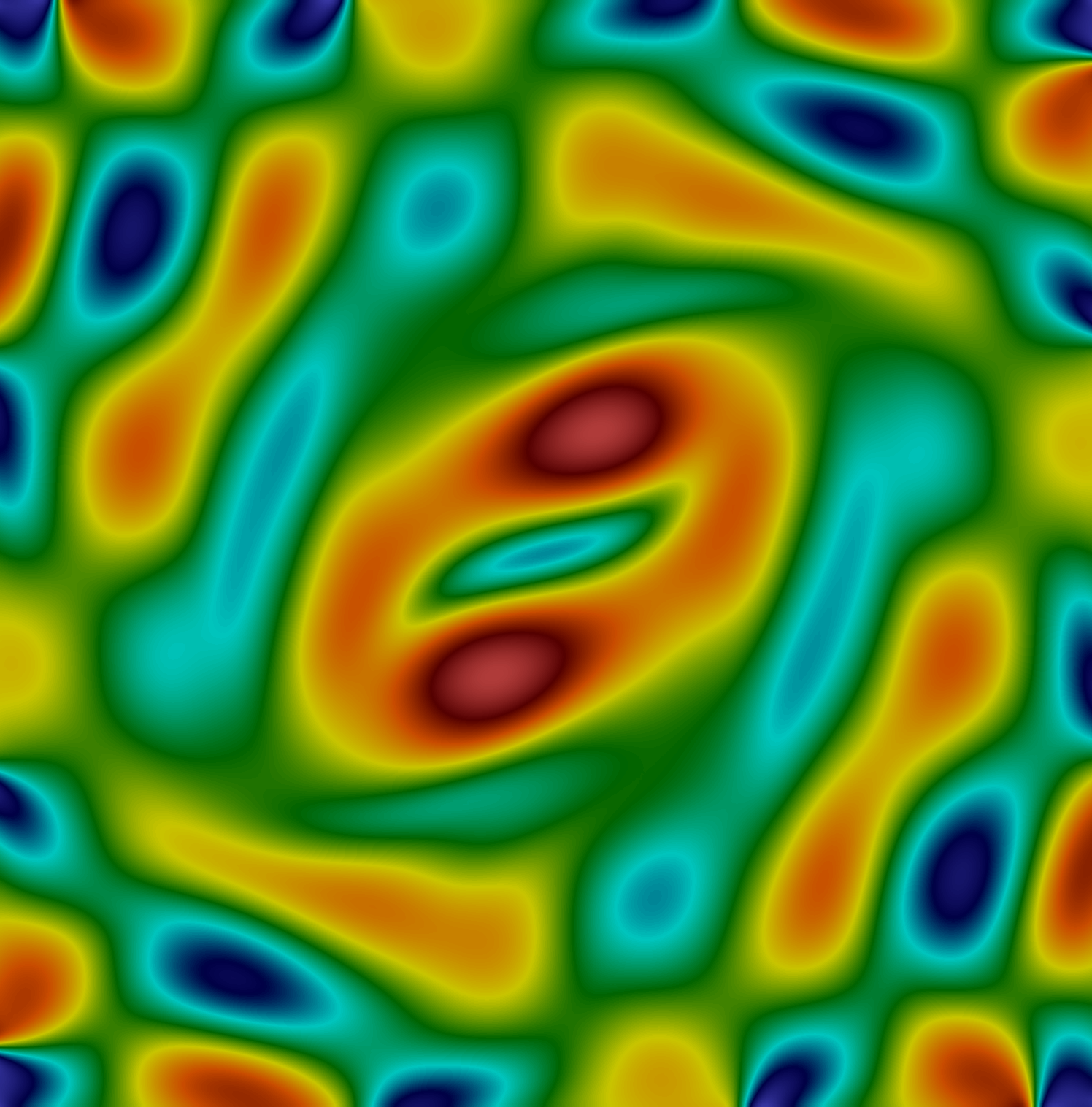}
      \end{overpic}
 \begin{overpic}[width=0.07\textwidth]{img/legenda_800_0_1_zeta.png}
      \end{overpic}\\
            \vspace{0.5cm}
\begin{overpic}[width=0.175\textwidth]{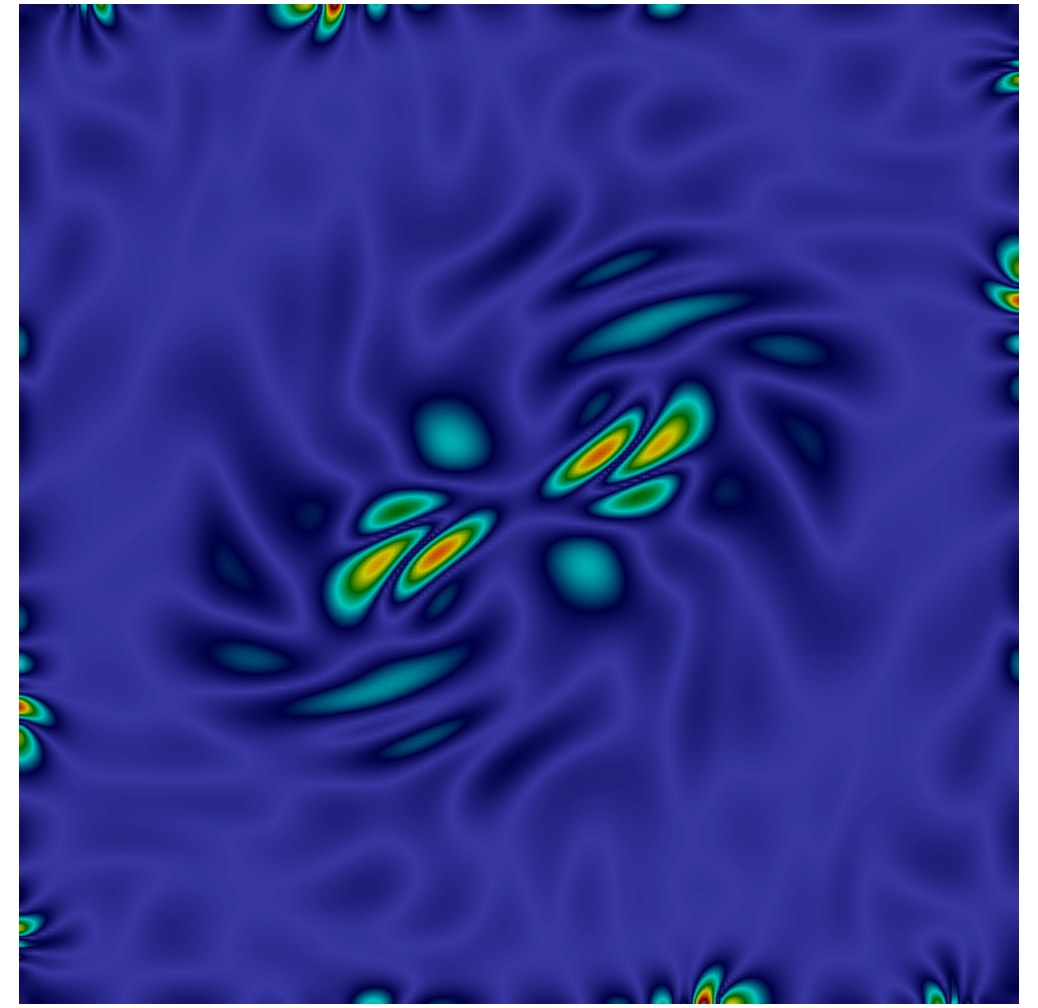}
        \put(-32,48){Diff.}
      \end{overpic}
 \begin{overpic}[width=0.07\textwidth]{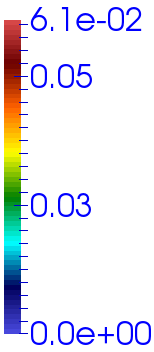}
      \end{overpic}
   \begin{overpic}[width=0.175\textwidth]{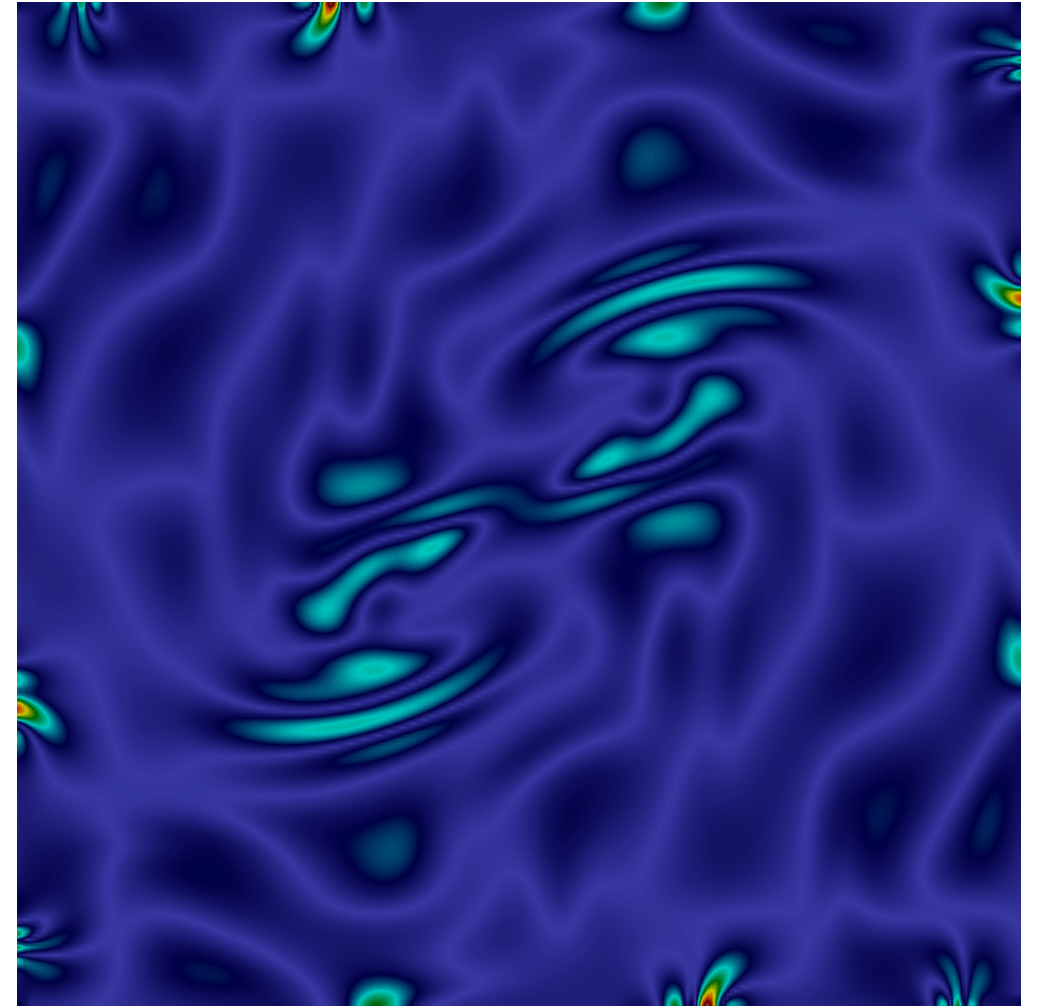}
      \end{overpic}
 \begin{overpic}[width=0.07\textwidth]{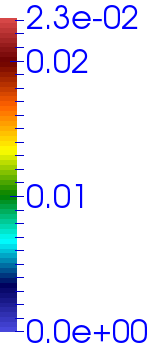}
      \end{overpic}
      \begin{overpic}[width=0.175\textwidth]{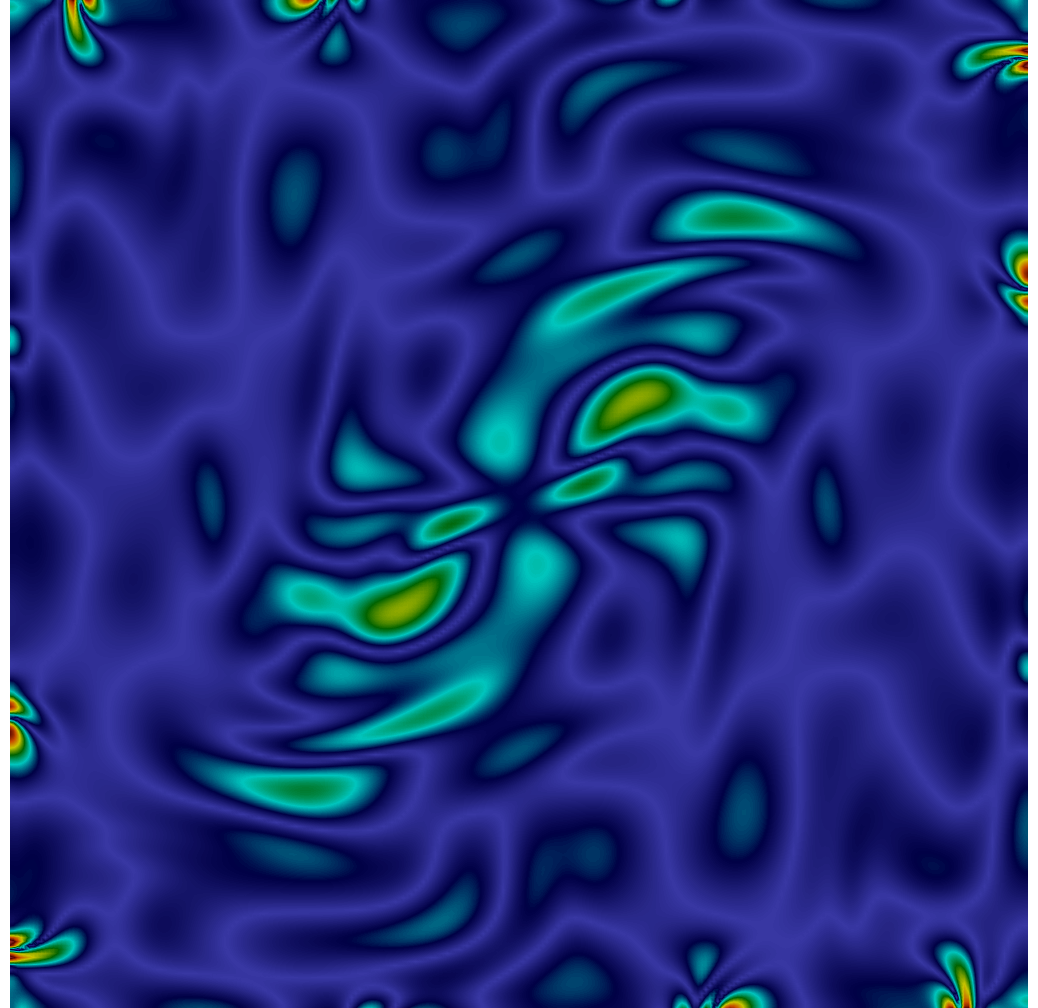}
      \end{overpic}
 \begin{overpic}[width=0.07\textwidth]{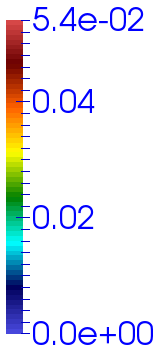}
      \end{overpic}\\
\caption{ROM validation - $\gamma$ parameterization:
vorticity  $\omega$ computed by the FOM (first row) and the ROM (second row), and difference between the two fields in absolute value (third row) for $\gamma = 0.05$ (first column), $\gamma = 0.075$ (second column) and $\gamma = 0.1$ (third one) at time $t = 10$. Twelve modes for $\omega$ were considered.}
\label{fig:errors_zeta_absolute_gamma}
\end{figure}

To conclude, we qualitatively compare the solutions computed by FOM and ROM for the three test values of $\gamma$ at $t = 10$ in 
Fig.~\ref{fig:errors_psi_absolute_gamma}
and \ref{fig:errors_zeta_absolute_gamma}. Once again, we see that our ROM approach provides a good global reconstruction of both stream function and vorticity. 
In fact, the maximum relative difference in absolute value does not exceed $6.3e-3$ for $\psi$ and $6.1e-2$ for $\omega$. 
%\subsubsection{Parametrization with respect to the Reynolds number}

%\michele{(Anna, fai un check di tutti i numeri scritti nel testo, i valori degli errori ecc..grazie! :) )}

%\subsubsection{Parametrization with respect to forcing term}

\begin{comment}

figura 11: errori relativi per Re (500 e 1000)
FIgura 12: confronto fom/rom per psi per Re (uno dei due)
fifura 13: confronto fom/rom per zeta per Re (uno dei due)
figura 14, 15 e 16: stesse cose per gamma...
gli autovalori forse no??? perchè piu' o meno sono gli stessi...li enunciamo solo!!!
\end{comment}
\section{Conclusions and future developments}\label{sec:conc}
This work presents a POD-Galerkin based Reduced Order Method 
for the Navier--Stokes equations in stream function-vorticity  formulation within a Finite Volume environment. % For the numerical employment, we choose a Finite Volume method because of its computational efficiency. T
The main novelties of the proposed ROM approach are: (i) the use of different coefficients to approximate the stream function and vorticity fields and (ii) the use of a global POD basis space for parametric studies. 
We assessed our ROM approach with the vortex merger problem,
a classical benchmark used for the validation of numerical methods for the stream function-vorticity formulation of the Navier--Stokes equations. 
The numerical results show that our ROM is able to capture the flow features with a good accuracy both in the reconstruction of the flow field evolution and in a physical parametric setting. 
In addition, for the simple vortex merger problem we observed that our ROM enables substantial computational time savings. 
%an accuracy comparable to other ROMs
%applied to similar benchmarks in [53, 55, 57]. 

Next, we are going to extend the ROM approach presented here to the Quasi-Geostrophic Equations. In particular, we intend to couple such equations with a differential filter \cite{Girfoglio2021a, Girfoglio2021b,Girfoglio2019,Girfoglio_JCP,Girfoglio2021c} in order to simulate two-dimensional turbulent geophysical flows on under-refined meshes in the spirit of \cite{Holm2003, Monteiro2015}.

\section*{Acknowledgements}
We acknowledge the support provided by the European Research Council Executive Agency by the Consolidator Grant project AROMA-CFD ``Advanced Reduced Order Methods with Applications in Computational Fluid Dynamics" - GA 681447, H2020-ERC CoG 2015 AROMA-CFD, PI G. Rozza, and INdAM-GNCS 2019-2020 projects.
This work was also partially supported by US National Science Foundation through grant DMS-1953535. A.~Quaini acknowledges 
support from the Radcliffe Institute for Advanced Study at Harvard University where she has been a 2021-2022 William and Flora Hewlett Foundation Fellow.

\bibliographystyle{plain}
\bibliography{biblio.bib}

\end{document}